\documentclass[10pt,a4paper]{book}
\usepackage{latexsym}
\usepackage{amsbsy}
\usepackage{amssymb}
\usepackage{amsmath}
\usepackage{a4}
\usepackage{makeidx}
\usepackage[all]{xy}
\usepackage[latin1]{inputenc} 
\usepackage{cancel}
\usepackage[mathscr]{euscript}
\usepackage[hidelinks]{hyperref}
\usepackage{verbatim}

\usepackage{graphicx}

\usepackage{amscd}

\numberwithin{equation}{section}

\newcommand{\beqs}{\begin{displaymath}}
\newcommand{\eeqs}{\end{displaymath}}

\newcommand{\ep}{\hspace*{\fill}$\Box$\medskip\\} 
\newcommand{\eps}{\varepsilon}
\newcommand{\pr}{\noindent{\bf Proof. }} 

\newcommand{\R}{\mathbb R}
\newcommand{\N}{\mathbb N}

\newcommand{\Q}{\mathbb Q}

\newcommand{\comp}{\Subset} 
\newtheorem{thr}{\hspace*{-3mm} \bf}[section] 
\newcommand{\bt}{\begin{thr}{\bf Theorem. }} 
\newcommand{\et}{\end{thr}} 
\newcommand{\bp}{\begin{thr}{\bf Proposition. }} 
\newcommand{\bc}{\begin{thr} {\bf Corollary. }}
\newcommand{\blem}{\begin{thr} {\bf Lemma. }} 
\newcommand{\bex}{\begin{thr}{\bf Example. }\rm}
\newcommand{\bexs}{\begin{thr} {\bf Examples. }\rm}
\newcommand{\bd}{\begin{thr} {\bf Definition. }}
\newcommand{\brem}{\begin{thr} {\bf Remark. }\rm}
\newcommand{\beast}{\begin{eqnarray*}} 
\newcommand{\eeast}{\end{eqnarray*}}
\newcommand{\vphi}{\varphi}
\newcommand{\mc}{\mathcal}

\newcommand{\mcu}{{\mathcal U}}

\newcommand{\mcn}{{\mathcal N}}

\newcommand{\la}{\langle}
\newcommand{\ra}{\rangle}
\newcommand{\pa}{\partial}
\newcommand{\rk}{\mathrm{rk}}

\newcommand{\id}{\mathrm{id}}

\newcommand{\cinfty}{{\mathcal C}^\infty}
\newcommand{\proj}{\mathrm{pr}}

\newcommand{\supp}{\mathrm{supp}}

\newcommand{\X}{{\mathfrak X}}

\newcommand{\Fl}{\mathrm{Fl}}

\newcommand{\grad}{\mathrm{grad}}

\newcommand{\g}{{\mathfrak g}}
\newcommand{\Right}{{\mathscr R}}
\newcommand{\ddt}{\left.\frac{d}{d t}\right|_0}
\newcommand{\dde}{\left.\frac{d}{d \eps}\right|_0}
\newcommand{\pdl}[3]{\left.\frac{\partial #1}{\partial #2}\right|_{#3}}
\newcommand{\pd}[2]{\frac{\partial #1}{\partial #2}}
\newcommand{\sse}{\subseteq}
\newcommand{\Xloc}{{\mathfrak X}_{\mathrm{loc}}(M)}
\newcommand{\cU}{{\mathcal U}}
\newcommand{\cD}{\mathcal D}
\newcommand{\set}{\text{set}}
\newcommand{\mcS}{{\mathcal S}}
\newcommand{\xtil}{{\tilde x}}
\newcommand{\util}{{\tilde u}}
\newcommand{\ftil}{{\tilde f}}
\newcommand{\prol}{\mathrm{pr}^{(n)}}
\newcommand{\prolm}{\mathrm{pr}^{(n-1)}}
\newcommand{\prolp}{\mathrm{pr}^{(n+1)}}
\newcommand{\prolo}{\mathrm{pr}^{(1)}}
\newcommand{\prolt}{\mathrm{pr}^{(2)}}
\newcommand{\prolmm}{\mathrm{pr}^{(m)}}
\newcommand{\prolmmp}{\mathrm{pr}^{(m+1)}}
\newcommand{\un}{u^{(n)}}
\newcommand{\um}{u^{(m)}}

\newcommand{\unp}{u^{(n+1)}}
\newcommand{\hn}{^{(n)}}
\newcommand{\hm}{^{(m)}}
\newcommand{\hk}{^{(k)}}
\newcommand{\hnp}{^{(n+1)}}
\newcommand{\hnm}{^{(n-1)}}
\newcommand{\heins}{^{(1)}}
\newcommand{\hzwei}{^{(2)}}
\newcommand{\al}{\alpha}
\newcommand{\mcl}{{\mathcal L}}
\newcommand{\D}{{\mathcal D}}
\newcommand{\Om}{\Omega}
\newcommand{\tr}{\mathrm{tr}}
\newcommand{\ddee}{\frac{d}{d \eps}}
\newcommand{\Div}{\mathrm{Div}\,}
\newcommand{\be}{\beta}
\newcommand{\sump}{\sum_{j=1}^p}
\newcommand{\sumq}{\sum_{\al=1}^q}
\newcommand{\tof}{{\mathcal T}_\Om[f]}
\newcommand{\bfx}{\mathbf{x}}
\newcommand{\bfa}{\mathbf{a}}
\newcommand{\lam}{\lambda}
\newcommand{\od}[2]{\frac{d #1}{d #2}}

\newcommand{\odk}[2]{\frac{d^k #1}{d #2^k}}

\makeindex

\begin{document}

\begin{titlepage}
\begin{center}

{\bf \huge Lie Transformation Groups 
} \\[1em] 
{\Large \bf An Introduction to Symmetry Group Analysis of Differential Equations}\\[2em] 
{\Large Summer Term 2015}\\[2em]
{\Large Michael Kunzinger} \\[1em]
{\tt michael.kunzinger@univie.ac.at}\\
Universit\"at Wien \\
Fakult\"at f\"ur Mathematik\\
Oskar-Morgenstern-Platz 1\\
A-1090 Wien

\end{center}
  
\end{titlepage}

\newpage
\
\newpage

\pagenumbering{roman}
\setcounter{page}{1}

\chapter*{Preface}
\addcontentsline{toc}{part}{Preface}
These lecture notes cover the material of a four hour course on Lie
transformation groups. The goal of this course is to provide the
foundations of symmetry group analysis of differential equations,
a vast and very active research area with connections to 
Lie group theory, differential geometry, differential equations,
calculus of variations, integrable systems, and mathematical
physics, among others. The standard reference in this field is 
P.\ J.\ Olver's seminal work \cite{O}, and to a large part this
set of notes will closely follow his exposition. We do, however,
take a slightly different approach to the foundations of the
field and try to give a reasonably complete exposition of the 
underlying theory of local transformation groups, which in 
turn is based on H.\ Sussmann's study \cite{Sussmann} of the integrability of
distributions of non-constant rank. For connecting this material
with Olver's approach, we have found V.\ Schmidt's doctoral 
thesis \cite{S} very helpful.

I have tried to make these notes self-contained, pre-supposing only
that the reader be familiar with the material covered in basic
courses on Differential Geometry and Lie Groups. More precisely,
the presentation is a direct continuation of my courses \cite{KAoM}
and \cite{KLie_new} (which in turn is mostly based on \cite{BC}).
In particular, these notes supply proofs for many results that are
left to the reader in \cite{O}. Clearly what we can provide here
is only a first step into the field, but it should allow the 
reader to delve into the subject, e.g., by
further studying \cite{O}.  

I am greatly indebted to Roman Popovych for carefully reading the entire
manuscript and for suggesting numerous improvements.
\vskip1cm

\hspace*{\fill} Michael Kunzinger, summer term 2015

\vskip1cm
P.S. (Dec.\ 2024): Since the first version of these lecture notes, Lavau's article \cite{Lavau} has appeared, pointing out 
some wrong claims in the original literature on the integrability of singular distributions, one of which
I also had overlooked. This has now been rectified (cf., in particular, \ref{sussmain2}, \ref{th:local_lie_generated}).
I thank Theresa P\"oll, whose Master's thesis \cite{Poell} was of great help in clarifying the situation.

\
\newpage

\tableofcontents

\newpage


\pagenumbering{arabic}
\setcounter{page}{1}
\pagestyle{plain}

\chapter{Lie transformation groups} \label{ltg_chapter}
\section{Basic concepts}
Throughout this course, by a manifold we mean a smooth ($C^\infty$) manifold.
However, we do not a priori assume further properties of the natural manifold topology like Hausdorff
or second countable. In this chapter we will follow \cite{BC}.

We begin by recalling \cite[Def.\ 16.1]{KLie_new}:
\bd\label{transgroup} A {\em transformation}\index{transformation} of a manifold $M$ is a diffeomorphism from $M$ onto $M$.
A group $G$ acts on $M$ as a transformation group\index{transformation group} (on the left) if there exists a map
$\Phi\colon G\times M \to M$ satisfying: 
\begin{itemize}
\item[(i)] For every $g\in G$ the map $\Phi_g\colon m\mapsto \Phi(g,m)$ is a transformation of $M$.
\item[(ii)] For all $g$, $h\in G$, $\Phi_g\circ \Phi_h = \Phi_{gh}$.
\end{itemize}
In particular, $\Phi_e= \mathrm{id}_M$.

$G$ acts {\em effectively} on $M$ if $\Phi_g(m) = m$ for all $m\in M$ implies $g=e$. $G$ acts
{\em freely} on $M$ if it has no fixed points, i.e., if $\Phi_g(m) = m$ for some $m\in M$ implies 
that $g=e$.
\et
Our main object of study is introduced next.
\bd\label{ltgdef} A Lie group $G$ acts as a {\em Lie transformation group}\index{Lie transformation group} 
on a manifold $M$ if there exists 
a smooth surjection $\Phi\colon G\times M \to M$ such that for all $g$, $h\in G$, $\Phi_g\circ 
\Phi_h = \Phi_{gh}$.
\et
We first note that any Lie transformation group is a transformation group in the sense of 
\ref{transgroup}. To see this, note first that $\Phi$ being surjective implies that
$\Phi_e=\mathrm{id}_M$. In fact, given $m\in M$ there exist $m'\in M$, $h\in G$ with $\Phi_h(m')=m$.
Therefore, $\Phi_e(m) = \Phi_e\circ \Phi_h(m') = \Phi_h(m')=m$. Consequently, $\Phi_g^{-1} = 
\Phi_{g^{-1}}$ for all $g\in G$, so each $\Phi_g$ is a transformation of $M$, and we have
verified all conditions from Def.\ \ref{transgroup}. We will often briefly write $gm$ instead
of $\Phi_g(m)$.
\bex\label{complvf} Let $X$ be a complete vector field on a Hausdorff manifold $M$. Then the flow of $X$ defines
a smooth map $\Psi:=\Fl^X:\R \times M \to M$ with $\Psi(s,\Psi(t,m))$ $=$ $\Psi(s+t,m)$ for all $s$, $t\in \R$
and all $m\in M$. $\Psi$ therefore induces a Lie transformation group of the Lie group $(\R,+)$
on $M$. Each $\Psi_t := m\mapsto \Psi(t,m)$ is a transformation on $M$ and $t\mapsto \Psi_t$ is
a homomorphism of $\R$ into the group of transformations on $M$.
\et
\blem\label{subgr} Let $G$ be a Lie group acting as a Lie transformation group on $M$ and let $H$ be a Lie
subgroup of $G$. Then also $H$ acts as a Lie transformation group on $M$.
\et
\pr Let $j: H \hookrightarrow G$ be the natural injection and let $\Phi\colon G\times M \to M$ be as
in \ref{ltgdef}. Then $\Phi_H: H\times M\to M$, $\Phi_H:=\Phi\circ(j\times \id_M)$ does the trick.
\ep
If $G$ acts as a transformation group on $M$ via $\Phi$ then a subset $A$ of $M$ is called 
{\em invariant}\index{invariant} under $G$ if $\Phi(G\times A) \subseteq A$.
\bp \label{regsub} If a regular submanifold $M'$ of $M$ is invariant under the action of a Lie transformation
group $G$ on $M$ then $G$ acts naturally on $M'$ as a Lie transformation group.
\et
\pr Let $j:G\times M' \hookrightarrow G\times M$ be the natural injection. Then since $M'$ is regular, the map
$\Phi':G\times M' \to M'$ induced by $\Phi\circ j$ is smooth (see \cite[3.3.14]{KAoM}), hence it defines an action of $G$ on $M'$.
\hspace*{\fill}$\Box$
\bex\label{onex}
$GL(n,\R)$ acts on $\R^n$ as a Lie transformation group via $(A,x)\mapsto Ax$. By \ref{subgr},
so does $O(n,\R)$. Since $S^{n-1}$ is a regular submanifold of $\R^n$ invariant under $O(n,\R)$,
\ref{regsub} implies that $O(n,\R)$ acts on $S^{n-1}$ as a Lie transformation group.
\et
For the notion of quotient manifold, cf.\ \cite[Sec.\ 15]{KLie_new}. We say that a transformation group
$G$ {\em preserves}\index{preserving an equivalence relation} an equivalence relation $\rho$ on $M$ 
if $(m_1,m_2)\in \rho$ and $g\in G$ imply
$(gm_1,gm_2)\in \rho$.
\bp If $M/\rho$ is a quotient manifold of $M$ and if the equivalence relation $\rho$ is preserved 
by a Lie transformation group $G$ on $M$ then $G$ acts naturally on $M/\rho$ as a Lie transformation
group.
\et
\pr Denote by $\pi: M \to M/\rho$ the natural surjection. Then $\id\times \pi: G\times M \to G\times (M/\rho)$
is a submersion, so $G\times (M/\rho)$ is a quotient manifold of $G\times M$. Moreover, the smooth map
$\pi\circ \Phi\colon G\times M \to M/\rho$ is an invariant of the corresponding equivalence relation on $G\times M$
since $(g,m_1)\sim (g,m_2)$ if and only if $(gm_1,gm_2)\in \rho$, which implies $\pi\circ \Phi(g,m_1)
= \pi(gm_1) = \pi(gm_2) = \pi\circ \Phi(g,m_2)$. By \cite[15.13]{KLie_new} it therefore projects to a smooth
map
$$
\Phi_\rho: G\times (M/\rho) \to M/\rho, \qquad (g,\pi(m)) \mapsto \pi(gm)
$$
which defines the required action of $G$ on $M/\rho$ since
\begin{equation*}
\begin{split}
\Phi_\rho(g_1,\Phi_\rho(g_2,\pi(m))) &= \Phi_\rho(g_1,\pi\circ\Phi(g_2,m)) = 
\pi\circ \Phi(g_1,\Phi(g_2,m)) \\
&= \pi\circ \Phi(g_1g_2,m) = \Phi_\rho(g_1g_2,\pi(m)).
\end{split}
\end{equation*}
\hspace*{\fill}$\Box$
\bex \label{gmodh} Any Lie group $G$ acts on itself as a Lie transformation group via multiplication $\mu\colon G\times G\to G$.
A subgroup $H$ of $G$ determines an equivalence relation $\rho$ on $G$ by 
$$
(g_1,g_2)\in \rho :\Leftrightarrow \exists h\in H: g_1 = g_2h
$$
and this relation is preserved by the transformation group. If the set of orbits $G/H$ 
(i.e., the set of left cosets) is a quotient manifold of $G$ (cf.\ \cite[Sec.\ 15, 20]{KLie_new})
then $G$ acts naturally on $G/H$ as a Lie transformation group with action $\Phi$ given 
by $(g,aH) \mapsto (ga)H$. A sufficient condition for $G/H$ to be a quotient manifold of $G$
is that $H$ is closed, not open, and connected (see \cite[20.5]{KLie_new}).
\et
\bp If $H$ is a closed subgroup of a Lie group $G$ such that $G/H$ is a quotient manifold of $G$, then 
$G/H$ is a Hausdorff manifold.
\et
\pr Recall from \cite[15.9]{KLie_new} that the natural manifold topology of $G/H$ is the quotient topology with
respect to the natural projection $\pi: G \to G/H$. Suppose first that $gH\in G/H$ is distinct from $H$.
Then $gH$ is a closed subset of $G$ that does not contain $e$. Hence there exists a neighborhood $U$ of $e$
in $G$ which is disjoint from $gH$. Pick a neighborhood $V$ of $e$ in $G$ such that $V^{-1}V\subseteq U$.
We show that $\pi(Vg)$ and $\pi(V)$ are disjoint neighborhoods of $gH$ and $H$, respectively. Indeed,
if $\pi(Vg)\cap \pi(V)\not=\emptyset$ then there exist $a$, $b\in V$ such that $agH=bH$, so $a^{-1}b\in gH$.
But $a^{-1}b\in V^{-1}V\subseteq U$, so we obtain a contradiction.

In general, if $gH$ and $g'H$ are distinct elements of $G/H$ then so are $H$ and $(g^{-1}g')H$ and by
what we have shown above they possess disjoint neighborhoods $W$ and $W'$. Denoting by $\phi$ the action
of $G$ on $G/H$ from \ref{gmodh}, it follows that $\Phi_g(W)$ and $\Phi_g(W')$ are disjoint neighborhoods
of $gH$ and $g'H$.
\ep
Recall from \cite[Sec.\ 4, 6]{KLie_new} that a vector field $X\in \X(G)$ is called left- (resp.\ right)-invariant if 
$L_g^*X=X$ (resp.\ $R_g^*X=X$) for each $g\in G$. We denote by $\X_L(G)$ (resp.\ $\X_R(G)$) the space of 
left- (resp.\ right)-invariant vector fields on $G$. Also, we identify $\X_L(G)$ with the Lie algebra $\g$ 
via $T_eG \ni v \mapsto L^v\in \X_L(G)$, where $L^v(g):=T_eL_g(v)$. The inverse of this map is $X\mapsto X(e)$.

Similarly, $\g$ is isomorphic to 
$\X_R(G)$ via $T_eG \ni v \mapsto R^v\in \X_R(G)$ with $R^v(g):=T_eR_g(v)$. We usually equip $\g$
with the Lie algebra structure induced by $[v,w]:=[L^v,L^w](e)$. Alternatively, we may set $[v,w]_R:=[R^v,R^w](e)$.
The resulting Lie algebra will be denoted by $\g_R$. We have:
\bp\label{leftright} Let $G$ be a Lie group. Then
\begin{itemize}
\item[(i)] The map $F: L^v \mapsto - R^v$ is a Lie algebra isomorphism from $\X_L(G)$ onto $\X_R(G)$.
\item[(ii)] The map $f\colon  v\mapsto -v$ is a Lie algebra isomorphism from $\g$ onto $\g_R$. 
\end{itemize}
\et
\pr (i) By what was said above, $F$
is a linear isomorphism, since both $L^v \mapsto v\in \g$ and $v\mapsto R^v$ are. Furthermore,
by \cite[6.1]{KLie_new} we have:
\begin{equation*}
[F(L^v),F(L^w)] = [R^v,R^w] = R^{-[v,w]} = -R^{[v,w]} = F(L^{[v,w]}) = F([L^v,L^w]),
\end{equation*}
so $F$ is a Lie algebra isomorphism.\medskip

(ii) We have $[v,w]_R=[R^v,R^w](e)=R^{-[v,w]}(e)=-[v,w]$. Hence 
$$
[f(v),f(w)]_R = [-v,-w]_R = [v,w]_R = - [v,w] = f([v,w]).
$$
\ep
\brem By \cite[6.1]{KLie_new} The map $F$ from \ref{leftright} can also be written as $X\mapsto \nu^*X$. 
In fact,
if $X\in \X_L(G)$ then since $\nu\circ R_g = L_{g^{-1}}\circ \nu$ we have
$$
R_g^*(\nu^*X) = (\nu\circ R_g)^*X = (L_{g^{-1}}\circ \nu)^*X = \nu^*(L_{g^{-1}}^*X) = \nu^*X,
$$
so $\nu^*X\in \X_R(G)$, and by symmetry $\nu^*:\X_L(G)\to \X_R(G)$ is an isomorphism. Moreover,
$\nu^*(L^v)(e) = (T_e\nu)^{-1} \circ L^v\circ \nu(e) = - v$, so $\nu^*(L^v)=R^{-v}=F(L^v)$.
\et 
Let now $\Phi\colon G\times M \to M$ be a Lie transformation group and let $m\in M$. Then the map
\begin{equation}\label{phigdef}
\begin{split}
\Phi_m: G &\to M\\
\Phi_m(g) &:= \Phi(g,m)
\end{split}
\end{equation}
is smooth. Given any $v\in T_eG = \g$, let 
\begin{equation}\label{phidef}
\begin{split}
\Phi(v): M &\to TM\\
\Phi(v)(m) &:= T_e\Phi_m(v)
\end{split}
\end{equation}
Then $\Phi(v)(m) = T\Phi(v,0_m)$, so $\Phi(v)$ is smooth. Also, with $\pi:TM\to M$ and
$\tilde \pi:T(G\times M)\to G\times M$ the canonical projections,
$\pi\circ \Phi(v)(m) = \Phi(\tilde\pi(v,0_m)) = \Phi(e,m) = m$, so $\Phi(v)$ is a smooth section of $TM$, i.e.,
$\Phi(v)\in \X(M)$. We set
\begin{equation}\label{rightvf}
\Right(G,M) := \{\Phi(v) \mid v\in \g\}.
\end{equation}
$\Right(G,M)$ is called the {\em Killing algebra}\index{Killing algebra} of $\Phi$.
Since $v\mapsto \Phi(v)$ is linear, it is a finite-dimensional vector space.
In fact, we even have:
\bp\label{psihomo} The map $\Phi\colon v\mapsto \Phi(v)$ is a Lie algebra homomorphism from $\g_R$ onto $\Right(G,M)$.
\et
\pr Given any $v\in \g$, we first show that for any $m\in M$, the right-invariant vector field $Y$ on 
$G$ with $Y(e)=v$ is $\Phi_m$-related to $\Phi(v)$. In fact, noting that $\Phi_m\circ R_g(h) = \Phi_m(hg)
= \Phi(h,\Phi_m(g)) = \Phi_{m'}(h)$ with $m':=\Phi_m(g)$, we obtain:
$$
T_g\Phi_m(Y(g)) = T_g\Phi_m(T_eR_g(v)) = T_e(\Phi_m\circ R_g)(v) = T_e\Phi_{m'}(v) = \Phi(v)(\Phi_m(g)).
$$
Thus for any $v\in \g$,
\begin{equation}\label{olver1.46.2}
T_g\Phi_m(Y(g)) = \Phi(v)(\Phi_m(g)).
\end{equation}
From this, by \cite[4.4]{KLie_new} we conclude that for $Z=R^w$ another right-invariant vector field we have
$$
[\Phi(v),\Phi(w)]\circ \Phi_m = T\Phi_m\circ[Y,Z].
$$
Inserting $e$, this gives 
$$
[\Phi(v),\Phi(w)]_m = T_e\Phi_m([Y,Z]_e) = T_e\Phi_m([v,w]_R)
$$
(where we used that $[R^v,R^w]=R^{[v,w]_R}$, cf.\ \ref{leftright}). Thus, finally,
$[\Phi(v),\Phi(w)] = \Phi([v,w]_R)$, as claimed.
\ep
In the following result, we will start denoting the time variable in the flow map of a
vector field with $\eps$, for the sake of compatibility with the case of symmetry groups
of differential equations later on.
\bp\label{flowaction} If $M$ is $T_2$, then every vector field $\Phi(v)\in \Right(G,M)$ is complete, with flow 
$$
\Fl^{\Phi(v)}_\eps(m) = \Phi(\exp(\eps v),m) = \exp(\eps v)\cdot m. 
$$
In particular,
\begin{equation}\label{olver1}
\Phi(v)(m) = \dde \Phi(\exp(\eps v),m).
\end{equation}
\et
\pr In the proof of \ref{psihomo} we have seen that $R^v$ is $\Phi_m$-related to $\Phi(v)$. Thus by
\cite[7.2]{KLie_new} it follows that $\Phi_m\circ \Fl^{R^v}_\eps = \Fl^{\Phi(v)}_\eps\circ \Phi_m$. Moreover, by
\cite[8.2]{KLie_new}, $\Fl^{R^v}_\eps(g) = \exp(\eps v)\cdot g$. Therefore,
$$
\Fl^{\Phi(v)}_\eps(m) = \Fl^{\Phi(v)}_\eps\circ \Phi_m(e) = \Phi_m(\exp(\eps v)\cdot e) = \Phi(\exp(\eps v),m).
$$
Completeness of $\Phi(v)$ is immediate from this description.
\ep
\bex\label{rflowex} Let $\Phi\colon \R\times M \to M$ be a Lie transformation group with Lie group $(\R,+)$ and $T_2$-manifold $M$.
Then $\Right(\R,M)$ is generated by $\Phi(\dde)$.
Moreover, the flow of $\Phi(\dde)$ is given by $\Phi$ itself:
$\Fl^{\Phi(\dde)}_\eps(m) = \Phi(\eps,m)$: this is immediate from \ref{flowaction}
since $\exp(\eps\dde) = \eps$.

\et
Recall from \cite[Sec.\ 16]{KLie_new} that, given a Lie transformation group $\Phi\colon G\times M \to M$, we set
$$
K:=\{k\in G \mid \Phi_k = \id_M\}.
$$
Then $K$ is a normal subgroup of $G$ and \cite[16.3]{KLie_new} demonstrates that the quotient group
$G/K$ acts effectively on $M$. Moreover, $K$ is closed in $G$ since we can write $K=\bigcap_{m\in M} \Phi_m^{-1}(m)$.
Hence by \cite[21.6]{KLie_new}, if $K$ is not open then the quotient manifold $G/K$ is a Lie group. The
exceptional case here is not very interesting though:
\brem If $K$ is open, then $\Right(G,M)=0$. Indeed, if $K$ is open (and closed) then since $e\in K$ we must have $G_e\subseteq K$.
Then because for all $v\in \g$, $\exp(\eps v)\in G_e$ it follows that $\Phi_m(\exp(\eps v)) = m$ for all $m\in M$, and
differentiating with respect to $\eps$ at $\eps=0$ implies that $\Phi(v)(m)=T_e\Phi_m(v) = 0$ for all $m$. Thus
$\Phi(v)=0$ and thereby $\Right(G,M)=0$.
\et
\bp If $K$ is not open then $G/K$ acts effectively on $M$ as a Lie transformation group and $\Right (G/K,M)=
\Right(G,M)$.
\et
\pr Let 
\begin{equation*}
\begin{split}
\Phi': (G/K) \times M &\to M\\
(gK,m) &\mapsto \Phi(g,m)
\end{split}
\end{equation*}
To see that $\Phi'$ is well-defined and smooth let $\pi: G\to G/K$ be the natural projection. 
Then $\pi\times \id_M: G\times M \to (G/K)\times M$ is a submersion (cf.\ \cite[20.5]{KLie_new})
and $\Phi$ is an invariant of the corresponding equivalence relation on $G\times M$: if
$(g_1,m_1)\sim (g_2,m_2)$ then $g_1K=g_2K$ and $m_1=m_2$, so $\Phi(g_1,m_1) = \Phi'(g_1K,m_1)
= \Phi'(g_2K,m_2) = \Phi(g_2,m_2)$. Thus by \cite[15.13]{KLie_new}, the map $\Phi'$, being the
projection of $\Phi$, is itself smooth. That $\Phi'$ gives an effective action has been shown in 
\cite[16.3]{KLie_new}.

Finally, for any $m\in M$ we have $\Phi_m = \Phi'_m\circ \pi$, so for any $v\in \g$ we get
$T_e\Phi_m(v) = T_e\Phi'_m \circ T_e\pi(v)$, i.e., $\Phi(v)=\Phi'(T_e\pi(v))$. Since $T_e\pi$
is a surjection, this gives $\Right (G/K,M)=\Right(G,M)$.
\ep
\bp\label{phiiso} Let $\Phi$ be a transformation group that acts effectively on a $T_2$-manifold $M$. Then the map
$$
\Phi\colon \g_R \to \Right(G,M)
$$
defined by \eqref{phidef} is a Lie algebra isomorphism.
\et
\pr Using \ref{psihomo} it only remains to show that $\Phi$ is injective. Thus let $\Phi(v)=0$. Then by 
\ref{flowaction} it follows that $\exp(\eps v)\cdot m = m$ for all $m\in M$ and all $\eps\in \R$. Since $G$ acts effectively we 
conclude that $\exp(\eps v)= e$ for all $\eps$, giving $v=0$.
\ep
\bex Let $G = GL(n,\R)$ and $M=\R^n$ with $\Phi\colon (A,x)\mapsto A\cdot x$. To determine the map \eqref{phidef}
we let $A^{ij}$ be the standard coordinates on $GL(n,\R)$ and $x^k$ those on $\R^n$. Then
\begin{equation*}
\begin{split}
T_e\Phi_x\Big(\frac{\partial}{\partial A^{ij}}\Big) &= \sum_k \Big(\frac{\partial}{\partial A^{ij}}(x^k\circ \Phi_x)\Big)
\Big(\left. \frac{\partial}{\partial x^k}\right|_x\Big) \\
& = \sum_{k,r}  \Big(\pd{}{A^{ij}}(A^{kr}x_r)\Big)\pdl{}{x^k}{x} = x_j \pdl{}{x^i}{x},
\end{split}
\end{equation*}
so $\Phi = [B^{ij}] \mapsto \sum_{i,j} B^{ij}x^j \pd{}{x^i}$. Since $G$ act effectively, $\Phi$ is an
isomorphism.
\et
\brem There is a strong converse of \ref{phiiso}, due to R.\ Palais (see \cite[p.\ 95]{Palais}): {\em If 
$\mathscr A$ is a finite-dimensional Lie algebra of complete vector fields on a $T_2$-manifold $M$, then
there exists a connected Lie group $G$ which acts effectively on $M$ as a Lie transformation group and
such that $\mathscr{A} = \Right(G,M)$.}
\et
\section{Orbits under a Lie transformation group}
From \cite[Sec.\ 16]{KLie_new} we know that any transformation group $\Phi$ acting on a manifold $M$ induces an 
equivalence relation on $M$: $m\sim m'$ iff there exists some $g\in G$ with $m'=gm$. The equivalence class
of any $x\in M$ is called the {\em orbit}\index{orbit} of $m$ under $G$. It is the range of the map
$\Phi_m = g\mapsto \Phi(g,m)$, $G\to M$.
\bd Let $\Phi\colon G\times M \to M$ be a transformation group on a manifold $M$. For any $m\in M$, the subgroup 
$G_m:=\Phi_m^{-1}(m)=\{g\in G\mid g\cdot m = m\}$ is called the {\em isotropy group}\index{isotropy group} 
of $m$.
\et
\brem\label{isotropyrem} 
(i) Isotropy groups at equivalent points in $M$ are conjugate subgroups of $G$. In fact, suppose that
$m'=gm$. Then $(gG_mg^{-1})m' = m'$, so $gG_mg^{-1}\subseteq G_{m'}$, and analogously $g^{-1}G_{m'}g\subseteq G_m$,
so $gG_mg^{-1} = G_{m'}$.

(ii) For any $m\in M$, the map $\Phi_m: G\to M$ projects to a map $\Psi_m:G/G_m \to M$ defined by
$gG_m \mapsto gm$. The range of $\Psi_m$ is the orbit of $m$. Also, $\Psi_m$ is injective: if
$\Psi_m(g_1G_m) = \Psi_m(g_2G_m)$ then $g_1m=g_2m$, so $g_1^{-1}g_2\in G_m$, i.e., $g_1G_m = g_2G_m$.
\et
\brem\label{genltg}
Let $\Phi\colon G\times M\to M$ be a Lie transformation group on a $T_2$-manifold $M$. Then any $G_m$ is a
closed subgroup of $G$. Hence by \cite[21.7]{KLie_new}, $G_m$ is either discrete or it admits a unique
structure as a (regular) submanifold of $G$. In the latter case it is also a Lie subgroup of $G$.
If $m'=gm$ then by \ref{isotropyrem} (i), $G_m$ is mapped onto $G_{m'}$ by the diffeomorphism $L_g\circ R_{g^{-1}}$.
Thus the isotropy groups at points of an orbit are either all discrete or are regular submanifolds and
Lie subgroups of $G$ that are pairwise diffeomorphic (since $G_m$, $G_{m'}$ are regular submanifolds, the restriction
of $L_g\circ R_{g^{-1}}$ is also a diffeomorphism from $G_m$ onto $G_{m'}$).
\et
If $G_m$ is open then it must be a union of connected components of $G$, which, by \cite[2.4]{KLie_new}
are precisely the cosets of the normal subgroup $G_e$. If, for example, $G=g_1G_e\cup\dots\cup g_kG_e$
 and $G_m = g_1G_e\cup\dots\cup g_lG_e$, then 
$G/G_m = \{g_{1}G_m,\dots,g_kG_m\}$, since for $g\in g_jG_e$ we have
$$
\pi(g)=g\cdot \bigcup_{i=1}^l g_iG_e = \bigcup_{i=1}^l g G_e g_i = \bigcup_{i=1}^l g_jG_e g_i = g_j G_m.
$$
Also, the orbit of $m$ only consists of finitely many points
(namely $\Psi_m(G/G_m)=\{g_{1}\cdot m,\dots,g_k\cdot m\}$).
Otherwise, we have:
\bt\label{ltgiso} Let $\Phi\colon G\times M \to M$ be a transformation group on a $T_2$-manifold $M$ and let $m\in M$. 
If the isotropy group $G_m$ of $m$ is not open in $G$ then the map $\Psi_m$ from \ref{isotropyrem} (ii)
is an injective immersion of the quotient manifold $G/G_m$ into $M$.
\et
\pr Since $\Phi_m = \Psi_m\circ \pi$ (with $\pi:G\to G/G_m$), $\Psi_m$ is the projection of the smooth map 
$\Phi_m$, hence is itself smooth
(see \cite[3.3.9]{KAoM}). Also, $G/G_m$ is a quotient manifold of $G$ by \cite[21.5]{KLie_new}. Since 
$\Psi_m$ is injective by \ref{isotropyrem} (ii), it remains to show that its rank in any point equals 
the dimension of $G/G_m$. Since $\pi$ is a submersion, this is the case if and only if the rank of
$\Phi_m$ is everywhere equal to $\dim(G/G_m)$. We begin by showing that this is true at $e$.

Let $v\in T_eG$ such that $T_e\Phi_m(v)=0$. By \eqref{phidef} this means that the vector field
$\Phi(v)$ has a zero at $m$. Therefore \ref{flowaction} implies that $\Fl_\eps^{\Phi(v)}(m)=\exp(\eps v)m =m$
for all $\eps$, i.e., $\exp(\eps v)\in G_m$ for all $\eps$. Now if $G_m$ is discrete then the image of $\eps\to 
\exp(\eps v)$, being connected, must consists solely of $e\in G$, so $v = \dde\exp(\eps v)=0$. In this
case, then, $T_e\Phi_m$ is injective, and so the rank of $T_e\Phi_m$ equals the dimension of $G$, and
thereby the dimension of $G/G_m$.

If $G_m$ is non-discrete then by \ref{genltg} it is a regular submanifold of $G$. Hence
$\eps\mapsto \exp(\eps v)$ is smooth as a map into $G_m$ (see \cite[3.3.14]{KAoM}), and so 
$v = \dde\exp(\eps v) \in T_eG_m$. Altogether, we obtain that $\ker(T_e\Phi_m) \subseteq T_eG_m$.
Conversely, $\Phi_m$ is constant on $G_m$, so $T_e\Phi_m|_{T_eG_m} \equiv 0$, hence
in fact $\ker(T_e\Phi_m) = T_eG_m$. Consequently, using \cite[21.7]{KLie_new} we obtain
$$
\rk(T_e\Phi_m) = \dim G - \dim G_m = \dim G/G_m.
$$
Finally, if $g$ is an arbitrary point in $G$ then $\Phi_m\circ R_g = \Phi_{m'}$, where $m'=gm$. 
Then since $R_g$ is a diffeomorphism we have
$$
\rk_g(\Phi_m) = \rk_e(\Phi_{m'}) = \dim G/G_{m'}.
$$
Now by \ref{genltg} $G_m$ and $G_{m'}$ are either both discrete or they are diffeomorphic,
so we conclude that the rank of $\Phi_m$ equals $\dim G/G_{m}$ for every $m\in M$.
\ep
\bc\label{ltgisocor} Under the assumptions of \ref{ltgiso}, the orbit $G\cdot m$ of any $m\in M$ 
can be endowed with the structure of
an immersive submanifold of $M$ diffeomorphic to $G/G_m$.
\et
\pr For clarity, we write $\tilde \Psi_m$ for $\Psi_m$, viewed as a (bijective) map from $G/G_m$ to $G\cdot m$.
Declaring $\tilde \Psi_m$ to be a diffeomorphism provides $G\cdot m$ with a differentiable structure
with respect to which the inclusion map $j:G\cdot m \hookrightarrow M$ is an immersion since $j\circ \tilde\Psi_m = 
\Psi_m: G/G_m \to M$ is an immersion.
\ep
\brem\label{conrem} Suppose that $G$ is connected and let $m\in M$. 
By \ref{ltgiso} the orbit $G\cdot m$ can be discrete only if $G_m$ is open. In this case, $G_m$ is
open and closed in $G$, so $G_m=G$ and therefore $G\cdot m = \{m\}$.
\et
\bex\label{rotex} Let $M=\R^2$, $G=(\R,+)$, and $X=-y\partial_x+x\partial_y$ the rotation vector field on $M$.
Then we obtain a Lie transformation group $\Phi$ by
$$
\Phi\colon (\eps,(x,y))\mapsto \Fl^X_\eps(x,y) = (x\cos\eps - y\sin\eps, x\sin\eps + y\cos\eps).
$$
The orbits of $\Phi$ are the regular submanifolds $\{x^2+y^2=\text{const}\}$ and the singleton $\{(0,0)\}$.
\et
Recall that $G$ is said to act transitively\index{transitive group action} on $M$ if for any $m$, $m'\in M$ 
there exists some $g\in G$ with $gm=m'$. Such a group action possesses only a single orbit, namely the
manifold $M$ itself. By \ref{isotropyrem} (i) this means that any $\Psi_m$ is a bijection from $G/G_m$ onto $M$.
To further elaborate on this, we will need the following auxiliary result:
\blem \label{immerlem} Let $M^m$ and $N^n$ be manifolds and suppose that $M$ is second countable and that $m<n$.
Then an immersion $\Psi: M\to N$ cannot be onto any open subset $W$ of $N$. 
\et
\pr It suffices to take for $W$ the domain of 
a chart $\chi$ in $N$. Let $p$ be any point in $\Psi^{-1}(W)$ and choose charts
$(U,\vphi=(x^1,\dots,x^m))$ around $p$ in $M$ and $(V,\eta=(y^1,\dots,y^n))$ around 
$\Psi(p)$ such that $\Psi(U)\subseteq V\subseteq W$ and such that
$\eta\circ\Psi\circ \vphi^{-1} = (x^1,\dots,x^m) \to (x^1,\dots,x^m,0,\dots,0)$.
Then $\eta(\Psi(U))\subseteq \R^m\times \{0\} \subseteq \R^n$, so it has Lebesgue
measure $0$. It follows that also the image $\chi(\Psi(U))$ of this set under
the smooth map $\chi\circ \eta^{-1}$ has Lebesgue measure $0$.

As $p$ varies in $\Psi^{-1}(W)$, the domains $U$ cover the set $\Psi^{-1}(W)$. As $M$ is
second countable we may extract a countable subcover $\{U_k\mid k\in \N\}$ from this collection.
Then the sets $\Psi(U_k)$ cover $W\cap \Psi(M)$, entailing
$$
\chi(W\cap \Psi(M)) = \bigcup_{k\in \N} \chi(\Psi(U_k)).
$$
But then $\chi(W\cap \Psi(M))$ has Lebesgue measure $0$ and so it cannot be all of $\chi(W)$.
It follows that $\Psi(M)$ cannot contain all of $W$.
\ep
Using this, we can prove:
\bp\label{diffeotrans} Let $G$ be a second countable Lie group that acts transitively 
as a Lie transformation group on the $T_2$-manifold 
$M$. Then $M$ is diffeomorphic to $G/G_m$, for any $m\in M$.
\et
\pr Fix any $m\in G$. If $G_m$ is open, then so is any $gG_m$, hence any point in $G/G_m$ (because
$\pi^{-1}(\pi(g)) = gG_m$), which is
therefore discrete (cf.\ \cite[Rem.\ 20.1]{KLie_new}). Otherwise, by \cite[21.5]{KLie_new} $G/G_m$ possesses
a differentiable structure as a quotient manifold of $G$. In both cases, the quotient map $\pi:
G\to G/G_m$ is open and continuous, so also the topology of $G/G_m$ is second countable. Hence if
$G/G_m$ were discrete it would be countable. But $\Psi_m$ is bijective, so this would imply that
$M$ was countable, which is impossible. Hence $G/G_m$ is not discrete, and so it is a quotient 
manifold of $G$ with a countable basis for its topology. Also, by \ref{ltgiso} $\Psi_m$ is an
injective immersion of $G/G_m$ into $M$. Hence $\dim G/G_m \le \dim M$. Since $\Psi_m$ is onto $M$, \ref{immerlem}
implies that the dimensions in fact are equal. As $\Psi_m:G/G_m\to M$ is an
immersion, it follows that its tangent map is bijective at any point. Thus by the inverse function theorem it is
a local diffeomorphism, hence a global diffeomorphism since it is bijective.
\ep
\bc Let $G$ be a second countable Lie group that acts transitively and freely as a Lie transformation group
on the $T_2$-manifold $M$. Then $M$ is diffeomorphic to $G$.
\et
\pr If $G$ acts freely then $G_m=\{e\}$ for every $m\in M$. Hence $G/G_m=G$ and the result follows from
\ref{diffeotrans}.
\ep
\bc Let $G$ be a compact Lie group that acts transitively as a Lie transformation group on the $T_2$-manifold 
$M$. Then $M$ is compact.
\et
\pr Since $G$ is compact, it is second countable. By \ref{diffeotrans}, $M$ is diffeomorphic to $G/G_m$, for
any $m\in M$. Let $\pi:G\to G/G_m$ be the quotient map. Then $\pi(G)=G/G_m$ is compact, hence so is $M$.
\ep
\bex We continue our study of the action of $O(n,\R)$ on $S^{n-1}$ as a Lie transformation group from
\ref{onex}. This action is given by
\begin{equation*}
\begin{split}
\Phi\colon O(n,\R)\times S^{n-1} &\to S^{n-1}\\
(A,x) &\mapsto Ax.
\end{split}
\end{equation*}
For any vector $v\in S^{n-1}$ we can find an orthogonal matrix $T$ with $v$ as its first column, so
$Te_1=v$. Hence $\Phi$ has only a single orbit, $S^{n-1}$, i.e., it is transitive. The isotropy group
$G_{e_1}$ of $e_1$, i.e., those elements of $O(n,\R)$ that leave $e_1$ unchanged, are the matrices of the form
$$
A = 
\begin{pmatrix}
1 & 0\\
0 & D	
\end{pmatrix}
$$
with $D\in O(n-1,\R)$. By \ref{diffeotrans} we conclude that $S^{n-1}$ is diffeomorphic to the 
quotient manifold $O(n,\R)/G_{e_1}$, i.e., $S^{n-1}\cong O(n,\R)/O(n-1,\R)$.
\et
\section{Groups of transformations on a $T_2$-manifold}
\bd A set $G$ of transformations on a $T_2$-manifold $M$ which is a group under composition
$(g_1,g_2)\mapsto g_1\circ g_2$ $(g_1,\,g_2\in G)$ is called a group of transformations\index{group of transformations} 
of $M$.
\et
According to \ref{transgroup}, the map $\Phi\colon (g,m)\mapsto g(m)\equiv gm$ then defines a 
transformation group on $M$. Moreover, this action obviously is effective. In this section we
want to analyze whether $G$ can be endowed with a Lie group structure such that $\Phi$ is smooth.
The following example shows that such a structure may not be unique in general:
\bex Let $G$ be the group of translations on $\R^2$ and let $\phi_a:= z\mapsto z+a$. Then $G$ is
bijectively mapped onto $\R^2$ by $\phi_a\mapsto a$, and this defines a $\cinfty$-structure on $G$,
for any given $\cinfty$-structure on $\R^2$. But on $\R^2$ there are different such structures,
e.g., the standard one and the one from \cite[Ex.\ 14.3]{KLie_new}. Both of these induce corresponding
structures on $G$ such that $\Phi$ is smooth.
\et
Recall from \ref{complvf} that a complete vector field $X$ on $M$ induces an action of $(\R,+)$ on
$M$ as a Lie transformation group via the flow of $X$, $\Psi_\eps:=m\mapsto \Fl^X_\eps(m)$. 
We call the vector field $X$ 
{\em tangent}\index{vector field!tangent to group of transformations}
to $G$ if all the resulting transformations $\Psi_\eps$ belong to $G$.
\bd \label{unicond} A group of transformations on a $T_2$-manifold $M$ is called a {\em Lie group of 
transformations}\index{Lie group of transformations} if it admits a Lie group structure such that
\begin{itemize}
	\item[(i)] The map $\Phi\colon(g,m)\mapsto gm$ is smooth.
	\item[(ii)] If $X\in \X(M)$ is complete and tangent to $G$ then the group homomorphism $\eps\mapsto \Psi_\eps$
	is a one-parameter subgroup of $G$.
\end{itemize}
\et
The above is the desired property entailing uniqueness:
\bp A group of transformations of a $T_2$-manifold $M$ admits at most one structure as a Lie
group of transformations of $M$.
\et
\pr Let $G$ and $G'$ be two Lie group structures on a group of transformations that both satisfy \ref{unicond}.
We have to show that $\id: G\to G'$ is a diffeomorphism. We first note that by \ref{unicond} (i), $G$ is a 
Lie transformation group on $M$, denoted by $\Phi\colon G\times M \to M$. Hence
by \ref{flowaction} any $v\in T_eG$ defines a complete vector field $\Phi(v)\in \X(M)$ with
$\Fl_\eps^{\Phi(v)}(m)= \exp(\eps v)m$. It follows that $\Phi(v)$ is tangent to $G$. Now
\ref{unicond} (ii) implies that $\eps\mapsto \id\circ \exp(\eps v):\R\to G'$ is a one-parameter subgroup of $G'$.
By \cite[8.3]{KLie_new} there is a unique $v'\in T_eG'$ such that $\id\circ \exp(\eps v) = \exp'(\eps v')$
(with $\exp'$ the exponential map of $G'$). Thereby we obtain a well-defined map $v\mapsto v'$ from
$T_eG$ to $T_eG'$.

Fixing a basis $(v_1,\dots,v_n)$ of $T_eG$, let $\vphi=(x^1,\dots,x^n)$ be the corresponding canonical
chart of the second kind, cf.\ \cite[Rem.\ 8.5]{KLie_new}. Then since $\id$ is a group homomorphism, 
$$
\id\circ(\exp(x^1 v_1)\dots\exp(x^n v_n)) = \exp'(x^1 v_1')\dots\exp(x^n v_n'),
$$
implying that $\id$ is smooth on the domain of $\vphi$. Smoothness at an arbitrary point of
$G$ then follows by writing $\id = L'_g\circ \id \circ L_{g^{-1}}$.
Interchanging $G$ and $G'$ in the above argument shows that $\id$ in fact is a diffeomorphism.
\ep
\brem By \ref{unicond} (i), any Lie group $G$ of transformations is also a Lie transformation group
that acts effectively on $M$. Therefore, \ref{phiiso} implies that $\Right(G,M)$ is Lie algebra-isomorphic
to $\g_R$.
\et
\bp\label{gtransnec} If a group $G$ of transformations of a $T_2$-manifold $M$ admits the structure of a Lie group of
transformations then the set $\X_t(M)$ of complete vector fields on $M$ tangent to $G$ is a finite-dimensional
Lie algebra (namely $\Right(G,M)$).
\et
\pr Denote the action of $G$ on $M$ by $\Phi\colon G\times M\to M$. By \ref{flowaction}, every $\Phi(v)\in\Right(G,M)$
is complete with corresponding transformations $\Psi_\eps: m\mapsto \Fl_\eps^{\Phi(v)}(m)= \exp(\eps v)m$, hence
belonging to $G$. Thus $\Phi(v)$ is a complete vector field on $M$ that is tangent to $G$, i.e., $\Right(G,M)
\sse \X_t(M)$.

Conversely, let $X\in \X_t(M)$ and let $\Psi_\eps:= m\mapsto \Fl^X_\eps(m)$ be the corresponding transformation
of $M$. By \ref{unicond} (ii), $\eps\mapsto \Psi_\eps$ is a one-parameter subgroup of $G$, so by 
\cite[8.3]{KLie_new} there is a unique $v\in T_eG$ with $\Psi_\eps(m) = \exp(\eps v)m$ for all $m\in M$.
Thus the vector fields $X$ and $\Phi(v)$ have the same maximal integral curves and hence coincide,
so $X\in \Right(G,M)$.
\ep
\brem There is in fact also a converse to \ref{gtransnec}, see \cite[p.\ 103]{Palais}: Let $G$ be a group of
transformations of a $T_2$-manifold $M$. If $\X_t(M)$ generates a finite-dimensional Lie algebra
$\mathcal A$ of vector fields on $M$ then $G$ admits the structure of a Lie group of transformations
of $M$ and ${\mathcal A}=\Right(G,M)$.
\et
\chapter{Integrability of distributions of non-constant rank}\label{nonconstchap}
\section{Distributions of non-constant rank}
In \cite[17.32]{KLie_new} we proved the classical Frobenius theorem\index{Frobenius Theorem} on the integrability of 
distributions of constant rank:
\bt\label{classfrob} Let $M$ be an $n$-dimensional $T_2$-manifold and let $\Delta$ be a $k$-dimensional distribution
on $M$. Then the following are equivalent:
\begin{itemize}
	\item[(i)] $\Delta$ is involutive.
	\item[(ii)] Every point in $M$ lies in the domain of a flat chart\index{flat chart} $\vphi=(x^1,\dots,x^n)$ for $\Delta$,
	i.e., such that $\pa_{x^1}|_m,\dots,\pa_{x^k}|_m$ forms a basis of $\Delta(m)$ for each $m$ in the domain of $\vphi$.
	\item[(iii)] Every point of $M$ is contained in an integral manifold of $\Delta$.
	\item[(iv)] Every point $m$ of $M$ lies in a cubic chart $(\vphi=(x^1,\dots,x^n),U)$, $\vphi(U)=[-c,c]^n$ 
	centered around $m$ such that the slices $U_a=\vphi^{-1}(\R^k\times \{a\})$ are integral manifolds 
	of $\Delta$. If $M'$ is a connected integral-manifold of $\Delta$ with $M'\sse U$ then $M'$ is contained
	in one such slice.
\end{itemize}
\et 
For the application to orbits of (local) transformation groups we have in mind, \ref{classfrob} is too
restrictive in that it requires $\Delta$ to have constant rank (i.e., dimension) $k$ everywhere. In this 
chapter we mainly follow H.\ Sussmann's article \cite{Sussmann}, as well as
P.\ Stefan's approach \cite{Stefan1,Stefan2} to develop a theory of integrability for distributions
of non-constant rank. Additional very helpful sources were the overview article \cite{Lavau}, 
P.\ Michor's exposition in \cite{M}, and the Master's thesis of T.\ P\"oll \cite{Poell}.

Throughout this chapter, all manifolds are supposed to be $T_2$ and paracompact.
\bd We call
$$
\Xloc := \bigcup \{\X(U)\mid U\sse M \text{ open}\}
$$
the space of local vector fields\index{local vector field} on $M$.
\et
If $X$, $Y$ are local vector fields on $M$, then so is $[X,Y]$, defined on the intersection
of the domains of $X$ and $Y$. We agree to consider the `empty vector field' an element of
$\Xloc$, to avoid having to add formulations like `if the domains of $X$ and $Y$ intersect',
and similar for further local notions to be introduced below.

For any $X\in \Xloc$ we denote the maximal domain of $\Fl^X$ in $\R\times M$ by $U_X$.
$U_X$ is an open subset of $\R\times M$ (see \cite[2.3.3]{KAoM}). For each $t\in \R$, $\Fl^X_t$ is a diffeomorphism
of some maximal open set $U_t(X)$ (which may be empty) onto some open set $U_t'(X)$. Note that
$U_t'(X) = U_{-t}(X)$.

For any $n\ge 1$ in $\N$, any $\xi=(X_1,\dots,X_n)\in \Xloc^n$, any $T=(t_1,\dots,t_n)\in \R^n$,
and any $m\in M$ we set 
\begin{equation}\label{xitdef}
\xi_T(m) := \Fl^{X_1}_{t_1}(\Fl^{X_2}_{t_2}(\dots \Fl^{X_n}_{t_n}(m)\dots)).
\end{equation}
The maximal domain of $(T,m)\mapsto \xi_T(m)$ then is an open subset of $\R^n\times M$
that will be denoted by $U_\xi$, and we let $U_T(\xi)$ be the set of all $m\in M$ such
that $\xi_T(m)$ is defined, i.e., $U_T(\xi) = \{m\in M\mid (T,m)\in U_\xi\}$.
\bd A diffeomorphism $f\colon U\to U'$ between open subsets of $M$ is called a local 
diffeomorphism.\index{local diffeomorphism}
\et
If $f_i: U_i\to U_i'$ ($i=1,2$) are local diffeomorphisms, then so is $f_1\circ f_2$, with domain
$f_2^{-1}(U_1)$ and range $f_1(U_2'\cap U_1)$. Moreover, $f_1^{-1}:U_1'\to U_1$ is a local diffeomorphism 
as well. 
A group of local diffeomorphisms\index{local diffeomorphisms!group of} is a set $G$ of local 
diffeomorphisms that is closed under composition and inverses. 

Our main examples of local diffeomorphisms are flows of local vector fields. For any $X\in \Xloc$
the set of all $\Fl^X_t$ ($t\in \R$) is called the group of local diffeomorphisms generated by $X$,
and is denoted by $G_X$. More generally, if $D$ is any subset of $\Xloc$ then there exists 
a smallest group of local diffeomorphisms containing $\bigcup\{G_X\mid X\in D\}$. This group
(more precisely, pseudogroup) will be denoted by $G_D$. It is called the group of local diffeomorphisms
generated by $D$.\index{group of local diffeomorphisms} By definition, we have
\begin{equation}\label{gddef}
G_D = \{\xi_T\mid \exists n:\, \xi\in D^n,\, T\in \R^n\}.
\end{equation}
Given finite tuples $\lambda=(\lambda_1,\dots,\lambda_m)$, $\mu=(\mu_1,\dots,\mu_k)$, by
$\lambda \mu$ we denote their concatenation $(\lambda_1,\dots,\lambda_m,\mu_1,\dots,\mu_k)$,
and by $\hat\lambda$ we denote the reverse tuple $(\lambda_m,\dots,\lambda_1)$. With these notations
we have
$$
\xi_T\eta_{T'}:=\xi_T\circ\eta_{T'} = (\xi\eta)_{TT'}\quad \text{and}\quad (\xi_T)^{-1} = \xi_{-\hat T}.
$$
A subset $D$ of $\Xloc$ is called {\em everywhere defined}\index{everywhere defined} if the union of
all domains of elements of $D$ covers $M$. An analogous definition applies to groups of local
diffeomorphisms.
\bd\label{orb} Let $G$ be an everywhere defined group of local diffeomorphisms. Two elements $m_1$, $m_2$ of $M$
are called $G$-equivalent if there exists some $g\in G$ such that $g(m_1)=m_2$. The equivalence classes
of the resulting equivalence relation on $M$ are called the orbits\index{orbit} of $G$. If $D\sse \Xloc$
is everywhere defined then the $G_D$-orbits are also called $D$-orbits. The orbit of $m\in M$ 
is called {\em trivial}\index{orbit!trivial} if it equals $\{m\}$.
\et
Thus $m_1$, $m_2$ belong to the same $D$-orbit if and only if there exists some $n\ge 1$, some
$\xi\in D^n$ and some $T\in \R^n$ such that $\xi_T(m_1)=m_2$. From this we immediately conclude:
\blem Two points $m_1$, $m_2$ belong to the same $D$-orbit if and only if there exists a 
piecewise
smooth curve $\gamma:[a,b]\to M$ such that $\gamma(a)=m_1$, $\gamma(b)=m_2$, with the following property:
There exist $a=t_0<t_1<\dots < t_r=b$ and $X_i\in D$ ($i=1,\dots,r$) such that, for each $i$, 
$\gamma|_{[t_{i-1},t_i]}$ is an integral curve of $X_i$ or of $-X_i$.
\et
Such a curve $\gamma$ will be called a piecewise integral curve\index{piecewise integral curve} of $D$.

\brem\label{orbittop}
We now want to endow the orbits of $D$ with a natural topology. To this end, for any $m\in M$ and $\xi\in D^n$,
let $\rho_{\xi,m}:= T \mapsto \xi_T(m)$, and let $U_{\xi,m}\sse\R^n$ be its domain. Then the $D$-orbit
$S_m$ of $m$ is the union of all the images of the mappings $\rho_{\xi,m}$. We equip $S:=S_m$ with the finest
topology such that each $\rho_{\xi,m}$ ($n\ge 1$, $\xi\in D^n$) is continuous. Since each $\rho_{\xi,m}$
is continuous for the trace topology of $M$ on $S$ it follows that the topology of $S$ is finer than
the trace topology, i.e., $S\hookrightarrow M$ is continuous. In particular, $S$ is $T_2$. In general,
the topology on $S$ will be strictly finer than the trace topology. Since all $U_{\xi,m}$ are connected
and their images all contain $m$, $S$ is connected.
\et
\brem The topology on $S$ does not depend on the choice of $m\in S$. To see this, denote by $S_m$ the set 
$S$ equipped with the topology induced by the maps $\rho_{\xi,m}$. By symmetry, it will suffice to show that
$\id:S_m\to S_{m'}$ is continuous for all $m$, $m'\in S$. Pick $\eta$, $T_0$ such that $\eta_{T_0}(m')=m$.
Then $\rho_{\xi,m} = T\mapsto \xi_T\eta_{T_0}(m')$, which is the composition of $T\to TT_0$ and 
$\rho_{\xi\eta,m'}$. Now $\rho_{\xi\eta,m'}$ is continuous into $S_{m'}$, so $\rho_{\xi,m}:U_{\xi,m}\to S_{m'}$
is continuous. By the universal property of the finest topology, the claim follows.
\et
We next generalize the definition (\cite[Def.\ 17.1]{KLie_new}) of distribution\index{distribution} to the
variable rank setting: 
\bd A distribution on a manifold $M$ is a mapping $\Delta$ that assigns to every $m\in M$ a linear subspace
$\Delta(m)$ of $T_mM$. A set $D$ of local vector fields is said to span $\Delta$ if, for every $m\in M$,
$\Delta(m) = \text{span}\{X(m)\mid X\in D\}$.
\et
If $D\sse \Xloc$ is everywhere defined then it spans a unique distribution, which will be denoted by $\Delta_D$.
Any distribution which is of the form $\Delta_D$ for some everywhere defined family $D\sse \Xloc$ is called
smooth.\index{distribution!smooth} 

We say that a local vector field $X$ on $M$ belongs to a distribution $\Delta$ if $X(m)\in \Delta(m)$ for every
$m$ in the domain of $X$. Let
$$
D_\Delta := \{X\in \Xloc \mid X \text{ belongs to } \Delta\}.
$$
Then $\Delta$ is $\cinfty$ if and only if $\Delta$ is spanned by $D_\Delta$. Also, we always have
$\Delta_{D_{\Delta}}=\Delta$.

A distribution\index{distribution!invariant} $\Delta$ is called invariant under a group of local
diffeomorphisms $G$ if
\begin{equation}\label{ginvdef}
\forall m\in M: \quad T_mg(\Delta(m))\sse \Delta(g(m))
\end{equation}
for all $g\in G$ such that $m$ is in the domain of $g$. In this case also $T_{gm}g^{-1}$ maps 
$\Delta(g(m))$ into $\Delta(m)$, i.e., $T_mg(\Delta(m)) = \Delta(g(m))$. It follows that the
dimension of $\Delta(m)$ is the same for all points $m$ in the same $G$-orbit.

A distribution $\Delta_1$ is said to be contained in a distribution $\Delta_2$ if $\Delta_1(m)
\sse\Delta_2(m)$ for all $m\in M$. If $\Delta$ is a distribution and $G$ is a group of local
diffeomorphisms on $M$ then there is a smallest distribution $\Delta^G$ which contains 
$\Delta$ and is $G$-invariant (namely the intersection of all such distributions). More precisely,
$\Delta^G(m)$ is the linear span of all vectors $v\in T_mM$ such that $v\in \Delta(m)$ or
$v=T_{m'}g(w)$ for some $g\in G$ and $m'\in M$ with $w\in \Delta(m')$ and $m=gm'$.

\brem\label{deltagcinfty} Let $\Delta$ be spanned by $D\sse \Xloc$. Then $\Delta^G$ is spanned by
the family
$$
D \cup \{g_*X \mid X\in D,\, g\in G \text{ s.t. } g_*X \text{ is defined }\}.
$$
It follows that if $\Delta$ is $\cinfty$ then so is $\Delta^G$. 
\et
Note that if $G=G_D$ then in fact $\Delta^{G_D}$ is spanned by
\begin{equation}\label{gxd}
\D:=\{g_*X \mid X\in D,\, g\in G \text{ s.t. } g_*X \text{ is defined }\}
\end{equation}
(i.e., the union with $D$ is not required). To see this, note that any $X\in D$
with domain, say, $U$ can be written as $X=g_*X$ for $g=\Fl^X_0=\id_U\in G_D$.
If $D\sse \Xloc$ then a $G_D$-invariant distribution is called $D$-invariant. The smallest
$D$-invariant distribution which contains $\Delta$ is denoted by $\Delta^D$, i.e.,
$$
\Delta^D := \Delta^{G_D}.
$$
 
Next, let $D$ be an everywhere defined subset of $\Xloc$. Then we set 
\begin{equation}\label{pddef}
P_D := \Delta_D^D := (\Delta_D)^D.
\end{equation}
Thus $P_D$ is the smallest distribution that is $D$-invariant and contains $\Delta_D$. By 
\ref{deltagcinfty}, $P_D$ is smooth and the dimension of $P_D(m)$ depends only on the
$D$-orbit of $m$. 
\brem\label{imponote}
It is important to note that the $D$-orbits are precisely the
$P_D$-orbits: to see this, by \eqref{gddef} and the remark following \ref{orb} it suffices to
note that if $g = \xi_T\in G_D$ and $X\in D$, then $X\sim_{\xi_T} g_*X$, and therefore 
(by \cite[17.8]{KLie_new}),
$$
\Fl^{g_*X}_t = \xi_T\circ\Fl^X_t\circ\xi_{-\hat T}.
$$
Hence the flows of the $g_*X$ do not alter the $D$-orbits. 
\et
Thus the following definition makes sense:
\bd\label{orbitrank} Let $S$ be an orbit of an everywhere defined subset $D$ of $\Xloc$. For any $m\in S$,
the dimension of $P_D(m)$ is called the rank\index{orbit!rank of} of $S$.
\et
Trivial orbits are characterized by the following result:
\blem\label{trivorb} Let $D\sse \Xloc$ be everywhere defined and let $S$ be the orbit of $m\in M$. The
following are equivalent:
\begin{itemize}
	\item[(i)] $P_D(m)=\{0\}$.
	\item[(ii)] The orbit of $m$ is trivial\index{orbit!trivial}, i.e.,  $S=\{m\}$. 
\end{itemize}
\et
\pr (i)$\Rightarrow$(ii): By \eqref{gxd}, for any $X\in D$ and 
$g\in G_D$ such that $g_*X$ is defined, $g_*X(m)=0$. 
Therefore, the flow of any element of $P_D$ leaves $m$ unchanged.

(ii)$\Rightarrow$(i): Suppose that $P_D(m)\not=\{0\}$, then by \eqref{gxd} there would exist some $g\in G_D$
and some $X\in D$ with $g_*X(m)\not=0$. But then the flow of $g_*X$ would leave $\{m\}$,
contradicting (ii).
\ep
\bex As in \ref{rotex}, let $M=\R^2$ and $X:=-y\partial_x + x\partial_y$. Now set $D:=\{X\}$. Then the
$D$-orbit of any $(x_0,y_0)$ is a circle through $(x_0,y_0)$ with center $(0,0)$, while the orbit at
$(0,0)$ is trivial. In this example, $P_D = \Delta_D$: In fact, this is always true when $D$ consists of
only one vector field $X$ because in that case  
\begin{equation*}
\Fl^{(\Fl^X_t)_*X}_s = \Fl^X_t\circ \Fl^X_s\circ \Fl^X_{-t} = \Fl^X_s,
\end{equation*}
so $X$ and $(\Fl^X_t)_*X$ have the same flow and therefore coincide.
\et
\bd\label{dstar} A set $D\sse\Xloc$ is called involutive\index{involutive} if for any $X,\, Y\in \Xloc$ that belong to
$D$ also $[X,Y]$ belongs to $D$. If $D$ is any subset of $\Xloc$ then the smallest involutive subset of
$\Xloc$ that contains $D$ will be denoted by $D^*$. A smooth distribution $\Delta$ is called involutive
if the corresponding set $D_\Delta$ is involutive.
\et
\blem\label{incldist} 
Let $D$ be an everywhere defined subset of $\Xloc$. Then
$$
\Delta_D \sse \Delta_{D^*} \sse P_D.
$$
\et
\pr The first inclusion is immediate since $D\sse D^*$. The second inclusion
follows since $P_D$ is involutive by \ref{sussmain1} (iv) below.
\ep
\bex\label{rankdist} Clearly the first inclusion in \ref{incldist} can be proper. This example shows that the same may 
happen for the second inclusion. Let $M=\R^2$, $X_1:=\frac{\partial}{\partial x}$, $X_2=\phi\frac{\partial}{\partial y}$,
where $\phi(x,y)=\psi(x)$ and $\psi$ is a smooth function with $\psi(x)=0$ for $x\le 0$ and $\psi(x)>0$
for $x>0$. Let $D:=\{X_1,\,X_2\}$. Then since $D(x,y)$ has dimension $2$ for $x>0$, the same is true
for $P_D\supseteq \Delta_D$. Moreover, any point in $\R^2$ can be joined to a point $(x,y)$ with $x>0$ by
a piecewise integral curve of $D$, so 
in fact $P_D$ has dimension $2$ everywhere. However, for $x\le 0$ the distribution $\Delta_{D^*}$ has
dimension $1$.
\et
\bd\label{intedef} An immersive submanifold $S$ of $M$ is called an integral mani\-fold\index{integral manifold} of a distribution
$\Delta$ on $M$ if, for all $s\in S$, we have $T_sS = \Delta(s)$. 
A $C^\infty$-distribution $\Delta$ is said
to be integrable\index{integrable}, or to
have the integral manifold property\index{integral manifold property} if, for every $m\in M$,
either the orbit of $m$ is trivial or there exists
an integral manifold $S$ of $\Delta$ such that $m\in S$.
\et
Note that we notationally suppress the inclusion map $j: S\hookrightarrow M$ here. If a $\cinfty$-distribution
$\Delta$ is integrable then a smooth vector field $X$ belongs to $\Delta$ if and only if
$X$ is tangent to every integral manifold of $\Delta$. By \cite[17.14]{KLie_new}, any such $X$ can locally 
be viewed as a smooth vector field on any given integral manifold. Moreover, by \cite[17.22]{KLie_new}
and \ref{trivorb} we have:
\blem\label{intinv} Any integrable $\cinfty$-distribution is involutive.
\et
Note, however, that the converse is
not true in the present situation (contrary to constant rank distributions!). In fact, \ref{rankdist} provides
an example of a distribution $\Delta_D$ that is involutive:
\begin{equation*}
[X_1,X_2](x,y) = \psi'(x)\partial_y = \left\{ 
\begin{array}{rl}
\frac{\psi'(x)}{\psi(x)}X_2 & x>0\\
0 & x\le 0
\end{array}
\right.
\end{equation*}
but cannot have the integral manifold property (the dimension of the integral manifolds at $x=0$ 
would have to be $1$ and $2$).
\bd Let $\Delta$ be a smooth distribution on $M$. A maximal integral manifold\index{integral manifold!maximal} of $\Delta$
is a connected immersive submanifold $S$ of $M$ such that
\begin{itemize}
\item[(i)] $S$ is an integral manifold of $\Delta$, and
\item[(ii)] every connected integral submanifold of $\Delta$ that intersects $S$ is an open submanifold
of $\Delta$. 
\end{itemize}
$\Delta$ is said to have the maximal integral manifold property if through each point of $M$ 
with $P_{D_\Delta}(m)\not=0$ (cf.\ \ref{trivorb}) there passes a maximal integral manifold of $\Delta$.
\et
In particular, any two maximal integral submanifolds through the same point $m$ must coincide.

For the discussion below we will need the following auxiliary result:
\blem\label{eduard} Let $X$, $Y$ be smooth local vector fields on $M$. Then for any $m$ in the intersection of
the domains of $X$, $Y$ we have:
$$
[X,Y](m)=\ddt \Fl^Y_{-\sqrt{t}}(\Fl^X_{-\sqrt{t}}(\Fl^Y_{\sqrt{t}}(\Fl^X_{\sqrt{t}}(m))))
$$
\et
\pr For any local smooth function $f$ we have 
\begin{equation}\label{floff}
\frac{d}{dt} (\Fl^X_t)^*f = X(f)\circ \Fl^X_t = (\Fl^X_t)^*(L_Xf) = (\Fl^X_t)^*(Tf(X)).
\end{equation}
Since 
$$
\ddt f(\Fl^Y_{-\sqrt{t}}(\Fl^X_{-\sqrt{t}}(\Fl^Y_{\sqrt{t}}(\Fl^X_{\sqrt{t}}(m)))))
= Tf\left(\ddt \Fl^Y_{-\sqrt{t}}(\Fl^X_{-\sqrt{t}}(\Fl^Y_{\sqrt{t}}(\Fl^X_{\sqrt{t}}(m))))\right)
$$
and $Tf([X,Y])(m) = ([X,Y](f))(m)$, the claim will follow if we can show that for any such $f$ we have
$$
\ddt\left((\Fl^X_{\sqrt{t}})^*(\Fl^Y_{\sqrt{t}})^*(\Fl^X_{-\sqrt{t}})^*
(\Fl^Y_{-\sqrt{t}})^*f\right)(m) = ([X,Y](f))(m).
$$
Using \eqref{floff} we obtain 
\begin{equation*}
\begin{split}
\frac{d}{dt} \Big((\Fl^X_{\sqrt{t}})^*(\Fl^Y_{\sqrt{t}})^*&(\Fl^X_{-\sqrt{t}})^*
(\Fl^Y_{-\sqrt{t}})^*f\Big)\\ 
= & \big((\Fl^X_{\sqrt{t}})^*L_X((\Fl^Y_{\sqrt{t}})^*(\Fl^X_{-\sqrt{t}})^*
(\Fl^Y_{-\sqrt{t}})^*f) \\
&+ (\Fl^X_{\sqrt{t}})^*(\Fl^Y_{\sqrt{t}})^*L_Y((\Fl^X_{-\sqrt{t}})^*
(\Fl^Y_{-\sqrt{t}})^*f)\\
&- (\Fl^X_{\sqrt{t}})^*(\Fl^Y_{\sqrt{t}})^*(\Fl^X_{-\sqrt{t}})^*L_X((\Fl^Y_{-\sqrt{t}})^*f)\\
&- (\Fl^X_{\sqrt{t}})^*(\Fl^Y_{\sqrt{t}})^*(\Fl^X_{-\sqrt{t}})^*
(\Fl^Y_{-\sqrt{t}})^*(L_yf)
\big)\cdot \frac{1}{2\sqrt{t}} 
=: \frac{g(\sqrt{t})}{2\sqrt{t}}
\end{split}
\end{equation*}
We need to calculate the limit as $t\searrow 0$ of this expression. Now since $g(0)=0$ it 
follows that 
$$
\lim_{t\searrow 0} \frac{g(\sqrt{t})}{2\sqrt{t}} = \frac{1}{2}g'(0).
$$
Again using \eqref{floff} we calculate:
\begin{equation*}
\begin{split}
g'(0) =&
\big(L_XL_X + L_X(L_Y-L_X-L_Y) + (L_X+L_Y)L_Y + L_Y(-L_X-L_Y) \\
&- (L_X+L_Y-L_X)L_X + L_XL_Y -L_XL_Y-L_YL_Y+L_XL_Y+L_YL_Y\big)f \\
&= 2(L_XL_Y-L_YL_X)f = 2[X,Y]f,
\end{split}
\end{equation*}
giving the claim.
\ep
\brem To clarify the geometric meaning of $\Delta_{D^*}$ and $P_D$, and also to motivate the structure
of the main results in the following section, suppose that a subset $D$ of $\Xloc$ is given and that
we want to find a distribution $\Delta$
with the property that the orbits of $D$ are precisely the maximal integral manifolds of $\Delta$. 
It is then geometrically natural to define $\Delta(m)$ as
the set of all tangent vectors of smooth curves that pass through $m$ and lie entirely 
in the $D$-orbit of $m$. Call this set of curves $\Gamma_m$. 

For any $X\in D$, the integral curve $t\mapsto \Fl^X_t(m)$ belongs to $\Gamma_m$. Consequently, $\Delta(m)$
must contain $\ddt\Fl^X_t(m) = X(m)$. Moreover,  for $X$, $Y\in D$ the curve 
$$
t \mapsto \Fl^Y_{-\sqrt{t}}(\Fl^X_{-\sqrt{t}}(\Fl^Y_{\sqrt{t}}(\Fl^X_{\sqrt{t}}(m))))
$$
belongs to $\Gamma_m$. By \ref{eduard},  the derivative of this curve at $t=0$
is $[X,Y](m)$, which therefore must also lie in $\Delta(m)$. Iterating this procedure it follows that
$\Delta_{D^*}$ must be contained in $\Delta$.

There may, however, be further vectors beside those in $\Delta_{D^*}(m)$ that have to be contained in $\Delta(m)$:
Let $X\in D$, fix $t'\in \R$ and set $m':=\Fl^X_{-{t'}}(m)$. If $\gamma$ is a smooth curve with $\gamma(0)=m'$
and $\gamma\in \Gamma_{m'}$ then the curve $\sigma:=t\mapsto \Fl^X_{t'}(\gamma(t))$ belongs to $\Gamma_m$.
Setting $v:=\gamma'(0)$ it follows that $\sigma'(0)= T_{m'}\Fl^X_{t'}(v)$. We conclude that for any $v\in \Delta(m')$
we must have $T_{m'}\Fl^X_{t'}(v)\in \Delta(m)$, i.e., $\Delta$ must be $D$-invariant (i.e., $G_D$-invariant).

These considerations suggest to define $\Delta$ as the smallest $D$-invariant distribution that contains $\Delta_{D^*}$.
We shall see below that this distribution coincides with the smallest $D$-invariant distribution that contains
$\Delta_{D}$, i.e., with $P_D$.

The above also explains why $\Delta_{D^*}$ by itself may be too small to serve our purpose: it may not contain
sufficiently many directions: one may move within the orbit of $m$ along an integral curve of some $X\in D$,
catch up a new direction there and come back. Only $P_D$ is large enough to also contain these directions.
\et

\section{Orbit structure and integrability}
Throughout this section, let $D$ be an everywhere defined subset of $\Xloc$, and let $S$ be an orbit of $D$
(cf.\ \ref{trivorb}).
We equip $S$ with the natural topology introduced in \ref{orbittop}, and we will use the notations introduced
there. We set
$$
D^\infty := \bigcup_{n\in \N_{>0}} D^n,
$$
the set of all finite tuples of elements of $D$. If $\xi\in D^n$ then $U_{\xi,m}$ is open in $\R^n$ and
$\rho_{\xi,m}: U_{\xi,m} \to M$, $T \mapsto \xi_T(m)$ is $\cinfty$. Also, for any $m\in S$, 
$\rho_{\xi,m}: U_{\xi,m} \to S$ is continuous. 

Given $\xi\in D^\infty$, $m\in M$, and $T\in U_{\xi,m}$, we set
$$
V(\xi,m,T) := T_T\rho_{\xi,m}(T_TU_{\xi,m}),
$$
the image of the tangent space of $U_{\xi,m}$ at $T$ under the tangent map of $\rho_{\xi,m}$. Setting 
$m_0:=\xi_T(m)$, $V(\xi,m,T)$ is a linear subspace of $T_{m_0}M$. 
\blem\label{susslem1} Let $\xi\in D^\infty$, $m\in S$, $T\in U_{\xi,m}$, and $m_0:=\xi_T(m)$. Then
$$
V(\xi,m,T)\sse P_D(m_0).
$$
\et
\pr To begin with, let $n=1$, $\xi=X$ and $T=t_0$. Then $U_{\xi,m}=\{t\in \R\mid \exists \Fl^X_t(m)\}$
and $T_{t_0}U_{\xi,m}=\R$. Also, $\rho_{\xi,m}(t)=\Fl^X_t(m)$, so setting 
$m':=\Fl^X_{t_0}(m)$ we have
$T_{t_0}\rho_{\xi,m}=X(m')$ and therefore
$$
V(X,m,t_0) = \text{span}(X(m')) \sse \Delta_D(m')\sse P_D(m').
$$
Suppose now that the claim is already true for $n-1$. Let $\xi\in D^n$ and $T\in U_{\xi,m}$.
Then we may write $\xi=X\eta$, $T=t_0T'$ for suitable $X\in D$ and $\eta\in D^{n-1}$, $T'\in U_{\eta,m}$
and $t_0\in \R$. By definition, $V(\xi,m,T)=\text{im}(T_T\rho_{\xi,m})$. Here, we have
$\rho_{\xi,m}(t_0,T')=\Fl^X_{t_0}(\rho_{\eta,m}(T'))$, so that
\begin{equation}\label{tangsplit}
\begin{split}
T_{(t_0,T')}\rho_{\xi,m} &= X(\Fl^X_{t_0}(\rho_{\eta,m}(T'))) \oplus T_{\rho_{\eta,m}(T')}\Fl^X_{t_0}(T_{T'}\rho_{\eta,m})\\
&=X(\xi_T(m)) \oplus T_{\rho_{\eta,m}(T')}\Fl^X_{t_0}(T_{T'}\rho_{\eta,m}).
\end{split}
\end{equation}
Therefore, $V(\xi,m,T)\sse \text{span}(X(\xi_T(m))) + T_{\rho_{\eta,m}(T')}\Fl^X_{t_0}(V(\eta,m,T'))$.
By our induction assumption, $V(\eta,m,T')\sse P_D(\eta_{T'}(m))$, and since $P_D$ is $D$-invariant we 
obtain
$$
T_{\rho_{\eta,m}(T')}\Fl^X_{t_0}(V(\eta,m,T')) \sse P_D(\xi_T(m)).
$$
Moreover, $X(\xi_T(m))\in \Delta_D(\xi_T(m))\sse P_D(\xi_T(m))$, so altogether $V(\xi,m,T)\sse P_D(\xi_T(m))$,
concluding the proof for $n$.
\ep
\blem\label{susslem2} Let $m_0\in S$. Then there exist $\xi\in D^\infty$, $m\in S$, and $T\in U_{\xi,m}$
such that $\xi_T(m)=m_0$ and $V(\xi,m,T)=P_D(m_0)$.
\et
\pr We will see that the claim follows once we establish the following two statements:
\begin{itemize}
\item[(i)] If $\xi,\,\eta\in D^\infty$, $m,\,m'\in S$, $T\in U_{\xi,m}$, and $T'\in U_{\eta,m'}$
are such that $\xi_T(m)=\eta_{T'}(m')=m_0$ then there exist $\sigma\in D^\infty$, $m''\in S$, 
and $T''\in U_{\sigma,m''}$ such that 
$$
V(\xi,m,T)\cup V(\eta,m',T') \sse V(\sigma,m'',T'').
$$
\item[(ii)] There exists a finite subset $A$ of $P_D(m_0)$ that spans $P_D(m_0)$ and satisfies: for every $v\in A$
there exist $\xi\in D^\infty$, $m\in S$, and $T\in U_{\xi,m}$ such that $\xi_T(m)=m_0$ and $v\in V(\xi,m,T)$.
\end{itemize}
In fact, let us suppose that (i) and (ii) have already been established. By (ii), any element of $P_D(m_0)$
is a linear combination of the elements of $A=:\{v_1,\dots,v_k\}$ and there exist
$\xi_i$, $m_i$, $T_i$ ($1\le i\le k$) such that $v_i\in V(\xi_i,m_i,T_i)$ and $\xi_{T_i}(m_i)=m_0$ for all $i$.
By (i), then, there exists one $V(\xi,m,T)$ containing all $V(\xi_i,m_i,T_i)$, so $P_D(m_0)\sse V(\xi,m,T)$.
Together with \ref{susslem1} this finishes the proof. It therefore only remains to show (i) and (ii). 

To see (i), set $m'':=m'$, $\sigma:=\xi\hat\xi\eta$, and $T'':=T(-\hat T)T'$. Then
$$
\rho_{\sigma,m''}(T'') = \sigma_{T''}(m'') = \xi_T(\hat\xi_{-\hat T}(\eta_{T'}(m'))) = \eta_{T'}(m') = m_0.
$$
Analogously to \eqref{tangsplit} we may therefore split the tangent map of $\rho_{\sigma,m''}$ at $T''$
into a direct sum of three maps, corresponding to differentiation with respect to $T$, $-\hat T$, and $T'$,
respectively. Therefore $V(\sigma,m'',T'')=\text{im}(T_{T''}\rho_{\sigma,m''})$ contains the sum 
of the images of these partial maps. Since $\hat\xi_{-\hat T}\eta_{T'}(m'')=m$, the first of these maps 
is $T_T\rho_{\xi,m}$, so $V(\xi,m,T)=\text{im}(T_T\rho_{\xi,m})\sse V(\sigma,m'',T'')$. Moreover,
since $\xi_T\hat\xi_{-T}$ is the identity, the third of these maps has image $\text{im}(T_{T'}\rho_{\eta,m'})
=V(\eta,m',T')$, so also $V(\eta,m',T')\sse V(\sigma,m'',T'')$. This proves (i).

To see (ii), let $\tilde A$ be the set of all vectors $Y(m_0)$, where $Y$ is of the form $g_*X$ for some 
$X\in D$ and some $g\in G_D$. By \eqref{gxd}, $\tilde A$ spans
$P_D(m_0)$. Since $\dim P_D(m_0)\le \dim(M)$ we may extract a finite subset $A$ from $\tilde A$ that still
spans $P_D(m_0)$. Given any $v\in A$ we have $v=T_mg(w)$, where $g\in G_D$, $m\in S$, $g(m)=m_0$ and
$w\in T_mM$ is of the form $w=X(m)$ for some $X\in D$. Also, since $g\in G_D$ we have $g=\xi_T$ for some
$\xi\in D^\infty$ and some $T\in U_{\xi,m}$. 

Now set $\eta:=\xi X$ and $T':=(T,0)$. Then $\eta_{T'}(m)=\xi_T(m)=m_0$, 
so $\rho_{\eta,m}(T') = \eta_{T'}(m) = \xi_T\circ \Fl^X_0(m)$, and splitting the tangent map
as above we find
$$
T_{T'}\rho_{\eta,m} = T_T\rho_{\xi,m} \oplus T_m\xi_T(X(m)).
$$
In particular, $v = T_m\xi_T(X(m))\in \text{im}(T_{T'}\rho_{\eta,m}) = V(\eta,m,T')$, which gives (ii).
\ep
\blem\label{susslem3} If a connected integral manifold $N$ of $P_D$ intersects $S$ then $N$ is an 
open subset of $S$ (in the topology introduced in \ref{orbittop}).
\et
\pr Note first that if $N$ intersects $S$ in, say, $m$ then $S$ must be non-trivial by \ref{trivorb} 
because $1\le \dim T_mN$ and $T_mN = P_D(m)$. Let $\D$ as in \eqref{gxd},
then from \ref{imponote} we know that the $\D$-orbits are precisely the $D$-orbits.

By \eqref{gxd}, for any $m\in N$ there exist $X_1,\dots,X_p\in \cD$ such that 
$\{X_1(m),\dots,X_p(m)\}$ form a basis for $P_D(m)=T_mN$. In particular, $p=\dim(N)$. Let
\begin{equation*}
\begin{split}
\Phi\colon \R^p &\to N \\
(t_1,\dots,t_p) &\mapsto \Fl^{X_1}_{t_1}\circ\Fl^{X_2}_{t_2}\circ\dots\circ \Fl^{X_p}_{t_p}(m).
\end{split}
\end{equation*}
Then $\Phi$ is a diffeomorphism from an open neighborhood of $0$ in $\R^p$ onto
a neighborhood of $m$ in $N$: In fact, since $X_1,\dots,X_p$ are tangential to $N$ it follows 
from \cite[17.14]{KLie_new} that their restrictions to $N$ can be viewed as vector fields on $N$, 
and since $\partial_i\Phi(0) = X_i(m)$,  
the inverse function theorem gives the claim.

From what was said above it follows that every point in the image of $\Phi$ lies in the same
$D$-orbit as $m$. Thus any point in $N$ has a neighborhood that is contained in one orbit of $D$.
Since $N$ is connected, it follows that $N$ is contained in a single $D$-orbit $S$ (given $n_1$, $n_2\in N$
we may connect them by a smooth curve $\gamma$ in $N$. Now covering $\gamma$ by neighborhoods as
above it follows that the entire curve lies in the same orbit). Consequently, if $N\cap S\not=\emptyset$
for some $D$-orbit $S$, then $N\sse S$.

It remains to show that $N$ (as a set) is open in $S$. By \ref{orbittop} we need to see that
for any $m\in S$ and $\xi\in D^n$ ($n\in \N$), $\rho_{\xi,m}^{-1}(N)$ is open in $\R^n$.
Thus let $T=(t_1,\dots,t_n)\in U_{\xi,m}\sse \R^n$ be such that $\rho_{\xi,m}(T)\in N$.
Since $\xi = (X_1,\dots,X_n)\in D^n$, it follows that 
$$
\rho_{\xi,m}(T) = \Fl^{X_1}_{t_1}\circ\Fl^{X_2}_{t_2}\circ\dots\circ \Fl^{X_n}_{t_n}(m).
$$ 
Note that $X_i\in D\sse \cD$ for all $i$, so as in our considerations concerning $\Phi$ above we may view the
restrictions to $N$ of the $X_i$ as
local vector fields on $N$. Thus for $1\le i\le n$, $X_i|_N\in \X_{\mathrm{loc}}(N)$.
As in \eqref{xitdef} it follows from this that the maximal domain of the map $\R^n\times N\to N$, 
$(T,m')\mapsto \rho_{\xi,m'}(T)$ is open in $\R^n\times N$, so in particular
there exists an open neighborhood of $T$ in $\R^n$
which is mapped by $\rho_{\xi,m}$ into $N$.
We conclude that $\rho_{\xi,m}^{-1}(N)$ is open in $\R^n$, as claimed.
\ep
After these preparations we are now ready to prove the first main result of this section:
\bt\label{sussmain1} Let $D\sse \Xloc$ be everywhere defined. Then
\begin{itemize}
\item[(i)] If $S$ is a non-trivial orbit of $D$, then $S$, equipped with the topology from \ref{orbittop},
admits a unique $\cinfty$-structure such that $S$ becomes an immersive submanifold of $M$.
The dimension of $S$ equals the rank of $S$ (see \ref{orbitrank}).
\item[(ii)] With the topology and differentiable structure from (i), each non-trivial orbit of $D$ is a
maximal integral submanifold of $P_D$. In fact, the non-trivial orbits of $D$ are exactly the
maximal integral submanifolds of $P_D$.
\item[(iii)] $P_D$ has the maximal integral manifold property.
\item[(iv)] $P_D$ is involutive.
\end{itemize}
\et
\pr We note first that (iii) is immediate from (ii), and that (iv) follows from (iii) via \ref{intinv}.

(i) Let $m_0\in S$. Then by \ref{susslem2} there exist $m\in S$, $\xi\in D^\infty$, and $T\in U_{\xi,m}$
such that $\xi_T(m)=m_0$ and $V(\xi,m,T)=P_D(m_0)$. It follows that 
$$
k:=\rk(S)=\dim(P_D(m_0))=\dim(V(\xi,m,T)) = \dim(\text{im}(T_T\rho_{\xi,m})),
$$
so $\rk_T(\rho_{\xi,m})=k$. As we noted before \ref{orbitrank}, $k=\dim(P_D(p))$ for any $p\in S$.
Therefore, if $T'\in U_{\xi,m}$ then since $\text{im}(T_{T'}\rho_{\xi,m})=V(\xi,m,T')\sse P_D(\rho_{\xi,m}(T'))$
by \ref{susslem1},
it follows that $\rk(T_{T'}\rho_{\xi,m})\le k$, i.e., the rank of $\rho_{\xi,m}$ cannot exceed $k$ anywhere
in $U_{\xi,m}$.

On the other hand, the rank of $\rho_{\xi,m}$ locally can only increase, so 
by the rank theorem (\cite[3.3.3]{KAoM})
there exist open neighborhoods $U$
of $T$ in $\R^n$ and $V$ of $m_0$ in $M$ and diffeomorphisms $\vphi: U \to \vphi(U)\sse \R^n$
(where $n$ is such that $\xi\in D^n$), $\vphi(T)=0$ 
$\psi: V\to \psi(V)\sse \R^l$ (with $l=\dim(M)$), $\psi(m_0)=0$, such that the following diagram commutes:
\begin{equation}\label{rhodia}
\begin{CD}
U @>\rho_{\xi,m}>> V  \\
@V\vphi VV @VV \psi V \\
(-1,1)^n @>i_{n,l,k}>> (-1,1)^l 
\end{CD}
\end{equation}
Here, $i_{n,l,k}$ is the map 
$$
(x_1,\dots,x_n)\mapsto (x_1,\dots,x_k,\underbrace{0,\dots,0}_{l-k}).
$$
Denote by $N$ the (regular) submanifold $\psi^{-1}(i_{n,l,k}((-1,1)^n))$ of $M$. Then, as a set,
$N = \rho_{\xi,m}(U)$. If $T'\in U$ and $m':=\rho_{\xi,m}(T')$, then 
since \eqref{rhodia} commutes we obtain
$$
T_{m'}N = T\psi^{-1}(i_{n,l,k}(T_{T'}\vphi(\R^n))) = T_{T'}\rho_{\xi,m}(\R^n) = V(\xi,m,T').
$$
By \ref{susslem1}, $V(\xi,m,T')\sse P_D(m')$ and since both spaces have dimension $k$ they must 
be equal. Thus $N$ is an integral manifold of $P_D$. By definition of $\rho_{\xi,m}$, $N\sse S$,
so \ref{susslem3} shows that, as a set, $N$ is open in $S$. In the same way, \ref{susslem3}
applies to any open connected (in the natural manifold topology of $N$) subset of $N$.
As these sets form a basis of the manifold topology of $N$, the natural inclusion map 
$I:N\hookrightarrow S$ is an open map. Moreover, due to \eqref{rhodia} $I$ can be decomposed as
\[
\begin{CD}
U @>\rho_{\xi,m}>> \rho_{\xi,m}(U)=N @>I'>> V\\
@V\vphi VV  &   & @VV \psi V \\
(-1,1)^n &  @< {\hphantom{\text{xxxxxx}}} i_{l,n,k} {\hphantom{\text{xxxxxx}}} << & (-1,1)^l 
\end{CD}
\]
Here, $I'$ is the inclusion map and (reversing the roles of $l$ and $n$ above)
$$
i_{l,n,k} =  (x_1,\dots,x_l) \mapsto (x_1,\dots,x_k,\underbrace{0,\dots,0}_{n-k})
$$
(Note that $k\le \min(l,n)$ by the above). Hence 
$$
I = \rho_{\xi,m}\circ \vphi^{-1}\circ i_{l,n,k} \circ \psi \circ I'
$$ 
By definition of the topology on $S$, $\rho_{\xi,m}$ is continuous as a map into $S$, so also $I$ is.
We conclude that $I$ is a homeomorphism onto its image and that this image is open in $S$.

Denote by $\mcn$ the family of all manifolds $N$ constructed as above. Also, let $\set(N)$ denote
the underlying set of the manifold $N$. By the above, $\{\set(N)\mid N\in \mcn\}$ is an open cover
of $S$ and each $I:N\hookrightarrow S$ is a homeomorphism onto its image. This provides a family
of differentiable structures on the elements of an open cover of $S$. Our aim is to define a 
differentiable structure on $S$ such that each $N\in \mcn$ becomes an open submanifold of $S$.
To this end it suffices to show that for any $N_1,N_2\in \mcn$ the differentiable structures of
$\set(N_1)\cap \set(N_2)$ as an open submanifold of $N_1$ resp.\ $N_2$ coincide. Call these manifolds
$W_1$, $W_2$ and let $j:W_1\to W_2$ be the identity. By symmetry we only have to show that 
$j$ is $\cinfty$. Since both $i_1: W_1\hookrightarrow M$ and $i_2: W_2\hookrightarrow M$ are immersions,
by \cite[3.3.8]{KAoM} to see this it suffices to show that $j$ is continuous. This, however,
is immediate since both $i_1$ and $i_2$ are homeomorphisms onto the same open subset of $S$ and
$j=i_2^{-1}\circ i_1$.

It follows that $S$ possesses a structure of a smooth manifold whose natural manifold topology is
precisely the topology from \ref{orbittop}. Also, since each $N$ as above is an open submanifold of 
$M$ it follows that $S$ itself is an immersive submanifold of $M$ which is an integral manifold of $P_D$.
Suppose that $S'$ is another such smooth structure on $S$.
Then $i:S\to S'$ is a homeomorphism and both $S\hookrightarrow M$ and $S'\hookrightarrow M$ are
immersions, so again by \cite[3.3.8]{KAoM} $i$ is a diffeomorphism. Hence the smooth structure on
$S$ is unique. 

(ii) That $S$ is an integral manifold of $P_D$ was already shown in (i). In particular,
$\dim(S)=k=\dim(P_D(m_0))$. Also, $S$ is connected by \ref{orbittop}.
Now let $R$ be any connected integral manifold of $P_D$ with $R\cap S\not = \emptyset$.
By \ref{susslem3}, $\set(R)$ is an open subset of $S$. As we did in the proof of (i) for
$S$, we may apply \ref{susslem3} to each open connected subset of $R$ to see that 
the inclusion map $R\hookrightarrow S$ is an open map. Denote by $R'$ the open submanifold
of $S$ with underlying set $\set(R)$ and let $I:R'\to R$ be the identity map. Then by 
what we have just shown, $I$ is continuous. Note that both $i_R:R\hookrightarrow M$ and 
$i_{R'}:R'\hookrightarrow M$ are immersions, and $i_{R'}=i_R\circ I$. Therefore, again by 
\cite[3.3.8]{KAoM} it follows that $I$ is smooth, and is in fact an immersion. Since both
$R$ and $R'$, being integral manifolds of $P_D$, have the same dimension $k$, $I$ is 
even a local diffeomorphism and, due to its injectivity, a diffeomorphism. We conclude
that $R=R'$ as a manifold, i.e., $R$ is an open submanifold of $S$. Thus $S$ is indeed a 
maximal integral manifold of $P_D$, as claimed. Finally, let $m\in M$ and let $R$ be
a maximal integral manifold of $P_D$ containing $m$. Then by what we have just shown,
$R$ is contained in the $D$-orbit $S$ of $m$. Since, conversely, $S$ is a connected integral
manifold of $P_D$ and $R$ is maximal, we in fact have $R=S$. Thus the maximal integral
manifolds of $P_D$ are precisely the orbits of $D$.  
\ep

The second main result is as follows:
\bt\label{sussmain2} Let  $\Delta$ be a smooth distribution 
on $M$ and let $D\sse \Xloc$ span $\Delta$ (in particular, $D$ is everywhere defined). Then the 
following statements are equivalent: 
\begin{itemize}
\item[(i)] $\Delta$ is integrable.
\item[(ii)] $\Delta$ has the maximal integral manifold property.
\item[(iii)] $\Delta$ is $D$-invariant.
\item[(iv)] For every $X\in D_\Delta$, $t\in \R$, and $m\in M$ such that  $\Fl^X_t(m)$ is defined,
$$
T_m \Fl^X_t(\Delta(m)) \sse \Delta(\Fl^X_t(m)).
$$
\item[(v)] $\Delta = P_D$.
\end{itemize}
\et
\pr Clearly, (iv)$\Rightarrow$(iii)$\Rightarrow$(v), and (ii)$\Rightarrow$(i).
Also, (v)$\Rightarrow$(ii) follows from \ref{sussmain1}.

(i)$\Rightarrow$(iv): Let $X\in D_\Delta$ and let $S$ be the integral manifold of $\Delta$ that contains $m\in M$.
Then $X$ is tangential to $S$. Hence (cf.\ \cite[17.14]{KLie_new}) there exists a vector field $X'\in \X(S)$ that is $j$-related to $X$, where $j$ denotes the inclusion $j:S\hookrightarrow M$. So $Tj\circ X' = X\circ j$, implying that $j\circ \Fl^{X'}_t = \Fl^X_t\circ j$ (cf.\ \cite[2.3.17]{KAoM}). 
Consequently,
\begin{equation}\label{eq:flow_Tj}
Tj\circ T\Fl^{X'}_t = T\Fl^X_t\circ Tj. 
\end{equation}
Since $S$ is an integral manifold of $\Delta$ we obtain
\begin{align*}
T_m\Fl^X_t(\Delta(m)) &= T_m\Fl^X_t(T_mS) = T_m\Fl^X_t(Tj(T_mS)) \stackrel{\eqref{eq:flow_Tj}}{=}
Tj(T\Fl^{X'}_t(T_mS))\\ 
&\subseteq T_{\Fl^X_t(m)}S = \Delta(\Fl^X_t(m)). 
\end{align*}
\ep
Next, we want to connect the above results to the approach taken in \cite{Stefan1,Stefan2}, following \cite{M,Poell}. We begin by introducing several concepts that will be required below.
\bd $D\sse \Xloc$ is called stable\index{stable} if for all $X,\, Y\in D$ we have 
$(\Fl^X_t)^*Y\in D$, for all $t$ such that this expression is defined.
A local vector field $X$ on $M$ is called an {\em infinitesimal automorphism}\index{infinitesimal automorphism} 
of a distribution $\Delta$ if $T_m\Fl^X_t(\Delta(m))\sse \Delta(\Fl^X_t(m))$ whenever defined.
The set of infinitesimal automorphisms of $\Delta$ is denoted by $\text{\rm aut}(\Delta)$.
\et
For any $D\sse \Xloc$, let $\D\sse \Xloc$ be as in \eqref{gxd} (for $G=G_D$).
\blem \label{michcomp} Let $D\sse \Xloc$. Then
\begin{itemize}
\item[(i)] $\D$ is the smallest stable subset of $\Xloc$ containing $D$.
\item[(ii)] $D$ is stable if and only if $\D=D$.
\item[(iii)] If $D$ is stable then $P_D=\Delta_D$.
\end{itemize}
\et 
\pr (i) By \eqref{gxd}, $D\sse\D$. Moreover, $\D$ is stable: if $g_*X$, $h_*Y\in \D$ (for $g=\xi_T, h=\eta_{T'}\in G_D$ and 
$X, Y\in D$) then by \ref{imponote} we obtain
$$
(\Fl^{g_*X}_t)^*(h_*Y) = (\xi_T\circ \Fl^X_t\circ \xi_{-\hat T})^*(h_*Y)=
(\eta_{-\hat T'}\circ \xi_T\circ \Fl^X_t\circ \xi_{-\hat T})^*Y \in \D.
$$
On the other hand, if $D'\supseteq D$ is stable then given $X\in D$ and $g=\xi_T \in G_D$ 
for $T=(t_1,\dots,t_n)$ then $(\Fl^{X_n}_{-t_n})^*X\in D'$, so $(\Fl^{X_{n-1}}_{-t_{n-1}})^*(\Fl^{X_n}_{-t_n})^*X\in 
D'$, etc., leading to 
$
g_*X = (\Fl^{X_1}_{-t_1})^*\dots (\Fl^{X_n}_{-t_n})^*X\in D',
$
i.e., $\D\sse D'$.

(ii) is immediate from (i), and (iii) follows since $\Delta_\D=P_D$ by \eqref{gxd}.
\ep
\bc\label{cor:stable-implies-integrable} Let $D\sse \Xloc$ be everywhere defined and stable. Then $\Delta_D$
is integrable.
\et
\pr By \ref{michcomp}, $\Delta_D = P_D$, so the claim follows from \ref{sussmain2}.
\ep
After these preparations, we have (see \cite[3.24]{M}):
\bt\label{th:integrable_is_stable} Let $M$ be a smooth manifold and let $\Delta$ be a smooth distribution 
on $M$. Then the following statements are equivalent:
\begin{itemize}
	\item[(i)] $\Delta$ is integrable.
	\item[(ii)] $D_\Delta$ is stable.
	\item[(iii)] There exists a subset $D\sse \Xloc$ such that $\D$ 
	spans $\Delta$ (i.e., s.t.\ $P_D=\Delta$). 
	\item[(iv)]  $\text{\rm aut}(\Delta)\cap D_\Delta$ spans $\Delta$.
\end{itemize}
\et
\pr (i)$\Leftrightarrow$(iii): this is \ref{sussmain2}, (i)$\Leftrightarrow$(v), applied to $D:=D_\Delta$.

(i)$\Rightarrow$(ii): Let $X, Y\in D_\Delta$ and $m\in M$. Then $Y(\Fl^X_t(m))\in \Delta(\Fl^X_t(m))$,
and so \ref{sussmain2} (i)$\Rightarrow$(iv), applied to $D=D_\Delta$ gives
$$
(\Fl^X_t)^*Y(m) = T\Fl^X_{-t}(Y(\Fl^X_t(m)))\sse \Delta(m).
$$
Thus $(\Fl^X_t)^*Y\in D_\Delta$.

(ii)$\Rightarrow$(iii): Set $D:=D_\Delta$, then the claim follows from \ref{michcomp} (ii).

(i)$\Rightarrow$(iv): By \ref{sussmain2}, (i)$\Rightarrow$(iv), applied to $D_\Delta$, it follows
that $D_\Delta\sse \mathrm{aut}(\Delta)$. Thus $\mathrm{aut}(\Delta)\cap D_\Delta = D_\Delta$, which
spans $\Delta$ by definition.

(iv)$\Rightarrow$(iii): Set $D:=\text{\rm aut}(\Delta)\cap D_\Delta$. Then $D$ spans $\Delta$ by assumption,
and for $X\in D$ and $Y\in D_\Delta$ we have $T\Fl^X_{-t}(Y(\Fl^X_t(m)))\in \Delta(m)$ since $X\in \mathrm{aut}(\Delta)$,
so $(\Fl^X_t)^*Y\in D_\Delta$. Iterating this argument it follows that for $n\in \N$, $X_i, Y\in D$ 
($1\le i\le n$), $\xi=(X_1,\dots,X_n)$ and $T\in \R^n$, $(\xi_T)_*Y\in D_\Delta$. Therefore,
$D\sse \D\sse D_\Delta$, so $\D$ spans $\Delta$ and the claim follows.
\ep
Next, following \cite[3.25]{M}, we want to analyze the local structure of the integral manifolds
of an integrable distribution. For this we first introduce an important class of immersive
submanifolds:
\bd\label{inimandef} Let $A$ be any subset of $M^l$, 
and for $m\in A$ denote by $C_m(A)$ the set of all points in 
$A$ that can be joined to $m$ by a smooth curve\footnote{By a smooth curve here we simply mean a $\cinfty$-map
from some interval into $M$.} in $M$ lying in $A$. A subset $N$ of $M$
is called an {\em initial submanifold}\index{initial submanifold} of $M$ of dimension $n$, 
if for each $m\in N$ there exists a chart $(U_m,\vphi_m)$ of $M$ centered at $m$ such that
\begin{equation}\label{initman}
\vphi_m(C_m(U_m\cap N)) = \vphi_m(U_m)\cap (\R^n\times \{0\}) \sse \R^n\times \R^{l-n}.
\end{equation}
\et
In order to see that calling such subsets submanifolds is justified, we need an auxiliary result:
\blem\label{walschap}
Any piecewise smooth curve $c:[a,b]\to M$ admits a reparametrization as a $C^\infty$-curve.
\et
\pr Let $a=t_0<t_1<\dots <t_k=b$ be such that each curve $c|_{[t_{i-1},t_i]}$ is smooth.
For each $i$ pick a smooth map $\phi_i:[a,b]\to \R$ such that $\phi_i(t)=0$ for $t\le t_{i-1}$,
$\phi_i(t)=1$ for $t\ge t_{i}$, and $\phi$ is strictly increasing on $[t_{i-1},t_i]$. Then the map
$$
\phi := t_0 + \sum_{i=1}^k (t_i-t_{i-1})\phi_i
$$
is smooth, strictly increasing on $[a,b]$, $\phi([a,b])=[a,b]$, $\phi(t_i)=t_i$, and $\phi^{(m)}(t_i)=0$
for all $i$ and all $m\ge 1$. Thus $c\circ\phi$ is the desired reparametrization.
\ep
\blem\label{openlem} Under the assumptions of \ref{inimandef},  let $m_1,m_2\in N$ and $U_{m_1}$, $U_{m_2}$ 
be as in \eqref{initman}. Then $\vphi_{m_1}(C_{m_1}(U_{m_1}\cap N)\cap 
C_{m_2}(U_{m_2}\cap N))$ is open in $\R^n\times \{0\}$.
\et
\pr We may suppose that the intersection is nonempty, so let $p\in C_{m_1}(U_{m_1}\cap N)\cap 
C_{m_2}(U_{m_2}\cap N)$. Then 
$$
W:=C_p(U_{m_1}\cap U_{m_2}\cap N)\sse C_{m_1}(U_{m_1}\cap N)\cap 
C_{m_2}(U_{m_2}\cap N):
$$
In fact, there is a smooth curve $c_1$ from $m_1$ to $p$ in $U_{m_1}\cap N$ and
if $q\in W$ then there is also a smooth curve $c_2$ from $p$ to $q$ in $U_{m_1}\cap U_{m_2}\cap N$.
The concatenation of $c_1$ and $c_2$ can be reparametrized smoothly by \ref{walschap}, 
so $q\in C_{m_1}(U_{m_1}\cap N)$, and analogously for $m_2$. We claim that
\begin{equation}\label{dings}
\vphi_{m_1}(W) = C_{\vphi_{m_1}(p)}(\vphi_{m_1}(U_{m_1}\cap U_{m_2})\cap (\R^n\times \{0\})).
\end{equation}
$\sse$: $\vphi_{m_1}(W)\sse \vphi_{m_1}(C_{m_1}(U_{m_1}\cap N))= \vphi_{m_1}(U_{m_1})\cap 
(\R^n\times \{0\})$ and $\vphi_{m_1}(W)\sse \vphi_{m_1}(U_{m_1}\cap U_{m_2})$. Moreover, 
any $q\in W$ is connected to $p$ within $W$ by a smooth curve $c$,
and so $\vphi_{m_1}\circ c$ connects the corresponding images.

$\supseteq$: Let $\tilde c$ be a smooth curve from $\vphi_{m_1}(p)$ to $\vphi_{m_1}(q)$
in $\vphi_{m_1}(U_{m_1}\cap U_{m_2})\cap (\R^n\times \{0\})$. Then $c:=\vphi_{m_1}^{-1}\circ \tilde c$
is a smooth curve in $U_{m_1}\cap U_{m_2}\cap N$ from $p$ to $q$.

The right hand side of \eqref{dings} is precisely the connected component of $\vphi_{m_1}(p)$ in
$\vphi_{m_1}(U_{m_1}\cap U_{m_2})\cap (\R^n\times \{0\})$ in $\R^n\times \{0\}$, hence is open. Since
$p$ was arbitrary, this shows that $\vphi_{m_1}(C_{m_1}(U_{m_1}\cap N)\cap 
C_{m_2}(U_{m_1}\cap N))$ is open, as claimed.
\ep
\bt\label{initmanuniv} Let $N$ be an initial submanifold of dimension $n$ of $M^l$. Then there is a unique $\cinfty$-structure on 
$N$ such that the inclusion map $i:N\hookrightarrow M$ becomes an injective immersion and
such that the following universal property holds:
For any manifold $R^k$ and any map $f\colon R\to N$, $f$ is smooth if and only if $i\circ f\colon R\to M$ is smooth.
Moreover, $N$ is paracompact.
\et
\pr For any $m\in N$, by \ref{inimandef} there exists a chart $(U_m,\vphi_m)$ of $M$ around $m$
such that $\vphi_m(C_m(U_m\cap N)) = \vphi_m(U_m)\cap (\R^n\times \{0\})$. We define a chart for $N$
at $m$ by $(C_m(U_m\cap N), \psi_m:=\vphi_m|_{C_m(U_m\cap N)})$. Then the chart transition functions 
$\psi_{m_1}\circ\psi_{m_2}^{-1}$ are the restrictions of the smooth transition functions 
$\vphi_{m_1}\circ\vphi_{m_2}^{-1}$ of $M$ to 
subsets of $\R^n\times \{0\}$ that are open by \ref{openlem}, hence we obtain a differentiable structure on $N$.
From this choice of charts it immediately follows that $i$ is an immersion. Thus $N$ is 
an immersive submanifold of $M$. In particular, the natural manifold topology of $N$ is finer than 
the trace topology of $M$ on $N$.

To show (the non-trivial part of) the universal property, let $f\colon R\to N$ be such that
$i\circ f\colon R\to M$ is smooth, let $r\in R$ and choose a chart $(U,\vphi)$ at $f(r)$ in $M$
such that $\vphi(C_{f(r)}(U\cap N)) = \vphi(U)\cap (\R^n\times \{0\})$. Then since $f^{-1}(U)$ is
open in $R$ we may pick a chart $(V,\psi)$ in $R$ at $r$ with $V\sse f^{-1}(U)$ such that $\psi(V)$ is a ball in $\R^k$. Then since $V$ is $\cinfty$-contractible and $i\circ f$ is smooth it follows that any point in 
$f(V)$ can be connected by a smooth curve in $M$ that lies entirely in $U\cap N$ with $f(r)$.
Therefore, $f(V)\sse C_{f(r)}(U\cap N)$, and we can write
$$
(\vphi|_{C_{f(r)}(U\cap N)})\circ f\circ \psi^{-1} = \vphi\circ f\circ \psi^{-1}.
$$
Here, the right hand side is smooth by assumption and the left hand side is the chart representation
of $f$ around $r$, as a map from $R$ to $N$, so indeed $f\colon R\to N$ is $\cinfty$.

To see uniqueness, denote by $N'$ another differentiable structure on $N$ with the universal property.
Then $\id: N\to M$ and $\id: N'\to M$ are smooth, hence so are $\id: N\to N'$ and $\id: N'\to N$,
i.e., the smooth structures in fact coincide.

Concerning the paracompactness of $N$, note that any connected component of $N$ is contained in 
a connected component of $M$, hence is second countable by \cite[1.3.15]{KAoM}.
Alternatively, $M$ can be equipped with a Riemannian metric $g$ (simply by gluing local
Riemannian metrics on charts via a partition of unity) and
this induces a Riemannian metric $i^*g$ on $N$. The claim then follows from \cite[Satz 1.1.2]{H}.
\ep
Next we show that also a converse of the previous theorem holds. 
\bt\label{imfc} Let $f\colon  M^k \to N^l$ be an injective immersion between manifolds which has the
universal property from \ref{initmanuniv}, i.e.: if $h\colon P\to N$ is smooth and $h(P)\sse f(M)$,
then the induced map $\bar h\colon  P\to M$ with $f\circ \bar h = h$ is smooth.
Then $f(M)$ is an initial submanifold of $N$.
\et
\pr By \cite[3.3.3]{KAoM}, given any $m\in M$, there exist charts $(\vphi,W)$ centered at $m$ in $M$
and $(\tilde\psi,V)$ centered at $f(m)$ in $N$ such that 
$$
(\tilde\psi\circ f\circ\vphi^{-1}) = i:=(x^1,\dots,x^k)
\mapsto (x^1,\dots,x^k,0,\dots,0)\in \R^l.
$$ 
We may assume that $W = f^{-1}(V)$, so $i:\vphi(W)\to \tilde\psi(V)$.
Pick $r>0$ so small that $\{x\in \R^k\mid |x|<2r\}\sse\vphi(W)$
and $\{y\in \R^l\mid |y|<2r\}\sse\tilde\psi(V)$. Define $B^k_r(0)$ resp.\
$B^l_r(0)$ to be the open ball of radius $r$ in $\R^k$ resp.\ $\R^l$, and set 
\begin{equation*}
U := \tilde\psi^{-1}(B^l_r(0))\sse N,\quad W_1 := \vphi^{-1}(B^k_r(0))\sse M
\end{equation*}
We show that $(\psi:=\tilde\psi|_U,U)$ satisfies \eqref{initman}. In fact, since 
\begin{equation*}
\begin{split}
i(i^{-1}(\{(y^1,\dots,y^k,0,\dots,0)\mid |y|<r\})) &= \{(y^1,\dots,y^k,0,\dots,0)\mid |y|<r\}\cap i(\vphi(W))\\
&= \{(y^1,\dots,y^k,0,\dots,0)\mid |y|<r\},
\end{split}
\end{equation*}
we have
\begin{equation*}
\begin{split}
&\psi^{-1}(\psi(U)\cap (\R^k\times \{0\})) = \psi^{-1}(\{(y^1,\dots,y^k,0,\dots,0)\mid |y|<r\})\\
&\hphantom{xxx}=f\circ\vphi^{-1}( (\psi\circ f\circ\vphi^{-1})^{-1}(\{(y^1,\dots,y^k,0,\dots,0)\mid |y|<r\}))\\
&\hphantom{xxx}= f\circ\vphi^{-1}(B^k_r(0))=f(W_1).
\end{split}
\end{equation*}
Now $\tilde\psi\circ f(W_1)=i(B^k_r(0))\sse B^l_r(0)$, so $f(W_1)\sse U\cap f(M)$. Also, $f(W_1)$
is $\cinfty$-contractible, so altogether we get 
$$
\psi^{-1}(\psi(U)\cap (\R^k\times \{0\}))= f(W_1) \sse C_{f(m)}(U\cap f(M)).
$$
Conversely, let $n\in C_{f(m)}(U\cap f(M))$. This means that there exists a smooth curve 
$c:[0,1]\to N$ with $c(0)=f(m)$, $c(1)=n$, and $c([0,1])\sse U\cap f(M)$. By the universal 
property of $f$ it follows that the unique curve $\bar c:[0,1]\to M$ with 
$c=f\circ\bar c$ is smooth. 

We show that $\bar c([0,1])\sse W_1$.  Since $\bar c([0,1])\sse f^{-1}(U)\sse f^{-1}(V)=W$
and $\bar c(0) = m$ (since $f$ is injective), it follows that $\vphi(\bar c(0))=0$. Thus if
$\bar c([0,1])\not\sse W_1$, by continuity of $\bar c$ there must exist some $t\in(0,1]$
where $\vphi\circ \bar c$ intersects $\pa B^k_r(0)$, i.e., with $\bar c(t)\in \vphi^{-1}(\pa B^k_r(0))$.
But then
\begin{equation*}
\begin{split}
\psi\circ f(\bar c (t)) \in \psi\circ f\circ\vphi^{-1}(\pa B^k_r(0))=i(\pa B^k_r(0))\sse \pa B^l_r(0),
\end{split}
\end{equation*}
and so $\psi\circ c(t)=\psi\circ f\circ \bar c (t)\not\in B^l_r(0)$, i.e., $c(t)\not\in U$, a contradiction.
We conclude that $\bar c([0,1])\sse W_1$, and therefore $n=f(\bar c(1))\in f(W_1)$. Altogether,
$$
C_{f(m)}(U\cap f(M)) = f(W_1) = \psi^{-1}(\psi(U)\cap (\R^k\times \{0\})),
$$
which shows that $f(M)$ is an initial submanifold of $N$.
\ep
The previous result in particular applies to the situation where $f=j:M\hookrightarrow N$,
i.e., where $M$ is an immersive submanifold with the universal property.
In fact, the setup of \ref{imfc} is only seemingly more general: 
if $f\colon  M\to N$ is an injective immersion, then by transporting the
manifold structure of $M$ to $f(M)$ via $f$, i.e., by declaring $f\colon M\to f(M)$ to
be a diffeomorphism it follows that $j:f(M)\hookrightarrow N$ becomes an immersion,
hence $f(M)$ turns into an immersive submanifold of $N$, and the universal property
from \ref{imfc} translates into the one from \ref{initmanuniv}.

We now return to the study of integrable distributions.
\bt\label{locfol} Let $\Delta$ be an integrable distribution 
on $M^l$. Let $S$ be a non-trivial orbit of $D_\Delta$ and let $m\in S$.
Then there exists a cubic chart $(U,\vphi=(x_1,\dots,x_l))$ centered at $m$, 
$\vphi(U)=(-\eps,\eps)^l$ for
some $\eps>0$, some $k\ge 1$ and a countable set $A\sse \R^{l-k}$ such that 
$$
\vphi(U\cap S) = \{x\in \vphi(U)\mid (x^{k+1},\dots,x^l)\in A\}.
$$
If the distribution is of constant rank $k$ then the above holds for every non-trivial orbit intersecting $U$, with
the same $k$. Moreover, each non-trivial orbit is an initial submanifold of $M$.
\et
\pr 
Let $k:=\dim(S)$ and pick $X_1,\dots,X_k\in D_\Delta$ such that $\{X_1(m),\dots,X_k(m)\}$
is a basis of $\Delta(m)$. Next, choose a chart $(\chi=(y^1,\dots,y^l),W)$ around $m$ in $M$ such that
$X_1(m),\dots,X_k(m),\pdl{}{y^{k+1}}{m},\dots,\pdl{}{y^l}{m}$ is a basis of $T_mM$.
Let 
$$
f(t^1,\dots,t^l):=(\Fl^{X_1}_{t^1}\circ\dots \circ \Fl^{X_k}_{t^k})(\chi^{-1}(0,\dots,0,t^{k+1},\dots,t^l)).
$$ 
Then $f$ is a diffeomorphism from some neighborhood of $0\in \R^l$ onto a neighborhood of $m$ in $M$, and
we take $\vphi:=f^{-1}$ on a suitable neighborhood $U$ of $m$, for which we may suppose
that $\vphi(U)$ is a cube $(-\eps,\eps)^l$ with center $0=\vphi(m)$. Since $S$ is an orbit of $D_\Delta$, 
$$
m'\in S \Leftrightarrow \Fl^{X_1}_{t^1}\circ\dots \circ \Fl^{X_k}_{t^k}(m')\in S
$$
for all $m'$ and $t^1,\dots,t^k$ where the right hand side is defined. Therefore, 
for any $m'=f(t^1,\dots,t^l)\in U$ we have
\begin{equation}\label{locstru}
m'=f(t^1,\dots,t^l)\in S \Leftrightarrow f(0,\dots,0,t^{k+1},\dots,t^l)\in S.
\end{equation}
This means that $U\cap S$ is the disjoint union of connected sets of the form
$$
U_c:=\{m'\in U\mid x^{k+1}(m')=c_{k+1},\dots,x^l(m')=c_l\}
$$
where $c=(c_{k+1},\dots,c_l)$ is constant. By assumption, $\Delta$ is integrable, so \ref{sussmain2} and 
\ref{sussmain1} show that any orbit is a maximal integral manifold of $\Delta$.
Therefore, since
$$
U_c = \{(\Fl^{X_1}_{t^1}\circ\dots \circ \Fl^{X_k}_{t^k})(\chi^{-1}(0,c))\mid (t^1,\dots,t^k)\in (-\eps,\eps)^k\},
$$ 
and $S$ is an integral manifold of $\Delta$, the proof of \ref{susslem3} demonstrates that $U_c$ is an open 
(and connected) submanifold of $S$.
Now $S$, being a connected immersive submanifold, is contained
in a connected component $C$ of $M$, and $C$ is second countable since $M$ is paracompact.
Thus by \cite[14.7]{KLie_new}, $S$ is itself second countable. This shows that there
can at most be countably many $U_c$ as above. If $\Delta$ is of constant rank $k$ then clearly the above
construction works for this same $k$ for any orbit that intersects $U$. 

Finally, from \eqref{locstru} it follows that $m'\in C_m(U\cap S)$ if and only if 
$m'\in U_0$, i.e.,
$$
\vphi(C_m(U\cap S)) = \vphi(U)\cap (\R^k\times \{0\}),
$$ 
so $S$ is indeed an initial submanifold of $M$.
\ep
\brem\label{alt34} The previous result provides an alternative proof of \ref{classfrob}, (iii) $\Leftrightarrow$ (iv).
\et
\bd A chart as in \ref{locfol} is called a {\em distinguished chart}\index{distinguished chart} for $\Delta$.
The connected components of $U\cap S$ are called plaques\index{plaque} (or slices).
\et
\bd Let $D\sse \Xloc$ be everywhere defined. $D$ is said to satisfy the {\em reachability condition}\index{reachability 
condition} if the $D$-orbits are exactly the connected components of $M$.
\et
A necessary and sufficient condition for reachability is given in the following result:
\bt\label{reachth} Let $D\sse \Xloc$ be everywhere defined. Then the following are equivalent:
\begin{itemize}
	\item[(i)] $D$ satisfies the reachability condition.
	\item[(ii)] For every $m\in M$ we have $\dim P_D(m) = \dim(M)$. 
\end{itemize}
\et
\pr Let $n:=\dim(M)$.

(i)$\Rightarrow$(ii) Suppose that for some $m\in M$ we have $k:=\dim P_D(m) < n$. 
Then (by \ref{sussmain1} (i)) the orbit $S$ of
$D$ through $m$ is a $k$-dimensional connected immersive submanifold of $M$. It follows that the interior
of $S$ in the topology of $M$ is empty. But any connected component of $M$ is open, so $S$ cannot be 
such a connected component, contradicting our assumption (i).

(ii)$\Rightarrow$(i) By \cite[14.1]{KLie_new}, every maximal integral manifold of $P_D$ is open
(being an immersive submanifold of the same dimension as $M$), and connected. 
Also, $M$ is the disjoint union of
the orbits of $D$ which, by \ref{sussmain1} (ii), are exactly the maximal integral manifolds of $P_D$. 
Thus these orbits are the connected components of $M$. 
\ep
Recall from \ref{dstar} that if $D$ is any subset of $\Xloc$ then the smallest involutive subset of
$\Xloc$ that contains $D$ is denoted by $D^*$.
\bc With $M$, $n$, $D$ as in \ref{reachth}, if $\Delta_{D^*}$ has dimension $n$ for every $m\in M$ then
$D$ satisfies the reachability condition.
\et
\pr By \ref{incldist}, $\Delta_{D^*}\sse P_D$, so $\dim(P_D(m))=n$ for
all $m\in M$, and the claim follows from \ref{reachth}.
\ep
The following result provides a practically useful sufficient condition for the integrability of a singular distribution.
\bt\label{th:local_lie_generated} Let $M$ be a smooth manifold and let $D \subseteq \Xloc$ be everywhere defined. 
Suppose that each $m\in M$ has a neighborhood $U$ on which there exists a finite-dimensional Lie algebra $\mathcal{L}_U$ of vector fields which generates $D|_U$ in the sense that each element of $D|_U$ is a $C^\infty$-linear combination of elements of $\mathcal{L}_U$. 
Then the induced distribution $\Delta_D$ on $M$ is integrable.
\et
\pr Due to \ref{cor:stable-implies-integrable} it suffices to show that $D$ is stable. Let $\mathcal{L}_U$ be the $\R$-linear span of $X_1,\dots,X_r\in \X(U)$ (i.e., $\{X_1,\dots,X_r\}$ is a basis for $\mathcal{L}_U$). Then by assumption we have $D|_U = \mathrm{span}_{C^\infty(U)}(X_1,\dots,X_r)$ and it suffices to show that
\[
\forall\ X, Y \in D|_U:  \ (\Fl^X_t)^*Y|_U \in \D|_U
\]
for all $t$ where this is defined. 

Suppose first that $Y=X_j$ for some $j\in \{1,\dots,r\}$ and set $X_j(t):=(\Fl^X_t)^*X_j$. Then by \cite[2.3.13]{KAoM} we have
\begin{align}\label{eq:ddtXj}
\frac{d}{dt}X_j(t) = (\Fl^X_t)^*L_XX_j = (\Fl^X_t)^*([X,X_j]).
\end{align}
By assumption there are $a_i\in \cinfty(U)$ such that $X=\sum_{i=1}^r a_i X_i$, and using \cite[2.2.17]{KAoM} we obtain
\[
[X,X_j] = \sum_{i=1}^r \left( a_i [X_i,X_j] - X_j(a_i)X_i \right).
\]
Since $\mathcal{L}_U$ is a Lie algebra, this shows that $[X,X_j]\in \mathrm{span}_{C^\infty(U)}(X_1,\dots,X_r)$, i.e., there are $f_j^i\in\cinfty(U)$ with $[X,X_j] = \sum_{i=1}^n f_j^i X_i$. Combining this with
\eqref{eq:ddtXj} we arrive at
\begin{equation}\label{eq:ode_for_X_i}
\frac{d}{dt}X_j(t) = \sum_{i=1}^r (f_j^i\circ\Fl^X_t) (\Fl^X_t)^*X_i =  \sum_{i=1}^r (f_j^i\circ\Fl^X_t) X_i(t).
\end{equation}
This is a linear system of ODEs for $(X_1(t),\dots,X_n(t))$ with initial condition $X_i(0)=X_i$ for $i=1,\dots,r$, hence possesses a unique global (in $t$) solution, and thereby uniquely determines all $X_j(t)$. To show that $X_j(t)\in \mathrm{span}_{C^\infty(U)}(X_1,\dots,X_r)$ for all $t$ it therefore suffices to show that there are $\cinfty(U)$-linear combinations $\sum_{i=1}^r c_j^i(t)X_i$ ($j=1,\dots,r$) which solve the initial value problem \eqref{eq:ddtXj}, because it will then follow that $X_j(t)=\sum_{i=1}^r c_j^i(t)X_i$ for $j=1,\dots,r$. Inserting $\sum_{i=1}^r c_j^i(t)X_i$ into \eqref{eq:ddtXj} gives
\begin{equation*}
\sum_{i=1}^r \frac{d}{dt}[c^i_j(t)] X_i = \sum_{k=1}^r\sum_{i=1}^r  (f_j^k\circ\Fl^X_t) c_k^i(t) X_i
\end{equation*}
$(j=1,\dots,r)$. Since the $X_i$ are linearly independent, this system of equations holds if and only if
\[
\frac{d}{dt}[c^i_j(t)] = \sum_{k=1}^r  (f_j^k\circ\Fl^X_t) c^i_k(t) \qquad (i,j=1,\dots,r),
\]
which is a system of linear ODEs for the coefficient functions $c^i_j$ and thereby uniquely solvable given the initial conditions $c^i_j(0)=\delta^i_j$. Altogether, we get $X_j(t)\in \mathrm{span}_{C^\infty(U)}(X_1,\dots,X_r)$ for each $j=1,\dots,r$ and each $t$.

Finally, if $Y=\sum_{i=j}^r h_j X_j$ is a general $\cinfty(U)$-linear combination, then
\[
(\Fl^X_t)^*Y = \sum_{j=1}^r (h_j\circ\Fl^X_t) (\Fl^X_t)^*X_j \in \mathrm{span}_{C^\infty(U)}(X_1,\dots,X_r)
\]
for each $t$ by what we have just shown in the special case $Y=X_j$.
\ep
To conclude this chapter we show how the classical results on the integrability of constant rank distributions
can be derived from the results established above.
\bt\label{frob2} (Frobenius) Let $\Delta$ be a smooth distribution on $M$ of constant rank $k$. Then the following are
equivalent:
\begin{itemize}
	\item[(i)] $\Delta$ has the maximal integral manifold property.
	\item[(ii)] $\Delta$ is involutive.
\end{itemize}
\et
\pr (i)$\Rightarrow$(ii) This is \ref{intinv}.

(ii)$\Rightarrow$(i) Let $m\in M$ and pick $X^1,\dots,X^k\in \Delta$ such that $\{X^1(m),\dots,X^k(m)\}$
is a basis of $\Delta(m)$. Then also $\{X^1(m'),\dots,X^k(m')\}$ is linearly independent for all 
$m'$ in some neighborhood $U$ of $m$ in $M$. Since $\dim(\Delta(m'))=k$ for all $m'$, any local vector
field on $U$ that belongs to $\Delta$ is a linear combination of $X^1,\dots,X^k$ with smooth coefficients.
If $X$ belongs to $\Delta$ then since $\Delta$ is involutive, $[X,X^i]\in \Delta$.  Inspection of the proof
of \ref{th:local_lie_generated} shows that these properties suffice to show that $D_\Delta$ is stable and
thereby that $\Delta$ is integrable.
\ep
\bc Let $D\sse \Xloc$ be everywhere defined, involutive and of constant rank. Then $\Delta_D$ has the 
maximal integral property.
\et
\pr By \ref{frob2} it suffices to show that $\Delta_D$ is involutive. 
Let $m\in M$ and pick $X^1,\dots,X^k\in \Delta$ such that $\{X^1(m'),\dots,X^k(m')\}$ 
is linearly independent for all $m'$ in some neighborhood $U$ of $m$ in $M$. If $X,\,X'\in \Delta_D$
on $U$ then both $X$ and $X'$ are linear combinations of $X^1,\dots,X^k$ with smooth coefficients.
But then also $[X,X']$ is such a linear combination of the $X^j$ and of brackets of $X^j$ and $X^l$,
which also belong to $\Delta_D$ due to the involutivity of $D$. Thus $[X,X']\in \Delta_D$.
\ep
\brem Collecting some of the results proved above we obtain an independent proof of the 
classical Frobenius theorem \ref{classfrob}.
In fact, (iii)$\Leftrightarrow$(iv) follows from \ref{alt34}. \ref{frob2} shows that 
(i)$\Rightarrow$(iii), and (iii)$\Rightarrow$(i) is clear. (ii)$\Rightarrow$(i)
is the easy part of \cite[17.11]{KLie_new}.

Finally, to see that (iv)$\Rightarrow$(ii), let $(\vphi,U)$ be a chart as in (iv). Then for any
point $m$ in $U$, $\pa_{x^1}|_m,\dots,\pa_{x^k}|_m$ span the tangent space at $m$ of the slice $U_a$ containing
$m$. But $U_a$ is an integral manifold of $\Delta$, so this tangent space equals $\Delta(m)$.
\et
\chapter{Symmetry groups of differential equations}
\section{Local transformation groups}
In Chapter \ref{ltg_chapter} we studied Lie transformation groups on differentiable manifolds. Such
group actions are always defined globally ($\Phi\colon G\times M\to M$). For the applications to symmetry
groups of differential equations we have in mind, the natural actions will typically not be 
defined globally, however. In this section, following \cite{S,O}, we therefore study local transformation groups by means 
of the tools developed in Chapter \ref{nonconstchap}. Throughout, we will assume $M$ to 
be a connected paracompact (hence Hausdorff and second countable) $\cinfty$-manifold.
\bd\label{ltgdef2} A local transformation group\label{local transformation group} (or local Lie transformation group) 
on $M$ consists of a Lie group $G$,
an open subset $\cU$ of $G\times M$ with $\{e\}\times M \sse \cU$, and a smooth map $\Phi\colon\cU\to M$
such that
\begin{itemize}
	\item[(i)] If $(h,m)\in \cU$, $(g,\Phi(h,m))\in\cU$ and $(gh,m)\in \cU$ then
	$$
	\Phi(g,\Phi(h,m)) = \Phi(gh,m).
	$$
	\item[(ii)] For all $m\in M$, $\Phi(e,m)=m$.
	\item[(iii)] If $(g,m)\in \cU$ then $(g^{-1},\Phi(g,m))\in \cU$ and (by (i), (ii))
	$$
	\Phi(g^{-1},\Phi(g,m)) = m.
	$$
\end{itemize}
We will often abbreviate $\Phi(g,m)$ by $g\cdot m$. We set 
\begin{equation*}
\begin{split}
\cU_g &:= \{m\in M\mid (g,m)\in \cU\}\quad (g\in G)\\
\cU_m &:= \{g\in G\mid (g,m)\in \cU\}\quad (m\in M)
\end{split}
\end{equation*}
and
\begin{equation*}
\begin{split}
\Phi_g\colon  \cU_g\to M,\ m\mapsto \Phi(g,m)\\
\Phi_m: \cU_m\to M,\ g\mapsto \Phi(g,m)
\end{split}
\end{equation*}
\et
For $\cU=G\times M$ we obtain a global transformation group as in \ref{ltgdef}.
\bd\label{olvorb} A subset $\emptyset\not=S\sse M$ is called an orbit\index{orbit} of the local transformation
group $\Phi\colon\cU\sse G\times M \to M$, or a $\Phi$-orbit, if it is a minimal $G$-invariant subset of $M$, i.e.,
\begin{itemize}
	\item[(i)] If $m\in S$, $g\in G$ and $(g,m)\in \cU$ then $\Phi(g,m)\in S$.
	\item[(ii)] If $S'\sse S$ is another subset of $M$ satisfying (i) then either $S'=\emptyset$ or $S' = S$.
\end{itemize}
For any $m\in M$ we denote by $S_{\Phi,m}$ the orbit of $m$ under $\Phi$.
\et
\brem If $\Phi$ is a global transformation group then 
$$
S_{\Phi,m}=\Phi_m(G)=G\cdot m = \{\Phi(g,m)\mid g\in G\}.
$$
For a local transformation group $\Phi\colon\cU\to M$ we obtain
\begin{equation}\label{orbitform}
\begin{split}
S_{\Phi,m}=\{& p\in M\mid \exists k\in \N,\, \exists g_i\in G\ (1\le i\le k): g_k\in \cU_m,\,\\
&g_i\in \cU_{g_{i+1}\cdot\dots\cdot g_k\cdot m}\ (1\le i\le k-1)\ \text{and}\ g_1\cdot\dots\cdot g_k\cdot m = p\}
\end{split}
\end{equation}
\et
In what follows, we will illustrate many of the concepts we consider in the following
example:
\bex\label{osex} Let $M=\R^2$, $G=(\R,+)$, and 
$$
\Phi(\eps,(x,y)):=\left(\frac{x}{1-\eps x},\frac{y}{1-\eps x}\right).
$$
The natural domain of $\Phi$ is 
$$
\cU=\{(\eps,(x,y))\mid \eps<\frac{1}{x}\ \text{for}\ x>0,\ \eps>\frac{1}{x}\ \text{for}\ x<0,\ \eps\in \R
\ \text{for}\ x=0\}\sse \R\times \R^2,
$$
which is an open subset of $G\times M$, and $\Phi$ is $\cinfty$ on $\cU$. One easily checks that
$$
\Phi(\eps_1,\Phi(\eps_2,(x,y))) = \Phi(\eps_1+\eps_2,(x,y)),
$$
whenever both sides are defined. Thus we obtain a local transformation group. Its orbits are the points
on the $y$-axis and the straight half-rays emanating from the origin (except for the positive and negative
$y$-axis). Thus they are either single points or (regular) submanifolds of $M$.

Note that $\Phi$ cannot be realized as the restriction to $\cU$ of some global Lie transformation group 
on $\R^2$: in fact, for any $x\not=0$ we have $|\Phi(\eps,(x,y))|\to \infty$ as $\eps\to \frac{1}{x}$.
\et
\bd\label{blanket} A local transformation group $\Phi\colon G\times M\supseteq \cU \to M$ is called
\begin{itemize}
	\item[(i)] {\em connected}\index{local transformation group!connected}, if
	\begin{itemize}
	\item[(a)] $M$ and $G$ are connected.
	\item[(b)] $\cU$ is connected.
	\item[(c)] $\cU_m$ is connected for each $m\in M$. 
  \end{itemize} 
  \item[(ii)] {\em semi-regular},\index{local transformation group!semi-regular} if
   all orbits can be endowed with a smooth structure as immersive submanifolds of $M$ of
	the same dimension.
	\item[(iii)] {\em regular},\index{local transformation group!regular} if it is semi-regular and
	 every $m\in M$ possesses a neighborhood basis of open sets $U$ such that for 
	 every orbit $S$ of $G$ the set $U\cap S$ is connected in $S$. 
\end{itemize}
\et
{\bf Blanket assumption:} From now on we will assume all local transformation groups to be connected in the
above sense.
\brem\label{regmanrem}
\begin{itemize}
	\item[(i)] Since $S$ (being a manifold) is locally pathwise connected, a semi-regular transformation group
   is regular if and only if every $m\in M$ possesses a neighborhood basis of open sets $U$ such that for 
	 every orbit $S$ of $G$ the set $U\cap S$ is pathwise connected in $S$.
	 \item[(ii)] By \cite[3.3.12]{KAoM}, every orbit of a regular transformation group is a regular submanifold 
   of $M$.
\end{itemize}
\et
\bex 
(i) The map $\Phi$ from \ref{osex} defines a regular transformation group on $\R^2\setminus\{(0,y)\mid y\in \R\}$.

(ii) (Cf.\ \cite[Ex.\ 18.7]{KLie_new})
Let $M=T^2=S^1\times S^1$ be the two-dimensional torus and $G=(\R,+)$. Fixing $\omega\in \R$ and using 
angular coordinates $(\theta,\rho)$ on $M$ we set
$$
\Phi(\eps,(\theta, \rho)) := (\theta+\eps,\rho+\omega\eps)\quad \text{mod}\ 2\pi.
$$
Then the orbits of $\Phi$ are immersive submanifolds of dimension $1$, so $G$ acts semi-regularly on
$M$. If $\omega\in \Q$ then the orbits are closed curves and $\Phi$ acts regularly. However, if
$\omega$ is irrational then the orbits are dense in $M$ and therefore cannot be regular submanifolds
by \cite[14.1]{KLie_new}. So in this case $\Phi$ does not act regularly on $M$.
\et 
Our next aim is to show that the (non-discrete) orbits of any local transformation group can naturally be
endowed with the structure of an immersive (indeed even initial) submanifold of $M$.
\brem\label{globorbrem} 
If $\Phi\colon G\times M\to M$ is a global (connected) Lie transformation group on $M$ then we have already shown
this in \ref{ltgiso} and \ref{ltgisocor}: Denoting by $G_m$ the isotropy group of $m$ in $G$ we have
by \ref{isotropyrem} that the map $\Psi_m:G/G_m\to G\cdot m$, $gG_m\mapsto gm$ is a bijection.
If the isotropy group $G_m$ of $m$ is open then by \ref{conrem} the orbit of $m$ is the singleton $\{m\}$.
Otherwise, declaring $\Psi_m$ to be a diffeomorphism we may endow $G\cdot m$ with a smooth structure as
an immersive submanifold of $M$. The fact that $G\cdot m$ is in fact an initial submanifold of $M$ will follow 
from \ref{orbitequ} and \ref{orbcoin} below. 
\et
Turning now to the case of local transformation groups, we need to come up with a different construction
since in general $G_m$ will no longer be a subgroup of $G$ in this case: in fact, if $g_1$, $g_2\in G_m$
then in general we will not have that $(g_1^{-1}\cdot g_2,m)\in \cU$. But this fact was used in \ref{isotropyrem}
to obtain injectivity of $\Psi_m$.

The route we will take to finding a smooth structure on the orbits of $G$ goes via the infinitesimal generators
of the action of $\Phi$. We first note that also for a local transformation group we may introduce
the definitions of $\Phi(v)$ and $\Right(G,M)$ exactly as in \eqref{phidef}. 
Also, as in \ref{psihomo} it follows that the map $\Phi\colon v\mapsto \Phi(v)$ is a Lie algebra homomorphism 
from $\g_R$ onto $\Right(G,M)$, the Killing algebra\index{Killing algebra} of $\Phi$. Concerning 
\ref{flowaction} we have to be careful about the domain $\cU_m$ and note that for a local
transformation group the vector field $\Phi(v)$ need no longer be complete: 
\bp\label{flowaction2} 
Let $\Phi\colon\cU\to M$ be a local transformation group, and let $m\in M$ and $v\in \g$. Then
for all $\eps$ such that $\exp(\eps v)\in \cU_m$ we have
$$
\Fl^{\Phi(v)}_\eps(m) = \Phi(\exp(\eps v),m) = \exp(\eps v)\cdot m. 
$$
In particular,
\begin{equation*}
\Phi(v)(m) = \dde \Phi(\exp(\eps v),m).
\end{equation*}
\et
\pr For $\eps$ in a sufficiently small interval around $0$, $\exp(\eps v)\in \cU_m$. On any such
interval (hence on the maximal such interval), we can argue as in the proof of \ref{flowaction}.
\ep
For any $m\in M$ we set
\begin{equation*}
\Right_m(G,M):= \{\Phi(v)(m)\mid v\in \g\}.
\end{equation*}
Any $\Right_m(G,M)$ is a linear subspace of $T_mM$.
\bex\label{osex2} For the example from \ref{osex} we have 
\begin{equation*}
\begin{split}
\Phi(\partial_\eps|_0)(x,y)=T_0\Phi_{(x,y)}(\partial_\eps|_0)&=
\partial_\eps|_0\left(\frac{x}{1-\eps x}\right)\partial_x + \partial_\eps|_0\left(\frac{y}{1-\eps x}\right)\partial_y\\ 
&= x^2 \partial_x + xy\partial_y =:X(x,y).
\end{split}
\end{equation*}
Therefore, $\Right_{(x,y)}(G,M)$ is one-dimensional if $x\not=0$ and equals $\{0\}$ for $x=0$.

To verify \ref{flowaction} in this example (cf.\ \ref{rflowex}) we first calculate the flow of $X$. The integral curve
$c(\eps)=(x(\eps),y(\eps))$ of $X$ through $(x_0,y_0)$ has to satisfy the following initial value problem:
\begin{equation*}
\begin{split}
\partial_\eps x(\eps) &= x(\eps)^2\\
\partial_\eps y(\eps) &= x(\eps)y(\eps)\\
(x(0),y(0)) &= (x_0,y_0).
\end{split}
\end{equation*}
Indeed we obtain 
$$
(x(\eps),y(\eps))=\left(\frac{x_0}{1-\eps x_0},\frac{y_0}{1-\eps x_0}\right)=\Phi(\eps,(x_0,y_0))
$$
with maximal domain $(-\infty,\frac{1}{x_0})$ if $x_0>0$, $(\frac{1}{x_0},\infty)$ if $x_0<0$,
and $\R$, respectively, if $x_0=0$, i.e., the maximal domain of $c$ is precisely $\cU_{(x_0,y_0)}$.
In general, one can only expect that $\cU_{(x_0,y_0)}$ is contained in the domain of the corresponding
maximal integral curve. Also, the image of $c$ is the orbit of $(x_0,y_0)$.
\et
\bp\label{killaldist} Let $\Phi\colon G\times M\supseteq \cU \to M$ be a local transformation group. Then 
the Killing algebra $\Right(G,M)$ of $\Phi$ spans an integrable distribution 
$$
\Delta_\Phi := \Delta_{\Right(G,M)}
$$ 
on $M$.
\et
\pr By definition (see \eqref{rightvf}), $\Right(G,M)$ consists of the smooth local vector 
fields $\Phi(v)$ for $v\in \g$, hence $\Delta_\Phi$ is a smooth distribution on $M$. 
Pick any basis $\{v_1,\dots,v_k\}$ of $\g$ and set $X_i:=\Phi(v_i)\in \Right(G,M)$ for $1\le i \le k$.
Then by \ref{psihomo}, $\{X_1,\dots, X_k\}$ is a Lie algebra of local vector fields 
spanning $\Delta_\Phi$. Together with \ref{th:local_lie_generated} it follows that $\Delta_\Phi$ is integrable. 
\ep
Note that the dimension of $\Delta_\Phi$ may vary from point to point, so we really need the
theory of distributions of non-constant rank from Chapter \ref{nonconstchap}.
From \ref{killaldist} it follows by \ref{sussmain1} and \ref{locfol} that the non-trivial orbits $S_{\Delta_\Phi,m}$ of 
$\Delta_\Phi$ are initial submanifolds of $M$. It therefore remains to show that these orbits in
fact coincide with those introduced in \ref{olvorb}, i.e., that $S_{\Phi,m}=S_{\Delta_\Phi,m}$. 
For this we need some auxiliary results.
\blem\label{schmidtlem1} Let $\Phi\colon\cU\to M$ be a local transformation group and as in \ref{ltgdef2}, for $m\in M$ let
$\cU_m = \{g\in G\mid (g,m)\in \cU\}\quad (m\in M)$. Then for any $g\in \cU_m$ there exist
$n\in \N$ and $v_1,\dots,v_n\in\g$, $t_1,\dots,t_n\in \R_{\ge 0}$ such that
\begin{itemize}
	\item[(i)] $g=\exp(t_1v_1)\cdot\dots\cdot\exp(t_nv_n)$
	\item[(ii)] $\forall s\in [0,t_n]:\ \exp(sv_n)\in \cU_m$ 
	\item[(iii)] $\forall i=1,\dots,n-1\ \forall s\in [0,t_i]:\ \exp(sv_i)\exp(t_{i+1}v_{i+1})
	\dots \exp(t_nv_n)\in \cU_m$.
\end{itemize}
\et
\pr Denote by $\mc W$ the set of all $g\in \cU_m$ that satisfy (i)--(iii). Then $\mc W\not=\emptyset$ since
$e\in \mc W$. To show that $\mc W$ is open, let $g\in \mc W$ and pick an absolutely convex neighborhood $V$
of $0$ in $\g$ with $\exp(V)g\sse \cU_m$. Then $\exp(V)g\sse\mc W$: given $h\in \exp(V)g$ we can write
$$
h = \exp(v)\exp(t_1v_1)\cdot\dots\cdot\exp(t_nv_n),
$$
and since $V$ is absolutely convex it follows that also $\exp(tv)\exp(t_1v_1)\cdot\dots\cdot\exp(t_nv_n)
\in \cU_m$
for all $t\in [0,1]$. But also $\cU_m\setminus \mc W$ is open: suppose that $g\in \cU_m\setminus \mc W$
and again choose an absolutely convex neighborhood $V$ of $0$ in $\g$ with $\exp(V)g\sse \cU_m$.
Suppose that $(\exp(V)g)\cap \mc W\not=\emptyset$, then there exists some $v\in V$ such that $h:=\exp(v)g\in \mc W$.
Since $\exp(V)g\sse \cU_m$ we also have $\exp(tv)g\in \cU_m$ for all $|t|\le 1$. 
It follows that
$$
g_t:=\exp(t(-v))h = \exp((1-t)v)g\in \cU_m \quad (t\in [0,1]).
$$
This, however, shows that $g=g_1$ is an element of $\mc W$, contradicting our assumption.
Therefore, $\exp(V)g\sse \cU_m\setminus \mc W$, implying that $\cU_m\setminus \mc W$ is open. Since $\cU_m$
is connected by our blanket assumption following \ref{blanket}, ${\mc W}=\cU_m$.
\ep
\blem\label{schmidtlem2} Let $\Phi\colon\cU\to M$ be a local transformation group, let $m\in M$, $g\in \cU_m$ and
pick $n\in \N$, $v_1,\dots,v_n\in\g$, $t_1,\dots,t_n\in \R_{\ge 0}$ such that
(i)--(iii) of \ref{schmidtlem1} are satisfied. Define curves
$$
\gamma_l:[0,\sum_{i=1}^n t_i] \to G \quad (l=1,2)
$$
for $s\in [0,t_n]$ by
$$
\gamma_1(s) := \exp(sv_n)m \ \text{ and } \ \gamma_2(s) := \Fl^{\Phi(v_n)}_s(m),
$$
and for $s\in [\sum_{j=i+1}^n t_j,\sum_{j=i}^n t_j]$ $(i=1,\dots,n-1)$ by
\begin{equation*}
\begin{split}
\gamma_1(s) &:= \Big(\exp\Big(\Big(s-\sum_{j=i+1}^n t_j\Big) v_i\Big)\exp(t_{i+1}v_{i+1})\dots \exp(t_n v_n)\Big)\cdot m\\
\gamma_2(s) &:= \Fl^{\Phi(v_i)}_{s-\sum_{j=i+1}^n t_j}(\Fl^{\Phi(v_{i+1})}_{t_{i+1}}(\dots \Fl^{\Phi(v_n)}_{t_n}(m))\dots).
\end{split}
\end{equation*}
Then $\gamma_1(s)=\gamma_2(s)$ for all $s\in [0,\sum_{i=1}^n t_i]$.
\et
\pr For $s\in [0,t_n]$ this is immediate from \ref{schmidtlem1} and \ref{flowaction2}. 
Proceeding by induction, suppose that we already know that $\gamma_1=\gamma_2$ on 
$[0,\sum_{j=i+1}^{n}t_j]$. Now set 
\begin{equation*}
\begin{split}
\rho_1(s) &:= (\exp(sv_i)\exp(t_{i+1}v_{i+1})\dots \exp(t_n v_n))\cdot m\\
\rho_2(s) &:= \Fl^{\Phi(v_i)}_{s}(\Fl^{\Phi(v_{i+1})}_{t_{i+1}}(\dots \Fl^{\Phi(v_n)}_{t_n}(m))\dots).
\end{split}
\end{equation*}
Then it suffices to show that $\rho_1=\rho_2$ on $[0,t_i]$. We do this by showing that
$\rho_1$ satisfies the same initial value problem as $\rho_2$ on $[0,t_i]$.
Set $\exp(t_{i+1}v_{i+1})\dots \exp(t_n v_n) =: h$. By \cite[8.2]{KLie_new},
$R_h(\exp(sv_i))=\Fl^{R^{v_i}}_s(h)$, so 
$$
\frac{d}{ds}(R_h(\exp(sv_i))) = R^{v_i}_s(R_h(\exp(sv_i))).
$$
Using this, together with \eqref{olver1.46.2}, we calculate
\begin{equation*}
\begin{split}
\frac{d}{ds}\rho_1(s) &= \frac{d}{ds}\Phi_m(R_h(\exp(sv_i))) = T\Phi_m(R^{v_i}_s(R_h(\exp(sv_i))))\\
&=\Phi(v_i)(\Phi_m(R_h(\exp(sv_i)))) = \Phi(v_i)(\rho_1(s)),
\end{split}
\end{equation*}
which is precisely the defining ODE for $\rho_2$. Also, the initial values coincide since
$\rho_1(0)=\rho_2(0)$ by our inductive assumption.
\ep
After these preparations we may now show the equality of the two kinds of orbits we introduced above:
\bp\label{orbitequ} Let $\Phi\colon G\times M\supseteq \cU \to M$ be a local transformation group and let $m\in M$. Then
$$
S_{\Phi,m}=S_{\Delta_\Phi,m}.
$$
In particular, each non-trivial orbit is an initial submanifold of $M$.
\et
\pr Since $G$ is connected, once we choose a basis $\{v_1,\dots,v_k\}$ of $\g$, any element of
$G$ is a product of certain $\exp(\eps_iv_i)$. Thus by \ref{schmidtlem2} it follows that
$S_{\Phi,m}$ is trivial if and only if $\Phi(v_i)(m)=0$ for all $i$, i.e., if and only if $\Delta_\Phi(m)=0$.
Since $\Delta_\Phi$ is integrable by \ref{killaldist}, \ref{sussmain2} (v) gives that $\Delta_\Phi = P_{\Right(G,M)}$
(recall that $\Delta_\Phi=\Delta_{\Right(G,M)}$ by definition).
Consequently, \ref{trivorb} shows that the trivial orbits of $S_{\Phi,m}$ and $S_{\Delta_\Phi,m}$ coincide.

Turning now to the case of non-trivial orbits, let $m'\in S_{\Delta_\Phi,m}$ and let $g\in \cU_{m'}$. By \ref{schmidtlem1}
there exist $n\in \N$, $v_1,\dots,v_n\in \g$ and $t_1,\dots,t_n\in \R_{\ge 0}$ such that
$g=\exp(t_1v_1)\dots\exp(t_nv_n)$ and (i)-(iii) of \ref{schmidtlem1} are satisfied (with $m'$ instead of $m$).
Let $\xi:=(\Phi(v_1),\dots,\Phi(v_n))$ $\in$ $\Xloc^n$, and $T:=(t_1,\dots,t_n)$. Then using the notation
\eqref{xitdef}, \ref{schmidtlem2} shows that $g\cdot m' = \xi_T(m')$, and by definition $\xi_T(m')$ stays in 
the orbit $S_{\Delta_\Phi,m}$ of $m$, hence $g\cdot m'\in S_{\Delta_\Phi,m}$. From this, starting
with $m'=m$ and then continuing inductively, it follows that any $g_1\dots g_l\cdot m$ stays in $S_{\Delta_\Phi,m}$,
so by \eqref{orbitform} $S_{\Phi,m}\sse S_{\Delta_\Phi,m}$. Since in the above consideration we started out with
any $m'\in S_{\Delta_\Phi,m}$ and since the $\Phi$-orbits are disjoint by definition we have even shown that
$S_{\Delta_\Phi,m}$ is the disjoint union of certain $\Phi$-orbits.

To conclude the proof we show that any $\Phi$-orbit that is contained in $S_{\Delta_\Phi,m}$ is in fact
an open subset of $S_{\Delta_\Phi,m}$. Indeed, once we know this then due to $S_{\Delta_\Phi,m}$ being
a disjoint union of such sets it will follow that any $\Phi$-orbit contained in $S_{\Delta_\Phi,m}$ is
both open and closed in it, and since $S_{\Delta_\Phi,m}$ is connected the orbits must in fact coincide.

Thus let $m'\in S_{\Phi,m}\sse S_{\Delta_\Phi,m}$ (since we may move points within orbits
this is the only case we need to consider). Choose $X_i=\Phi(v_i)\in \Right(G,M)$ 
($1\le i\le n$, for a suitable $n$ and suitable elements $v_i$ of $\g$) 
such that $X_1,\dots,X_n$ is a basis for $\Delta_\Phi(m')$. Set $\xi:=(X_1,\dots,X_n)$. Then by 
the proof of \ref{sussmain1} (i) and \ref{sussmain2} (v) we have
$$
\rk(T_0\rho_{\xi,m'}) = \dim(\Delta_\Phi(m')) = n.
$$
As in \eqref{rhodia} it follows that there exists an open neighborhood $U$ of $0$ in $\R^n$ such that
$\rho_{\xi,m'}(U)$ is the domain of a coordinate chart of $S_{\Delta_\Phi,m}$. 

There exists an absolutely convex open neighborhood $V\sse U$ of $0$ in $\R^n$ such that
for any $T=(t_1,\dots,t_n)\in V$ we have 
$$
\exp(t_1v_1)\dots\exp(t_nv_n)\in \cU_{m'}.
$$ 
In particular, (i)-(iii) of \ref{schmidtlem1} are satisfied.
Since $\rho_{\xi,m'}(V)$ is open in $S_{\Delta_\Phi,m}$ it suffices to show that $\rho_{\xi,m'}(V)\sse S_{\Phi,m}$.
Also, because $V$ is absolutely convex we only need to show that $\rho_{\xi,m'}(T)\in S_{\Phi,m}$
for any $T=(t_1,\dots,t_n)\in V\cap \R_{\ge 0}^n$ (otherwise replace $v_i$ by $-v_i$). For such a $T$, \ref{schmidtlem2}
(together with \eqref{xitdef}) gives
\begin{equation}\label{fromproof}
\begin{split}
\rho_{\xi,m'}(T) &= \Fl^{X_1}_{t_1}(\Fl^{X_2}_{t_2}(\dots \Fl^{X_n}_{t_n}(m')\dots))\\
								 &= \exp(t_1v_1)\cdot\dots\cdot\exp(t_nv_n)\cdot m' \in S_{\Phi,m},
\end{split}
\end{equation}
as desired. The final claim was already shown in the remark following \ref{killaldist}.
\ep
If the transformation group $G$ is in fact global then at the moment we have two ways of
endowing the non-trivial orbits $S_\Phi$ of $\Phi$ with a differentiable structure:
the one from \ref{globorbrem} and the one from \ref{killaldist}. The following result shows
that these approaches in fact coincide.
\bp\label{orbcoin} Let $\Phi\colon G\times M \to M$ be a global Lie transformation group and let $m\in M$. Then 
the differentiable structures on any non-trivial orbit $S_{\Phi,m}$ introduced in \ref{globorbrem} and
in \ref{killaldist}, \ref{orbitequ} coincide.
\et
\pr For brevity, we set $S:=S_{\Phi,m}$, with the smooth structure from \ref{killaldist}, and we
write $S'$ for the smooth structure on $S_{\Phi,m}$ defined in \ref{globorbrem}.
We know that the inclusions $i:S\hookrightarrow M$ and $i':S'\hookrightarrow M$ are injective
immersions. Since $S$ is an initial submanifold of $M$ and $i' = i\circ \text{id}_{S'\to S}$ is smooth
it follows from the universal property in \ref{initmanuniv} that $\id:S'\to S$ is smooth.

Conversely, to see that $\id:S\to S'$ is smooth, 
let $\rho_{\xi,m'}(V)$ (with $m'\in S$) be as in the proof of \ref{orbitequ}. 
Then $\rho_{\xi,m'}: V\to \rho_{\xi,m'}(V)$ 
is a diffeomorphism onto an open submanifold (containing $m'$) of $S$, so we only need to show that
$\id_{S\to S'}\circ \rho_{\xi,m'}$ is smooth on $V$. By \ref{globorbrem}, this is the case if and only 
if $\Psi_{m'}^{-1}\circ \rho_{\xi,m'}$ is smooth on $V$. Now since $\Phi$ is global, 
it follows from \ref{flowaction} that for any $T\in V$ we have
\begin{equation*}
\begin{split}
\Psi_{m'}^{-1}\circ \rho_{\xi,m'}(T) &= \Psi_{m'}^{-1}(\exp(t_1v_1)\cdot\dots\cdot\exp(t_nv_n)\cdot m')\\
&= \exp(t_1v_1)\cdot\dots\cdot\exp(t_nv_n)\cdot G_{m'}
\end{split}
\end{equation*}
and this map is indeed smooth as the composition of the smooth maps 
$T\mapsto \exp(t_1v_1)\cdot\dots\cdot\exp(t_nv_n)$, $V\to G$,  and the quotient map $\pi: G\to G/G_{m'}$.
\ep
\brem\label{constrk} Suppose that $\Phi\colon G\times M \supseteq \cU \to M$ is a semi-regular local transformation group. Then 
since all orbits of $\Phi$ are immersive submanifolds of the same dimension, say $k$, it follows that the
integrable distribution $\Delta_\Phi$ is of constant rank $k\ge 1$. If $S$ is an orbit of $\Phi$
(i.e., by \ref{orbitequ}, an orbit of $D_{\Delta_\Phi}$) then by \ref{sussmain1} (ii) $S$ is a maximal
integral manifold of $\Delta_\Phi$ and thereby a leaf of $\Delta_\Phi$ (cf.\ \cite[Def.\ 17.23]{KLie_new}).
\et
Recall from \cite[Def.\ 17.27]{KLie_new} that a flat chart for a distribution $\Delta$ is called regular if
every leaf of $\Delta$ that intersects it does so in precisely one slice. Moreover, $\Delta$ is called regular if
every point of $M$ lies in the domain of a regular chart. The following result shows that 
for a regular local transformation group $\Phi$, $\Delta_\Phi$ is regular in this sense.
\blem\label{regslice} Let $\Phi\colon G\times M \supseteq \cU \to M$ be a regular local transformation 
group.\label{local transformation group!regular} Then
each point of $M$ lies in the domain of a regular chart for $\Delta_\Phi$, i.e., $\Delta_\Phi$ is 
a regular distribution.
\et
\pr Let $m\in M$ and pick a cubical chart $(U,\vphi)$ as in \ref{classfrob} centered at $m$, 
$\vphi(U)=[-c,c]^n$. Since $G$ acts regularly on $M$, there exists some neighborhood $V\sse U$ 
of $m$ such that the intersection of each orbit $S$ of $G$ (i.e., each maximal integral 
manifold of $\Delta_\Phi$) with $V$ is connected in $S$. This intersection
is therefore a connected integral manifold of $\Delta_\Phi$, hence by \ref{classfrob} it is 
contained in a single slice $U_a$ of $\vphi$. Now pick $c'\in (0,c)$
so that $U':=\vphi^{-1}((-c',c')^n)\sse V$. If $S$ is any orbit of $G$ that intersects $U'$ then
it also intersects $V$, hence $S\cap V\sse U_a$ for some $a$. Therefore, $S\cap U' \sse U'_a$. But since
$S$ is a maximal integral manifold we also have $U'_a\sse S$, so altogether we obtain $S\cap U'=U'_a$. 
\ep
\bex Returning once more to the example from \ref{osex}, let
\begin{equation*}
\begin{split}
U_+&:=\{(x,y)\in \R^2\mid x>0\}\\
U_-&:=\{(x,y)\in \R^2\mid x<0\}
\end{split}
\end{equation*}
and set
\begin{equation*}
\begin{split}
\vphi_\pm :U_\pm &\to \R\times\R_\pm\\
(x,y)&\mapsto \Big(x,\frac{y}{x}\Big)
\end{split}
\end{equation*}
We show that $(U_+,\vphi_+)$ is regular chart for the transformation group $\Phi$ from \ref{osex}
(and analogously for $(U_-,\vphi_-)$). Note that $\vphi^{-1}(r,s)=(r,rs)$.
From \ref{osex2} we know that $\Delta_\Phi$ is spanned by the vector field $X=x^2 \partial_x + 
xy\partial_y$. In the chart $(U_+,\vphi_+)$, with coordinates $(r,s):=\vphi_+(x,y)$
we have 
$$
\vphi_*(\partial_x) = \frac{\partial r}{\partial x}\partial_r + \frac{\partial s}{\partial x}\partial_s 
= \partial_r - \frac{s}{r}\partial_s,\ 
\vphi_*(\partial_y) = \frac{\partial r}{\partial y}\partial_r + \frac{\partial s}{\partial y}\partial_s
= \frac{1}{r}\partial_s.
$$
Therefore, $X$ has the representation
\begin{equation*}
\begin{split}
\vphi_*X &= (\vphi^{-1})_1^2 \vphi_*(\partial_x) + (\vphi^{-1})_1 (\vphi^{-1})_2 \vphi_*(\partial_y)\\
&= r^2\vphi_*(\partial_x) + rs\vphi_*(\partial_y) = r^2\partial_r.
\end{split}
\end{equation*}
By \ref{osex} the non-trivial orbits, i.e., the maximal integral manifolds of $\Phi$ are half-rays
emanating from $(0,0)$, except for those lying on the $y$-axis. Hence each such orbit $S$ is either contained
in $U_+$ or in $U_-$. If $S\sse U_+$ then 
$$
S=\{(x,y)\in \R^2\mid s=\frac{y}{x}=c \text{ and } x>0\},
$$
so $S$ is a slice of $\vphi_+$.  
This shows that $(U_\pm,\vphi_\pm)$ are regular charts for $\Delta_\Phi$. Nevertheless, $\Phi$ is not
regular because the points on the $y$-axis are trivial orbits, i.e., not every orbit of $\Phi$ has the
same dimension. Removing the $y$-axis from $M=\R^2$ we obtain a regular local transformation group.
\et
Finally, we have the following fundamental result on the space of orbits of a regular
local transformation group:
\bt\label{quotient} Let $\Phi\colon G\times M \supseteq \cU \to M$ be a regular local transformation 
group on an $n$-dimensional manifold $M$ with $k$-dimensional orbits. Then the set 
$M/G$ of orbits of $\Phi$ can be endowed with the structure of an $(n-k)$-dimensional
manifold with the following properties:
\begin{itemize}
	\item[(i)] The quotient map $\pi: M \to M/G$, $m\mapsto S_{\Phi,m}$ is a surjective submersion.
	\item[(ii)] $m$ and $m'$ belong to the same orbit if and only if $\pi(m)=\pi(m')$.
	\item[(iii)] For any $m\in M$, $\Delta_{\Phi}(m) = \ker T_m\pi$. 
\end{itemize}
\et
\pr Surjectivity of $\pi$ and (ii) are immediate from the definition. Since the orbits of $\Phi$
are precisely the leaves of the foliation induced by $\Delta_\Phi$, which is regular by \ref{regslice}, (i) is \cite[17.29]{KLie_new}.
For (iii), note that since $\pi$ is constant on any orbit $S$ of $\Phi$, $T_m\Phi$ must vanish
on the tangent space of $S$ at $m\in S$. But $S$ is an integral manifold of $\Delta_\Phi$, so $T_mS = \Delta_\Phi(m)$.
Thus $\Delta_\Phi(m)\sse \ker(T_m\pi)$. Finally, since $\pi$ is a submersion, $\dim\ker(T_m\pi)=k=\dim(\Delta_\Phi(m))$,
so we have equality. 
\ep
\brem\label{quotientconst} For later use we recall the construction of \cite[17.29]{KLie_new}: Given a regular chart
$(U,\vphi=(x^1,\dots,x^n))$, let $U':=\pi(U)$. Then the map $\vphi':U'\to \R^{n-k}$,
$m'\mapsto \proj_2(\vphi(m))$, where $m$ is any element of $\pi^{-1}(m')\cap U$ 
(recall that $\pi^{-1}(m')\cap U = U_a$ for some $a$ and $\proj_2\circ\vphi|_{U_a}\equiv a$) 
is a typical chart for $M/G$. The local
representation of $\pi$ with respect to the standard charts then is 
$$
\vphi'\circ\pi\circ \vphi^{-1}=\proj_2 = (x^1,\dots,x^n) \mapsto (x^{k+1},\dots,x^n):
$$
\begin{equation*}
\begin{CD}
M\supseteq U @>\pi>> U'\sse M/G  \\
@V\vphi VV @VV \vphi' V \\
\vphi(U) @> \proj_2>> \vphi'(U) 
\end{CD}
\end{equation*}
\et
\section{Symmetries of algebraic equations}
From this point of the course onwards we will closely follow Olver's work \cite{O}. 
Concerning notations, we will henceforth
typically denote points in a smooth manifold $M$ by $x,\, y,\dots$, since
the manifolds we are interested in will mainly be subsets of spaces of independent and dependent
variables of differential equations. As already announced before \ref{flowaction}, we will from now
on denote the group parameter with $\eps$ instead of $t$ since we will often need $t$ as a
variable in a differential equation. Moreover, given a local Lie group action $\Phi$ on $M$
we will often notationally suppress the Lie algebra homomorphism $v\mapsto \Phi(v)$ from \eqref{phidef},
i.e., we will often simply write $v$ instead of $\Phi(v)\in \Xloc$. Also, we will usually only write $G$
instead of $\Phi$.

By a system of algebraic equations\index{algebraic equation} 
we mean any system of equations
\begin{equation}\label{algsys}
F_\nu(x) = 0,\qquad \nu=1,\dots,l,
\end{equation}
where $F_1,\dots,F_l$ are smooth real-valued functions on $M$. The term `algebraic' is used to 
distinguish this situation from the case of differential equations to be considered later on.
It does {\em not}, however, restrict the form of the $F_\nu$ (e.g., to polynomials).
A {\em solution}\index{algebraic equation!solution}
of \eqref{algsys} is any point $x\in M$ such that $F_\nu(x) = 0$ for $\nu=1,\dots,l$.
A {\em symmetry group}\index{symmetry group} of \eqref{algsys} is any local transformation group $G$ 
on $M$ that transforms any solution of \eqref{algsys} into another solution.

More generally, we define:
\bd Let $G$ be a local transformation group on a smooth manifold $M$. A subset $S\sse M$
is called $G$-invariant, and $G$ is called a symmetry group of $S$ if, whenever $x\in S$
and $g\in G$ is such that $g\cdot x$ is defined, then also $g\cdot x\in S$.
\et
\bex\label{ase1} (i) Let $M=\R^2$ and consider the one-parameter group of translations
$$
G_c: (x,y) \mapsto (x+c\eps,y+\eps) \qquad (c\in \R).
$$
\et
Then any line $\{x=cy+d\}\sse \R^2$ is an orbit of $G_c$, hence is invariant under
$G_c$. Any $G_c$-invariant subset of $M$ is a union of such orbits.

(ii) Again let $M=\R^2$ and for $\alpha\in \R$ let
$$
G^\alpha: (x,y) \mapsto (\lambda x, \lambda^\alpha y) \qquad (\lambda>0).
$$
Then $\{(0,0)\}$, as well as the positive and negative $x$- and $y$-axes are $G^\alpha$-invariant.
Hence also the entire axes are invariant, being unions of invariant sets.
Moreover, $\{xy=0\}$ and $\{y=k|x|^\alpha\}$ (for $x<0$ or $x>0$) are invariant.

(iii) Most of the time we will be interested in invariant sets $S$ that are subvarieties,\index{subvariety}
given by the common zero set of smooth functions $F=(F_1,\dots,F_l)$:
$$
S=S_F = \{x \mid F_\nu(x)=0,\, \nu=1,\dots, l\}.
$$

(iv) If $S_1$, $S_2$ are invariant subsets, so are $S_1\cup S_2$ and $S_1\cap S_2$.

\bd Let $G$ be a local transformation group acting on $M$ and let $F: M\to N$ be a
map into a manifold $N$. $F$ is called $G$-invariant\index{invariant function} if
for all $x\in M$ and all $g\in G$ such that $g\cdot x$ is defined we have
$$
F(g\cdot x) = F(x).
$$
If $N=\R$ then $F$ is simply called an invariant\index{invariant} of $G$. 
\et
Obviously, $F=(F_1,\dots,F_l):M\to \R^l$ is $G$-invariant if and only if each $F_\nu$ is an
invariant of $G$.
\bex\label{ase2} (i) For $G_c$ the group of translations from \ref{ase1} (i), the function
$$
f(x,y) := x -cy
$$
is an invariant: $f(x+c\eps,y+\eps) = f(x,y)$. Moreover, any invariant of $G_c$ must be of 
the form $g(x-cy)$ for some smooth $g$.

(ii) As in \ref{ase1} (ii), let $G^1:(x,y)\mapsto (\lambda x,\lambda y)$ ($\lambda>0$). Then
$f(x,y):=x/y$ is an invariant for $G^1$, defined on $\{y\not=0\}$. Another invariant is
$(x,y)\mapsto xy/(x^2+y^2)$, defined on $\R^2\setminus \{(0,0)\}$. There is no smooth nonconstant
invariant defined on all of $\R^2$.
\et
\brem\label{123} If $F:M\to \R^l$ is a $G$-invariant function then every level set of $F$ is a $G$-invariant
subset of $M$: if $F(x)=c$ and $g\cdot x$ is defined then also $F(g\cdot x)=F(x)=c$.

However, if the zero-set $\{F(x)=0\}$ of a smooth map $F$ is $G$-invariant then $F$ itself need not
be invariant. For example, $\{(x,y)\mid xy=0\}$ is invariant under $G^1$ from \ref{ase1} (ii),
but $F(x,y)=xy$ is not $G$-invariant since $F(\lambda x,\lambda y)=\lambda^2xy\not=F(x,y)$.
To obtain a true statement we need to take all level sets into consideration:
\et
\blem Let $G$ be a local transformation group acting on $M$ and let $F\in \cinfty(M,\R^l)$.
The following are equivalent:
\begin{itemize}
	\item[(i)] $F$ is invariant under $G$.
	\item[(ii)] Every level set $\{F(x)=c\}$ ($c\in \R^l$) of $F$ is invariant under $G$.
\end{itemize}
\et
\pr (i)$\Rightarrow$(ii): See \ref{123}.

(ii)$\Rightarrow$(i): Let $x$, $g$ be such that $g\cdot x$ is defined and set $c:= F(x)$.
Then also $gx\in \{y\mid F(y)=c\}$, so $F(gx)=F(x)$.
\ep
The following result gives a simple, linear, criterion for a function to be an invariant
of a group action. It is a first typical example of Lie theoretic methods in symmetry analysis.
Recall from \ref{blanket} that we always assume transformation groups to be connected.
\bt \label{finvcrit} Let $G$ be a local transformation group acting on $M$ and let $f\in \cinfty(M,\R)$.
The following are equivalent:
\begin{itemize}
	\item[(i)] $f$ is an invariant of $G$.
	\item[(ii)] For every infinitesimal generator $v$ of $G$ we have $v(f)=0$.
\end{itemize}
\et
\pr (i)$\Rightarrow$(ii): For clarity, we de-identify again and write $\Phi(g,x)=gx$ for the action of $G$.
Let $v\in \g$. Then since $f$ is an invariant, given any $x\in M$,
by \ref{flowaction2} we obtain for $\eps$ small:
$$
f(x)=f(\exp(\eps v)\cdot x) \equiv f(\Phi(\exp(\eps v),x)) = f(\Fl^{\Phi(v)}_\eps(x))\equiv f(\Fl^{v}_\eps(x)).
$$
Differentiating this expression with respect to $\eps$ at $\eps=0$ we get 
$$
v(f)\equiv \Phi(v)(f) = 0.
$$
(ii)$\Rightarrow$(i): Since $v(f)$ vanishes identically on $M$ we have
$$
\frac{d}{d\eps} f(\exp(\eps v)\cdot x) = \frac{d}{d\eps} f(\Fl^{v}_\eps(x)) = v(f)(\Fl^{v}_\eps(x)) = 0,
$$
wherever defined, so $f(\exp(\eps v)\cdot x) = f(x)$ wherever defined.
Since $G$ is connected, by \ref{schmidtlem1} any $g\in {\mathcal U}_x$ can be written as 
a product of certain $\exp(\eps_iv_i)$ ($v_i\in \g$) satisfying (i)--(iii) from that result.
Then \ref{schmidtlem2} shows that we can iterate the above argument to conclude that 
$f(g\cdot x)= f(x)$ for all $g\in {\mathcal U}_x$.
\ep
It follows that, if $\{v_1,\dots,v_r\}$ is a basis of the local Lie algebra of infinitesimal generators of $G$
(i.e., if $\{\Phi(v_1),\dots,\Phi(v_r)\}$ is a local basis of $\Delta_\Phi$)
then $f$ is an invariant of $G$ (on the open set where $\{v_1,\dots,v_r\}$ is a basis) 
if and only if $v_k(f)=0$ for $k=1,\dots,r$. If $G$ acts effectively, then by \ref{phiiso}
we may take for $\{v_1,\dots,v_r\}$ any basis of $\g$. Writing $v_k=\sum_{i=1}^n \xi^i_k 
\partial_{x_i}$ in local coordinates, this means that $f$ has to satisfy the homogeneous system of linear
PDEs of first order
\begin{equation}\label{invhomsys}
v_k(f)(x) = \sum_{i=1}^n \xi^i_k(x) \frac{\partial f}{\partial x_i} = 0 \quad (k=1,\dots,r).
\end{equation}
\bex We return to the translation group $G_c$ from \ref{ase1} (i). Its infinitesimal generator is 
$$
v = \dde (x+c\eps,y+\eps) = (c,1)\equiv c\partial_x + \partial_y.
$$
We already know that $f(x,y)=x-cy$ is an invariant, and indeed $v(f)=0$.

In the case of the scale group $G^\alpha$ from \ref{ase1} (ii) we have an action of the multiplicative
group $(\R^+,\cdot)$ on $\R^2$ whose generator therefore is
$v=\partial_\lambda|_{\lambda=1}(\lambda x,\lambda^\alpha y) = x\partial_x + \alpha y\partial_y$ and it
is easily checked that the infinitesimal criterion is satisfied for the invariants of $G^\alpha$ given
in \ref{ase2} (ii).
\et
Turning now to symmetries of systems of algebraic equations, we begin by deriving a general criterion
for the local invariance of submanifolds under local group actions.
\bd Let $G$ be a local transformation group acting on $M$. A subset $S\sse M$ is called
locally $G$-invariant\index{locally $G$-invariant} if for every $x\in S$ there exists a
neighborhood $\tilde\mcu_x\sse \mcu_x$ (cf.\ \ref{ltgdef2}) of the identity in $G$
such that $g\cdot x\in S$ for all $g\in \tilde\mcu_x$. A smooth map $F:U\to N$ ($U$ open in $M$)
is called locally $G$-invariant if for each $x\in U$ there exists some neighborhood
$\tilde\mcu_x\sse \mcu_x$ of $e$ in $G$ such that $F(g\cdot x)=F(x)$ for all $g\in \tilde\mcu_x$ with $g\cdot x\in U$.
$F$ is called globally $G$-invariant if $F(g\cdot x) = F(x)$ for all $x\in U$ and $g\in G$ such that
$g\cdot x \in U$.
\et
\bex Let $G$ be the group of translations $(x,y)\mapsto (x+\eps,y)$ on $M=\R^2$. 
Then the set $S:=\{(x,0)\mid -1<x<1\}$ is locally $G$-invariant, but not $G$-invariant.

Let 
\begin{equation*}
f(x,y) := \left\{ 
\begin{array}{rl}
0 & y\le 0 \text{ or } y>0 \text{ and } x>0\\
e^{-1/y} & y>0 \text{ and } x< 0.
\end{array}
\right.
\end{equation*}
Then $f$ is smooth and locally $G$-invariant on $U:=\R^2\setminus\{(0,y)\mid y\ge 0\}$: in fact,
$f(x+\eps,y)=f(x,y)$ for $|\eps| < |x|$. But $f$ is not globally $G$-invariant.
\et
\bt\label{locinvth} Let $N$ be an initial submanifold of $M$. The following are equivalent:
\begin{itemize}
\item[(i)] $N$ is locally $G$-invariant.
\item[(ii)] The infinitesimal generators of $G$ are everywhere tangent to $N$, i.e.: 
$$
\forall x\in N\ \forall v\in \g\colon  v|_x \in T_xN
$$
\end{itemize}
\et
\pr If $\Phi\colon G\times M\supseteq \mcu \to M$ is the action of $G$ then the above condition, 
written out, means (cf.\ \eqref{phidef}):
$$
\forall x\in N\ \forall v\in \g\colon  \Phi(v)|_x \in T_xN
$$
(i)$\Rightarrow$(ii): Let $v\in \g$, $x\in N$. Since $N$ is locally $G$-invariant, 
for $|\eps|$ small we have $\Phi(\exp(\eps v),x)\in N$. As $N$ is initial, $\eps\mapsto \Phi(\exp(\eps v),x)$
is smooth as a map into $N$ as well. The derivative of this curve at $\eps=0$ therefore is an
element of $T_xN$. This implies (ii) since by \ref{flowaction2} we have
$$
\Phi(v)(x) = \dde \Phi(\exp(\eps v),x).
$$
(ii)$\Rightarrow$(i): By assumption, any $\Phi(v)$ is tangent to $N$, hence by \cite[17.14]{KLie_new}
it can be viewed as a vector field on $N$, and so the restriction of its flow to $N$ is a local 
diffeomorphism on $N$. 
Given $x\in N$, $\Fl^{\Phi(v)}_\eps(x)$ 
exists for $|\eps|$ small. Using \ref{flowaction2} we therefore have for $|\eps|$ small
\begin{equation}\label{singleeps}
\exp(\eps v)\cdot x = \Phi(\exp(\eps v),x) = \Fl^{\Phi(v)}_\eps(x)
\in N.
\end{equation}
Next, choose a basis $v_1,\dots,v_k$ of $\g$. 
Then as in the last part of the
proof of \ref{susslem3} it follows that there exists some $c>0$ such that 
$\Fl^{\Phi(v_1)}_{\eps_1}\circ\dots\circ \Fl^{\Phi(v_k)}_{\eps_k}(x)$ 
exists for all $|\eps_i|<c$. In addition, we may assume $c$ so small that
$\exp(\eps_1 v_1)\cdot \dots\cdot \exp(\eps_k v_k)\in \cU_x$ for $|\eps_i|<c$. Then (i)--(iii) of
\ref{schmidtlem1} hold and so, by \ref{schmidtlem2}, iterating \eqref{singleeps} gives
$$
\exp(\eps_1 v_1)\cdot \dots\cdot \exp(\eps_k v_k)\cdot x = 
\Fl^{\Phi(v_1)}_{\eps_1}\circ\dots\circ \Fl^{\Phi(v_k)}_{\eps_k}(x)\in N
$$
for these $\eps_i$. Therefore, defining the neighborhood $\tilde\mcu_x := 
\{\exp(\eps_1v_1)\dots\exp(\eps_kv_k)\mid |\eps_i|<c\}$ 
of $e\in G$ (shrinking $c$ if necessary, cf.\ \cite[8.5]{KLie_new}) we obtain that $g\cdot x\in N$ for all $g\in \tilde\mcu_x$. 
\ep
\brem As the proof of \ref{locinvth} shows, the implication (ii)$\Rightarrow$(i) remains correct
even for $N$ an immersive submanifold of $M$. However, (i)$\Rightarrow$(ii) is not true
in this generality: To see this, equip $M=\R^2$ with a new manifold structure $N$ whose
charts are $\vphi_a:(x,a)\mapsto x$ ($a\in \R$), see \cite[Ex.\ 14.3]{KLie_new}. Then $N$ is 
a one-dimensional immersive submanifold of $M$ with underlying set $\R^2$, 
hence is (even globally) invariant
under {\em any} local group action on $M$. But clearly not every group action has generators
tangential to $N$ (i.e., horizontal). 
\et
As the most important special case of \ref{locinvth} we consider zero sets of smooth maps:
\bc\label{algsym} Let $G$ be a local transformation group acting on $M$, $\dim(M)=m$. Let 
$F:M\to \R^l$ ($l\le m$) be smooth and of maximal rank ($=l$) at every solution $x$ of
the system $F_\nu(x)=0$ ($1\le \nu \le l$). Then $G$ is a symmetry group of this system
if and only if
\begin{equation}\label{algsymcrit}
\forall v\in \g\ \forall \nu=1,\dots l\ \forall x\in M \text{ with } F(x)=0:\ v(F_\nu)|_x = 0
\end{equation}
\et
\pr By \cite[3.3.23]{KAoM}, $N:=\{x\in M\mid F(x)=0\}$ is a regular (hence in particular initial)
submanifold of $M$. Thus by \ref{locinvth}, $N$ is locally invariant under $G$ if and
only if $v|_x\in T_xN$ for all $x\in N$ and $v\in \g$. By \cite[3.3.25]{KAoM}, 
$$
T_xN = \ker(T_xF)=\bigcap_{\nu=1}^l \ker(T_xF_\nu)
$$ 
for all $x\in N$. This implies that $N$ is locally $G$-invariant if and only if for each $x\in N$,
i.e., for each $x\in M$ with $F(x)=0$, each $v\in \g$, and each $\nu=1,\dots,l$ we have that 
$$
0 = T_xF_\nu(v|_x) = v(F_\nu)|_x.
$$
It remains to show that local invariance implies invariance of $N$ under $G$.
Thus let $x\in N$ and $g\in \cU_x$. By \ref{schmidtlem1}, $g= \exp(t_1v_1)\cdot \dots \cdot \exp(t_nv_n)$
with properties (i)--(iii) from that result. To show that $g\cdot x\in N$, we proceed by induction.
Suppose we already know that $g_i\cdot x\in N$, where $g_i:=\exp(t_{i+1}v_{i+1})\cdot\dots\cdot \exp(t_nv_n)$
and set $A_i:=\{t\in [0,t_i]\mid (\exp(tv_i)\cdot g_i)\cdot x \in N\}$. Then $0\in A_i$ and $A_i$
is closed since $N$ is. If $t\in A_i$ then by local invariance of $N$ there exists some $\alpha>0$
such that for all $|s|<\alpha$ we have $\exp(sv_i)\cdot [(\exp(tv_i)\cdot g_i)\cdot x]\in N$. Here,
$h_i:=\exp(tv_i)\cdot g_i\in \cU_x$ and $\exp(sv_i)\in \cU_{h_i\cdot x}$. Also, for $|s|<\alpha$ and
$t+s\in [0,t_i]$, $\exp(sv_i)\cdot\exp(tv_i)\cdot g_i=\exp((s+t)v_i)\cdot g_i\in \cU_x$ by \ref{schmidtlem1} (iii).
Hence \ref{ltgdef2} (i) gives $[(\exp((s+t)v_i)\cdot g_i)]\cdot x=
[\exp(sv_i)\cdot (\exp(tv_i)\cdot g_i)]\cdot x =
\exp(sv_i)\cdot [(\exp(tv_i)\cdot g_i)\cdot x] \in N $ for these
$s$, which establishes that $A_i$ is also open in $[0,t_i]$. Altogether, $A_i=[0,t_i]$, and 
our claim follows.
\ep
\bex (i) Let $G=SO(2)$ be the rotation group acting on $M=\R^2$, with infinitesimal generator 
$$
v= \dde(x\cos\eps - y\sin\eps, x\sin\eps + y\cos\eps)=-y\partial_x+x\partial_y
$$
Then $S^1=\{x^2+y^2=1\}$ is clearly invariant under $G$. We can also see this using \ref{algsym}:
let $F(x,y):=x^2+y^2-1$. Then $F$ is of maximal rank ($=1$) on its zero set $S^1$ and
$v(F)= -2xy + 2xy = 0$.

(ii) Let $H(x,y)=y^2-2y+1$. Then the zero set of $H$ is the horizontal line $\{y=1\}$, which
is manifestly not invariant under $G$. Nevertheless, $v(H)=2xy-2x=2x(y-1)$ vanishes on this
zero set. However, also $TH=(0,2y-2)$ vanishes there, so \ref{algsym} does not apply.
This demonstrates that the maximal rank condition in \ref{algsym} cannot be dropped.
\et
The following criterion will turn out to be useful many times later on:
\bt\label{hadamard} (Hadamard's Lemma)\index{Hadamard's Lemma} Let $F:M^m\to \R^l$ ($l\le m$) be 
of maximal rank on the subvariety $S_F:=\{x\in M\mid F(x)=0\}$. Then a smooth real-valued
function $f\colon M\to \R$ vanishes on $S_F$ if and only if there exist $Q_1,\dots,Q_l\in \cinfty(M,\R)$
such that
$$
\forall x\in M:\ f(x) = Q_1(x)F_1(x)+\dots+Q_l(x)F_l(x).
$$
\et
\pr The condition is clearly sufficient. Conversely, suppose that $f$ vanishes on $S_F$ and
let $x\in S_F$. Since $F$ is a submersion, by \cite[3.3.7]{KAoM} there exists a chart 
$(\vphi=(y^1,\dots,y^m),U)$ centered at $x$ such that $F_\vphi(y):=F\circ \vphi^{-1}(y^1,\dots,y^m) = (y^1,\dots,y^l)$.
It follows that $f_\vphi:=f\circ\vphi^{-1}$ vanishes on $\vphi(U)\cap (\{0\}\times \R^{m-l})$.
We can suppose that $\vphi(U)$ is a ball around $0$ in $\R^m$. Then for $y=(y',y''):=(y^1,\dots,y^l,y^{l+1},\dots,y^m)$
we have
\begin{equation*}
\begin{split}
f_\vphi(y',y'')&=f_\vphi(y',y'')-f_\vphi(0,y'')=\int_0^1 \partial_t f_\vphi(ty',y'')\,dt\\
&= \sum_{i=1}^l \int_0^1 \frac{\partial f_\vphi}{\partial y_i}(ty',y'')\cdot y^i\,dt
=: \sum_{i=1}^l Q_\vphi^i(y)\cdot y^i = \sum_{i=1}^l Q_\vphi^i(y)\cdot (F_i)_\vphi(y).
\end{split}
\end{equation*}
Therefore, setting $Q^i_U:=Q_\vphi^i\circ \vphi\in \cinfty(U)$ we obtain $f|_U=\sum_{i=1}^l Q_U^i\cdot F_i|_U$.

On the other hand, if $x\not\in S_F$ then there exists some open neighborhood $U$ of $x$ and some $i\in \{1,\dots,l\}$
such that $F_i(y)\not=0$ for all $y\in U$. Then we set $Q^i_U:=f/F_i$ and $Q^j_U:=0$ for $j\not=i$ to again
obtain $f|_U=\sum_{i=1}^l Q_U^i\cdot F_i|_U$.
Now pick a locally finite open cover $(U_\alpha)_{\alpha\in A}$ of neighborhoods as above and
a subordinate partition of unity $(\chi_\alpha)_{\alpha\in A}$ with $\supp(\chi_\alpha)\sse U_\alpha$ for all $\alpha$.
Then
$$
f = \sum_\alpha \chi_\alpha f = \sum_\alpha \chi_\alpha \sum_{i=1}^l Q^i_{U_\alpha}F_i =: \sum_{i=1}^l Q_i F_i.
$$
\ep
It follows that an equivalent reformulation of \eqref{algsymcrit} is that for any $v\in \g$ 
there exist smooth functions $Q_{\nu\mu}: M\to \R$ ($\mu,\nu=1,\dots,l$) such that
\begin{equation}\label{algsymcrit2}
v(F_\nu)(x) = \sum_{\mu=1}^l Q_{\nu\mu}(x) F_\mu(x) \quad \nu=1,\dots,l,\quad x\in M.
\end{equation}
\bex The functions $Q_\nu$ in \ref{hadamard} are in general not unique: Let $M=\R^3$ and
$F(x,y,z):=(x,y)$. Then the function $f(x,y,z)=xz+y^2$ vanishes on $S_F$ (which is the $z$-axis) and
we have
$$
f=zF_1 + yF_2 = (z-y)F_1 + (x+y)F_2.
$$
\et
If $f = \sum_\nu Q_\nu F_\nu = \sum_\nu \tilde Q_\nu F_\nu$, then the differences $R_\nu:=Q_\nu-\tilde Q_\nu$
satisfy the homogeneous system
\begin{equation}\label{homsys}
\sum_{\nu=1}^l R_\nu(x) F_\nu(x) = 0.
\end{equation}
For such functions we have:
\bp\label{hadex} Let $F:M^m\to \R^l$ be of maximal rank on $S_F=\{x\in M\mid F(x)=0\}$. If $R_1,\dots,R_l\in \cinfty(M,\R)$
satisfy \eqref{homsys}, then each $R_\nu$ vanishes identically on $S_F$. Equivalently, there exist
$S^\mu_\nu\in \cinfty(M,\R)$ ($\nu,\mu = 1,\dots,l$) such that, for all $x\in M$ and all $\nu=1,\dots,l$
\begin{equation}\label{homsysmn}
R_{\nu} = \sum_{\mu=1}^l S^\mu_\nu F_\mu.
\end{equation}
In addition, the $S^\mu_\nu$ can be chosen to be anti-symmetric: $S^\mu_\nu = -S^\nu_\mu$. In this case,
\eqref{homsysmn} is necessary and sufficient for \eqref{homsys}.
\et
\pr As in the proof of \ref{hadamard}, by using a partition of unity and suitable charts we may 
reduce the proof to neighborhoods of points with either $x\in S_F$ or $x\not\in S_F$.

1) If $x\in S_F$, using a chart as in \ref{hadamard}, we may suppose that
$M$ is a ball in $\R^m$ and $F=(x^1,\dots,x^m)\mapsto (x^1,\dots,x^l)$. Then 
\eqref{homsys} reads
\begin{equation}\label{homsys2}
\sum_{\nu=1}^l R_\nu(x) x^\nu = 0.
\end{equation}
We show by induction that \eqref{homsys2} implies the existence of a skew-symmetric matrix 
$S^\mu_\nu$ of smooth functions satisfying \eqref{homsysmn}.

For $l=1$, \eqref{homsys2} reduces to $R_1(x)x^1=0$ for all $x$, so by continuity $R_1(x)=0$ for all
$x$, and we can choose $S^1_1=0$.

$l-1 \rightarrow l$: Setting $x':=(x^{l+1},\dots,x^m)$, \eqref{homsys2} implies that 
$R_l(0,\dots,0,x^l,x')\cdot x^l=0$, so in fact
$R_l(0,\dots,0,x^l,x')=0$ for all $(x^l,x')$. Applying \ref{hadamard} to $F:=(x^1,\dots,x^m) \mapsto (x^1,\dots,x^{l-1})$,
it follows that there exist smooth functions $S_l^\mu$ ($\mu=1,\dots,l-1$) of $x$ with
$R_l(x) = \sum_{\mu=1}^{l-1} S_l^\mu(x)x^\mu$. By \eqref{homsys2}, therefore,
$$
0 = \sum_{\mu=1}^{l-1} (R_\mu(x)x^\mu + S_l^\mu(x)x^\mu x^l) = \sum_{\mu=1}^{l-1} (R_\mu(x)+S_l^\mu(x)x^l)x^\mu
=:\sum_{\mu=1}^{l-1} \tilde R_\mu x^\mu.
$$
By our induction hypothesis there exist smooth functions $\tilde S^\mu_\nu$ ($1\le \mu,\nu\le l-1$) such that $\tilde R_\nu = 
\sum_{\mu=1}^{l-1}\tilde S^\mu_\nu x^\mu$ for $1\le \nu\le l-1$ and such that $\tilde S^\mu_\nu$
is skew-symmetric. By definition, $R_\mu = \tilde R_\mu - S_l^\mu x^l$ for $1\le \mu\le l-1$.
Therefore,
$$
\begin{pmatrix}
R_1\\
\vdots \\
R_{l-1}\\
R_l	
\end{pmatrix}
=
\left(
\begin{array}{cccc}
 \tilde S_1^1 &\dots &\tilde S_1^{l-1} &-S^l_1\\
\vdots & &  &\vdots \\
\tilde S_{l-1}^1 &\dots &\tilde S_{l-1}^{l-1} &-S^l_{l-1}	\\
 S_l^1 &\dots &S_{l}^{l-1} & 0 
\end{array}
\right)
\cdot
\begin{pmatrix}
x^1\\
\vdots \\
x^{l-1}\\
x^l	
\end{pmatrix}
=: S\cdot x,
$$
giving the desired skew-symmetric $l\times l$ matrix $S$.

2) Now suppose that $x\not\in S_F$. Again we proceed by induction. For $l=1$, $F_1(x)\not=0$
and $R_1\cdot F_1=0$ imply that we can write $R_1(x) = 0 = 0\cdot F_1=:S^1_1\cdot F_1$ on a neighborhood of 
$x$ where $F_1\not=0$.

Next, assume that the result is already proved for $l-1$. Since $F(x)\not=0$, some component of
$F(x)$ must be non-zero, and without loss of generality we may suppose that $F_{l-1}(x)\not=0$
(the proof below works the same way for any other component as well). Then on a neighborhood $U$ of $x$
where $F_{l-1}(x)\not=0$ we have
$$
R_1F_1+\dots+R_{l-2}F_{l-2} + \Big(R_{l-1} + R_l\frac{F_{l}}{F_{l-1}}\Big)F_{l-1} = 0,
$$
and so by our induction assumption (and shrinking $U$ if necessary) there exists a skew-symmetric $(l-1)\times (l-1)$ 
matrix $\hat S^\mu_\nu$ such that
$$
\begin{pmatrix}
R_1\\
\vdots \\
R_{l-2}\\
R_{l-1} + R_l\frac{F_{l}}{F_{l-1}}
\end{pmatrix}
=
\left(
\begin{array}{ccc}
  \hat S_1^1 &\dots & \hat S_1^{l-1} \\
\vdots &   & \vdots \\
\hat S_{l-2}^1 &\dots &\hat S_{l-2}^{l-1} 	\\
 \hat S_{l-1}^1 &\dots &\hat S_{l-1}^{l-1}  
\end{array}
\right)
\cdot
\begin{pmatrix}
F_1\\
\vdots \\
F_{l-2}\\
F_l	
\end{pmatrix}.
$$
Consequently, 
$$
\begin{pmatrix}
R_1\\
\vdots \\
R_{l-2}\\
R_{l-1}\\
R_l
\end{pmatrix}
=
\left(
\begin{array}{cccc}
  \hat S_1^1 &\dots  & \hat S^{l-1}_1 & 0 \\
\vdots &  & \vdots  & \vdots \\
\hat S_{l-2}^1 &\dots  & \hat S^{l-1}_{l-2} & 0	\\
\hat S^{1}_{l-1} & \dots & \hat S^{l-1}_{l-1} &  -R_{l}/F_{l-1} \\
0 &\dots  & 0\ \  R_{l}/F_{l-1} & 0
\end{array}
\right)
\cdot
\begin{pmatrix}
F_1\\
\vdots \\
F_{l-2}\\
F_{l-1}\\
F_l	
\end{pmatrix},
$$
which gives the claim also in this case.

Conversely, if \eqref{homsysmn} is satisfied, then denoting by $\langle \,,\,\rangle$ the standard
scalar product on $\R^l$ we have 
$$
\langle R(x),F(x)\rangle = \langle S(x)F(x),F(x)\rangle = 0
$$
since $S^\top=-S$. This is \eqref{homsys}.
\ep
\section{Invariance and functional dependence}
Our goal in this section is to determine `how many' invariants a given group action possesses.
We first observe that if $\zeta^1,\dots,\zeta^k$ are invariants (either local or global) of a
group $G$ and $F=F(z^1,\dots,z^k)$ is any smooth function then also 
$$
\zeta(x):=F(\zeta^1(x),\dots,\zeta^k(x))
$$
is an invariant. However, $\zeta$ is completely determined by $\zeta^1,\dots,\zeta^k$ and
so adds no additional information. We therefore want to determine invariants up to this
kind of relation.
\bd\label{funindef} Let $\zeta^1,\dots,\zeta^k\in \cinfty(M,\R)$. Then
\begin{itemize}
\item[(i)] $\zeta^1,\dots,\zeta^k$ are called functionally dependent\index{functionally dependent}
if for each $x\in M$ there exists a neighborhood $U$ of $x$ and a smooth map $F:\R^k\to \R$
that does not vanish identically on any open subset of $\R^k$ such that
\begin{equation}\label{funcind}
F(\zeta^1(x),\dots,\zeta^k(x)) = 0
\end{equation}
for all $x\in U$.
\item[(ii)] $\zeta^1,\dots,\zeta^k$ are called functionally independent\index{functionally independent}
if they are not functionally dependent when restricted to any open subset $U$ of $M$.
Equivalently, if $F$ as in (i) is such that \eqref{funcind} holds for all $x$ in some open subset
$U$ of $M$ then $F(z_1,\dots,z_k)=0$ for all $z$ in some open subset of $\R^k$ (which contains
$(\zeta^1,\dots,\zeta^k)(U')$ for some $U'\sse U$).
\end{itemize}
\et
\bex (i) The functions $x/y$ and $xy/(x^2+y^2)$ are functionally dependent on $\{(x,y)\mid y\not=0\}$
because
$$
\frac{xy}{x^2+y^2} = \frac{x/y}{1+(x/y)^2} = f(x/y).
$$
On the other hand, $x/y$ and $x+y$ are functionally independent where defined since if
$F(x+y,x/y)\equiv 0$ for $(x,y)$ in any open subset of $\R^2$ then since $(x,y)\mapsto (x+y,x/y)$
is a local diffeomorphism, $F$ has to vanish identically on some open subset of the image of this set.

(ii) The functions
$$
\eta(x,y) = x, \quad \zeta(x,y) = \left\{ \begin{array}{cl}
x & y\le 0\\
x+e^{-1/y} & y>0 
\end{array}
\right.
$$
are functionally dependent on $\{y<0\}$, independent on $\{y>0\}$, but {\em neither} on
the entire space $\R^2$.
\et
For a characterization of functional (in)dependence we need two auxiliary results, both of
which are of independent interest as well.
\bp\label{ccc} Let $M$ be a second countable smooth manifold and let $S$ be a closed subset of $M$.
Then there exists a smooth function $f\colon M\to \R$ such that $S=\{x\in M\mid f(x)=0\}$.
\et
\pr Since $M\setminus S$ is open and $M$ is second countable, $M\setminus S$ can be written as a 
locally finite union of countably many sets $V_j$ ($j\in \N$) and there exist 
functions $\chi_j\in \cinfty(M)$ with $\chi_j\ge 0$ and $\chi_j>0$ precisely on $V_j$
(see the proof of \cite[1.3.10]{KAoM}). Then set
$$
f:= \sum_{j\in\N} \chi_j.
$$
Since $(V_j)_j$ is locally finite, locally around any point in $M$ only finitely many
summands are nonzero, so $f\in \cinfty(M)$. Moreover, for any $x\in M$ we have
$$
f(x)=0\Leftrightarrow \forall j: \chi_j(x)=0 
\Leftrightarrow x\in \bigcap_{j\in\N}(M\setminus V_j) = M\setminus \bigcup_{j\in \N}V_j = S.
$$
\ep
In the next theorem, by a {\em critical point}\index{critical point} of a smooth map $f\colon M\to N$
we mean a point $x\in M$ where $T_xf$ is not surjective. A point $y\in N$ is called a 
{\em critical value}\index{critical value} of $f$ if there exists a critical point $x\in M$
of $f$ such that $y=f(x)$.
\bt\label{sard}\index{Sard's Theorem}{\em (Sard's Theorem)} Let $f\colon M\to N$ be a smooth map between smooth manifolds. Then
the set of critical values of $f$ in $N$ has measure zero, in the sense that in any chart of $N$
the image of the set of critical values in that chart domain has Lebesgue measure zero.
\et
A proof of this result would take us too far afield, so we refer to \cite[3.2]{Kahn}.

Based on the previous results we now have:
\bt\label{funcindchar} Let $\zeta = (\zeta^1,\dots,\zeta^k)$ be a smooth map from $M$ to $\R^k$. Then the 
following are equivalent:
\begin{itemize}
\item[(i)] $\zeta^1,\dots,\zeta^k$ are functionally dependent on $M$.
\item[(ii)] $\forall x\in M: \ \rk(T_x\zeta)<k$.
\end{itemize}
\et
\pr (i)$\Rightarrow$(ii): Suppose that for some $x\in M$ we had $\rk(T_x\zeta)=k$. Since
the rank can't fall locally, there would then exist some open set $U$ around $x$
with $\rk(T_{x'}\zeta)=k$ for all $x'\in U$. We may assume $U$ small enough that
there is some $F$ as in \ref{funindef} with $F(\zeta^1(x'),\dots,\zeta^k(x')))\equiv 0$ on $U$.
Since $\zeta$ is a submersion on $U$, it is an open map. Thus $F$ vanishes on the open set
$\zeta(U)$, a contradiction. 

(ii)$\Rightarrow$(i): 
Let $U$ be open and relatively compact in $M$. By \ref{ccc} there exists some $F\in \cinfty(\R^k)$ 
such that $F(z)=0$ if and only if $z\in \zeta(\overline U)$. Since $\zeta(\overline U) \sse 
\{\zeta(x)\mid \rk(T_x\zeta)<k\}$ and the latter set has measure zero by Sard's 
theorem \ref{sard}, $\zeta(\overline U)$ does not contain any non-empty
open set. Thus $\zeta^1,\dots,\zeta^k$ are functionally dependent.
\ep
\bt\label{maininv} Let $G$ act semi-regularly on the $m$-dimensional manifold $M$ with $k$-dimensional orbits and let 
$x_0\in M$. Then there exists some open neighborhood $U$ of $x_0$ such that there are precisely 
$m-k$ functionally independent local invariants $\zeta^1,\dots,\zeta^{m-k}$ of $G$ on $U$. Any other
invariant of $G$ on $U$ is of the form 
$$
f(x)=F(\zeta^1(x),\dots,\zeta^{m-k}(x))
$$ 
for some smooth function $F$. If the action of $G$ is regular then the invariants can be taken to be
globally invariant on a neighborhood of $x_0$.
\et
\pr By \ref{constrk}, any orbit of $G$ is a leaf of the $k$-dimensional integrable distribution $\Delta_\Phi$. 
Thus the classical Frobenius theorem \ref{classfrob} yields a cubic chart $(\vphi=(x^1,\dots,x^m),U)$
centered at $x_0$, $\vphi(U)=(-c,c)^m$, such that any orbit of $G$ that intersects $U$ does so in a 
(by \ref{locfol} at most countable) union of slices $U_a=\vphi^{-1}(\R^k\times \{a\})$. 
For $1\le i\le m-k$, let $\zeta^i:U\to \R$, $\zeta^i:=x^{k+i} (=\vphi^{k+i})$. Then by definition, each 
$\zeta^i$ is constant on each slice, hence is locally constant on any orbit that intersects $U$,
hence is a local invariant. Any other local invariant $f$ must be constant on each slice. In terms
of the chart $\vphi$ this means that $f$ does in fact not depend on the variables $x^1,\dots,x^k$,
and therefore is a function of the remaining variables $x^{k+i}$ ($1\le i\le k$) only, i.e., of the $\zeta^i$.

Finally, if the action of $G$ is regular, then by \ref{regslice} the chart $(\vphi,U)$ can be chosen such that
each $G$-orbit intersects $U$ in at most one slice. Then the $\zeta^i$ are constant on this slice, i.e., on 
$U\cap S$, so they are global invariants.
\ep
Since $\zeta=(\zeta^1,\dots,\zeta^{m-k})$ has maximal rank, the $\zeta^i$ are functionally independent 
by \ref{funcindchar}. Such a family of invariants (i.e., a functionally independent family such that any
other invariant is a function of the members of the family) is called a 
{\em complete set of functionally independent invariants}.
\bp\label{lociprop} Let $G$ act semi-regularly on the $m$-dimensional manifold $M$ with $k$-dimensional orbits.
Then: 
\begin{itemize}
\item[(i)] Any complete set of functionally independent local invariants of $G$ has $m-k$ elements.
\item[(ii)] If $\eta^1,\dots,\eta^{m-k}$ is a set of functionally independent local invariants, 
then locally around any point in its domain it is complete.
\end{itemize}
\et
\pr (i) Let $\eta=(\eta^1,\dots,\eta^l)$ and $(\theta^1,\dots,\theta^p)$ be two 
sets of functionally independent invariants of $G$ and assume that $l>p$. Since
$\theta$ is complete, there exists some smooth map $F:\Omega\to \R^l$ with $\Om\sse\R^p$
open such that $\eta = F\circ \theta$. But then for any $x$, $\rk(T_x\eta)\le \rk T_{\theta(x)}(F)\le p<l$,
a contradiction to \ref{funcindchar}. By symmetry, $l=p$, and by \ref{maininv} it follows that in fact $l=m-k$.

(ii) Suppose that $f$ is some invariant defined on the domain $\Omega$ of $\eta=(\eta^1,\dots,\eta^{m-k})$
and let $x\in \Om$. Also, pick $\zeta = (\zeta^1,\dots,\zeta^{m-k})$ around $x$ as in \ref{maininv}.
Then there exists some $F:\R^{m-k}\to\R^{m-k}$, such that $\eta = F\circ \zeta$ near $x$. Since
$\rk(\eta)=\rk(\zeta)=m-k$, we also must have $\rk(F)=m-k$, i.e., $F$ is a local diffeomorphism.
Since $f$ is a smooth function of $\zeta$ near $x$ it is therefore also a smooth function of $\eta$. 
\ep
\bp\label{invvar} 
Let $G$ act semi-regularly on the $m$-dimensional manifold $M$ with $k$-dimensional orbits
and let $\zeta^1,\dots,\zeta^{m-k}$ be a complete set of functionally independent invariants defined
on an open subset $U$ of $M$. If a subvariety $S_F=\{x\in M\mid F(x)=0\}$ (with $F\in\cinfty(M,\R^l)$) 
is $G$-invariant, then for each
solution $x_0\in S_F\cap U$ there is a neighborhood $\tilde U\sse U$ of $x_0$ and a smooth function $f$
such that the solution set of the corresponding invariant $\tilde F = x\mapsto f(\zeta^1(x),\dots,\zeta^{m-k}(x))$
on $\tilde U$ coincides with that of $F$, i.e.,
$$
S_F\cap \tilde U = S_{\tilde F}\cap \tilde U = \{x\in \tilde U\mid f(\zeta^1(x),\dots,\zeta^{m-k}(x))=0\}.
$$
\et
\pr Since (by \ref{funcindchar}) the rank of $(\zeta^1,\dots,\zeta^{m-k})$ is $m-k$ on $U$ we may find
coordinates $(y^1,\dots,y^m)$ on a neighborhood $\tilde U$ of $x_0$ such that $y^i=\zeta^i$ for
$1\le i\le m-k$. In fact, the remaining $k$ coordinates can be determined by picking suitable
$x^{i_j}$ (called {\em parametric variables}) from any given coordinates $(x^1,\dots,x^m)$. 
The change of coordinates is then of the form $y=\psi(x)=(\zeta(x),\hat x)$, with $\hat x$
the parametric variables.

We show that these coordinates are then flat for $\Delta_\Phi$: Each $\zeta^j$ is constant
on the orbits of $\Delta_\Phi$, hence if $v_1,\dots,v_k\in \g$ are such that 
$\Delta_\Phi(x)=\text{span}(v_1|_x,\dots,v_k|_x)$ then 
$$
T_x{\zeta^j}(v_i|_x)=v_i(\zeta^j)|_x = \ddt\zeta^j(\Fl^{v_i}_t(x))=0
$$
for all $x\in \tilde U$. Therefore, $v_i|_x \in \bigcap_{j=1}^k \ker(T_x\zeta^j) = \ker T_x\zeta$, and
since $\dim(\ker(T_x\zeta))=k$, $\Delta_\Phi(x)=\text{span}(v_1|_x,\dots,v_k|_x)=\ker T_x\zeta$.
On the other hand, $\ker(T_x\zeta)$ is spanned by $\partial_{y^{m-k+1}},\dots,\partial_{y^{m}}$ because
$T\zeta^i(\partial_{y^{j}})= \frac{\partial y^i}{\partial y^j}=0$ for $i\in \{1,\dots,m-k\}$ and
$j\in \{m-k+1,\dots,m\}$.

In terms of the new variables we can write $F(x) = F^*(y) = F^*(\zeta(x),\hat x)$, where $F^*:=F\circ\psi^{-1}$.
Denoting by $\hat x_0$ the value of the parametric variables at $x_0$ we now set
$f(z):=F^*(z,\hat x_0)$, and
$$
\tilde F(x) := f(\zeta(x)) = f(\zeta^1(x),\dots,\zeta^k(x)).
$$
By assumption, $S_F$ is $G$-invariant. Also, the orbits of $G$ intersect $U$ in the slices $\{\zeta(x)=c\}$.
Therefore, $F(x)=F^*(\zeta(x),\hat x)=0$ if and only if $\tilde F(x) = F^*(\zeta(x),\hat x_0)=0$
since $(\zeta(x),\hat x)$ and $(\zeta(x),\hat x_0)$ lie in the same slice.
\ep
We now turn to the problem of actually calculating invariants of a local transformation group.
To begin with, we consider the case of a one-parameter group $G$ with infinitesimal generator $v\in \Xloc(M)$.
In terms of local coordinates $x^1,\dots,x^n$ we can write
$$
v = \xi^1(x)\partial_{x^1}+\dots+\xi^m(x)\partial_{x^m}
$$
for certain smooth functions $\xi^i$. Then by \ref{finvcrit}, a local invariant $f$ of $G$ is a solution
of the linear homogeneous first order PDE
\begin{equation}\label{invpde}
v(f) = \xi^1(x)\frac{\partial f}{\partial x^1} +\dots+\xi^m(x)\frac{\partial f}{\partial x^m} = 0.
\end{equation}
By \ref{maininv} we know that if $v|_x\not=0$ then locally around $x$ there exist $m-1$ functionally 
independent invariants, i.e., $m-1$ functionally independent solutions of \eqref{invpde}. In PDE
theory, the method of choice for solving \eqref{invpde} is the 
{\em method of characteristics}.\index{method of characteristics}, cf.\ \cite[Sec.\ 3.2]{Ev}.
This basically consists in solving
the corresponding system of ODEs given by
\begin{equation}\label{invode}
\frac{dx^i}{dt} = \xi^i(x(t)) \quad (1\le i\le m).
\end{equation}
Solutions of \eqref{invpde} are then functions that are constant along solutions of \eqref{invode}.
Such functions are also called {\em first integrals}\index{first integral} of \eqref{invode}. In this terminology, \ref{maininv}
says that \eqref{invode} possesses $m-1$ functionally independent first integrals locally around any point.

\brem Using the straightening-out theorem (cf.\ \cite[17.12]{KLie_new}) it is very easy to see all 
of the above directly. In fact, in an appropriate coordinate system $y^1,\dots,y^m$ we have $v=\partial_{y^1}$
near the point $x$. It is then immediate that the functions $y^2,\dots,y^m$ form a complete set of
functionally independent invariants of $\Fl^v$. Furthermore, in these coordinates \eqref{invode}
reads $\frac{dy^1}{dt} = 1$, and $\frac{dy^i}{dt} = 0$ for $i=2,\dots,m$, with general solution
$y^1(t) = t + c_1$, $y^i=c_i$ for $i=2,\dots,m$. A complete set of functionally independent first 
integrals then of course is also given by $y^2,\dots,y^m$, because these functions are constant along this general solution
of the ODE.
\et
\bex Consider the rotation group $SO(2)$ on $\R^2\setminus \{(0,0)\}$ with infinitesimal generator $v=-y\partial_x +x\partial_y$.
Then the characteristic system reads
$$
\frac{dx}{dt} = -y \quad \frac{dy}{dt} =x.
$$
A first integral is $f(x,y)=x^2+y^2$, since 
\begin{equation*}
\begin{split}
\frac{d}{dt}f(x(t),y(t)) &= \partial_x f(x(t),y(t))(-y(t)) + \partial_y f(x(t),y(t))(x(t))\\ 
&= -2x(t)y(t)+2x(t)y(t)=0.
\end{split}
\end{equation*}
Any other first integral, i.e., any other invariant of $v$, is a function of $f$.
\et
In the remainder of this section we want to explore how invariance of functions or subvarieties
under a local transformation group can be expressed in terms of quotient manifolds, using
\ref{quotient}. The main result is as follows:
\bt\label{quotob} Let $G$ be a local transformation group on $M$ that acts regularly on the $m$-dimensional
manifold $M$ with $k$-dimensional orbits. 
\begin{itemize}
\item[(i)] A smooth map $F: M\to \R^l$ is $G$-invariant if and only if there exists a smooth 
function $\tilde F: M/G \to \R^l$ such that $F(x)=\tilde F(\pi(x))$ for all $x\in M$. 
\item[(ii)] Let $F\in \cinfty(M,\R^l)$. Then the corresponding subvariety $\mcS_F=\{x\in M\mid F(x)=0\}$
is $G$-invariant if and only if there exists a smooth map $\tilde F:M/G\to \R^l$ such that,
for all $x\in M$, $F(x)=0$ if and only if $\tilde F(\pi(x))=0$.
\item[(iii)] An $n$-dimensional regular submanifold $N$ of $M$ is $G$-invariant if and only
if there exists a smooth $(n-k)$-dimensional regular submanifold $\tilde N = N/G$ of $M/G$
such that $N=\pi^{-1}(\tilde N)$ (and therefore $\tilde N=\pi(N)$).
\end{itemize}
\et
\pr (i) By definition, $\pi$ is invariant under $G$, so any $\tilde F\circ\pi$ is $G$-invariant as well.
Conversely, suppose that $F\in \cinfty(M,\R^l)$ is $G$-invariant. Using a regular chart $(U,\vphi)$,
the proof of \ref{maininv}, together with \ref{quotientconst} shows that on $\pi(U)$ there exists
a smooth map $\tilde F_U$ such that $\tilde F_U\circ\pi= F$ on $U$. The sets $\pi(U)$ 
form a covering of $M/G$ by chart neighborhoods (since $\pi(U)=U'$), so we may find a partition
of unity $(\chi_j)$ on $M/G$ subordinate to it, say with $\supp \chi_j\comp U_j'$ for all $j$. 
Then setting $\tilde F:= \sum_j \chi_j\cdot \tilde F_{U_j}$
it follows that, for any $x\in M$,
$$
\tilde F(\pi(x)) =  \sum_j \chi_j(\pi(x))\cdot \tilde F_{U_j}(\pi(x)) =
\sum_j \chi_j(\pi(x))\cdot F(x) = F(x).
$$
(ii) Suppose first that there exists some $\tilde F\in \cinfty(M/G,\R^l)$ such that $F(x)=0$ if
and only if $\tilde F(\pi(x))=0$. Then if $y$ is in the same orbit as $x$, $\pi(x)=\pi(y)$.
Thus, if $F(x)=0$ then $0=\tilde F(\pi(x))=\tilde F(\pi(y))$, which in turn is equivalent to
$F(y)=0$. It follows that $\mcS_F$ is $G$-invariant.

Conversely, if $\mcS_F$ is $G$-invariant then by \ref{invvar} (together with \ref{quotientconst}) it 
follows that any point in $M$ has a neighborhood $U$ such that there exists some smooth map
$\tilde F_U$ on $U'=\pi(U)$ with the property that, for all $x\in U$, $F(x)=0$ if and
only if $\tilde F_U(\pi(x))=0$. Let $\tilde F_U(z)=(\tilde F_1(z),\dots,\tilde F_{l}(z))$. Then 
we may replace $\tilde F_i$ by $(\tilde F_i)^2$ for all $i\in \{1,\dots,l\}$ while retaining
the property stated above. Hence without loss of generality we may assume that each $\tilde F_i$
is non-negative. As in (i), pick a partition
of unity $(\chi_j)$ on $M/G$ subordinate to a covering by such neighborhoods $U_j'$ ($j\in \N$) and set 
$\tilde F:= \sum_j \chi_j\cdot \tilde F_{U_j}$. Let $x\in M$ with $F(x)=0$ and suppose that
$\pi(x)\in U_j'=\pi(U_j)$. Then for some $x_j\in U_j$ we have $\pi(x)=\pi(x_j)$, so since $S_F$
is $G$-invariant we obtain $x_j\in S_F$ as well. Hence $0=\tilde F_{U_j}(\pi(x_j))=
\tilde F_{U_j}(\pi(x))$ for any such $j$, implying $\tilde F(\pi(x))=0$.
Conversely, suppose that
$\tilde F(\pi(x))=0$ and pick $k$ such that $\chi_k(\pi(x))>0$. Then since
each $\tilde F_{U_j}^ i$ is non-negative, it follows
that $\tilde F_{U_k}(\pi(x))=0$. Since $U_k'=\pi(U_k)$, there exists some $x_k\in U_k$
with $\pi(x_k)=\pi(x)$. By construction of $\tilde F_{U_k}$, therefore, $F(x_k)=0$, i.e., $x_k\in S_F$.
But then also $F(x)=0$ since $S_F$ is $G$-invariant.

(iii) Suppose that $N = \pi^{-1}(\tilde N)$ and let $x\in N$ and $g\in G$ such that $y:=g\cdot x$
exists. Then $\pi(y)=\pi(x)\in \tilde N$, so $y\in \pi^{-1}(\tilde N)=N$. Thus $N$ is $G$-invariant.

Conversely, if $N$ is a $G$-invariant regular submanifold of $M$ and $\Phi\colon\cU \to M$ denotes the
action of $G$ on $M$ (cf.\ \ref{ltgdef2}), then $\cU':=\cU\cap (G\times N)$ is open in $G\times N$
and contains $\{e\}\times N$. Also, $\Phi\colon\cU'\to M$ is smooth. Moreover, $\Phi(\cU')\sse N$
since $N$ is $G$-invariant and because $N$ is regular it follows that $\Phi\colon\cU'\to N$ is smooth,
hence is a local transformation group. Let $x\in N$ and let $S$ be the $G$-orbit of $x$. Then 
$S$ is a regular submanifold of $M$ (by \ref{regmanrem}) and is contained in $N$. Thus by 
\cite[13.6]{KLie_new} it is an immersive submanifold of $N$. Also, it carries the trace topology of
$M$, as does $N$, so in fact it carries the trace topology of $N$, i.e., it is a regular submanifold
of $N$. Furthermore, if $\mathcal B$ is a neighborhood basis of $x$ in $M$ such that $S'\cap U$
is connected in $S'$ for every orbit $S'$ and every $U\in \mathcal B$ then $\{U\cap N\mid U\in \mathcal B\}$
has the same property in $N$ (since it carries the trace topology of $M$), 
so $G$ acts regularly on $N$ as well. We are therefore in a position
to apply \ref{quotient} to $\Phi\colon\cU'\to N$ and obtain that $\pi|_N: N\to N/G$ is a surjective 
submersion and $N/G=\pi(N)=:\tilde N$ is an $(n-k)$-dimensional manifold. Since $N$ is $G$-invariant,
if $y\in \pi^{-1}(\tilde N)$ then $\pi(y)=\pi(x)$ for some $x\in N$, so $y\in N$, i.e., $N=\pi^{-1}(\tilde N)$.

It only remains to prove that $\tilde N$ is a regular submanifold of $M/G$. To see this, we
will proceed analogously to the proof of \ref{locfol}.
Let $x\in N$, let $S$ be the orbit of $x$ (so $S\sse N$),
and pick $v_1,\dots,v_k\in\g$ such that $\{v_1(x),\dots,v_k(x)\}$
is a basis of $\Delta_\Phi(x)$. Next, choose an adapted chart $(\chi=(y^1,\dots,y^m),W)$ 
around $x$ in $M$ such that $((y^1,\dots,y^n),W\cap N)$ is a chart of $N$, $\chi(x)=0\in \R^m$, and
$\chi(W\cap N)=\{y^{n+1}=\dots=y^m=0\}$. Also, we may suppose that
$v_1(x),\dots,v_k(x),\pdl{}{y^{k+1}}{x},\dots,\pdl{}{y^n}{x}$ is a basis of $T_xN$, and
$v_1(x),\dots,v_k(x),\pdl{}{y^{k+1}}{x},$ $\dots,\pdl{}{y^m}{x}$ is a basis of $T_xM$.
Let 
$$
f(t^1,\dots,t^m):=(\Fl^{v_1}_{t^1}\circ\dots \circ \Fl^{v_k}_{t^k})(\chi^{-1}(0,\dots,0,t^{k+1},\dots,t^m)).
$$ 
Then $f$ is a diffeomorphism from some neighborhood $(-d,d)^m$
of $0\in \R^m$ onto a neighborhood of $x$ in $M$, and
we take $\vphi:=f^{-1}$ on a suitable neighborhood $U$ of $x$ to obtain that
$(U,\vphi)$ is a chart of $M$ around $x$. Furthermore, since $\chi$ is an
adapted chart and $N$ is $G$-invariant,
$$
(t^1,\dots,t^n)\mapsto f(t^1,\dots,t^n,0,\dots,0)
$$
maps into $N$ and can (by making $d$ smaller if necessary) also be assumed to 
be a diffeomorphism. This means that $(\vphi,U)$ is also an adapted chart
for $N$.
 
Since $S$ is an orbit of $G$, 
$$
y\in S \Leftrightarrow \Fl^{v_1}_{t^1}\circ\dots \circ \Fl^{v_k}_{t^k}(y)\in S
$$
for all $y$ and $t^1,\dots,t^k$ where the right hand side is defined. Therefore, 
for any $y=f(t^1,\dots,t^m)\in U$ we have
\begin{equation}\label{locstru2}
y=f(t^1,\dots,t^m)\in S \Leftrightarrow f(0,\dots,0,t^{k+1},\dots,t^m)\in S.
\end{equation}
This means that $U\cap S$ is the disjoint union of connected sets of the form
$$
U_c:=\{y\in U\mid x^{k+1}(y)=c_{k+1},\dots,x^m(y)=c_m\}
$$
where $c=(c_{k+1},\dots,c_m)$ is constant. Since $G$ acts regularly, we may
suppose in addition that $U\cap S$ is connected. Then since $f(0,\dots,0)=x\in S$
it follows that the only non-empty $U_c$ is $U_0$, so $S\cap U =
\{y\in U\mid x^{k+1}(y)=\dots=x^m(y)=0\}$. 

Summing up, we have that $(\vphi,U)$ is a flat chart for $G$ on $M$ and $(\vphi,U\cap N)$
is a flat chart for $G$ on $N$. Since $G$ acts regularly, the proof of \ref{regslice}
shows that we may in addition assume that $U$ is so small that $(\vphi,U)$ is 
regular for $\Delta_\Phi$ on $M$ and $(\vphi,U\cap N)$ is regular for $\Delta_\Phi$ on $N$.
Then by \ref{quotientconst}, on $U'=\pi(U)$ we obtain a chart $\vphi'$ for $M/G$ around $\pi(x)$ by
$$
\vphi': y' \mapsto (x^{k+1}(y),\dots,x^m(y)), \quad U'\to \R^{m-k}
$$
with $y$ any element of $\pi^{-1}(y')\cap U$. By the same reasoning,
$$
y' \mapsto (x^{k+1}(y),\dots,x^n(y)), \quad U'\cap \tilde N \to \R^{n-k}
$$
(with $y$ any element of $\pi^{-1}(y')\cap U\cap N$) is a chart for $\tilde N =N/G$ around $\pi(x)$. 
This shows that $(\vphi',U')$ is a chart adapted to $\tilde N$, i.e., $\tilde N$ is a regular
submanifold of $M/G$.
\ep
\section{Dimensional analysis}\index{dimensional analysis}
In a typical physics problem, there are certain fundamental physical quantities, such as 
length, time, mass, \dots, which can all be scaled independently of each other.
Let $z^1,\dots,z^r$ denote these quantities, then this scaling can be described
by the action of a scaling group
$$
(z^1,\dots,z^r) \mapsto (\lam_1z^1,\dots,\lam_rz^r),
$$
with fixed scaling factors $\lam = (\lam_1,\dots,\lam_r)\in\R^r$. The underlying group therefore
is the $r$-th power of the multiplicative group $\R^+$.

Furthermore, there typically are certain derived quantities, such as velocity, acceleration, force, 
density, \dots, which also scale if the fundamental quantities are scaled. Call these quantities
$x=(x^1,\dots,x^m)$. Then the action of the scaling group on the derived quantities takes the form
\begin{equation}\label{derform}
\lam\cdot (x^1,\dots,x^m) = (\lam_1^{\al_{11}}\lam_2^{\al_{21}}\dots \lam_r^{\al_{r1}}x^1,\dots,
\lam_1^{\al_{1m}}\lam_2^{\al_{2m}}\lam_r^{\al_{rm}}x^m),
\end{equation}
where $\al_{ij}$, $i=1,\dots,r$, $j=1,\dots,m$ is a matrix of real numbers determined by the problem
at hand. For example, if $z^1$ is length, $z^2$ is time and $z^3$ is mass, then the action of
scaling on velocity $x^1$ and force $x^2$ are given by
$$
\lam\cdot (x^1,x^2) = (\lam_1\lam_2^{-1}x^1,\lam_1^1\lam_2^{-2}\lam_3 x^2).
$$
If a derived quantity is invariant under the corresponding scaling then it is called
{\em dimensionless}\index{dimensionless}. To describe some physical situation, one often
has functional relations of the form $F(x^1,\dots,x^m)=0$ for the derived quantities.
Such a relation is called {\em unit-free} if it remains unchanged under a rescaling
of the fundamental quantities. It turns out that such unit-free relations are often
of great physical significance. The Buckingham Pi-Theorem shows that any unit-free
relation can be written solely in terms of dimensionless quantities.
\bt\label{buckpi} (Buckingham Pi-theorem)\index{Buckingham Pi-theorem}
Let $z^1,\dots,z^r$ be fundamental physical quantities that scale independently according to
$z^i\mapsto \lambda_iz^i$. Let $x^1,\dots,x^m$ be derived quantities that scale according to
\eqref{derform} for some $(r\times m)$-matrix of constants $A=(\al_{ij})$ and 
let $s$ be the rank of $A$. Then 
\begin{itemize}
\item[(i)] There exist $m-s$ functionally independent dimensionless `power products'
\begin{equation}\label{power}
\pi^k = (x^1)^{\be_{1k}}(x^2)^{\be_{2k}}\dots(x^m)^{\be_{mk}},\quad k=1,\dots, m-s
\end{equation}
such that any other dimensionless quantity can be written as a function of 
$\pi^1,\dots,$ $\pi^{m-s}$. The columns of the matrix $B=(\beta_{jk})$ in \eqref{power} can be
taken to be any basis of $\ker A$. 
\item[(ii)] If $F(x^1,\dots,x^m)=0$ is any unit-free relation among the given derived quantities,
then there is an equivalent relation $\tilde F=0$ which can be expressed
solely in terms of the above dimensionless power products:
$$
F(x)=0 \Leftrightarrow \tilde F(\pi^1(x),\dots,\pi^{m-s}(x))=0.
$$
\end{itemize}
\et
\pr The proof will basically be an application of \ref{quotob}. As the underlying 
manifold $M$ we take the positive octant in $\R^m$, so
$$
M = \{x=(x^1,\dots,x^m)\mid x^i>0,\, i=1,\dots, m\}.
$$
Then the multiplicative group $G=(\R^+)^r$ acts globally on $M$ via \eqref{derform}.
Since the action is multiplicative, its generators are found by differentiating 
\eqref{derform} with respect to $\lam_i$ and then setting $\lam_1=\dots=\lam_r=1$:
$$
v_i = \al_{i1}x^1\pd{}{x^1} + \al_{i2}x^2\pd{}{x^2}+\dots+\al_{im}x^m\pd{}{x^m},
\quad i=1,\dots,r.
$$
It follows that the dimension of the span of $v_1,\dots,v_r$ is precisely the rank
of $A$, namely $s$, so $G$ has $s$-dimensional orbits. For a function $f$
to be a global invariant of $G$, by \ref{finvcrit} it has to satisfy
$$
v_i(f) = 0, \quad i=1,\dots,r.
$$
In particular, if $f=\pi_k$ is given by \eqref{power}, then this condition is 
satisfied if and only if the exponents $\beta_{jk}$ satisfy the system
of linear equations
\begin{equation}\label{linsys}
\sum_{j=1}^m \al_{ij}\be_{jk}=0,\quad i=1,\dots,r.
\end{equation}
Since $\rk(A)=s$, there are $m-s$ linearly independent solutions to this system.
Choosing these solutions for the $\beta_{jk}$, let $\pi:=(\pi^1,\dots,\pi^{m-s})$.
Then the Jacobian of $\pi$ at $x=(1,\dots,1)$
is the transpose of the matrix $(\be_{jk})$, $j=1,\dots,m$, $k=1,\dots,m-s$, and its rank
therefore is $m-s$. For any $a>0$ we have (as matrix-multiplication)
$$
\pi(ax^1,x^2,\dots,x^m) = \mathrm{diag}(a^{\be_{11}},\dots,a^{\be_{1(m-s)}})\cdot \pi(x^1,x^2,\dots,x^m),
$$
so since $D\pi(1,1,\dots,1)$ has a non-zero $(m-s)\times (m-s)$-sub-determinant, the same
is true of $D\pi(a,1,\dots,1)$. Analogously we can argue for varying the $x^2$-component,
and so on. Altogether, it follows that the rank of the Jacobian of $\pi$ is $(m-s)$
everywhere on $M$. Thus by (the easy direction of) \ref{funcindchar}, $(\pi^1,\dots,\pi^{m-s})$
is a functionally independent system of invariants of $G$ on $M$.

Next we show that the orbits of $G$ are precisely the level sets of the function $\pi$.
To see this, let $x$, $\xtil$ be any two points in $M$. Then by the definition of $M$ there 
exist exponents $t_j$, $j=1,\dots,m$ such that $x^j = e^{t_j}\cdot \xtil^j$ for all $j$.
Therefore
\begin{equation*}
\begin{split}
\pi^k(x) &= e^{\be_{1k}t_1}e^{\be_{2k}t_2}\dots e^{\be_{mk}t_m}(\xtil^1)^{\be_{1k}}
\dots (\xtil^m)^{\be_{mk}} = \pi^k(\xtil)\\
&\Leftrightarrow \sum_{j=1}^m\be_{jk}t_j = 0,
\end{split}
\end{equation*}
so $x$ and $\xtil$ lie in the same level set of $\pi$ if and only if the $t_j$ satisfy
\begin{equation}\label{orbiteq}
\sum_{j=1}^m\be_{jk}t_j = 0, \quad k=1,\dots,m-s.
\end{equation}
By construction, the columns of the $m\times (m-s)$-matrix $B:=(\beta_{jk})$ form a basis of
$\ker A$, so in particular $A\cdot B=0$. Moreover, \eqref{orbiteq} means that
$B^\top\cdot t = 0$, where $t=(t^1,\dots,t^m)^\top$, i.e., $t\in \ker(B^\top)$.
We claim that this, in turn, is equivalent to the existence of $s_1,\dots,s_r$
with 
\begin{equation}\label{orbiteq2}
t_j=\sum_{i=1}^r s_i\al_{ij}, \quad j=1,\dots,m
\end{equation}
i.e., to $t\in \text{im}(A^\top)$. This means we have to show that $\text{im}(A^\top)=\ker(B^\top)$.
Now if $w=A^\top u\in \text{im}(A^\top)$, then $B^\top w = B^\top A^\top u = (AB)^\top u =0$,
so $\text{im}(A^\top)\sse\ker(B^\top)$. On the other hand, $\dim \text{im}(A^\top) = \dim \text{im}(A)=s$,
and also $\dim\ker(B^\top)=m-\dim \text{im}(B^\top)= m-\dim \text{im}(B) = s$, so indeed we have equality.

Hence two points $x$ and $\xtil$ in $M$ have the same image under $\pi$ if and only if
$x=\lam \xtil$, where $\lam_i=e^{s_i}$, i.e., if and only if they lie in the same orbit under $G$.
Since $\pi$ has maximal rank $m-s$ everywhere, its level sets are regular submanifolds of dimension $s$, 
so $G$ has $s$-dimensional regular submanifolds as orbits.
Moreover, the proof of \ref{invvar} 
shows that, around any point, $\pi^1,\dots,\pi^{m-s}$ 
can be completed by a set of $s$ parametric variables to
form a flat chart for the group action of $G$. In the present situation, such
parametric variables can even be chosen globally on all of $M$: indeed, pick
numbers $\beta_{jk}$, $j=1,\dots m$, $k=m-s+1,\dots,m$ that supplement the matrix 
$B$ from above to obtain an invertible $m\times m$ matrix $\tilde B$ and set
$$
\tilde \pi(x) := \Big(\prod_{i=1}^m(x^i)^{\beta_{i1}},\dots,\prod_{i=1}^m(x^i)^{\beta_{im}}\Big).
$$ 
Then the first $k$ components of $\tilde \pi$ are precisely $\pi$, and the remaining $m-s$ components
are our new parametric variables. By the same reasoning as above for $\pi$ it follows that
$\tilde \pi$ has rank $m$ everywhere on $M$, hence is a local diffeomorphism. In addition,
$\tilde \pi$ is injective: to see this, note first that $\tilde \pi(x)=\tilde \pi(y)$ implies
$x=y$ if and only if $\tilde\pi^i(x)=1$ implies $x^i=1$ for all $i$. Now if $\tilde \pi^i(x)=1$
for all $i$ then taking logarithms it follows that $\tilde B^\top \cdot (\mathrm{ln}(x^1),\dots,\mathrm{ln}(x^m))^\top=0$,
which due to the invertibility of $\tilde B$ indeed implies $x^i=1$ for all $i$. Thus we obtain the
desired global chart $\tilde \pi$ for $M$.

Any orbit intersects this chart in
a single slice ($\pi=$ const.), so the action of $G$ is regular.
Moreover, by \ref{quotientconst}, $\pi$
induces global coordinates on $M/G$, hence can be identified with the quotient map $M\to M/G$
(and $M/G$ can be identified with $\R_+^{m-s}$).
We are therefore in the position to apply \ref{quotob}. Since being $G$-invariant by definition
is the same as being dimensionless, (i) follows from \ref{quotob} (i). Moreover, an equation
is unit-free if and only if it is $G$-invariant, so (ii) follows from \ref{quotob} (ii).
\ep
\bex The energy yield of the first atomic explosion. In 1947, when the amount of energy
released by the first atomic explosion\index{atomic explosion} was still classified, G.\ Taylor calculated this
energy using dimensional analysis (see \cite{Tay}, our presentation follows \cite{BCo}).

We model the explosion by the radius $R$ of the spherical fireball emanating from the point of 
explosion. We assume that $R$ is a function
\begin{equation}\label{nuke}
R = F(x^1,x^2,x^3,x^4),
\end{equation}
where 
\begin{description}
	\item[]\hspace*{2em} $x^1=E$, the energy released by the explosion,
	\item[]\hspace*{2em} $x^2 = t$, the time elapsed since the explosion,
	\item[]\hspace*{2em} $x^3=\rho_0$, the initial ambient air density,
	\item[]\hspace*{2em} $x^4=P_0$, the initial ambient air pressure.
\end{description}
$R$ itself also is a derived quantity, so we set $x^5=R$. It follows that for this problem $m=5$.
The fundamental physical quantities needed to describe the derived quantities are
$z^1=$ length, $z^2 = $ mass, and $z^3=$ time, so $r=3$. 
The matrix $A=(\al_{ij})_{i=1,\dots,3}^{j=1,\dots,5}$ has the form
$$
A = \begin{pmatrix}
2 & 0 & -3 & -1 & 1\\
1 & 0 & 1 & 1 & 0\\
-2 & 1 & 0 & -2	& 0
\end{pmatrix}
$$
E.g., energy has dimension $\text{length}^2\times \text{mass}/\text{time}^2$, hence the first column.
The rank of this matrix is $s=3$, so $m-s=5-3=2$. Therefore, by \ref{buckpi} there are
two functionally independent dimensionless power products $\pi^1$, $\pi^2$, in terms of which any dimensionless 
quantity can be described. By \eqref{power}, these power products are of the form
\begin{equation*}
\begin{split}
\pi^1(x) &= (x^1)^{\be_{11}}(x^2)^{\be_{21}}(x^3)^{\be_{31}}(x^4)^{\be_{41}}(x^5)^{\be_{51}}\\
\pi^2(x) &= (x^1)^{\be_{12}}(x^2)^{\be_{22}}(x^3)^{\be_{32}}(x^4)^{\be_{42}}(x^5)^{\be_{52}}.
\end{split}
\end{equation*}
From the proof of \ref{buckpi} we know that the $5\times 2$ matrix $B=(\be_{jk})$
can be constructed by finding a basis of $\ker A$. It follows that
$$
B=
\begin{pmatrix}
-2 & -1\\
6 & -2\\
-3 & 1\\
5 & 0\\
0 & 5	
\end{pmatrix},
$$
so 
\begin{equation*}
\begin{split}
\pi^1(x) &= (x^1)^{-2}(x^2)^{6}(x^3)^{-3}(x^4)^{5} = \frac{P_0^5 t^6}{E^2\rho_0^3}\\
\pi^2(x) &= (x^1)^{-1}(x^2)^{-2}(x^3)(x^5)^{5} = \frac{R^5\rho_0}{Et^2}.
\end{split}
\end{equation*}
By \ref{buckpi} (ii), the relation \eqref{nuke} can equivalently be expressed in the
form $\tilde F(\pi^1,\pi^2)=0$, and we assume that in fact we can solve for $\pi^2$.
Thus there exists a smooth function $\tilde g$ of $\pi^1$ such that \eqref{nuke} is equivalent
to $\pi^2=\tilde g(\pi^1)$, i.e., to
\begin{equation}\label{nuke2}
R = \left(\frac{Et^2}{\rho_0}\right)^{1/5}g(\pi^1).
\end{equation}
for some function $g$. By continuity, $g(\pi^1)\approx g(0)$ for $\pi^1$ small, so Taylor derived the approximative
formula
$$
R = At^{2/5}, \quad \text{ with }  A = \left(\frac{E}{\rho_0}\right)^{1/5}g(0).
$$
Here, $\rho_0$ is known, and to find $g(0)$ one can plot $\log_{10} R$ versus $\log_{10} t$ using
data from experiments with conventional explosives (whose $E$ is known). Taylor then used 
motion picture records of the first atomic explosion, due to J.E.\ Mack, to plot
$5/2 \log_{10} R$ versus $\log_{10} t$ (in c.g.s.-units, but as we know any other system
of units would have given the same result!), to obtain an accurate estimate for $E$:
$$
E \approx 7,14 \times 10^{20} \text{ ergs} = 7,14 \times 10^{13} J,
$$
which corresponds to the energy release of about $16 800$ tons of TNT.
\et
\section{Groups and differential equations}
We now want to start applying symmetry methods to differential equations. To this end we consider a system
$\mcS$ of differential equations involving $p$ independent variables $x=(x^1,\dots,x^p)$ and $q$
dependent variables $u=(u^1,\dots,u^q)$. Solutions of $\mcS$ will be functions $u=f(x)$, or,
in components,
$$
u^\alpha = f^\alpha(x^1,\dots,x^p) \quad (\alpha=1,\dots, q).
$$
Henceforth we will use the convention of using Latin indices for the independent variables and
Greek indices for the dependent ones. We let $X=\R^p$ be the space of independent variables and
$U=\R^q$ the space of dependent variables. Then basically a symmetry group of $\mcS$ will be
a local transformation group $G$ acting on $M:=X\times U$ in such a way that it `transforms
solutions of $\mcS$ into solutions of $\mcS$.

We first need to clarify how such a group of transformations is to act on a function $f$. The key 
idea here is to identify $f\colon X\supseteq \Omega \to U$ with its graph
$$
\Gamma_f = \{(x,f(x))\mid x\in\Omega\}.
$$
Then $\Gamma_f$ is a $p$-dimensional submanifold of $M$. If $g\in G$ and $\Gamma_f\sse \mcu_g$,
the domain of $g$, then $g$ acts on $\Gamma_f$ by
$$
g\cdot \Gamma_f := \{(\xtil,\util) = g\cdot (x,u)\mid (x,u)\in \Gamma_f\}.
$$
Note that $g\cdot \Gamma_f$ will in general no longer be the graph of a function. However,
since $G$ acts smoothly and $e\cdot \Gamma_f = \Gamma_f$, by shrinking $\Omega$ and restricting
to elements of $G$ near to $e$ we can always achieve that $g\cdot \Gamma_f = \Gamma_{\ftil}$,
where $\ftil$ (corresponding to $\util=\ftil(\xtil)$) is a well-defined function. In this case we write
$\ftil = g\cdot f$ and call $\ftil$ the transform of $f$ by $g$. 
\bex\label{rotact2} Let $X=\R$, $U=\R$, so $p=q=1$ and again consider the rotation group $G=SO(2)$ on $M=X\times U=\R^2$.
The action of $G$ is given by
\begin{equation}\label{rotact}
(\xtil, \util) = \eps\cdot (x,u) = (x\cos\eps -u\sin\eps,x\sin\eps + u\cos\eps)
\end{equation}
Now if $f\colon x\mapsto u=f(x)$ is a function then $G$ acts on $f$ by rotating its graph, and in general
this rotated graph will not be the graph of a well-defined function: think of a straight line or a 
parabola. The example of the parabola in particular demonstrates that it will in general
not suffice to choose $g$ near $e$ (i.e., $\eps$ near $0$) but that one has to also restrict
the domain of $f$. But for $|\eps|$ small and restricting $f$ to a sufficiently small domain we 
do obtain a well defined function $\eps\cdot f = \tilde f$.

To explicitly calculate the action of $G$ on a linear function, let $u=f(x)=ax+b$. Then any
point in $\Gamma_f$ is of the form $(x,ax+b)$, which is rotated by $\eps$ to the point
$$
(\xtil, \util) =  (x\cos\eps - (ax+b)\sin\eps,x\sin\eps + (ax+b)\cos\eps).
$$
To determine $\util =\ftil(\xtil)$ we have to eliminate $x$ from this equation, which is
possible for $\eps$ small (so that $\cot\eps \not=a$). Then
$$
x=\frac{\xtil + b\sin\eps}{\cos\eps - a\sin\eps},
$$
and therefore
$$
\util = \ftil(\xtil) = \frac{\sin \eps + a\cos\eps}{\cos\eps - a\sin\eps}\xtil + \frac{b}{\cos\eps - a\sin\eps},
$$
which, as expected, is again a linear function.
\et
The general procedure for calculating $\ftil=g\cdot f$ from $f$ is as follows. Write
\begin{equation}\label{xiphidef}
(\xtil,\util) = g\cdot (x,u) = (\Xi_g(x,u),\Phi_g(x,u)),
\end{equation}
with $\Xi_g,\Phi_g$ smooth. The graph $\Gamma_\ftil = g\cdot \Gamma_f$ of $\ftil = g\cdot f$ is then
given parametrically by the equations (for $x\in\Omega$)
\begin{equation*}
\begin{split}
\xtil = \Xi_g(x,f(x)) &= \Xi_g\circ (\id\times f)(x)\\
\util = \Phi_g(x,f(x)) &= \Phi_g\circ (\id\times f)(x)
\end{split}
\end{equation*}
To calculate $\ftil$ we have to eliminate $x$ from these equations. For $g=e$, 
$\Xi_e\circ (\id\times f)=\id$, hence for $g$ in a neighborhood of $e$
the Jacobian of $\Xi_g\circ (\id\times f)$ is nonsingular and so by the inverse
function theorem we can locally solve for $x$:
$$
x=[\Xi_g\circ (\id\times f)]^{-1}(\xtil).
$$
Consequently,
\begin{equation}\label{gdotf}
g\cdot f = [\Phi_g\circ (\id\times f)]\circ [\Xi_g\circ (\id\times f)]^{-1}.
\end{equation}
\blem\label{grcomp} Let $f$ be a local smooth map, $f\colon \Omega\to U$, and let $g$, $h\in G$. Then
$h\cdot (g\cdot f) = (h\cdot g)\cdot f$, wherever defined.
\et
\pr Since any function is uniquely determined by its graph it suffices to note that
$$
\Gamma_{(h\cdot g)\cdot f}  = (h\cdot g) \Gamma_{f} = h\cdot (g\cdot \Gamma_{f})
= h\cdot \Gamma_{g\cdot f} = \Gamma_{h\cdot (g\cdot f)}.
$$
\ep
An important special case where $\ftil$ is automatically defined unrestrictedly is that
of a {\em projectable}\index{group action!projectable} group action. For such actions,
$\Xi_g$ is a function of $x$ only, i.e.,
$$
(\xtil,\util) = g\cdot (x,u) = (\Xi_g(x),\Phi_g(x,u)),
$$
and so $\ftil(\xtil) = [\Phi_g\circ (\id\times f)]\circ \Xi_g^{-1}(\xtil)$.
Using this terminology, we can now define:
\bd\label{dgsgdef} Let $\mcS$ be a system of differential equations. A symmetry 
group\index{symmetry group!of differential equation}
of $\mcS$ is a local transformation group $G$ acting on an open subset of $M=X\times U$ such that
whenever $f$ is a solution of $\mcS$ and $g\cdot f$ is defined for some $g\in G$ then also $g\cdot f$
is a solution of $\mcS$. (Here {\em solution} means any smooth solution defined on any subdomain $\Omega$ of $X$.)
\et
\bex (i) Let $\mcS$ consist of the equation $u_{xx}=0$. The solutions are precisely the linear functions
on $\R$. Since $G=SO(2)$ transforms linear functions to linear functions, it is a symmetry group of $\mcS$.

(ii) Let $\mcS=\{u_t=u_{xx}\}$, the heat equation\index{heat equation}, and consider the group action
$$
(x,t,u) \mapsto (x+\eps a,t+\eps b,u) \quad (\eps\in \R).
$$
This is a symmetry group of $\mcS$ since $f(x-\eps a,t-\eps b)$ is a solution of the heat equation
whenever $f$ is.
\et
One immediate advantage of knowing the symmetry group of a differential equation is that it
gives a straightforward way of constructing new solutions from known ones (sometimes even from
trivial ones, e.g., constant solutions) simply by applying the group transformations.
Our aim ultimately is to derive infinitesimal criteria for a local transformation group to 
be a symmetry group of a differential equation, which can then be used to explicitly 
calculate the full symmetry group of any given differential equation. For this purpose we
first need to develop some more technical machinery.
\section{Prolongation}\index{prolongation}
In the previous section we identified any smooth function with its graph in order to define
the action of a transformation group on it. Since we are now interested in differential equations
we need to find a way of simultaneously considering a function and its derivatives up to a
certain order. The mechanism we are going to employ is that of prolongation, a simplified
version of the theory of jet bundles.

Given a smooth function $f\colon \R^p\to\R$, $x\mapsto f(x^1,\dots,x^p)$ and any $k\in\N_0$, we let
$J=(j_1,\dots,j_k)$ be an unordered $k$-tuple with $1\le j_l \le p$ for all $l$. The order of
this tuple is defined as $\sharp J = k$, and the corresponding partial derivative of $f$ is
$$
\partial_Jf(x) = \frac{\partial^k f}{\partial x^{j_1}\partial x^{j_2}\dots \partial x^{j_k}}.
$$
There are
$$
p_k := \begin{pmatrix}
p+k-1\\ k	
\end{pmatrix}
$$
different possible $k$-th order partial derivatives of $f$. To see this, note that since the order
of derivatives is irrelevant we may first bring $J$ in ascending order, and then uniquely describe
any possible choice by writing $*$ for each number that appears in $J$ and $|$ to signify that we increase the number.
For example, if $p=5$ and $k=2$, the string $**||||$ corresponds to the tuple $(1,1)$, $|*|*||$ to
$(2,3)$, and $(|||*|*)$ to $(4,5)$, etc. We need $k$ stars and $p-1$ vertical lines, so our problem 
is equivalent to determining the number of possible selections of $k$ elements from a $p-1+k$-element
set, hence the above formula.

If $f\colon X\to U$ with $X=\R^p$ and $U=\R^q$, so $u=f(x)=(f^1(x),\dots,f^q(x))$, then we need $q\cdot p_k$
numbers $u^\alpha_J = \partial_Jf^\alpha(x)$ to represent all possible $k$-th order derivatives of 
$f$ at any point $x$. Therefore, we set $U_k:=\R^{q\cdot p_k}$ and write $u^\alpha_J$, with $\alpha\in
\{1,\dots,q\}$ and $J$ any unordered multi-index as above. Moreover, we let
$$
U^{(n)} := U\times U_1\times \dots \times U_n.
$$
The coordinates of $U^{(n)}$ can be used to represent all the derivatives of any map $f\colon X\to U$
of orders $0$ to $n$. The dimension of $U^{(n)}$ is
$$
q+qp_1+\dots + qp_n = q \begin{pmatrix}
p+n\\ n	
\end{pmatrix}
=: qp^{(n)}.
$$
We will denote elements of $U^{(n)}$ by $u^{(n)}$, with $q\cdot p^{(n)}$ components $u^\alpha_J$,
with $\alpha\in \{1,\dots,q\}$ and $J=(j_1,\dots,j_k)$ an unordered multiindex, $1\le j_l\le p$, 
$0\le k \le n$. For $k=0$ there is only one such multiindex, denoted by $0$ and $u^\alpha_0$
is the component $u^\alpha$ of $u$.
\bex Let $p=2$, $q=1$. Then $X=\R^2$ has coordinates $(x^1,x^2)=(x,y)$, and $U=\R$ has the coordinate
$u$. $U_1$ equals $\R^2$, with coordinates $(u_x,u_y)$, representing all first order derivatives of $u$.
Also, $U_2=\R^3$ has coordinates $(u_{xx},u_{xy},u_{yy})$, representing all second order derivatives
of $u$. In general, $U_k=\R^{k+1}$ since there are $k+1$ derivatives of $u$ of order $k$, namely
$\frac{\partial^k u}{\partial x^i \partial y^{k-i}}$, $i=0,\dots,k$. Also, $U^{(2)} = U\times U_1\times U_2
= \R^6$, with coordinates $u^{(2)} = (u;u_x,u_y;u_{xx},u_{xy},u_{yy})$, representing all derivatives of
$u$ with respect to $x$ and $y$ of order at most $2$.
\et
Now for any smooth function $f\colon  X\to U$ we define its $n$-th 
prolongation\index{prolongation!of function}
$\prol f\colon  X\to U^{(n)}$, $u^{(n)}=\prol f(x)$ by
$$
u^\alpha_J = \partial_J f^\alpha(x).
$$
Thus for any $x$, $\prol f(x)$ is a vector with $q\cdot p^{(n)}$ entries representing the values of
$f$ and all its derivatives up to order $n$ at $x$. In this sense we might also identify
$\prol f(x)$ with the $n$-th Taylor polynomial of $f$ at $x$.
\bex For $p=2$, $q=1$, we have $u=f(x,y)$. Then the second prolongation of $f$, 
$u^{(2)}=\text{pr}^{(2)}f(x,y)$ is given by
$$
(u;u_x,u_y;u_{xx},u_{xy},u_{yy}) = \left(f;\partial_x f,\partial_y f,\partial_x^2 f,\partial_{xy} f,
\partial_y^2 f\right),
$$ 
evaluated at $(x,y)$. 
\et
The space $X\times U^{(n)}$ is also called the $n$-th order jet-space\index{jet space} of 
$X\times U$, and $\prol f(x)$ is called the $n$-jet of $f$ at $x$. If $M$ is an open subset of
$X\times U$, then we set
$$
M^{(n)} := M\times U_1\times \dots \times U_n.
$$
If $u=f(x)$ is a function whose graph lies in $M$ then the $n$-th prolongation
$\prol f$ is a function whose graph lies in $M^{(n)}$.
\section{Systems of differential equations}\label{syssect}
For a system $\mcS$ of $n$-th order differential equations in $p$ independent and $q$
dependent variables we write
$$
P_\nu(x,u^{(n)}) = 0 \quad \nu=1,\dots,l,
$$
where $x=(x^1,\dots,x^p)$, $u=(u^1,\dots,u^q)$. Here the functions 
$$
P(x,\un)=(P_1(x,\un),\dots,P_l(x,\un))
$$
are assumed to be smooth, i.e., $P\in \cinfty(X\times U^{(n)},\R^l)$.
The differential equations in $\mcS$ determine where this map vanishes on
$X\times U^{(n)}$, i.e., they determine a corresponding subvariety
$$
\mcS_P = \{(x,\un)\mid P(x,\un)=0\}\sse X\times U^{(n)}
$$
of the total jet space $X\times U^{(n)}$. We may therefore identify the
system $\mcS$ with this set $\mcS_P$, and we shall do so at many 
places below.

A solution of $\mcS$ is a smooth function $u=f(x)$ such that
$$
P_\nu(x,\prol f(x)) = 0,\quad  \nu =1,\dots,l,
$$
for all $x$ in the domain of $f$. Geometrically, this means that the graph of
the prolongation $\prol f$ is contained in $\mcS_P$:
$$
\Gamma_f^{(n)} \equiv \{(x,\prol f(x))\}\sse \mcS_P = \{P(x,\un)=0\}.
$$ 
\bex Consider the Laplace equation in two variables:
\begin{equation}\label{laplace}
u_{xx} + u_{yy} = 0.
\end{equation}
Here, $p=2$, $q=1$, and $n=2$ since the equation is of second order. The coordinates of 
$X\times U^{(2)}$ are $(x,y;u;u_x,u_y;u_{xx},u_{xy},u_{yy})$, and \eqref{laplace}
describes a hyperplane $\mcS_P$ in $X\times U^{(2)}$. A function $f$ is a solution
of \eqref{laplace} if the graph of $\text{pr}^{(2)}f$ is contained in $\mcS_P$.
For example, let $f(x,y)=x^3-3xy^2$. Then 
$$
(x,y,\text{pr}^{(2)}f(x,y)) = (x,y;x^3-3xy^2;3x^2-3y^2,-6xy;6x,-6y,-6x) \in \mcS_P.
$$
\et
\section{Prolongation of group actions}
In this section, given a local transformation
group $G$ on an open subset $M$ of $X\times U$, we will construct an induced action 
$\prol G$, the $n$-th prolongation of $G$, on $M^{(n)}$. The idea is to define
$\prol G$ in such a way that it transforms the derivatives of any smooth function
$u=f(x)$ into the corresponding derivatives of the transformed function $\util = \ftil(\xtil)$.

Thus let $(x_0,u^{(n)}_0)\in M^{(n)}$ and choose any smooth function $u=f(x)$ defined in 
a neighborhood of $x_0$ whose graph lies in $M$ and such that it has the given derivatives
at $x_0$, i.e., such that $\un_0=\prol f(x_0)$. To see that such an $f$ indeed exists,
one may take it to be an appropriate Taylor polynomial:
\begin{equation}\label{taylor}
f^\alpha(x) := \sum_J \frac{u^\alpha_{J0}}{\tilde J!}(x-x_0)^J, \quad \alpha=1,\dots,q.
\end{equation}
Here the sum is over all $J=(j_1,\dots,j_k)$ with $0\le k\le n$, and
$$
(x-x_0)^J:=(x^{j_1}-x_0^{j_1})(x^{j_2}-x_0^{j_2})\dots (x^{j_k}-x_0^{j_k}).
$$
Also, for given $J$ we set $\tilde J:=(\tilde j_1,\dots,\tilde j_p)$, where $\tilde j_i$
equals the number of $j_\kappa$'s which equal $i$. For example, if $J=(1,1,1,2,4,4)$, $p=4$,
$k=6$, then $\tilde J = (3,1,0,2)$. Finally, $\tilde J = \tilde j_1!\dots\tilde j_p!$.

If $g$ is close to $e$ then the transformed function $\ftil = g\cdot f$ (as in \eqref{gdotf})
is defined in a neighborhood of $(\xtil_0,\util_0)=g\cdot (x_0,u_0)$, where $u_0=f(x_0)$. 
We then define
$$
\prol g\cdot (x_0,\un_0) = (\xtil_0,\util_0^{(n)}),
$$
where
\begin{equation}\label{groupprol}
\tilde u^{(n)}_0 = \prol(g\cdot f)(\xtil_0).
\end{equation}
An important point to note here is that, due to the chain rule, the value of $\prol g\cdot (x_0,\un_0)$
depends exclusively on the derivatives of $f$ at the point $x_0$ up to order $n$, i.e., only on 
$(x_0,\un_0)$ itself. It is therefore independent of the choice of the auxiliary function $f$ from above
and therefore provides a well-defined operation on $M^{(n)}$.
\blem\label{prolloctrans} For any local transformation group $G$ on $M$ as above, 
$\mathrm{pr}^{(n)} G$ is a local transformation group on $M^{(n)}$.
\et
\pr We first note that $(g,(x_0,\un_0))\mapsto \prol g\cdot (x_0,\un_0)$ is smooth. Indeed, taking
the $f$ from \eqref{taylor} this follows directly from \eqref{gdotf}. To see this, it suffices to note that
the inverse function theorem yields an inverse that automatically depends smoothly on all parameters
that the original function depends upon (in our case, the $\un_0$): if $f=f(x,\alpha)$ ($\alpha$
representing any parameters), then $f^{-1}$ is the solution of the implicit equation $F(x,y,\alpha)
= f(y,\alpha) - x = 0$, so the implicit function theorem yields the claim.

It remains to verify (i) and (ii) from \ref{ltgdef2}. Here, (ii) is automatic because for $g=e$ we
get $g\cdot f = f$. Thus we are left with proving that
\begin{equation}\label{prolcomp}
\prol h (\prol g \cdot (x_0,\un_0)) = \prol(h\cdot g)(x_0,\un_0).
\end{equation}
To see this, let $(x_0,\un_0)\in M^{(n)}$ and pick a smooth function $f$ defined near $x_0$
with $\un_0=\prol f(x_0)$. Then 
$$
(\xtil_0,\util_0^{(n)}) = g \cdot (x_0,\un_0) = (\tilde x_0,\prol(g\cdot f)(\xtil_0)).
$$
To determine the left hand side of \eqref{prolcomp} we have to find any smooth function $F$
such that $\prol F(\xtil_0)=\prol(g\cdot f)(\xtil_0)$, calculate $h\cdot F$ and then
determine $\prol h\cdot F$ at the corresponding point. Naturally, we take $F:=g\cdot f$. Then
using the terminology from \eqref{xiphidef} we obtain
\begin{equation*}
\prol h (\prol g \cdot (x_0,\un_0)) = (\Xi_h (\xtil_0,(g\cdot f)(\xtil_0)),
\prol(h\cdot F)(\Xi_h (\xtil_0,(g\cdot f)(\xtil_0))))
\end{equation*}
Here,
\begin{equation*}
\begin{split}
\Xi_h (\xtil_0,(g\cdot f)(\xtil_0)) &= \Xi_h(g\cdot(x_0,u_0)) = \mathrm{proj_1}(h\cdot (g\cdot(x_0,u_0)))\\
&= \mathrm{proj_1}((h\cdot g)\cdot(x_0,u_0)) = \Xi_{h\cdot g}(x_0,u_0),
\end{split}
\end{equation*}
and, by \ref{grcomp}, $h\cdot F = h\cdot (g\cdot f) = (h\cdot g)\cdot f$,
so altogether we get
\begin{equation*}
\begin{split}
\prol h (\prol g \cdot (x_0,\un_0)) &= (\Xi_{h\cdot g}(x_0,u_0),
\prol((h\cdot g) f)(\Xi_{h\cdot g}(x_0,u_0)))\\
&= \prol(h\cdot g)(x_0,\un_0).
\end{split}
\end{equation*}
\ep
\bex\label{rotact3} Again let $G=SO(2)$ be the rotation group on $M=\R^2=X\times U$. We determine the first prolongation
of $G$. We have $p=q=1$, so $X\times U^{(1)}=\R^3$, with coordinates $(x,u,u_x)$. If $u=f(x)$ is a local
smooth function then $\prolo f(x) = (f(x),f'(x))$. Now let $(x^0,u^0,u_x^0)\in X\times U^{(1)}$ and let
$\eps$ be an angle of rotation. We wish to determine
$$
\prolo \eps \cdot (x^0,u^0,u_x^0) = (\tilde x^0,\tilde u^0,\tilde u_x^0).
$$
To this end, let $f$ be the Taylor polynomial
$$
f(x) = u^0 + u_x^0(x-x^0) = u^0_x\cdot x +(u^0-u^0_x x^0),
$$
which satisfies $f(x^0)=u^0$, $f'(x^0)=u^0_x$. Then by \ref{rotact2} we get
$$
\ftil(\xtil) = (\eps\cdot f)(\xtil) = \frac{\sin\eps + u^0_x\cos\eps}{\cos\eps - u^0_x\sin\eps}\xtil
+ \frac{u^0-u^0_x x^0}{\cos\eps - u^0_x\sin\eps},
$$
which is well-defined for $u^0_x\not=\cot\eps$. By \eqref{rotact} we have $\xtil^0 = x^0\cos\eps - u^0\sin\eps$,
so either again by \eqref{rotact} or by inserting we get
$$
\util^0 = \ftil(\xtil^0) = x^0\sin\eps + u^0\cos\eps.
$$
Finally, we get for the first derivative
$$
\util^0_x = \ftil'(\xtil^0) = \frac{\sin\eps + u^0_x\cos\eps}{\cos\eps - u^0_x\sin\eps}.
$$
Combining this and dropping the $0$-superscripts we obtain the action of $\prolo SO(2)$
on $X\times U^{(1)}$:
$$
\prolo\eps\cdot (x,u,u_x) = \left(x\cos\eps - u\sin\eps,x\sin\eps + u\cos\eps,
\frac{\sin\eps + u_x\cos\eps}{\cos\eps - u_x\sin\eps},
\right),
$$
defined for $|\eps|<|\mathrm{arccot}(u_x)|$. The important point to note is that although
the action of $SO(2)$ is linear and globally defined, its first prolongation is nonlinear
and only locally defined.
\et
In the previous example we see that the first prolongation of $G$ acts on the original
variables $(x,u)$ in exactly the same way as $G$ itself. This is actually a general
phenomenon. In fact, it follows directly from the definition of $\prol G$ via evaluating
derivatives of $\ftil$ that the action of $\prol g$ on $(x,u^{(k)})$ with $k\le n$
coincides with that of $\mathrm{pr}^{(k)}$. In particular, $\mathrm{pr}^{(0)}G=G$.
To give a concise formulation of this property we introduce the natural projections
\begin{equation*}
\begin{split}
\pi^n_k: M^{(n)} &\to M^{(k)}\\
\pi^n_k(x,\un) &= (x,u^{(k)}),
\end{split}
\end{equation*}
where (for $k\le n$), $u^{(k)}$ consists of the components $u^\alpha_J$, $\sharp J\le k$
of $\un$. For example, if $p=2$, $q=1$, then
\begin{equation*}
\begin{split}
\pi^2_0(x,y;u;u_x,u_y;u_{xx},u_{xy},u_{yy}) &= (x,y;u)\\
\pi^2_1(x,y;u;u_x,u_y;u_{xx},u_{xy},u_{yy}) &= (x,y;u;u_x,u_y).
\end{split}
\end{equation*}
Then for any $k\le n$ and any $g\in G$ we have
\begin{equation}\label{graprol}
\pi^n_k\circ \prol g = \mathrm{pr}^{(k)} g.
\end{equation}
In particular, this observation allows to calculate prolongations of group actions inductively.
\section{Invariance of differential equations}
In Section \ref{syssect} we have identified any given system $\mcS$ of differential
equations $P_\nu(x,\un)=0$ ($\nu=1,\dots,l$) with a subvariety $\mcS_P$ of $M^{(n)}$.
This identification will allow us to utilize the methods of symmetry group analysis
of algebraic equations in the context of differential equations. The essential 
translation result between the two realms is the following theorem
\bt\label{maintrans} Let $M$ be an open subset of $X\times U$ and let $\mcS$:
$$
P_\nu(x,\un)=0\quad  \nu=1,\dots,l
$$
be an $n$-th order system of differential equations defined on $M$, with corresponding
subvariety $\mcS_P$ of $M^{(n)}$. Suppose that $G$ is a local transformation group
acting on $M$ whose prolongation leaves $\mcS_P$ invariant, i.e.,
$$
(x,\un)\in \mcS_P \Rightarrow \prol g\cdot (x,\un)\in \mcS_P,
$$
whenever defined. Then $G$ is a symmetry group of the system $\mcS$ in the sense of
\ref{dgsgdef}.
\et
\pr Let $u=f(x)$, $f\colon \Omega\to U$ be a local solution of $P(x,\un)=0$. This means that the graph
$$
\Gamma\hn_f = \{(x,\prol f(x))\mid x\in \Omega\}
$$
of $\prol f$ is contained in $\mcS_P$. If $g\in G$ is such that $g\cdot f$ is well-defined then,
by the very definition of $\prol g$, the graph of $\prol (g\cdot f)$ is given by the image
of the graph of $\prol f$ under $\prol g$, i.e., 
$$
\Gamma\hn_{g\cdot f} = (\prol g)\cdot \Gamma\hn_f.
$$
Since $\mcS_P$ is invariant under $\prol g$, it follows from this that the graph of 
$\prol(g\cdot f)$ is contained in $\mcS_P$ as well. But this just means that $g\cdot f$
is a local solution of the system $\mcS$.
\ep
\section{Prolongation of vector fields}
To make the calculation of symmetry groups of differential equations accessible to the
infinitesimal methods we developed for algebraic equations we need to be able to
effectively calculate the infinitesimal generators of prolonged group actions.
\bd\label{vprolo} Let $M\sse X\times U$ be open and suppose that $v$ is the infinitesimal generator
of a local one-parameter transformation group $\eps\mapsto \Fl^v_\eps$. Then the $n$-th
prolongation of $v$, denoted by $\prol v$, is the local vector field on $M\hn$ which is the
infinitesimal generator of the corresponding prolonged one-parameter group $\prol \Fl^v_\eps$:
$$
\prol v|_{(x,\un)} = \dde \prol(\Fl^v_\eps)(x,\un),
$$
for any $(x,\un)\in M\hn$.
\et
Note that from this definition, together with \ref{flowaction2} it immediately follows that
\begin{equation}\label{flowstuff}
\prol(\Fl^v_\eps) = \Fl^{\prol v}_{\eps}.
\end{equation}
Any vector field $v$ on $M=X\times U$ can be written in the form
$$
v|_{(x,u)} = \sum_{i=1}^p \xi^i(x,u)\partial_{x^i} + \sum_{\alpha=1}^q \phi_\alpha(x,u)\partial_{u^\alpha}.
$$
Its $n$-th prolongation, being a vector field on $M\hn$, therefore has to be of the form
\begin{equation}\label{prolv}
\prol v|_{(x,\un)} = \sum_{i=1}^p \xi^i(x,u)\partial_{x^i} + 
\sum_{\alpha=1}^q \sum_J\phi_\alpha^J(x,\un)\partial_{u^\alpha_J},
\end{equation}
where the last sum is over all $J$ with $0\le \sharp J\le n$.
Note that due to \eqref{graprol} the prolonged group action agrees with 
the original one on $M^{(0)}=M$, so the coefficients $\xi^i$ and $\phi^0_\alpha=\phi_\alpha$
agree with those of $v$. Moreover, for the same reason, if $\sharp J = k$ then the 
coefficient $\phi^J_\alpha$ of $\partial_{u^\alpha_J}$ will only depend on $k$-th and lower
order derivatives of $u$, $\phi^J_\alpha = \phi^J_\alpha(x,u^{(k)})$, i.e.,
\begin{equation}\label{projproj}
T\pi^n_k(\prol v) = \mathrm{pr}^{(k)}v \quad n\ge k,
\end{equation}
where $\mathrm{pr}^{(0)}=v$. Again this reflects the possibility of calculating the coefficients
$\phi^J_\alpha$ recursively. Our main goal will be to derive a general formula for calculating
the $\phi^J_\alpha$ given $\phi_\alpha$.
\bex\label{rotact4} Returning once more to the rotation group $G=SO(2)$ on $\R^2$ from 
\ref{rotact2}, \ref{rotact3}, 
with infinitesimal generator $v = -u\partial_x + x\partial_u$ and group action
\begin{equation*}
\Fl^v_\eps(x,u) = \eps\cdot (x,u) = (x\cos\eps -u\sin\eps,x\sin\eps + u\cos\eps)
\end{equation*}
we have
$$
\prolo[\Fl^v_\eps] (x,u,u_x) = \left(x\cos\eps - u\sin\eps,x\sin\eps + u\cos\eps,
\frac{\sin\eps + u_x\cos\eps}{\cos\eps - u_x\sin\eps},
\right),
$$
By \ref{vprolo}, the first prolongation of $v$ is calculated by differentiating this
expression at $\eps=0$. This gives
\begin{equation}\label{rotvprol}
\prolo v = -u\partial_x + x\partial_u + (1+u_x^2)\partial_{u_x}
\end{equation}
As predicted by \eqref{projproj}, the first two terms coincide with those of $v$ itself.
\et
In order to be able to apply the criterion \ref{algsym} to the present situation we need
a suitable maximal rank condition for systems of differential equations.
\bd A system
$$
P_\nu(x,\un) = 0 \quad \nu=1,\dots,l
$$
of differential equations is said to be of maximal rank\index{maximal rank} if the $l\times(p+qp\hn)$
Jacobian matrix
$$
J_P(x,\un) = \left(\frac{\partial P_\nu}{\partial x^i},\frac{\partial P_\nu}{\partial u^\alpha_J}\right)
$$
of $P$ with respect to all its variables $(x,\un)$ is of rank $l$ whenever $P(x,\un)=0$.
\et 
\bex (i) The $2$-dimensional Laplace equation $P=u_{xx}+u_{yy}=0$ is of maximal rank since the 
Jacobian of $P$ with respect to $(x,y;u;u_x,u_y;u_{xx},u_{xy},u_{yy})\in X\times U^{(2)}$ is
$$
J_P=(0,0;0;0,0;1,0,1),
$$
which has rank $1$ everywhere.

(ii) The equation $P=(u_{xx}+u_{yy})^2=0$ is not of maximal rank, because
$$
J_P=(0,0;0;0,0;2(u_{xx}+u_{yy}),0,2(u_{xx}+u_{yy}))
$$
vanishes on $P=0$.
\et
\brem In practice, the maximal rank condition is not much of a restriction, i.e., `most' systems
automatically satisfy it. In fact, if $\mcS_P=\{(x,\un)\mid P(x,\un)=0\}$ is a regular ($l$-dimensional)
submanifold of $M\hn$ then locally around any point in $M\hn$ we can pick an adapted coordinate
system (cf.\ \cite[3.3.12]{KAoM}), and thereby (extracting the relevant coordinates) 
find a smooth map $\tilde P$ of rank $l$ such that
$\mcS_P=\mcS_{\tilde P} = \{(x,\un)\mid \tilde P(x,\un)=0\}$.
\et
Using the above condition we can now formulate the main result on the calculation of symmetry
groups of systems of differential equations:
\bt\label{mainpde} Let 
$$
P_\nu(x,\un) = 0 \quad \nu=1,\dots,l
$$
be a system of differential equations of maximal rank defined on an open subset $M$
of $X\times U$. If $G$ is a local transformation group acting on $M$ and
\begin{equation}\label{mainsymcrit}
\prol v [P_\nu(x,\un)] = 0, \quad \nu=1,\dots,l, \text{ whenever } P(x,\un)=0
\end{equation}
for every infinitesimal generator $v$ of $G$, then $G$ is a symmetry group of
the system.
\et
\pr By \ref{maintrans} it suffices to show that $\mcS_P$ remains invariant under 
$\prol G$. By \ref{prolloctrans}, $\prol G$ is a local group of
transformations on $M\hn$, whose infinitesimal generators are exactly the $\prol v$,
for $v$ the infinitesimal generators of $G$ (by \ref{vprolo}).
Since $P$ is of maximal rank, the result now follows from \ref{algsym}.
\ep
We shall see in \ref{nesu} that if $P$ satisfies a certain local solvability condition
then \eqref{mainsymcrit} is in fact necessary and sufficient for $G$ to be a symmetry group
of the system.

The Hadamard Lemma \ref{hadamard} allows the following equivalent reformulation of \eqref{mainsymcrit}:
there exist smooth functions $Q_{\nu\mu}=Q_{\nu\mu}(x,\un)$ such that
\begin{equation}\label{mainsymhadamard}
\prol v [P_\nu(x,\un)] = \sum_{\mu=1}^l Q_{\nu\mu}(x,\un) P_\mu(x,\un)
\end{equation}
holds identically on $M\hn$.
\bex\label{rotact5} Continuing our analysis of the rotation group $SO(2)$ on $\R^2$ from \ref{rotact4},
consider the first order ODE
\begin{equation}\label{rotode}
P(x,u,u_x) = (u-x)u_x + u +x =0.
\end{equation}
The Jacobian of $P$ is 
$$
J_P=(\partial_x P,\partial_u P,\partial_{u_x} P)=(1-u_x,1+u_x,u-x),
$$
hence has rank $1$ everywhere.
\et
We now apply the infinitesimal generator $\prolo v$ from \eqref{rotvprol} to $P$:
\begin{equation*}
\begin{split}
\prolo v(P) &= -u\partial_xP + x\partial_u P + (1+u_x^2)\partial_{u_x}P\\
&= -u(1-u_x) + x(1+u_x) + (1+u_x)^2(u-x)\\ 
&= u_x[(u-x)u_x + u +x] = u_x P.
\end{split}
\end{equation*}
Therefore, \eqref{mainsymcrit} (or also \eqref{mainsymhadamard}) is satisfied, and we conclude
that the rotation group transforms solutions of \eqref{rotode} to other solutions.

The calculation of even the first prolongation of the infinitesimal generator of the rotation
group in \ref{rotact3} demonstrates that the direct method of first determining the prolonged
group action and then calculating the corresponding generators is not feasible in practice.
We need to find an algorithmic way of directly calculating prolongations of vector fields,
avoiding the detour of determining the prolonged group action along the way.

Before deriving a general formula, let us start out by considering some special cases first.
Let 
$$
v = \sum_{i=1}^p \xi^i(x)\partial_{x^i}
$$
on $M\sse X\times U$, with $U=\R$. According to \eqref{xiphidef}, 
the corresponding group action $g_\eps = \Fl^v_\eps$ then is of the
form
$$
(\xtil,\util) = g_\eps\cdot (x,u) = (\Xi_\eps(x),u),
$$
where 
\begin{equation}\label{xiabl}
\dde \Xi^i_\eps(x) = \xi^i(x).
\end{equation}
On $M^{(1)}$ we have coordinates $(x,u^{(1)}) =(x^i,u,u_j)$ with $u_j:=\partial_{x^j}u$.
By definition, to calculate $\prolo v$ in $(x,u\heins)$, we may take any function $f$ with 
$u=f(x)$ and $u_j=\partial_{x^j}f(x)$ for $j=1,\dots,p$. Then 
$$
\prolo g_\eps\cdot (x,u\heins) = (\xtil,\util\heins),
$$
where $\xtil = \Xi_\eps(x)$, $\util = u$, and $\util_j$ are the derivatives of the transformed
function $\ftil_\eps = g_\eps\cdot f$, which, by \eqref{gdotf}, is given by
$$
\util = \ftil_\eps(\xtil) = f(\Xi_\eps^{-1}(\xtil)) = f(\Xi_{-\eps}(\xtil)).
$$
Therefore,
\begin{equation}\label{utildiff}
\util_j = \pd{\ftil_\eps}{\xtil^j}(\xtil) = \sum_{k=1}^p\pd{f}{x^k}(\Xi_{-\eps}(\xtil))\cdot
\pd{\Xi^k_{-\eps}}{\xtil^j}(\xtil).
\end{equation}
Since $\Xi_{-\eps}(\xtil)=x$, this simplifies to 
$$
\util_j =  \sum_{k=1}^p
\pd{\Xi^k_{-\eps}}{\xtil^j}(\Xi_\eps(x)) u_k.
$$
This is the explicit formula for the action of $g_\eps$ on the $u_j$-variables. To find
the infinitesimal generator of $\prolo g_\eps$, we have to differentiate this expression
with respect to $\eps$ at $\eps=0$. Also, since $\prolo g_\eps$ acts on $(x,u)$ as $g_\eps$,
the coefficients of $\partial_{x^i}$ and $\partial_{u}$ in $\prolo v$ stay the same. Hence
\begin{equation}\label{pr1v}
\prolo v = \sum_{i=1}^p\xi^i(x)\partial_{x^i} + \sum_{j=1}^p \phi^j(x,u\heins)\partial_{u_j},
\end{equation}
where
$$
\phi^j(x,u\heins) = \dde \sum_{k=1}^p \pd{\Xi^k_{-\eps}}{\xtil^j}(\Xi_\eps(x)) u_k.
$$
Here, using \eqref{xiabl}, we have
\begin{equation*}
\dde\left[\pd{\Xi^k_{-\eps}}{\xtil^j}(\Xi_\eps(x))\right] = -\pd{\xi^k}{x^j}(x) + 
\sum_{l=1}^p \left[\frac{\partial^2 \Xi^k_{-\eps}}{\partial \xtil^j \partial \xtil^l}(\Xi_\eps(x))
\frac{d}{d\eps}\Xi^l_\eps(x)\right]_{\eps=0},
\end{equation*}
and the second term vanishes since $\Xi_0=\id$. Therefore,
\begin{equation}\label{phij}
\phi^j(x,u,u_x) = - \sum_{k=1}^p \pd{\xi^k}{x^j} u_k.
\end{equation}
The second special case we look at is a local group action as above, only this time acting
exclusively on the dependent variable:
$$
(\xtil, \util) =g_\eps\cdot(x,u) = (x,\Phi_\eps(x,u)).
$$
In this case, $v=\phi(x,u)\partial_u$, where
$$
\phi(x,u) = \dde \Phi_\eps(x,u).
$$
Let $f$ be a local smooth function with $f(x)=u$. Then $\ftil_\eps=g_\eps\cdot f$ is given by
\begin{equation}\label{fepss}
\util = \ftil_\eps(x)=\Phi_\eps(x,f(x)).
\end{equation}
To determine the prolonged group action, we need to differentiate $\ftil_\eps$:
$$
\util_j=\pd{\ftil_\eps}{x^j}(x) = \pd{\Phi_\eps}{x^j}(x,f(x)) + \pd{\Phi_\eps}{u}(x,f(x))\pd{f}{x^j}(x).
$$
Consequently, $\prolo g_\eps\cdot(x,u\heins) =(x,\util\heins)$, where
\begin{equation}\label{utilj}
\util_j = \pd{\Phi_\eps}{x^j} + \pd{\Phi_\eps}{u}u_j
\end{equation}
As before, to determine $\phi^j$ in
$$
\prolo v = v + \sum_{j=1}^p \phi^j(x,u\heins)\partial_{u_j},
$$
we have to differentiate \eqref{utilj} at $\eps=0$:
\begin{equation}\label{pr1spec}
\phi^j(x,u\heins) = \dde \util_j = \pd{\phi}{x^j} + u_j\pd{\phi}{u}.
\end{equation}
We have here the first instance of a generally useful operation, namely that of 
total derivative\index{total derivative}:
$$
\phi^j(x,u\heins) = \pd{}{x^j}[\phi(x,f(x))].
$$
Thus $\phi^j(x,u\heins)$ is calculated from $\phi(x,u)$ by differentiating with respect
to $x^j$ while treating $u$ as a function of $x$. Accordingly, the total derivative
of $\phi$ with respect to $x^j$, $D_j$ is defined as 
$$
D_j\phi := \pd{\phi}{x^j} + u_j\pd{\phi}{u}.
$$
More generally, we define:
\bd Let $(x,\un)\mapsto P(x,\un)$ be a smooth function of $x$, $u$, and derivatives
of $u$ up to order $n$, defined on an open subset $M\hn$ of $X\times U\hn$. The total
derivative\index{total derivative} of $P$ with respect to $x^i$ is the unique smooth
function $D_iP(x,\unp)$ defined on $M\hnp$ with the property that if $u=f(x)$ is any
smooth function of $x$ then
$$
D_iP(x,\prolp f(x)) = \pd{}{x^i}[P(x,\prol f(x))].
$$
\et
Thus $D_iP$ is calculated by differentiating $P$ with respect to $x^i$ while treating
all $u^\alpha$ and their derivatives as functions of $x$.
\blem For any $P(x,\un)$ we have
\begin{equation}\label{totder1}
D_iP = \pd{P}{x^i} + \sum_{\alpha=1}^q \sum_J u^\alpha_{J,i}\pd{P}{u^\alpha_J},
\end{equation}
where, for $J=(j_1,\dots,j_k)$,
\begin{equation}\label{totder2}
u^\alpha_{J,i} = \pd{u^\alpha_J}{x^i} = \pd{^{k+1} u^\alpha}{x^i \partial x^{j_1}\dots \partial x^{j_k}},
\end{equation}
and the sum in \eqref{totder1} extends over all $0\le \sharp J \le n$.
\et
For example, with $X=\R^2$, $U=\R$ we have
\begin{equation*}
\begin{split}
D_xP &= \pd{P}{x} + u_x \pd{P}{u} + u_{xx} \pd{P}{u_x} + u_{xy} \pd{P}{u_y} + u_{xxx} \pd{P}{u_{xx}} + \dots\\
D_yP &= \pd{P}{y} + u_y \pd{P}{u} + u_{xy} \pd{P}{u_x} + u_{yy} \pd{P}{u_y} + u_{xxy} \pd{P}{u_{xx}} + \dots
\end{split}
\end{equation*}
Higher order total derivatives are defined inductively: let $J=(j_1,\dots,j_k)$ be a $k$-th order
multiindex with $1\le j_\kappa\le p$ for each $\kappa$, then we set
$$
D_J := D_{j_1}D_{j_2}\dots D_{j_k}.
$$
Also, note that the order of differentiation does not matter for total derivatives: $D_jD_l = D_lD_j$
for all $j$, $l$.

In the proof of the general prolongation formula we will need the following auxiliary result:
\blem\label{detdiff} Let $\eps\mapsto M(\eps)$ be a smooth map that takes values in the space of invertible 
$n\times n$-matrices. Then
$$
\frac{d}{d\eps}[M(\eps)^{-1}] = -M(\eps)^{-1}\frac{dM(\eps)}{d\eps}M(\eps)^{-1}.
$$
\et
\pr It suffices to differentiate the identity $M(\eps)M(\eps)^{-1} = I$ by the product rule.
\ep
\brem\label{jetsub} As a final prerequisite for the proof of \ref{gpf} below we need the 
observation that given any open subset $M$ of $X\times U$, one can view
the $(n+1)$-st jet space $M\hnp$ as a subspace of the first jet space $(M\hn)\heins$ of the
$n$-th jet space. In fact, any $(n+1)$-st derivative $u^\alpha_J$ can be viewed as a 
first order derivative of an $n$-th order derivative (which, in general, can be done
in many different ways). For example, let $p=2$, $q=1$, so that the coordinates on $M\heins$
can be written as $(x,y;u;u_x,u_y)$. Then we may view $(u_x,u_y)$ as new dependent variables,
$u_x=v$, $u_y=w$. This turns $M\heins$ into a subset of $X\times \tilde U$, where $X$ is
still $\R^2$, but $\tilde U$ has three dependent variables $u,\,v,\,w$. Consequently,
the first jet space $(M\heins)\heins$ of $M\heins$ turns into an open subset of
$X\times \tilde U\heins$, with coordinates $(x,y;u;v;w;u_x,u_y,v_x,v_y,w_x,w_y)$.
Since we defined $v=u_x$, $w=u_y$ it follows that $M\hzwei\sse (M\heins)\heins$ is
the subspace defined by the relations
$$
v=u_x,\quad w=u_y,\quad v_y = w_x
$$
in $X\times \tilde U\heins$. Here, the third relation is a consequence of $u_{xy}=u_{yx}$.
\et
\bt\label{gpf} (The general prolongation formula)\index{prolongation formula} Let
\begin{equation}\label{vfform}
v = \sum_{i=1}^p \xi^i(x,u)\pd{}{x^i} + \sum_{\alpha=1}^q \phi_\alpha(x,u)\pd{}{u^\alpha}
\end{equation}
be a vector field on an open subset $M\sse X\times U$. The $n$-th prolongation of $v$ is the
vector field
\begin{equation}\label{theprol}
\prol v = v + \sum_{\alpha=1}^q\sum_J \phi^J_\alpha(x,\un)\pd{}{u^\alpha_J}
\end{equation}
defined on $M\hn \sse X\times U\hn$, where $J$ runs over all multi-indices $J=(j_1,\dots,j_k)$
with $1\le j_\kappa\le p$, $1\le k\le n$. The coefficient functions $\phi_\alpha^J$ are given by
\begin{equation}\label{phiform}
\phi_\alpha^J(x,\un) = D_J\left(\phi_\alpha - \sum_{i=1}^p \xi^iu^\alpha_i\right)
+ \sum_{i=1}^p \xi^iu^\alpha_{J,i},
\end{equation}
where $u^\alpha_i = \partial_{x^i}u^\alpha$ and $u^\alpha_{J,i} = \partial_{x^i}u^\alpha_J$.
\et
\pr We start out with the case $n=1$. Let $g_\eps = \Fl^v_\eps$ and 
$$
(\xtil,\util) = g_\eps\cdot (x,u) = (\Xi_\eps(x,u),\Phi_\eps(x,u)).
$$
Then
\begin{equation}\label{xiphi}
\begin{split}
\xi^i(x,u) &= \dde \Xi^i_\eps(x,u), \quad i=1,\dots,p,\\
\phi_\alpha(x,u) &= \dde \Phi^\alpha_\eps(x,u), \quad \alpha=1,\dots,q.
\end{split}
\end{equation}
Given any $(x,u\heins)\in M\heins$, pick a local smooth map $f$ with $\prolo f(x)=u\heins$,
i.e., with $u^\alpha = f^\alpha(x)$, $u^\alpha_i =\partial_{x^i} f^\alpha(x)$.
By \eqref{gdotf}, for $\eps$ small and restricting the domain of $f$ sufficiently, we have
$$
\util = \ftil_\eps(\xtil) = (g_\eps\cdot f)(\xtil) = 
[\Phi_\eps\circ (\id\times f)]\circ [\Xi_\eps\circ (\id\times f)]^{-1}(\xtil).
$$
By the chain rule, we conclude from this that for the Jacobian $J\ftil_\eps = 
(\partial \ftil_\eps^\alpha/\partial\xtil^i)$
we have (using that $x=[\Xi_\eps\circ (\id\times f)]^{-1}(\xtil)$)
\begin{equation}\label{jac}
J\ftil_\eps(\xtil) = J[\Phi_\eps\circ (\id\times f)](x)\cdot[J[\Xi_\eps\circ (\id\times f)](x)]^{-1}
\end{equation}
The matrix entries of $J\ftil_\eps(\xtil)$ then give formulae for $\prolo g_\eps$.

To calculate $\prolo v$ we have to differentiate \eqref{jac} at $\eps=0$. In doing this,
note that since $g_0=\id$, we have
\begin{equation}\label{xiphi0}
\Xi_0(x,f(x)) = x, \quad \Phi_0(x,f(x)) = f(x),
\end{equation}
so with $I$ the $p\times p$ unit matrix we get
$$
J[\Xi_0\circ (\id\times f)](x) = I,\quad J[\Phi_0\circ (\id\times f)](x) = Jf(x).
$$
Therefore, using \ref{detdiff},
\begin{equation*}
\begin{split}
\dde J\ftil_\eps(\xtil) &= \dde J[\Phi_\eps\circ (\id\times f)](x)
- Jf(x)\dde J [\Xi_\eps\circ (\id\times f)](x)\\
&= J[\phi\circ (\id\times f)](x) - Jf(x)\cdot J[\xi\circ (\id\times f)](x).
\end{split}
\end{equation*}
The matrix entries of this expression give the coefficient functions $\phi^k_\alpha$
of $\partial_{u^\alpha_k}$ in $\prolo v$. The $(\alpha,k)$-th entry is
$$
\phi^k_\alpha(x,\prolo f(x)) = \pd{}{x^k}[\phi_\alpha(x,f(x))]
-\sum_{i=1}^p \pd{f^\alpha}{x^i}\cdot \pd{}{x^k}[\xi^i(x,f(x))],
$$
or, in terms of total derivatives,
\begin{equation}\label{prol1form}
\begin{split}
\phi^k_\alpha(x,u\heins) &=  D_k\phi_\alpha(x,u^{(1)}) -\sum_{i=1}^p D_k\xi^i(x,u^{(1)})u^\alpha_i\\
&= D_k\left[\phi_\alpha -\sum_{i=1}^p \xi^i u^\alpha_i\right] + \sum_{i=1}^p \xi^i u^\alpha_{ki},
\end{split}
\end{equation}
with $u^\alpha_{ki}=\partial^2 u^\alpha/\partial x^k\partial x^i$. This is \eqref{phiform}
for $n=1$.

In the general case we proceed by induction. Here we make use of the key observation \ref{jetsub},
which allows us to view
the $(n+1)$-st jet space $M\hnp$ as a subspace of the first jet space $(M\hn)\heins$ of the
$n$-th jet space. The strategy is to regard $\prolm v$ as a vector field $M\hnm$ and then
use the case $n=1$ above to prolong it to $(M\hnm)\heins$. Then we restrict the resulting
vector field to the subspace $M\hn$ of $(M\hnm)\heins$, and this will give the $n$-th
prolongation of $v$. The fact that this is possible and gives a well-defined result will
follow from the explicit formula below.

The coordinates on $(M\hnm)\heins$ are given by $u^\alpha_{J,k}=\pd{u^\alpha_J}{x^k}$,
where $J=(j_1,\dots,j_{n-1})$, $1\le k\le p$, and $1\le \alpha\le q$. By \eqref{projproj},
the only new coefficients we need to calculate are those of the highest order derivatives
$\partial/\partial u^\alpha_{J,k}$. By \eqref{prol1form}, they are given by
\begin{equation}\label{recursion}
\phi^{J,k}_\alpha = D_k\phi^J_\alpha - \sum_{i=1}^p D_k\xi^i \cdot u^\alpha_{J,i}.
\end{equation}
To finish the proof, it suffices to show that the explicit formula \eqref{phiform} solves
the recursion relation \eqref{recursion} in closed form. This follows by induction:
supposing that $\phi^J_\alpha$ has already been shown to satisfy \eqref{phiform}, 
it follows from \eqref{recursion} that
\begin{equation*}
\begin{split}
\phi^{J,k}_\alpha &= D_k\left[D_J\left(\phi_\alpha-\sum_{i=1}^p\xi^i u^\alpha_i\right) + 
\sum_{i=1}^p\xi^i u^\alpha_{J,i}\right] - \sum_{i=1}^p D_k\xi^i\cdot u^\alpha_{J,i}\\
&= D_kD_J\left(\phi_\alpha-\sum_{i=1}^p\xi^i u^\alpha_i\right) + 
\sum_{i=1}^p(D_k\xi^i\cdot u^\alpha_{J,i} + \xi^i u^\alpha_{J,ik}) - 
\sum_{i=1}^p D_k\xi^i\cdot u^\alpha_{J,i}\\
&=D_kD_J\left(\phi_\alpha-\sum_{i=1}^p\xi^i u^\alpha_i\right) + \sum_{i=1}^p \xi^i u^\alpha_{J,ik},
\end{split}
\end{equation*}
where $u^\alpha_{J,ik}=\partial^2 u^\alpha_J/\partial x^i\partial x^k$. Thus 
$\phi^{J,k}_\alpha$ is of the form \eqref{phiform} as well, which completes the
induction step.
\ep
\brem Note that \eqref{recursion} provides a useful way of calculating prolongations
of vector fields recursively.
\et
\bex As an illustration of the relative ease with which one may calculate prolongations
of vector fields using the general prolongation formula, we once more return to the
rotation group $G=SO(2)$ acting on $\R^2$, cf.\ \ref{rotact5}. The infinitesimal generator
is $v=-u\partial_x + x\partial_u$, so $\phi=x$ and $\xi=-u$. By \ref{gpf} we have
$$
\prolo v = v + \phi^x\partial_{u_x},
$$
where 
$$
\phi^x = D_x(\phi-\xi u_x) + \xi u_{xx} = D_x(x+uu_x)-uu_{xx} = 1 + u_x^2,
$$
which confirms \eqref{rotvprol}. To calculate the coefficient $\phi^{xx}$
of $\partial_{u_{xx}}$ in $\prolt v$, we can either use \eqref{phiform}:
$$
\phi^{xx} = D_x^2(\phi-\xi u_x) + \xi u_{xxx} = D_x^2(x+uu_x)-uu_{xxx} = 3u_xu_{xx},
$$
or the recursion formula \eqref{recursion}:
$$
\phi^{xx} = D_x\phi^x-u_{xx}D_x\xi = D_x(1+u_x^2)+u_xu_{xx} = 3u_xu_{xx}.
$$
Consequently,
$$
\prolt v = -u\partial_x + x\partial_u + (1 + u_x^2)\partial_{u_x} + 3u_xu_{xx}\partial_{u_{xx}}.
$$
Note that to derive this result via the second prolongation of the group action (as in \ref{rotact3}) 
would be much more involved.

Based on this, and using the invariance criterion in \ref{mainpde}, it follows that the ODE 
$u_{xx}=0$ has $SO(2)$ as a symmetry group. In fact,
$$
\prolt v(u_{xx}) = 3u_xu_{xx} = 0 \quad\text{whenever }\ u_{xx}=0.
$$
Since the solutions of $u_{xx}=0$ are linear functions, this is basically just the statement
that rotations take straight lines to straight lines.

Next, consider the function
\begin{equation}\label{curv}
\kappa(x,u\hzwei) = \frac{u_{xx}}{(1+u_x^2)^{3/2}},
\end{equation}
describing the Frenet-curvature of a planar curve. Then one easily checks that $\prolt v(\kappa)=0$,
so \ref{finvcrit} shows that $\kappa$ is an invariant of $\prolt SO(2)$. Geometrically, this means
that the curvature of a planar curve is invariant under rotations.
\et
\bt Let $v$, $w$ be smooth vector fields on the open subset $M$ of $X\times U$. Then
\begin{itemize}
\item[(i)] $\prol(cv+w) = c \cdot\prol v + \prol w \quad (c\in \R)$.
\item[(ii)] $\prol[v,w] = [\prol v,\prol w]$.
\end{itemize}
\et
\pr (i) This is immediate from \ref{gpf}.

(ii) This also could be verified using \ref{gpf}, but we take a different route.
First, by \ref{prolloctrans} we have 
$$
\prol(g\cdot h) = (\prol g)\cdot (\prol h).
$$
Also, defining $(g+h)\cdot x := g\cdot x + h\cdot x$ we obtain 
$$
\prol(g + h) = \prol g + \prol h
$$
directly from the definition \eqref{groupprol} (using sums of representing functions
to represent sums of points in the jet bundle). Also, $\prol \id_M = \id_{M\hn}$, and 
$v\mapsto \prol v$ is clearly continuous. Finally, recall from \eqref{flowstuff} that
$$
\prol(\Fl^v_\eps) = \Fl^{\prol v}_{\eps},
$$
so we can use \ref{eduard} to calculate as follows:
\begin{equation*}
\begin{split}
[\prol v,\prol w] &= \lim_{\eps\to 0} \frac{
\Fl^{\prol w}_{-\sqrt{\eps}}\circ\Fl^{\prol v}_{-\sqrt{\eps}}\circ\Fl^{\prol w}_{\sqrt{\eps}}
\circ\Fl^{\prol v}_{\sqrt{\eps}} - \id_{M\hn}}{\eps}\\
&=\lim_{\eps\to 0} \frac{\prol[
\Fl^{w}_{-\sqrt{\eps}}\circ\Fl^{v}_{-\sqrt{\eps}}\circ\Fl^{w}_{\sqrt{\eps}}
\circ\Fl^{v}_{\sqrt{\eps}} - \id_{M}]}{\eps}\\
&=\prol\left[\lim_{\eps\to 0} \frac{
\Fl^{w}_{-\sqrt{\eps}}\circ\Fl^{v}_{-\sqrt{\eps}}\circ\Fl^{w}_{\sqrt{\eps}}
\circ\Fl^{v}_{\sqrt{\eps}} - \id_{M}}{\eps}\right]\\
&=\prol [v,w].
\end{split}
\end{equation*}
\ep
If we call a vector field $v$ an {\em infinitesimal symmetry}\index{infinitesimal symmetry} 
of a system $\mcS$ of differential
equations if it satisfies \eqref{mainsymcrit}, then we obtain from the previous theorem:
\bc\label{infliesym} Let $\mcS$ be a system of differential equations of maximal rank defined on $M\sse X\times U$.
Then the set of all infinitesimal symmetries of $\mcS$ forms a Lie algebra of vector fields on $M$.
\et
To conclude this section we note the following equivalent way of writing the general prolongation
formula \ref{gpf}. For $v$ a vector field on $M\sse X\times U$, let
\begin{equation}\label{chardef}
Q_\alpha(x,u\heins) := \phi_\alpha(x,u) - \sum_{i=1}^p \xi^i(x,u)u^\alpha_i, \quad \alpha=1,\dots,q.
\end{equation}
The $q$-tuple $Q(x,u\heins)=(Q_1,\dots,Q_q)$ is called the 
{\em characteristic}\index{characteristic!of vector field}
of the vector field $v$. Using this, \eqref{phiform} becomes
\begin{equation}\label{phiform2}
\phi^J_\alpha = D_JQ_\alpha + \sum_{i=1}^p \xi^i u^\alpha_{J,i}
\end{equation}
Substituting this into \eqref{theprol}, a brief calculation gives
\begin{equation}\label{doch}
\prol v = \sum_{\alpha=1}^q\sum_J D_JQ_\alpha \pd{}{u^\alpha_J} + 
\sum_{i=1}^p \xi^i\left(\pd{}{x^i} + \sum_{\alpha=1}^q\sum_J u^\alpha_{J,i} \pd{}{u^\alpha_J}\right),
\end{equation}
which, by \eqref{totder1}, means that
\begin{equation}\label{charprol}
\prol v = \prol v_Q + \sum_{i=1}^p \xi^iD_i,
\end{equation}
where 
\begin{equation}\label{charparts}
v_Q:=\sum_{\alpha=1}^qQ_\alpha(x,u\heins)\pd{}{u^\alpha},\quad \prol v_Q = 
\sum_{\alpha=1}^q\sum_J D_JQ_\alpha \pd{}{u^\alpha_J}.
\end{equation}
Here, the last equality comes from applying \eqref{doch} to $v_Q$. In all of the above
equations, the sums over $J$ extend over all $J$ with $0\le \sharp J\le n$. 
Note that the individual terms on the right hand side of \eqref{charprol} actually involve
$(n+1)$-st order derivatives of $u$, and only their combination gives a genuine vector
field on $M\hn$.
\section{Calculation of symmetry groups}\index{symmetry group!calculation of}
In this section we want to illustrate how to explicitly calculate the symmetry group of any given 
system $\mcS$ of differential equations. From the previous sections, we obtain the following
general strategy for tackling this problem:
\begin{itemize}
\item[(i)] Make an ansatz for an infinitesimal generator $v$ of a
prospective one-parameter symmetry group of $\mcS$ in the form \eqref{vfform}
with unknown functions $\xi_i$ and $\phi_\alpha$.
\item[(ii)] Calculate $\prol v$ according to (\ref{theprol}), \eqref{phiform}, where $n$ is the order of $\mcS$.
\item[(iii)] Insert the resulting expression into (\ref{mainsymcrit}). Since these
equations only have to hold on $\mcS_P$, eliminate any dependencies of
the derivatives of the $u$'s by using the equations in $\mcS$. The 
equations achieved in this way have to hold identically in $x$, $u$ and 
the remaining partial derivatives of the $u$'s.
\item[(iv)] Solving the resulting so-called {\em defining equations}\index{defining equations} 
produces a number of (usually elementary)
partial differential equations for the unknown functions $\xi_i$ and
$\phi_\alpha$.
\item[(v)] Compute the $\xi_i$ and $\phi_\alpha$ from these equations. 
By \ref{infliesym}, the vector fields gathered from this procedure
form a Lie algebra of infinitesimal symmetries.
\item[(vi)] The corresponding one-parameter symmetry groups are found by
calculating the flows of the respective infinitesimal symmetries.
\end{itemize}
Let us systematically work through this list in a concrete example:
\bex\label{heatsym} Consider the one-dimensional heat equation
\begin{equation}\label{heatequ}
u_t=u_{xx}
\end{equation}
\et
Thus we have $X=\R^2$ with variables $x$, $t$, and $U=\R$, and the equation 
$P(x,t,u\hzwei)=u_t-u_{xx}$ determines a subvariety in $X\times U\hzwei$.
For the infinitesimal generators of symmetries of $\mcS=\{P=0\}$ we make the
ansatz
\begin{equation}\label{ansatz}
v = \xi(x,t,u) \pd{}{x} + \tau(x,t,u)\pd{}{t} + \phi(x,t,u)\pd{}{u}
\end{equation}
Our aim is to determine all possible coefficient functions $\xi,\,\tau,\,\phi$
such that $\Fl^v_\eps$ is a symmetry group of the heat equation \eqref{heatequ}.
Since we wish to apply \ref{mainpde}, we will need the second prolongation of $v$,
$$
\prolt v = v + \phi^x\pd{}{u_{x}} + \phi^t\pd{}{u_{t}} + \phi^{xx}\pd{}{u_{xx}}
+ \phi^{xt}\pd{}{u_{xt}} + \phi^{tt}\pd{}{u_{tt}}.
$$
Inserting this into \eqref{mainsymcrit}, it follows that the only restriction on 
the coefficients of $v$ is that
\begin{equation}\label{onlyres}
\phi^t = \phi^{xx} \quad \text{ whenever }\quad u_t=u_{xx}.
\end{equation}
Using \eqref{phiform}, we calculate
\begin{equation}\label{phit}
\begin{split}
\phi^t &=  D_t(\phi-\xi u_x -\tau u_t) +\xi u_{xt} + \tau u_{tt}
= D_t\phi -u_xD_t\xi-u_tD_t\tau\\
&= \phi_t - \xi_t u_x + (\phi_u-\tau_t)u_t - \xi_uu_xu_t -\tau_u u_t^2
\end{split}
\end{equation}
and
\begin{equation}\label{phixx}
\begin{split}
\phi^{xx} &=  D_x^2(\phi-\xi u_x -\tau u_t) +\xi u_{xxx} + \tau u_{xxt}\\
&= D_x^2\phi - u_xD_x^2\xi- u_tD_x^2\tau -2 u_{xx}D_x\xi - 2u_{xt}D_x\tau\\
&= \phi_{xx} + (2\phi_{xu}-\xi_{xx})u_x -\tau_{xx}u_t + (\phi_{uu} -2\xi_{xu})u_x^2 -2\tau_{xu}u_xu_t\\
&\hphantom{=} -\xi_{uu} u_x^3 -\tau_{uu}u_x^2u_t + (\phi_u-2\xi_x)u_{xx} - 2\tau_xu_{xt} - 3\xi_uu_xu_{xx}\\
&\hphantom{=}-\tau_uu_tu_{xx}-2\tau_uu_xu_{xt}.
\end{split}
\end{equation}
If we insert these expressions into \eqref{onlyres} and replace $u_t$ by $u_{xx}$ wherever it occurs,
we obtain an equation on the jet space $X\times U\hzwei$ that has to hold identically in all
the remaining variables. We may therefore `split' the equation by equating the coefficient
functions of all monomials appearing in it. This gives the following table:
\vskip2em
\begin{center}
\begin{tabular}{ccccc}
\hline
Monomial & & Coefficient & & \\
\hline 
$u_xu_{xt}$ & & $0=-2\tau_u$ & & (a) \\
$u_{xt}$ & & $0=-2\tau_x$ & & (b) \\
$u_{xx}^2$ & & $-\tau_u=-\tau_u$ & & (c) \\
$u_x^2u_{xx}$ & & $0=-\tau_{uu}$ & & (d) \\
$u_xu_{xx}$ & & $-\xi_u=-2\tau_{xu}-3\xi_u$ & & (e) \\
$u_{xx}$ & & $\phi_u-\tau_t=-\tau_{xx} +\phi_u-2\xi_x$ & & (f) \\
$u_{x}^3$ & & $0=-\xi_{uu}$ & & (g) \\
$u_{x}^2$ & & $0=\phi_{uu} - 2\xi_{xu}$ & & (h) \\
$u_{x}$ & & $-\xi_t=2\phi_{xu} -\xi_{xx}$ & & (j) \\
$1$ & & $\phi_t=\phi_{xx}$ & & (k) \\
\hline
\end{tabular}
\end{center}
\vskip1em
We first note that (a), (b) imply that $\tau$ is a function of $t$ only: $\tau=\tau(t)$. Using this,
(e) shows that $\xi$ does not depend on $u$, and (f) gives $\tau_t=2\xi_x$. Therefore,
$$
\xi(x,t) = \frac{1}{2}\tau_tx + \sigma(t),
$$
where $\sigma$ is some function of $t$. (h) shows that $\phi$ is linear in $u$, so
$$
\phi(x,t,u) = \beta(x,t)u + \alpha(x,t)
$$
for some functions $\alpha$, $\beta$. By (j), $-2\beta_x = \xi_t = \frac{1}{2}\tau_{tt}x+\sigma_t$, so 
\begin{equation}\label{betae}
\beta(x,t) = -\frac{1}{8}\tau_{tt}x^2-\frac{1}{2}\sigma_tx + \rho(t).
\end{equation}
Finally, (k) implies that both $\alpha$ and $\beta$ must be solutions of the heat equation:
$$
\alpha_t = \alpha_{xx}, \qquad \beta_t=\beta_{xx}.
$$
Combined with \eqref{betae} we get
$$
\tau_{ttt} = 0, \qquad \sigma_{tt}=0, \qquad \rho_t = -\frac{1}{4}\tau_{tt}.
$$
Hence $\tau$ is quadratic in $t$, $\sigma$ is linear in $t$, so for suitable constants
$c_i$ we get:
\begin{equation*}
\begin{split}
\tau &= c_2 +2c_4 t + 4c_6 t^2\\
\sigma &= c_1 + 2c_5 t\\
\rho_t &= -\frac{1}{4}\tau_{tt} = -2c_6 \Rightarrow \rho = -2c_6 t + c_3\\
\Rightarrow\beta &=  -\frac{1}{8}\tau_{tt}x^2 -\frac{1}{2}\sigma_{t}x +\rho = -c_6x^2 -c_5x-2c_6t+c_3
\end{split}
\end{equation*}
Consequently,
\begin{equation*}
\begin{split}
\xi &= \frac{1}{2}\tau_t + \sigma = c_1 +c_4 x + 2c_5 t + 4c_6xt\\
\tau &= c_2 + 2c_4 t + 4c_6t^2\\
\phi &= \beta u+\alpha = (c_3 -c_5 x - 2c_6 t -c_6x^2)u + \alpha(x,t),
\end{split}
\end{equation*}
where $c_1,\dots,c_6$ are arbitrary constants and $\alpha$ is any solution of the heat equation.
It follows that the Lie algebra of infinitesimal symmetries of the heat equation is spanned
by the six vector fields
\begin{equation*}
\begin{split}
v_1 &= \pa_x\\
v_2 &= \pa_t\\
v_3 &= u\pa_u\\
v_4 &= x\pa_x + 2t\pa_t\\
v_5 &= 2t\pa_x -xu\pa_u\\
v_6 &= 4tx\pa_x + 4t^2\pa_t - (x^2+2t)u\pa_u
\end{split}
\end{equation*}
and the infinite-dimensional subalgebra
$$
v_\alpha = \alpha(x,t)\pa_u,
$$
where $\alpha$ is any solution of the heat equation. By \ref{infliesym} we know that
these generators span a Lie algebra of vector fields. Thus with any two generators 
also their Lie bracket is an infinitesimal symmetry of the heat equation.

The one-parameter symmetry group actions $\Fl^{v_i}_\eps$ corresponding to the generators can be 
calculated by solving the ODE 
$$
\frac{d}{d\eps}(x(\eps),t(\eps),u(\eps)) = v_i(x(\eps),t(\eps),u(\eps))
$$
with initial value $(x(0),t(0),u(0))=(x,t,u)$. This gives the groups
\begin{equation*}
\begin{split}
G_1:\ & (x+\eps,t,u)\\
G_2:\ & (x,t+\eps,u)\\
G_3:\ & (x,t,e^\eps u)\\
G_4:\ & (e^\eps x,e^{2\eps}t,u)\\
G_5:\ & (x+2\eps t,t,ue^{-\eps x -\eps^2 t})\\
G_6:\ & \left(\frac{x}{1-4\eps t},\frac{t}{1-4\eps t},u\sqrt{1-4\eps t}e^{\frac{-\eps x^2}{1-4\eps t}}\right)\\
G_\alpha:\ & (x,t,u+\eps\alpha(x,t))
\end{split}
\end{equation*}
Using \eqref{gdotf} we conclude that, given any solution $f$ of the heat equation, so are the following
functions (given by $g^{(i)}\cdot f$, for $g^{(i)}\in G_i$):
\begin{equation*}
\begin{split}
u^{(1)} &= f(x-\eps,t)\\ 
u^{(2)} &= f(x,t-\eps)\\
u^{(3)} &= e^\eps f(x,t)\\
u^{(4)} &= f(e^{-\eps}x,e^{-2\eps}t)\\
u^{(5)} &= e^{-\eps x +\eps^2 t}f(x-2\eps t,t)\\
u^{(6)} &= \frac{1}{\sqrt{1+4\eps t}}e^{\frac{-\eps x^2}{1+4\eps t}}f\left(\frac{x}{1+4\eps t},
\frac{t}{1+4\eps t}\right)\\
u^{(\alpha)} &= f(x,t) + \eps \alpha(x,t).
\end{split}
\end{equation*}
Here, $\eps$ is a real number and $\alpha$ is any solution of the heat equation.
$G_3$ and $G_\alpha$ are simply an expression of the linearity of the heat equation.
$G_1$ and $G_2$ reflect the space- and time-invariance of the equation. $G_4$ is 
a scaling symmetry, and $G_5$ represents a kind of Galilean boost to a moving
coordinate frame. $G_6$ is a genuinely local transformation group. Applying $G_6$
to a trivial solution $u=c$ ($c$ a constant) we obtain that also
$$
u = \frac{c}{\sqrt{1+4\eps t}}e^{\frac{-\eps x^2}{1+4\eps t}}
$$
is a solution. Setting $c=\sqrt{\eps/\pi}$ we obtain the fundamental solution to the
heat equation at the point $(x_0,t_0)=(0,-1/(4\eps))$. Translating in $t$ (using $G_2$)
we arrive at the fundamental solution at $(0,0)$,
$$
u = \frac{1}{\sqrt{4\pi t}}e^{\frac{-x^2}{4t}}.
$$
The most general one-parameter group of symmetries has a generator of the form $c_1v_1+\dots +c_6v_6 + v_\alpha$.
Also, an arbitrary element of the symmetry group of the heat equation (if sufficiently close
to the identity element) can be expressed in the form
$$
g = \Fl^{v_\alpha}_{\eps_\alpha}\circ \Fl^{v_6}_{\eps_6}\circ \dots \circ \Fl^{v_1}_{\eps_1},
$$
and so the most general solution obtainable from any given solution $u$ by such a transformation is
of the form
$$
u=\frac{1}{\sqrt{1+4\eps_6 t}}e^{\eps_3-\frac{\eps_5x + \eps_6x^2-\eps_5^2t}{1+4\eps_6t}}
f\left(\frac{e^{-\eps_4}(x-2\eps_5t)}{1+4\eps_6t} - \eps_1,\frac{e^{-2\eps_4}t}{1+4\eps_6t}-\eps_2\right) 
+\al(x,t),
$$
with $\eps_1,\dots,\eps_6$ real constants and $\alpha$ any solution of the heat equation.
\section{Nondegeneracy conditions}\index{nondegeneracy conditions}
The central result \ref{mainpde} only gives a sufficient condition for a local group
of transformations to be a symmetry group of a given system of differential equations.
In the present section we will see that under some mild condition on the system
the infinitesimal criterion from \ref{mainpde} is in fact also sufficient.

Let us first look at the difference between the criterion \ref{algsym} for the invariance
of a system of algebraic equations, which is necessary and sufficient, and \ref{mainpde}.
For a system of algebraic equations $F=0$, for each point $x_0$ on the subvariety $\mcS_F=\{x\mid F(x)=0\}$
there trivially exists a solution to the system, namely $x_0$ itself. In the case of a system
$P(x,\un)=0$ of differential equations, however, if $(x_0,\un_0)\in \mcS_P = \{(x,\un)\mid P(x,\un)=0\}$
there is no guarantee that there exists a solution $u=f(x)$ of the differential equation 
locally around $x_0$ such that $\un_0=\prol f(x_0)$. Our definition \ref{dgsgdef} of a symmetry group
of $P(x,\un)=0$ requires $\prol G$ to transform solutions into solutions, so we only can infer
that it will move points in $\mcS_P$ for which there exists a solution as above to other points
in $\mcS_P$. In general, these points need not constitute the entire set $\mcS_P$. Thus we define:
\bd A system $P(x,\un)=0$ of $n$-th order differential equations is called {\em locally solvable}\index{locally solvable}
at the point
$$
(x_0,\un_0)\in\mcS_P = \{(x,\un)\mid P(x,\un)=0\}
$$
if there exists a smooth function $u=f(x)$ defined near $x_0$ such that $\un_0=\prol f(x_0)$.
The system is called locally solvable if it is locally solvable at every point of $\mcS_P$.
It is called {\em nondegenerate}\index{nondegenerate} if it is locally solvable and
of maximal rank at every point $(x_0,\un_0)\in\mcS_P$.
\et
\brem\label{odenon} Consider a non-singular ODE of order $n$, i.e., of the form
\begin{equation}\label{oden1}
u_n = P(x,u,u_x,\dots,u_{n-1}),
\end{equation}
with $u_k=\pd{^ku}{x^k}$ for any $k$ and let $(x^0,u^0,\dots,u_n^0)\in \mcS_P$. Then we 
need to find a local smooth function $u=f(x)$ such that 
$$
u^0=f(x^0),\ u^0_x = f'(x^0),\dots, u^0_{n-1}=f^{(n-1)}(x^0),\ u^0_{n}=f^{(n)}(x^0).
$$
Here, the first $n$ conditions are the usual initial conditions for the ODE \eqref{oden1},
hence can always be satisfied by standard ODE solution theory. The last condition then
simply follows from \eqref{oden1} itself (i.e., from the fact that $(x^0,u^0,\dots,u_n^0)\in \mcS_P$). 
We conclude that non-singular (systems of) ODEs are always locally solvable, hence 
are in fact nondegenerate.
\et
Note that even for PDEs, local solvability is usually a very mild restriction -- e.g.,
in Cauchy problems one usually is looking for solutions with prescribed values on an
entire hypersurface.
\bex (i) The $2$-dimensional wave equation 
$$
P(x,t,u^{(2)}) = u_{tt}-u_{xx}=0
$$ 
is locally solvable.
In fact, let
$$
(x^0,t^0;u^0;u^0_x,u^0_t;u^0_{xx},u^0_{xt},u^0_{tt}) \in \mcS_P,
$$
i.e., $u^0_{tt}=u^0_{xx}$. Then we need to find a solution $u=f(x)$ of the wave equation
near $(x^0,t^0)$ with $\prolt f(x^0,t^0) = (u^0_x,u^0_t;u^0_{xx},u^0_{xt},u^0_{tt})$.
Since $u^0_{tt}=u^0_{xx}$, we can take for $f$ the polynomial solution
\begin{equation*}
\begin{split}
f(x,t) = u^0 + &u^0_x(x-x^0)+u^0_t(t-t^0) + \frac{1}{2}u^0_{xx}[(x-x^0)^2+(t-t^0)^2]\\
&+u^0_{xt}(x-x^0)(t-t^0).
\end{split}
\end{equation*}
(ii) The over-determined system
\begin{equation}\label{over}
u_x=yu,\quad u_y = 0
\end{equation}
is not locally solvable: given $(x_0,y_0)$, let
$$
u^0 = u(x_0,y_0)=1, \quad u^0_x = u_x(x_0,y_0)=y_0,\quad u^0_y = u_y(x_0,y_0)=0.
$$
Then there is no local solution with the prescribed jet at $(x_0,y_0)$. In fact,
by cross-differentiation we obtain
$$
0=u_{xy} = (yu)_y = yu_y + u = u,
$$
which is incompatible with the above prescription.

(iii) Apart from integrability conditions as in (ii), the second main reason why
a PDE may fail to be locally solvable is that there are no smooth solutions, even 
locally. The most famous example is due to H.\ Lewy (\cite{L}), who showed that
there exist smooth functions $h=h(x,y,z)$ such that the first order system
\begin{equation*}
\begin{split}
u_x-v_y+2yu_z+2xv_z &= h(x,y,z)\\
u_y+v_x-2xu_z+2yv_z &= 0
\end{split}
\end{equation*}
has no smooth solutions even locally around any point.

(iv) For analytic systems of PDEs (in Kovalevskaya-form), the Cauchy--Kovalevskaya theorem ensures local
solvability.
\et
For nondegenerate systems of differential equations we can now show that indeed
\eqref{mainsymcrit} is necessary and sufficient:
\bt\label{nesu} Let
$$
P_\nu(x,\un) = 0 \quad \nu=1,\dots,l
$$
be a nondegenerate system of differential equations defined on an open subset $M$
of $X\times U$. Let $G$ be a local transformation group acting on $M$.
Then $G$ is a symmetry group of the system if and only if
\begin{equation}\label{mainsymcrit1}
\prol v [P_\nu(x,\un)] = 0, \quad \nu=1,\dots,l, \text{ whenever } P(x,\un)=0
\end{equation}
for every infinitesimal generator $v$ of $G$.
\et
\pr The condition is sufficient by \ref{mainpde}. To see that it is also necessary,
note that by \ref{algsym}, condition \eqref{mainsymcrit1} is equivalent to 
$\prol G$ mapping $\mcS_P=\{(x,\un)\mid P(x,\un)=0\}$ into itself. Thus we have to
verify that any symmetry group $G$ of the system has this property.
Let $(x_0,\un_0)\in \mcS_P$. Then by local solvability we may find
a smooth function $u=f(x)$ defined near $x_0$ such that $\un_0=\prol f(x_0)$.
If $g\in G$ is such that $\prol g\cdot (x_0,\un_0)$ is defined, then by
shrinking the domain of $f$ if necessary, we can obtain that $\ftil = g\cdot f$
is a well-defined function near $\xtil_0$, where $(\xtil_0,\util_0)=g\cdot (x_0,u_0)$.
Then by \eqref{groupprol} we have
$$
\prol g \cdot (x_0,\un_0) = (\xtil_0,\prol (g\cdot f)(\xtil_0)) =: (\xtil_0,\util_0^{(n)}).
$$
Since $\ftil=g\cdot f$ is a solution, this shows that $(\xtil_0,\util_0^{(n)})\in \mcS_P$, 
as claimed.
\ep
\section{Integration of ordinary differential equations}\index{ordinary differential equation}
As an illustration of the power of the methods developed in this chapter, in the present section
we investigate the general problem of solving ordinary differential equations (ODEs). It turns out that 
all the standard solution techniques for special types of ODEs are instances of symmetry methods.
Moreover, knowledge of a sufficiently large group of symmetries of any given ODE allows to 
reduce the solution procedure to successive integrations (or quadratures,\index{quadrature} in classical terminology).

Consider a first-order ODE
\begin{equation}\label{firstode}
\pd{u}{x} = F(x,u).
\end{equation}
We will show that if a one-parameter symmetry group of \eqref{firstode} is known, then
it can be solved by integration. Let $G$ be a local one-parameter group of transformations acting
on an open subset $M$ of $X\times U=\R^2$, with infinitesimal generator
\begin{equation}\label{odeinf}
v = \xi(x,u)\pd{}{x} + \phi(x,u)\pd{}{u}.
\end{equation}
By \ref{gpf}, its first prolongation is given by
\begin{equation}\label{odeprol}
\prolo v = \xi(x,u)\pd{}{x} + \phi(x,u)\pd{}{u} + \phi^x(x,u,u_x)\pd{}{u_x},
\end{equation}
where
\begin{equation}\label{odeprol2}
\phi^x = D_x\phi - u_x D_x\xi = \phi_x + (\phi_u - \xi_x)u_x - \xi_u u_x^2.
\end{equation}
The infinitesimal criterion \eqref{mainsymcrit} for $v$ to generate a symmetry of \eqref{firstode} then
reads
\begin{equation}\label{odesymcrit}
\pd{\phi}{x} +\Big(\pd{\phi}{u} - \pd{\xi}{x}\Big)F - \pd{\xi}{u} F^2 = \xi \pd{F}{x} + \phi\pd{F}{u},
\end{equation}
for all $(x,u)$. Thus any solution to the PDE \eqref{odesymcrit} generates a symmetry of the ODE \eqref{firstode}.
Although at first sight, replacing the ODE \eqref{firstode} by the PDE \eqref{odesymcrit} may not seem 
to facilitate the problem, in practice one can often find symmetries directly, e.g.\ by geometric considerations.

Suppose that $v$ generates a symmetry of \eqref{firstode} and let $v|_{(x_0,u_0)}\not=0$. Then by straightening
out (see \cite[17.12]{KLie_new}) we may pick new coordinates 
\begin{equation}\label{newc}
y = \eta(x,u),\quad w= \zeta(x,u)
\end{equation}
near $(x_0,u_0)$ such that in the new coordinates we have $v=\pa_{w}$, and thereby
$$
\prolo v = v = \pd{}{w}
$$
The condition that \eqref{firstode} be invariant under $v$ then simply means that it must 
be independent of $w$, i.e., it can be written in the form
$$
\pd{w}{y} = H(y)
$$
for some smooth function $H$. This equation can then be solved trivially by integration:
$$
w = \int H(y)\,dy + \text{ const}.
$$
Transforming back to the $(x,u)$-coordinates gives the solution to \eqref{firstode}.

To find the required coordinate transform, note that in the new coordinates we require that
$$
v = v(\eta)\pd{}{y} + v(\zeta)\pd{}{w} \stackrel{!}{=} \pd{}{w},
$$
so in terms of $(x,u)$, we need to satisfy
\begin{equation}\label{etaz}
\begin{split}
v(\eta) &= \xi\pd{\eta}{x} + \phi\pd{\eta}{u} = 0\\
v(\zeta) &= \xi\pd{\zeta}{x} + \phi\pd{\zeta}{u} = 1.
\end{split}
\end{equation}
Here, the first equation means that $\eta$ is an invariant of $G$, hence can be found
by the method of characteristics indicated in \eqref{invode}. Once $\eta$ is known,
$\zeta$ can often be found by inspection. In any case, also the second equation can
be solved by the method of characteristics. As indicated above, there is no 
guarantee that solving \eqref{etaz} is indeed easier than solving the original 
equation \eqref{firstode}. For example, if 
\begin{equation}\label{umsonst}
\frac{\phi(x,u)}{\xi(x,u)} = F(x,u)
\end{equation}
then inserting $\phi=\xi\cdot F$ in \eqref{odesymcrit} gives an identity, i.e.,
finding a symmetry in this case is exactly the same problem as solving the original
equation. 
\bex\label{home} Consider the homogeneous equation\index{homogeneous ODE}
$$
\pd{u}{x} = F\left(\frac{u}{x}\right).
$$
This equation is invariant under the group of scaling transformations 
$$
G: (x,u)\mapsto (\lam x,\lam u),\quad \lam>0.
$$
Indeed, by \eqref{odeprol2}, the first prolongation of the generator $v=x\pa_x + u\pa_u$ is $v$ itself
and so one readily checks that \eqref{odesymcrit} is satisfied. New coordinates satisfying \eqref{etaz}
are given by
$$
y=\frac{u}{x}, \quad w=\log x.
$$
Then
$$
\pd{u}{x} = \frac{du/dy}{dx/dy} = \frac{x(1+yw_y)}{xw_y} = \frac{1+yw_y}{w_y},
$$
so the equation in the new coordinates becomes
$$
\frac{dw}{dy} = \frac{1}{F(y)-y},
$$
which can be solved by integration:
$$
w = \int \frac{dy}{F(y)-y} + c.
$$
\et
\bex Consider the equation $P\,dx + Q\,du = 0$. We show that if it possesses a one-parameter
symmetry with generator $v=\xi\pa_x + \phi\pa_u$, then the function
$$
R(x,u) = \frac{1}{\xi(x,u)P(x,u) + \phi(x,u)Q(x,u)}
$$
is an integrating factor.\index{integrating factor}

To see this, note first that the given ODE is of the form \eqref{firstode} with $F=-P/Q$.
Inserting this into \eqref{odesymcrit} and recalling \ref{nesu} and \ref{odenon} 
it follows that $v$ is a symmetry if and only if
\begin{equation}\label{intfac}
\begin{split}
\left(\xi\pd{P}{x} + \phi\pd{P}{u}\right)Q - &\left(\xi\pd{Q}{x} + \phi\pd{Q}{u}\right)P
+\pd{\phi}{x}Q^2  \\
&-\left(\pd{\phi}{u} 
- \pd{\xi}{x}\right)PQ  -\pd{\xi}{u}P^2 = 0.
\end{split}
\end{equation}
On the other hand, that $R$ is an integrating factor means that
$$
\pd{}{u}(RP) = \pd{}{x}(RQ),
$$
and inserting $R$ from above here gives
\begin{equation*}
\begin{split}
R^2 \Big( \phi \Big(Q\pd{P}{u} - P\pd{Q}{u}\Big) - &\pd{\xi}{u}P^2 - \pd{\phi}{u}PQ\Big) \\
&= R^2 \Big( \xi \Big(P\pd{Q}{x} - Q\pd{P}{x}\Big) - \pd{\xi}{x}PQ - \pd{\phi}{x}Q^2\Big),
\end{split}
\end{equation*}
which is equivalent to \eqref{intfac}, as claimed.
\et
Consider now a non-singular ODE of order $n$:
\begin{equation}\label{oden}
P(x,\un) = P(x,u,u_x,\dots,u_n),
\end{equation}
where $u_n=\frac{d^n u}{dx^n}$. We will show that knowledge of a one-parameter symmetry group $G$ of
\eqref{oden} with generator $v$ allows to reduce the order of the equation by one.

To see this, choose coordinates $y=\eta(x,u)$, $w=\zeta(x,u)$ as in \eqref{etaz} that straighten
out $v$, so that $v=\pa_w$. Then by the chain rule, we can re-write the derivatives of $u$ with
respect to $x$ in terms of $y$, $w$ and the derivatives of $w$ with respect to $y$:
$$
\frac{d^ku}{dx^k} = \delta_k\left(y,w,\od{w}{y},\dots,\odk{w}{y}\right),
$$
for suitable functions $\delta_k$. Substituting these expressions back into \eqref{oden}, we
obtain the equivalent $n$-th order equation
\begin{equation}\label{odene}
\tilde P(y,w\hn) = \tilde P(y,w,w_y,\dots,w_n) = 0 
\end{equation}
in the new coordinates. Now note that since \eqref{oden} has $G$ as a symmetry group, the same
is true of the transformed system \eqref{odene}: this follows from the definition of a symmetry
group and the fact that the prolongation of a group action according to \eqref{groupprol}
is equivariant under changing variables (again by the chain rule). In the new variables,
$\prol v = v = \pa_w$, and so the infinitesimal criterion \eqref{mainsymcrit} reduces to
$$
\prol v(\tilde P) = \pd{\tilde P}{w} = 0 \ \text{ whenever } \tilde P(y,w\hn)=0.
$$
Since $\tilde P$ is invariant under $G\hn$, by \ref{invvar} and \ref{nesu} (and provided that $G$ acts 
semi-regularly with orbits of constant dimension) there exists an equivalent
equation that depends exclusively on a complete set of functionally independent invariants
of $G\hn$. In the present situation, such a set of invariants is obviously given by
$y,w_y,w_{yy},\dots,w_n$. Consequently, there is an equivalent equation
$$
\hat P(y,w_y,\dots,w_n) = 0,
$$
which is independent of $w$. Thus, setting $z=w_y$ we obtain an equivalent ODE of order
$n-1$, as claimed:
$$
\hat P\Big(y,z,\dots,\frac{d^{n-1}z}{dy^{n-1}}\Big) = \hat P(y,z^{(n-1)}) = 0.
$$
Given any solution $z=h(y)$ of this ODE, $w=\int h(y)\,dy + c$ is a solution to
\eqref{odene}, and transforming back to $(x,u)$-coordinates gives a solution to \eqref{oden}.
\bex Any homogeneous linear ODE of second order
\begin{equation}\label{hom2}
u_{xx} + p(x)u_x + q(x)u = 0
\end{equation}
is invariant under the scaling group $G: (x,u)\mapsto (x,\lam u)$. In fact, the generator of
$G$ is $v=u\pa_u$, so $\prolt v = u\pa_u + u_x\pa_{u_x} + u_{xx}\pa_{u_{xx}}$, and so the 
infinitesimal criterion \eqref{mainsymcrit} is satisfied. Coordinates satisfying \eqref{etaz}
are given by
$$
y = x,\quad w=\log u,
$$
so $v=\pa_w$. To re-express the equation in the new coordinates, note that
$$
u=e^w,\quad u_x=w_x e^w,\quad u_{xx} = (w_{xx}+w_x^2)e^w,
$$
giving
$$
w_{xx} + w_x^2 + p(x)w_x + q(x) = 0,
$$
which, as expected, is independent of $w$. Setting $z=w_x=u_x/u$ we therefore obtain
the well-known transformation of \eqref{hom2} into the Riccati equation\index{Riccati equation}
$$
z_x = - z^2 - p(x)z -q(x).
$$
\et
\section{Differential invariants}\index{differential invariant}
In this section we look at a converse to the problem of determining the symmetry group of
a differential equation, namely: Given a local transformation group, what is the most general 
type of differential equation that admits it as a symmetry group? We will give an answer to 
this question in the setting of ordinary differential equations.

Recall from \ref{nesu} that a nondegenerate differential equation $P(x,\un)=0$ admits a local
transformation group as a symmetry group if and only if 
the corresponding subvariety $\mcS_P$ of $M\hn$ is invariant under $\prol G$.
Provided that $\prol G$ acts semi-regularly with orbits of constant dimension, \ref{invvar} then
implies that there is an equivalent equation $\tilde P=0$ describing the subvariety
$\mcS_P$ and depending only on the functionally independent invariants of $\prol G$.
We therefore give a name to the invariants of $\prol G$:
\bd Let $G$ be a local transformation group acting on $M\sse X\times U$ and let $n\ge 1$. An $n$-th
order {\em differential invariant}\index{differential invariant} of $G$ is a smooth function
$\eta:M\hn\to \R$ such that $\eta$ is an invariant of $\prol G$:
$$
\eta(\prol g \cdot (x,\un)) = \eta(x,\un)
$$
for all $(x,\un)\in M\hn$ and $g\in G$ such that $\prol g \cdot (x,\un)$ is defined.
\et
\bex\label{rotdi} For the rotation group $G=SO(2)$ on $X\times U=\R^2$ we have $v=-u\pa_x+x\pa_u$, and
by \eqref{rotvprol} we have
$$
\prolo v = -u\partial_x + x\partial_u + (1+u_x^2)\partial_{u_x}
$$
By definition, the first order differential invariants of $G$ are the invariants of $\prolo G$, 
which can be calculated by the method of characteristics, cf.\ \eqref{invode}. This gives 
the following complete set of first order invariants of $G$:
$$
y= \sqrt{x^2+u^2}, \quad w=\frac{xu_x-u}{x+uu_x}.
$$
Any other first order invariant must be a function of these (by \ref{lociprop}).
\et
Our next aim is to show that one can always calculate higher order differential invariants
inductively from known ones of lower order. To show this we need an auxiliary result first.
\blem Let $v=\xi\pa_x+\phi\pa_u$ be a vector field on $M\sse X\times U=\R^2$ and let
$\zeta = \zeta(x,\un)$ be smooth. Then
\begin{equation}\label{prnp1}
\prolp v (D_x\zeta) = D_x[\prol v(\zeta)] - D_x\xi\cdot D_x\zeta.
\end{equation}
\et
\pr By \eqref{charprol},
$$
\prolp v(D_x\zeta) = \prolp v_Q(D_x\zeta) + \xi D_x^2\zeta,
$$
and
$$
D_x[\prol v(\zeta)] = D_x[\prol v_Q(\zeta)] + D_x(\xi D_x\zeta).
$$
It therefore remains to show that 
$$
\prolp v_Q(D_x\zeta) = D_x[\prol v_Q(\zeta)].
$$
Using \eqref{charparts}, we calculate:
\begin{equation*}
\begin{split}
\prolp v_Q(D_x\zeta) &= \sum_{j=0}^{n+1} D_jQ \pd{}{u_j}\left(\pd{\zeta}{x}+
\sum_{k=0}^{n}u_{k+1}\pd{\zeta}{u_k}\right)\\
 &= \sum_{j=0}^{n+1} D_jQ \left(\pd{^2\zeta}{u_j\pa x} 
 +\sum_{k=0}^{n}\left(\delta^{k+1}_j\pd{\zeta}{u_k}+u_{k+1}\pd{^2\zeta}{u_j\pa u_k}\right)\right)\\
 &= \sum_{j=0}^n\left(D_{j+1}Q\pd{\zeta}{u_j} + D_jQ\left( \pd{^2\zeta}{u_j\pa x} +  
 \sum_{k=0}^{n}u_{k+1}\pd{^2\zeta}{u_j\pa u_k}\right)\right)\\
 &= D_x\left(\sum_{j=0}^{n} D_jQ \pd{\zeta}{u_j}\right) = D_x[\prol v_Q(\zeta)],
\end{split}
\end{equation*}
where in the penultimate line we used that $\zeta=\zeta(x,\un)$, hence its derivative with respect to
$u_{n+1}$ vanishes.
\ep
\bp\label{nextderp} Let $G$ be a local one-parameter group of transformations acting on $M\sse X\times U=\R^2$. 
Let $y=\eta(x,\un)$
and $w=\zeta(x,\un)$ be $n$-th order differential invariants of $G$ ($n\ge 1$). Then the derivative
\begin{equation}\label{nextder}
\od{w}{y} = \frac{dw/dx}{dy/dx} \equiv \frac{D_x\zeta}{D_x\eta}
\end{equation}
is a differential invariant of $G$ of order $n+1$.
\et
\pr Note first that for any vector field $Z$ and any smooth functions $f$, $g$ we have
the usual quotient rule $Z(f/g)=g^{-2}(Z(f)g-fZ(g))$, as follows from the derivation 
description of vector fields. Using this and \eqref{prnp1}, we obtain from \eqref{nextder}:
\begin{equation*}
\begin{split}
\prolp v \left(\od{w}{y}\right) &= \frac{1}{(D_x\eta)^2}\big(\prolp v(D_x\zeta)\cdot D_x\eta
-D_x\zeta\cdot \prolp v(D_x\eta)\big)\\
&= \frac{1}{(D_x\eta)^2}\big(D_x[\prol v(\zeta)]\cdot D_x\eta - D_x\xi\cdot D_x\zeta\cdot D_x\eta\\
&\hphantom{xxxxxxxxxx}-D_x\zeta\cdot D_x[\prol v(\eta)] + D_x\zeta\cdot D_x\xi\cdot D_x\eta\big)=0
\end{split}
\end{equation*}
since $\zeta$ and $\eta$ are $n$-th order invariants, so $\prol v(\zeta) = \prol(\eta) = 0$.
\ep
\bc\label{oi} Let $G$ be a local one-parameter group of transformations acting on $M\sse X\times U=\R^2$. 
Let $\{y=\eta(x,u),\, 
w=\zeta(x,u,u_x)\}$ be a complete set of functionally independent invariants of $\prolo G$. Then
(locally around any point) the derivatives
$$
y,w,\od{w}{y},\dots,\od{^{n-1}w}{y^{n-1}}
$$
form a complete set of functionally independent invariants for 
$\prol G$ for $n\ge 1$.
\et
\pr Successively applying \ref{nextderp} it follows that the derivatives up to order $k$
are invariants of $\mathrm{pr}^{(k)}G$, hence also of any higher prolongation of $G$
(by \eqref{graprol}). So it only remains to show that they are functionally independent.
But this is immediate because $\odk{w}{y}$ depends explicitly on $u_{k+1}$, hence
is independent of the previously constructed invariants
$y,w,\pd{w}{y},\dots,\od{^{k-1}w}{y^{k-1}}$, since those only depend on
$x,u,\dots,u_k$. By \ref{lociprop}, since $G$ has one-dimensional orbits 
and $\dim X\times U^{(n)} = n+2$, completeness follows.
\ep
\bex\label{rotsec} Applying the previous result to the first order invariants of the rotation group
given in \ref{rotdi} it follows that $y$, $w$, and
$$
\od{w}{y} = \frac{dw/dx}{dy/dx} = \frac{\sqrt{x^2+u^2}}{(x+uu_x)^3}[(x^2+u^2)u_{xx}
-(1+u_x^2)(xu_x-u)]
$$
together form a complete set of functionally independent invariants for $\prolt G$.
Thus any other second order differential invariant can be written as a smooth function
of these invariants. For example, the curvature $\kappa$ from \eqref{curv} is 
an invariant of $\prolt G$ (as one checks by verifying that $\prolt v(\kappa) = 0$),
and in fact
$$
\kappa =\frac{u_{xx}}{(1+u_x^2)^{3/2}} = \frac{w_y}{(1+w^2)^{3/2}} + \frac{w}{y(1+w^2)^{1/2}}.
$$
\et
The important point to note now is that once the differential invariants of a local
group of transformations $G$ on $M\sse X\times U=\R^2$ are known, we can determine 
all ODEs that admit $G$ as a symmetry group, i.e., all ODEs that can be integrated
using $G$:
\bt\label{it} Let $G$ be a connected local transformation group acting on $M\sse X\times U=\R^2$
and let $\eta^1(x,\un),\dots,\eta^k(x,\un)$ be a complete set of functionally independent
$n$-th order differential invariants of $G$. 
\begin{itemize}
\item[(i)] If $\prol G$ acts semi-regularly on $M\hn$, then
an $n$-th order nondegenerate differential equation $P(x,\un)=0$ admits $G$ as a symmetry group if and only if,
locally around any point in $\mcS_P=\{(x,\un)\mid P(x,\un)=0\}\sse M\hn$, there is an equivalent equation
$$
\tilde P(\eta^1(x,\un),\dots,\eta^k(x,\un))=0
$$
involving only the differential invariants of $G$.
\item[(ii)] If $G$ is a one-parameter group of transformations, any nondegenerate $n$-th
order differential equation having $G$ as a symmetry group is equivalent
to an $(n-1)$-st order equation
$$
\tilde P\Big(y,w,\od{w}{y},\dots,\od{^{n-1}w}{y^{n-1}}\Big) = 0,
$$
where $y=\eta(x,u),\, w=\zeta(x,u,u_x)$ form a complete set of functionally independent
invariants of $\prolo G$.
\end{itemize}
\et
\pr (i) Clearly if there is an equivalent equation depending only on invariants then
$\mcS_P$ itself is locally invariant. Since $P$ is non-degenerate, $\mcS$ is a closed
regular submanifold of $M\hn$, so the proof of \ref{algsym} shows that this implies
that $\mcS$ is invariant under $\prol G$, which in turn yields that $G$ is a symmetry group of $P$. 
Conversely, if $G$ is a symmetry group then the proof of \ref{nesu} shows that $\prol G$
leaves $\mcS_P$ invariant. The claim then follows from \ref{invvar}.

(ii) This is immediate from (i) and \ref{oi}.
\ep
\bex Continuing our analysis of the rotation group from \ref{rotsec}, we can use
\ref{it} to find the most general first and second order ODEs that admit the 
rotation group as a symmetry group. By \ref{it} (ii), any such first order equation
is equivalent to one involving only the first order invariants given in \ref{rotdi}.
If we solve for $w$ it follows that any such equation is of the form
$$
\frac{xu_x-u}{x+uu_x} = F(\sqrt{x^2+u^2}).
$$
Analogously, any second order equation invariant under rotations is equivalent
to one only involving $y$, $w$, and, say, the curvature $\kappa =\frac{u_{xx}}{(1+u_x^2)^{3/2}}$,
hence is of the general form
$$
u_{xx} = (1+u_x^2)^{3/2} F\Big(\sqrt{x^2+u^2},\frac{xu_x-u}{x+uu_x}\Big).
$$
\et
We may also employ \ref{it} to obtain an alternative way of reducing the order of a
given differential equation $P(x,\un)=0$ once we know a one-parameter symmetry group. 
In fact, we know from \ref{it} (ii) that the equation must be equivalent to 
one involving only $y,w,\od{w}{y},\dots,\od{^{n-1}w}{y^{n-1}}$. But this latter
equation automatically is of order $n-1$. Thus simply by re-expressing the equation
in terms of the differential invariants already reduces the order by $1$.
Furthermore, once the solution $w=h(y)$ of the reduced equation is known, the
solution of the original equation is found by solving the auxiliary first order equation
\begin{equation}\label{auxi}
\zeta(x,u,u_x) = h(\eta(x,u))
\end{equation}
obtained by substituting for $y$ and $w$ their expression in terms of $x$ and $u$.
As \eqref{auxi} depends only on the invariants $\eta$, $\zeta$ of $\prolo G$, it
has $G$ as a one-parameter symmetry group, hence can be integrated using the methods
described above for first order equations.
\bex The second order equation
\begin{equation}\label{secco}
x^2u_{xx} + xu_x^2 = uu_x
\end{equation}
is invariant under the scaling group $G: (x,u)\mapsto (\lam x,\lam u)$. In fact,
the infinitesimal generator of $G$ is $v=\left.\od{}{\lam}\right|_{\lam=1} (\lam x,\lam u)= 
x\pa_x + u\pa_u$. Hence \eqref{phiform} gives
$$
\prolt v = x\pa_x + u\pa_u - u_{xx}\pa_{u_{xx}},
$$
so 
$$
\prolt v(x^2u_{xx} + xu_x^2 - uu_x) = 2x^2u_{xx} + xu_x^2 -uu_x -x^2u_{xx} = 0
$$
for any $(x,u,u_x,u_{xx})\in \mcS_P$. By \ref{oi} and \ref{home}, a complete set of 
functionally independent invariants of second order is given by
$$
y=\frac{u}{x},\quad w=u_x,\quad \od{w}{y} = \frac{dw/dx}{dy/dx}=\frac{x^2u_{xx}}{xu_x-u}.
$$
Inserting this into \eqref{secco} gives the transformed equation
$$
(w-y)\od{w}{y} + w^2 = yw,
$$
which, as expected, is of first order. This has two families of solutions: either $w=y$
(which is a singular solution), or $\od{w}{y}=-w$, which implies $w=ce^{-y}$.
Re-substituting gives
$$
\od{u}{x}=\frac{u}{x}, \quad\text{ or }\quad  \od{u}{x}=ce^{-u/x}.
$$
This results in the one-parameter family of singular solutions $u=kx$, and the implicit
solutions
$$
\int\frac{dy}{ce^{-y}-y} = \log x + k,
$$
with $y=u/x$, which constitutes the general solution to the original equation \eqref{secco}.
\et
\bex We return to the two-dimensional heat equation $u_t=u_{xx}$, 
whose symmetry group we calculated in \ref{heatsym},
and use the methods of this section to determine the general form of {\em traveling wave
solutions}\index{traveling wave solution} to the heat equation.
In general, a traveling wave solution to a differential equation is one which is 
invariant under a translation group. Here we consider the symmetry of the heat equation
given by the translation group
$$
(x,t,u) \mapsto (x+c\eps,t+\eps,u),\quad \eps\in \R,
$$
with generator $v=\pa_t + c\pa_x$, for $c$ some fixed constant. Global invariants are
given by
\begin{equation}\label{heati}
y=x-ct,\quad v=u.
\end{equation}
Thus any group-invariant solution is of the form $v=h(y)$, or $u=h(x-ct)$, which is
a function that doesn't change its shape along the straight line $x-ct=$ const. If we
express the derivatives of $u$ in terms of the new variables we obtain
$$
u_t = -cv_y,\quad u_x=v_y,\quad u_{xx} = v_{yy}.
$$
Inserting into the heat equation, we arrive at the ODE determining all traveling 
wave solutions of the above form:
$$
-cv_y=v_{yy},
$$
with general solution
$$
v(y) = ke^{-cy} + l
$$
($k$, $l$ arbitrary constants). If we now substitute back into the heat equation
we obtain the general form of traveling wave solutions:
$$
u(x,t) = ke^{-c(x-ct)} + l.
$$
\et
Much more can be said about group-invariant solutions of PDEs. In fact, there is an entire
theory that allows to calculate solutions invariant under any subgroup of the full 
symmetry group of any given system of differential equations. While these methods also
basically 
rest on choosing new variables from invariants of the given sub-group, the case 
of higher (than one) dimensional symmetry groups is significantly more involved than what
we considered in this section. For a precise formulation, the theory of extended
jet bundles and group invariant prolongation is required. We refer to \cite[Ch.\ 3]{O}
for an in-depth study.

\chapter{Variational symmetries}
Conservation laws (like conservation of energy, momentum, etc.) play a central role in physics. In 
contemporary theoretical physics one usually formulates the governing principles of a theory in the
form of a variational principle -- roughly, this says that some relevant quantity is to take an optimal
value. It is a deep and far-reaching discovery of Emmy Noether (1918) that for such systems, every
conservation law comes from a corresponding symmetry property. In the present chapter we give an 
introduction to this field.
\section{Calculus of variations}
Let $X=\R^p$, with coordinates $(x^1,\dots,x^p)$, and $U=\R^q$ with coordinates $(u^1,\dots,$ $u^q)$.
Also, let $\Omega$ be a connected open subset of $X$ with smooth boundary $\pa \Omega$. By a 
{\em variational problem}\index{variational problem} we mean the problem of finding extrema 
(maxima, minima, or stationary points) of the functional
$$
\mcl[u] = \int_\Omega L(x,\un(x))\,dx
$$
in some class of admissible functions $u=f(x)$ on $\Omega$. The integrand $L:X\times U\hn\to \R$
is called the {\em Lagrangian}\index{Lagrangian} of the variational problem $\mcl$. The class 
of admissible functions varies with the problem under consideration (e.g., through the
imposition of various boundary or regularity conditions). We will confine ourselves
here to studying smooth variational problems. The calculus of variations is a vast subject and
we will here barely be able to scratch the surface, with a focus on symmetry methods. For a 
gentle introduction see, e.g., \cite{Sa}, a rather comprehensive treatise is \cite{GH1,GH2}.
\bex\label{first} (i) To find the shortest curve $u=f(x)$ connecting two points $(a,b)$,
$(c,d)$ in the plane, one needs to minimize the length of $u=f(x)$:
$$
\mcl[u] = \int_a^c \sqrt{1+u_x^2}\,dx.
$$
(ii) More generally, if $(M,g)$ is a smooth Riemannian manifold and $x$, $y\in M$ then the problem of
finding the shortest curve from $x$ to $y$ is a variational problem, namely that of minimizing
the arc-length of curves $u\colon [0,1]\to M$ connecting $x$ and $y$:
$$
\mcl[u] = \int_0^1 \sqrt{g(u'(t),u'(t))}\,dt.
$$
Solutions to this problem are called geodesics of $(M,g)$.
\et
The basic approach to finding extremals of variational problems is very similar to that
of extremizing smooth functions in basic real analysis. There, given a function $f\colon \Omega\to \R$, 
one notes that if $x$ is an extremum of $f$ then for any $y$, the one-dimensional function
$\eps\mapsto f(x+\eps y)$ must have an extremal. Therefore, we must have
$$
0= \dde f(x+\eps y) = \langle \nabla f(x),y\rangle \quad (y\in \R^p),
$$
and so the gradient $\nabla f$ of $f$ at $x$ vanishing is a necessary condition for
$x$ being an extremum of $f$. 

For functionals $\mcl[u]$, the role of the gradient is taken over by what is called
the {\em variational derivative}.\index{variational derivative} Moreover, the inner product
on $\R^n$ gets replaced by the $L^2$-inner product of functions $f,\,g\colon \Omega\to \R^q$:
$$
\la f,g\ra = \int_\Omega f(x)\cdot g(x)\,dx = \int_\Omega \sum_{\al=1}^q f^\al(x)g^\al(x)\,dx.
$$
With these notations we can define:
\bd Let $\mcl[u]$ be a variational problem. The variational derivative\index{variational derivative}
of $\mcl$ is the unique $q$-tuple 
$$
\delta\mcl[u] = (\delta_1\mcl[u],\dots,\delta_q\mcl[u]),
$$
$\delta\mcl:\cinfty(\Omega)\to \cinfty(\Omega)^q$ with the property that
\begin{equation}\label{varder}
\dde\mcl[f+\eps\eta] = \int_\Omega \delta\mcl[f](x)\cdot \eta(x)\,dx
\end{equation}
for every smooth function $f$ on $\Omega$ and any compactly supported smooth
function $\eta\in {\mathcal D}(\Omega)^q$ such that $f+\eps\eta$ is itself
admissible. $\delta_\alpha\mcl$ is called the variational derivative of $\mcl$
with respect to $u^\alpha$.
\et
Here, by ${\mathcal D}(\Omega)$ we denote the space of test functions on $\Omega$,
${\mathcal D}(\Omega)=\{f\in \cinfty(\Omega)\mid \supp(f)\comp \Omega\}$.
\bp If $u=f(x)$ is an extremal of $\mcl[u]$, then
\begin{equation}\label{cv0}
\delta\mcl[f](x)=0 \quad \forall x\in \Omega.
\end{equation}
\et
\pr For any $\eta\in \D(\Om)^q$ and for all $\eps$ small, $f+\eps\eta$ is admissible, and
by assumption, $\eps\mapsto \mcl[f+\eps\eta]$ must have an extremum at $\eps=0$. Hence
by classical real analysis, 
$$
0=\dde\mcl[f+\eps\eta] = \int_\Omega \delta\mcl[f](x)\cdot \eta(x)\,dx.
$$
As this has to hold for any $\eta$, the claim follows.
\ep
Using the notion of total derivative, cf. \eqref{totder1}, we now wish to derive an explicit
formula for the variational derivative. First, 
\begin{equation*}
\begin{split}
\dde\mcl[f+\eps\eta] &= \int_\Omega \dde L(x,\prol(f+\eps\eta)(x))\,dx\\
&= \int_\Om \Big( \sum_{\alpha,J} \pd{L}{u^\alpha_J}(x,\prol f(x))\cdot \pa_J\eta^\al(x)\Big)\,dx.
\end{split}
\end{equation*}
Since $\eta$ has compact support, we may apply the divergence theorem to this expression, with
all boundary terms vanishing. Note also that when a partial derivative $\pa_{x_j}$ is applied
to $\pa L/\pa u^\al_J$ then the result can be expressed by a total derivative, so
\begin{equation}
\dde\mcl[f+\eps\eta] = \int_\Om  \sum_{\alpha=1}^q \big[\sum_J (-D)_J\pd{L}{u^\alpha_J}
(x,\mathrm{pr}^{(n+\sharp J)} f(x)) \big]\eta^\al(x)\,dx.
\label{mcl}
\end{equation}
Here, for $J=(j_1,\dots,j_k)$ we set $(-D)_J:=(-D_{j_1})(-D_{j_2})\dots(-D_{j_k})$. This formula
is the first instance of the following operator that will play a central role below.
\bd For $1\le \al\le q$, the $\al$-th Euler operator\index{Euler operator} is defined by
\begin{equation}\label{eulerop}
E_\al = \sum_J (-D)_J \pd{}{u^\al_J},
\end{equation}
where the sum extends over all $J=(j_1,\dots,j_k)$ with $1\le j_\kappa\le p$, $k\ge 0$.
Moreover, we set $E(L):=(E_1(L),\dots,E_q(L))$.
\et
Although the above sum is formally infinite, whenever $E_\alpha$ is applied to some concrete
$L(x,\un)$ only finitely many terms are required. 

In terms of the Euler operator, \eqref{mcl} says that
$$
\delta\mcl = E(L)
$$
Together with \ref{cv0}, this gives the classical necessary condition for extremals of a
variational problem:
\bt If $u=f(x)$ is a smooth extremal of the variational problem $\mcl[u]=\int_\Om L(x,\un)\,dx$
then it must be a solution of the {\em Euler--Lagrange equations}\index{Euler--Lagrange equations}
$$
E_\alpha(L)=0, \quad \alpha=1,\dots,q.
$$
\hspace*{\fill}$\Box$
\et
\bex Let $p=q=1$, i.e., we consider a variational problem for a real-valued function of one variable.
In this case,
$$
E = \sum_{j=0}^\infty (-D_x)^j\pd{}{u_j} = \pd{}{u} - D_x\pd{}{u_x} + D_x^2\pd{}{u_{xx}}-\dots
$$ 
Thus the Euler--Lagrange equation for an $n$-th order variational problem 
$$
\mcl[u] = \int_a^b L(x,\un)\,dx
$$
is given by (setting $u_j=\pd{^j u}{x^j}$)
$$
0=E(L) = \pd{L}{u} - D_x\pd{L}{u_x} + D_x^2\pd{L}{u_{xx}} -\dots + (-1)^n D_x^n\pd{L}{u_n}.
$$
This is an ODE of order (at most) $2n$. The most common case of a first order variational problem
$L=L(x,u,u_x)$ gives
$$
0=E(L) = \pd{L}{u} - D_x\pd{L}{u_x} = \pd{L}{u} - \frac{\pa^2 L}{\pa x\pa u_x} - 
u_x \frac{\pa^2 L}{\pa u\pa u_x} - u_{xx} \pd{^2L}{u_x^2}.
$$
In particular, for the curve length problem \ref{first} the Euler--Lagrange equation reads
$$
-D_x\left(\frac{u_x}{\sqrt{1+u_x^2}}\right) = - \frac{u_{xx}}{(1+u_x^2)^{3/2}} = 0,
$$
i.e., $u_{xx}=0$. As geometrically expected, the solutions are straight lines $u=kx+d$.
\et
\bex The Dirichlet principle:\index{Dirichlet principle} Consider the variational problem
of minimizing the total energy (kinetic plus potential) of some system in the form
$$
\mcl[u] = \int_\Om \frac{1}{2}|\nabla u|^2 - u h\,dx 
$$
with some external potential $h$. Then the Euler--Lagrange equation reads
$$
0=E(L)=\pd{L}{u} - \sum_{i=1}^p D_i\pd{L}{u_i} = \pd{L}{u}-\sum_{i=1}^p D_i u_i
= -h -\Delta u,
$$
i.e., we obtain the Poisson equation\index{Poisson equation} $-\Delta u = h$.
\et
\section{Variational symmetries}
In this section we want to develop a notion of symmetry that applies to variational 
problems, similar to the symmetry groups of differential equations studied in the
previous chapter. Consider a variational problem
\begin{equation}\label{vp}
\mcl[u] = \int_{\Om_0} L(x,\un)\,dx.
\end{equation}
Let $G$ be a local transformation group on an open subset $M$
of $\Om_0\times U\sse X\times U$. If $u=f(x)$ is a smooth function on a 
sufficiently small subdomain $\Om\sse\Om_0$ such that the graph of $f$
lies in $M$ then any transformation $g\in G$ sufficiently close to the
identity will transform $f$ into another smooth function $\util=\ftil(\xtil)
=(g\cdot f) (\xtil)$ defined on some $\tilde\Om\sse\Om_0$. Therefore, we
define:
\bd A local transformation group $G$ acting on $M\sse \Om_0\times U$
is called a {\em variational symmetry}\index{variational symmetry} of the
functional \eqref{vp} if whenever $\Omega$ is a subdomain with $\overline{\Om}\sse\Om_0$,
$u=f(x)$ is a function defined over $\Om$ with graph contained in $M$, and
$g\in G$ is such that $\util=\ftil(\xtil) =(g\cdot f) (\xtil)$ is a well-defined
function defined on $\tilde\Om\sse \Omega_0$, then
\begin{equation}\label{varsym}
\int_{\tilde\Om} L(\xtil,\prol \ftil(\xtil))\,d\xtil = \int_\Om L(x,\prol f(x))\,dx.
\end{equation}
\et
Our first aim, in line with our previous considerations on symmetries of differential
equations, is to derive an infinitesimal criterion for variational symmetries. For this, 
we need a few preparations. 

To begin with, we need to study how a variational problem transforms under the action
of a local transformation group. Let
\begin{equation}\label{chvar}
\xtil = \Xi(x,u), \qquad \util = \Phi(x,u)
\end{equation}
be any change of variables. Then there is an induced change of variables
$$
\util\hn = \Phi\hn(x,\un)
$$
for the derivatives, given by prolongation (i.e., by differentiating \eqref{chvar} by the chain rule).
If the conditions of the inverse function theorem are satisfied then \eqref{chvar} determines
a transformed function $\util=\ftil(\xtil)$. Thereby, any functional
$$
\mcl[f] = \int_\Om L(x,\prol f(x))\,dx
$$
is transformed into a new functional
$$
\tilde \mcl[\ftil] = \int_{\tilde\Om} \tilde L(\xtil,\prol \ftil(\xtil))\,d\xtil.
$$
Here, the new domain $\tilde\Om = \{\xtil=\Xi(x,f(x))\mid x\in  \Om\}$ depends both on $\Om$
and on $f$. By the change of variables formula we get
\begin{equation}\label{trvar}
L(x,\prol f(x)) = \tilde L(\xtil,\prol \ftil(\xtil))\det J(x,\prolo f(x)),
\end{equation}
where $J$ is the Jacobian matrix with entries
\begin{equation}\label{jacob}
J^{ij}(x,\prolo f(x)) = \pd{}{x^j}\Xi^i(x,f(x)) = D_j\Xi^i(x,\prolo f(x)) \quad i,j=1,\dots,p,
\end{equation}
where, for simplicity, we assume that $\det J(x)>0$, otherwise we have to add absolute value signs.
\bd\label{totdiv} Let $x=(x^1,\dots,x^p)$ and let $P=P(x,\un)=(P_1(x,\un),\dots,$ $P_p(x,\un))$ be a $p$-tuple
of smooth functions. Then the {\em total divergence}\index{total divergence} of $P$ is
$$
\Div P = D_1P_1+\dots +D_pP_p,
$$
where $D_j$ denotes the total derivative with respect to $x^j$.
\et
With these notations, we have:
\blem Let  $v=\sum_i\xi^i\pa_{x^i}
+\sum_\al \phi^\al\pa_{u^\al}$ be a vector field on $M\sse \Omega_0\times U$
and set $g_\eps:=\Fl^{v}_\eps = (\Xi_{g_\eps},\Phi_{g_\eps})$. Denoting by $J_{g_\eps}$ the Jacobian 
\eqref{jacob} corresponding to $\Xi_{g_\eps}$, we have
\begin{equation}\label{jacobform}
\frac{d}{d\eps}(\det J_{g_\eps}(x,u\heins)) = (\Div\xi)(\prolo g_\eps\cdot (x,u\heins))\cdot\det J_{g_\eps}(x,u\heins).
\end{equation}
\et
\pr Recall first that for matrices $A$, $B$ with $A$ invertible, the derivative
of the determinant function at $A$ in the direction $B$ is given by $\det'(A)(B) = \det(A)\cdot\tr(A^{-1}B)$.
Consequently, for any fixed function $u$ we get
\begin{equation}\label{jumf}
\frac{d}{d\eps}\det J_{g_\eps} = \det{}'(J_{g_\eps})\Big(\frac{d}{d\eps}J_{g_\eps}\Big)
= \det J_{g_\eps}\cdot\tr(J_{g_\eps}^{-1}\frac{d}{d\eps}J_{g_\eps}).
\end{equation} 
By \eqref{jacob}, $J_{g_\eps}(x,u\heins(x)) = \pa_x(\Xi_{g_\eps}(x,u(x)))$ (with $\pa_x$ abbreviating the
$x$-gradient, and similar for $\pa_1$, $\pa_2$ below), so
\begin{equation*}
\begin{split}
&\frac{d}{d\eps}(J_{g_\eps}(x,u\heins)(x)) = \pa_x\Big(\ddee\Xi_{g_\eps}(x,u(x))\Big)
=\pa_x(\xi(\Xi_{g_\eps}(x,u(x)),\Phi_{g_\eps}(x,u(x))))\\
&=(\pa_1\xi)(\Xi_{g_\eps}(x,u(x)),\Phi_{g_\eps}(x,u(x)))\pa_x(\Xi_{g_\eps}(x,u(x)))\\
& \hphantom{=}+ (\pa_2\xi)(\Xi_{g_\eps}(x,u(x)),\Phi_{g_\eps}(x,u(x)))\pa_x(\Phi_{g_\eps}(x,u(x)))\\
&=[(\pa_1\xi)(\dots) + \pa_2\xi(\dots) \cdot \pa_x(\Phi_{g_\eps}(x,u(x)))\cdot \pa_x(\Xi_{g_\eps}(x,u(x)))^{-1}]
\pa_x(\Xi_{g_\eps}(x,u(x)))
\end{split}
\end{equation*}
Now let $\tilde u$ be the transformed function $g_\eps\cdot u$, so that with $\xtil(x)=\Xi_{g_\eps}(x,u(x))$
we have $\util(\xtil(x))=\Phi_{g_\eps}(x,u(x))$. Then
$$
\pa_x(\Phi_{g_\eps}(x,u(x))) = \pa_x(\util\circ\xtil)(x) = \pa_\xtil\util(\xtil(x))\cdot \pa_x\xtil(x)
=\pa_\xtil\util(\xtil(x))\cdot \pa_x(\Xi_{g_\eps}(x,u(x))),
$$
so, together with the above calculations we obtain
\begin{equation*}
\begin{split}
\frac{d}{d\eps}(J_{g_\eps}(x,u\heins)(x)) &= [(\pa_1\xi)(\xtil,\util(\xtil))+(\pa_2\xi)(\xtil,\util(\xtil))
\cdot \pa_\xtil\util(\xtil(x))]\cdot J_{g_\eps}(x,u\heins(x))\\
&= D_\xtil \xi(\prolo g_\eps(x,u\heins(x)))\cdot J_{g_\eps}(x,u\heins(x)).
\end{split}
\end{equation*}
Finally, using \eqref{jumf}, and the fact that $\tr(AB)=\tr(BA)$, we arrive at
\begin{equation*}
\begin{split}
\frac{d}{d\eps}(\det J_{g_\eps}(x,u\heins(x))) &= \det(J_{g_\eps}(x,u\heins(x)))\cdot \tr(D_\xtil 
\xi(\prolo g_\eps(x,u\heins(x))))\\
&= \det(J_{g_\eps}(x,u\heins(x)))\cdot (\Div \xi)(\prolo g_\eps(x,u\heins(x))).
\end{split}
\end{equation*}
\ep
The infinitesimal criterion for variational symmetries now is given by the following result.
\bt Let $G$ be a connected local transformation group acting on $M\sse \Om_0\times U$.
Then $G$ is a variational symmetry group of the variational problem \eqref{vp} if and only 
if 
\begin{equation}\label{varcrit}
\prol v(L) + L\cdot\Div\xi = 0
\end{equation}
for all $(x,\un)\in M\hn$ and every infinitesimal generator
$$
v = \sum_{i=1}^p \xi^i(x,u)\pd{}{x^i} + \sum_{\al=1}^q \phi_\al(x,u)\pd{}{u^\al}
$$
of $G$.
\et
\pr As usual, for $g\in G$ we write $(\xtil,\util)=g\cdot(x,u)=(\Xi_g(x,u),\Phi_g(x,u))$.
By definition of a variational symmetry, see \eqref{varsym}, we get
\begin{equation*}
\begin{split}
\int_\Om L(x,\prol f(x))\,dx &= \int_{\tilde\Om} L(\xtil,\prol \ftil(\xtil))\,d\xtil \\ 
&\stackrel{\eqref{trvar}}{=} \int_{\Om} L(\xtil(x),\prol \ftil(\xtil(x)))\det J_g(x,\prolo f(x))\,dx
\end{split}
\end{equation*}
Since this equality has to hold for all subdomains $\Omega$ and all functions $f$, the integrands
must agree pointwise, so $g$ is a variational symmetry if and only if
\begin{equation}\label{pointa}
L(\prol g \cdot(x,\un)) \det J_g(x,u\heins) = L(x,\un) \quad \forall (x,\un)\in M\hn
\end{equation}
for all $g$ such that both sides are defined.
We now insert $g=g_\eps=\Fl^v_\eps$ and differentiate with respect to $\eps$, 
which by \eqref{jacobform} and the product rule yields
\begin{equation}\label{jacobdiff}
(\prol v(L) + L\cdot \Div\xi)|_{(\prol g_\eps\cdot(x,\un))}\cdot \det J_{g_\eps}(x,u\heins) = 0
\end{equation}
Setting $\eps=0$, i.e., $g_\eps=\id$, this formula reduces to \eqref{varcrit}, proving necessity.

Conversely, if \eqref{varcrit} is satisfied then also \eqref{jacobdiff} holds for $\eps$ sufficiently
small. Therefore, the derivative of \eqref{pointa} (with $g=g_\eps$) with respect to $\eps$ vanishes
identically. Integrating from $0$ to $\eps$ then shows that \eqref{pointa} holds for any $g_\eps$. Since 
$G$ is connected the result then follows since any $g\in G$ can be written as a product of such
$g_\eps$.
\ep
\bex Returning to \eqref{first} (i), where $\mcl[u] = \int_a^c \sqrt{1+u_x^2}\,dx$, we 
expect that rotations should be variational symmetries, as they leave curve lengths invariant.
To verify this, let $v=-u\pa_x + x\pa_u$ be a generator of the rotation group on $\R^2$. Then
by \eqref{rotvprol},
$$
\prolo v = -u\pa_x + x\pa_u + (1+u_x^2)\pa_{u_x}.
$$
In particular $\xi=-u$, and we calculate
$$
\prolo v(L) +L\cdot D_x\xi = (1+u_x^2)\pd{}{u_x}\sqrt{1+u_x^2} - \sqrt{1+u_x^2}\cdot u_x =0,
$$
so indeed \eqref{varcrit} is satisfied.
\et
Next, we want to clarify the relationship between variational symmetries of some variational 
problem and symmetries of the corresponding Euler--Lagrange equations. To this end, we first 
need to understand how the Euler--Lagrange equations of variational problems transform under a 
change of variables.
\bp Let $L(x,\un)$, $\tilde L(\xtil,\util\hn)$ be two Lagrangians related by a change of variables
given by \eqref{chvar}, \eqref{trvar}, i.e.,
$$
\xtil = \Xi(x,u), \qquad \util = \Phi(x,u)
$$	
and
$$
L(x,\un) = \tilde L(\xtil,\util\hn)\det J(x,u\heins).
$$
Then
\begin{equation}\label{eulertrans}
E_{u^\al}(L)(x,u^{(2n)}) = \sum_{\beta=1}^q F_{\al\be}(x,u\heins)E_{\util^\be}(\tilde L)(\xtil,\util^{(2n)}), \quad \al=1,\dots,q,
\end{equation}
where $F_{\al\be}$ is the determinant of the following $(p+1)\times (p+1)$-matrix:
\begin{equation}\label{fmatrix}
F_{\al\be} = \det
\begin{pmatrix}
D_1\Xi^1 & \dots & D_p\Xi^1 & \pd{\Xi^1}{u^\al}\\
\vdots &          & \vdots & \vdots\\
D_1\Xi^p & \dots & D_p\Xi^p & \pd{\Xi^p}{u^\al}\\
D_1\Phi^\be & \dots & D_p\Phi^\be & \pd{\Phi^\be}{u^\al}\\
\end{pmatrix}.
\end{equation}
\et
\pr 
Let $(x_0,u_0^{(2n)})$ be any fixed point in $X\times U^{(2n)}$ and let $f$ be a smooth
function defined near $x_0$ with $(x_0,u_0^{(2n)})=(x_{0},\mathrm{pr}^{(2n)}f(x_0))$.
By the inverse function theorem, on an open ball $\Omega$ around $x_0$, $f$ is transformed by our change
of variables into a function $\util=\ftil(\xtil)$. 
If $\eta\in \D(\Om)^q$ and $\eps$ is sufficiently small, then
the perturbations $u_\eps = f(x,\eps) = f(x) +\eps\eta(x)$, required for calculating variational derivatives, 
are transformed into functions $\util = \ftil(\xtil,\eps)$, which are determined implicitly by
\begin{equation}\label{impldet}
\xtil = \Xi(x,f(x)+\eps\eta(x))=:\Xi_\eps(x),\quad \util=\Phi(x,f(x)+\eps\eta(x))=:\Phi_\eps(x).
\end{equation}
Again the inverse function theorem determines the domain of each $\ftil(\,.\,,\eps)$. More precisely,
by \cite[Suppl.\ 2.5 A]{AMR}, the minimal radius of a ball around $\tilde x_\eps:= \Xi(x,f(x_0)+\eps\eta(x_0))$ 
where $\ftil(\,.\,,\eps)$ is defined can be 
estimated from below in terms of the second derivatives of the maps $\Xi_\eps$ from \eqref{impldet}
in a neighborhood of $x_0$ as well as the norm of the inverse of their Jacobians at $x_0$.
Since $\eta$ is compactly supported, $\overline \Omega$ is compact, and $\xtil_\eps\to \xtil_0$ as $\eps\to 0$,
by choosing $\eps$ small 
we can therefore achieve that
all these balls contain a fixed neighborhood of $\tilde x_0$,
which we then take as a common domain $\tilde \Omega$ of all $\ftil(\,.\,,\eps)$, independently of $\eps$.
In fact, since the radius of the ball $\Omega$ around $x_0$ is itself determined
by the mentioned estimates (with $\eps=0$), we may in addition achieve that $\supp(\eta)$ 
is contained in $\Xi_\eps^{-1}(\tilde\Omega)$ for $\eps$ small. 

Now set $\tilde \eta(\xtil,\eps) := \ftil(\xtil,\eps)-\ftil(\xtil)$ and
note that $\Xi_\eps(x)=\Xi_0(x)$ as well as $\Phi_\eps(x)=\Phi_0(x)$, 
and thereby $\ftil(\Xi_\eps(x),\eps)$ $ = $ $ \ftil(\Xi_0(x))$ 
for $x\not\in \supp(\eta)$. This implies that $\tilde\eta(\tilde y,\eps)=0$ for $\tilde y
\not\in \Xi_\eps(\supp(\eta))$. Again due to the compactness of $\supp(\eta)$ this shows
that, for $\eps$ sufficiently small, all $\tilde\eta(\,.\,,\eps)$ have compact support contained in $\tilde \Omega$.
The same is therefore true of $\pd{\ftil}{\eps}\big|_{\eps=0} = \pd{\tilde\eta}{\eps}\big|_{\eps=0}$.
%

Using this, exactly as in \eqref{mcl} we obtain
\begin{equation}\label{ltifti}
\dde \tilde\mcl[\ftil] = \int_{\tilde\Om} E_\util(\tilde L)\cdot \left.\pd{\ftil}{\eps}\right|_{\eps=0}\,d\xtil,
\end{equation}
where $E_\util(\tilde L)$ is evaluated at $\util=\ftil$. To continue, we need to determine
$\pd{\ftil}{\eps}$.

To calculate $\ftil(\xtil,\eps)$, we need to solve the first equation in \eqref{impldet} for $x$ and
then insert into the second equation. Thus we have $x=x(\xtil,\eps)$ with 
$$
\xtil=\Xi(x(\xtil,\eps),f(x(\xtil,\eps))+\eps\eta(x(\xtil,\eps))).
$$
When computing variations of $\tilde L$, the base variables are not allowed to depend on $\eps$,
so for $i=1,\dots,p$ we must have
$$
\pd{\xtil^i}{\eps} = 0 = \sum_{j=1}^p D_j\Xi^i\pd{x^j}{\eps} + \sum_{\al=1}^q \pd{\Xi^i}{u^\al}\eta^\al.
$$
Using Cramer's rule we can solve for $\pd{x^j}{\eps}$ (and insert $\eps=0$), to obtain
$$
\left.\pd{x^j}{\eps}\right|_{\eps=0} = -\frac{1}{\det J}\sum_{i=1}^p K_{ij}\sum_{\al=1}^q 
\pd{\Xi^i}{u^\al}\eta^\al,
$$
where $K_{ij}$ is the $(i,j)$-cofactor (i.e., signed minor) of the Jacobian $J(x)$ from \eqref{jacob}.
Consequently,
\begin{equation*}
\begin{split}
\left.\pd{\ftil^\beta}{\eps}\right|_{\eps=0} &= \sumq \pd{\Phi^\be}{u^\al}\eta^\al + \sump D_j\Phi^\be 
\left.\pd{x^j}{\eps}\right|_{\eps=0}\\
&= \frac{1}{\det J} \sumq \Big[\pd{\Phi^\be}{u^\al}\det J - \sum_{i,j=1}^p D_j\Phi^\be\cdot 
K_{ij} \pd{\Xi^i}{u^\al}\Big]\eta^\al.
\end{split}
\end{equation*}
Here, the term in brackets is the expansion of the determinant \eqref{fmatrix} along the last column, wherein
$\sum_j D_j\Phi^\be\cdot K_{ij}$ is the row expansion of the $(i,p+1)$-st minor along its last row (recall
that the necessary signs are already contained in the $K_{ij}$). Thus we have shown that
$$
\left.\pd{\ftil^\beta}{\eps}\right|_{\eps=0} = \frac{1}{\det J}\sumq F_{\al\be}\eta^\al.
$$ 
If we insert this into \eqref{ltifti} and change variables we get
$$
\dde \tilde\mcl[\ftil] = \int_{\Om}\Big(\sum_{\al,\be=1}^q F_{\al\be}E_{\tilde u^\be}(\tilde L)\cdot \eta^\al\Big)\,dx,
$$
which, on the other hand, must equal
$$
\dde \mcl[f+\eps\eta] = \int_{\Om}\Big(\sum_{\al=1}^q E_{u^\al}(L)\cdot \eta^\al\Big)\,dx.
$$
Since $\eta$ was arbitrary, this proves that \eqref{eulertrans} holds at $(x_0,u_0^{(2n)})$.
\ep
Based on this result, we can now show:
\bt If $G$ is a variational symmetry group of the functional 
$$
\mcl[u]=\int_{\Om_0}L(x,\un)\,dx,
$$
then $G$ is also a symmetry group of the corresponding Euler--Lagrange equations $E(L)=0$.
\et
\pr Let $g\in G$ and let $f$ be a solution to the Euler--Lagrange equations. Then 
since $G$ is a variational symmetry group, \eqref{pointa} shows that
$L(\xtil,\prol \ftil(\xtil))$ and $L(x,\prol f(x))$ are related by the change of variables formula 
\eqref{trvar} with $L=\tilde L$, where $g^{-1}\cdot(x,u)=:(\Xi(x,u),\Phi(x,u))$. 
Interchanging the roles of $u$ and $\util$ in \eqref{eulertrans} we therefore have
$$
E_{\tilde u^\al}(L)(\xtil,\tilde u\hn) = \sum_{\beta=1}^q F_{\al\be}(\xtil,\tilde u\heins)
E_{u^\be}(L)(x,u\hn), \quad \al=1,\dots,q.
$$
Since $f$ is a solution to the Euler--Lagrange equations, $E_{u^\be}(L)(x,\prol f(x))=0$ for 
all $\be$ and all $x$. Thus also $E_{\tilde u^\al}(L)(\xtil,\prol \ftil(\xtil))=0$, so $g\cdot f$
is a solution as well, i.e., $G$ is a symmetry group of $E(L)=0$.
\ep
\section{Conservation laws}\index{conservation law}
\bd Let $P(x,\un)=0$ be a system of differential equations. A {\em conservation law} is a
divergence expression
\begin{equation}\label{cl}
\Div F = 0,
\end{equation}
which vanishes for all solutions $u=f(x)$ of the system. Here, for some $m$,
$$
F=(F_1(x,\um),\dots,F_p(x,\um))
$$
is a $p$-tuple of smooth functions and $\Div F$ is the total divergence from \ref{totdiv}.
$F$ then is called a conserved current.\index{conserved current}
\et
\bex (i) Let $P=\Delta(u)$, where $\Delta$ is the Laplace-operator. Then $\Delta$ itself is a conservation law because
$$
\Delta u = \Div(\grad u) = 0
$$
for every solution. Further conservation laws are obtained by multiplying the equation with
$u_i = \pd{u}{x^i}$:
$$
0 = u_i\Delta u = \sump D_j\Big(u_iu_j - \frac{1}{2}\delta^j_i\sum_{k=1}^p u_k^2\Big).
$$
(ii) If $P$ is any system of ODEs, with independent variable $x$, then any conservation law
is of the form $D_xF=0$ for all solutions $u=f(x)$ of the system, i.e., $F(x,\un)$ has to
be constant along solutions, so $F$ is a first integral\index{first integral} of the system
(cf.\ the remarks following \eqref{invode}). In light of this, \eqref{cl} can be viewed as
a generalization of the concept of a first integral to PDEs.
\et
In applications to physics, often there is a distinguished time-variable $t$, while the 
remaining variables $(x^1,\dots,x^p)$ are spatial variables. Then any conservation law
takes the form
$$
D_tT + \Div X = 0,
$$
where $\Div$ now is the spatial total divergence (with respect to $(x^1,\dots,x^p)$). 
Then $T$ is called the (conserved) {\em density}\index{density} and $X=(X_1,\dots,X_p)$ 
is called the {\em flux}\index{flux} of the conservation law. They are functions of
$x,t,u$ and the derivatives of $u$ with respect to $x$ and $t$.

Suppose now that $\Om\sse\R^p$ is a spatial domain and let $u=f(x,t)$ be a solution
defined for all $x\in \Om$ and $a\le t\le b$. Consider the functional
\begin{equation}\label{tofdef}
\tof(t) = \int_\Om T(x,t,\prolmm f(x,t))\,dx,
\end{equation}
which, for $f$, $\Om$ fixed, depends only on $t$. Then we have:
\bp $T$ and $X$ are the density and the flux of a conservation law of a system of differential equations
if and only if, for any bounded domain $\Om\sse\R^p$ with smooth boundary $\pa\Om$ and for any
solution $u=f(x,t)$ defined for $x\in\Om$ and $t\in [a,b]$ we have
\begin{equation}\label{tofcons}
\tof(t) - \tof(a) = -\int_a^t\int_{\pa\Om} X(x,\tau,\prolmm f(x,\tau))\cdot n\,dS\,d\tau,
\end{equation}
with $n$ the outward-pointing unit normal vector field on $\pa\Om$.
\et
\pr Differentiating \eqref{tofdef}, by the divergence theorem we have
\begin{equation}\label{divth}
\frac{d}{dt}\tof(t) = \int_\Om D_tT(x,t,\prolmmp f)\,dx = - \int_{\pa\Om} X(x,t,\mathrm{pr^{(m-1)}} f)\cdot n\,dS,
\end{equation}
so \eqref{tofcons} follows by integration.

Conversely, differentiating \eqref{tofcons} with respect to $t$ and applying the divergence theorem
we obtain
$$
\int_\Om (D_tT(x,t,\prolmmp f) + \Div X(x,t,\prolmmp f))\,dx = 0.
$$
As this is to hold for any $\Omega$ and $f$, the claim follows.
\ep
\bex Consider the motion of an incompressible, inviscid fluid. With $x\in \R^3$ representing
the spatial coordinates, let $u=u(x,t)\in\R^3$ denote the velocity of a fluid particle at
position $x$ and time $t$. Let $\rho(x,t)$ be the density and $p(x,t)$ be the pressure
of the fluid. We assume that the flow is isentropic (has constant entropy), then
$p$ depends only on $\rho$. The equation of continuity then takes the form
\begin{equation}\label{cons1}
\rho_t + \Div(\rho\cdot u) = 0,
\end{equation}
where $\Div(\rho u) = \sum_j \pd{(\rho u^j)}{x^j}$ is the spatial total divergence. Also, 
momentum balance gives the three equations
\begin{equation}\label{cons2}
\pd{u^i}{t} + \sum_{j=1}^3 u^j\pd{u^i}{x^j} = -\frac{1}{\rho}\pd{p}{x^i}, \quad i=1,2,3.
\end{equation}
Here, the equation of continuity is itself a conservation law with density $T=\rho$ and
flux $X=\rho u$. Using \eqref{divth} for the functional \eqref{tofdef} corresponding to the flux
we get
$$
\frac{d}{dt}\int_\Om \rho\,dx = - \int_{\pa\Om} \rho u\cdot n dS.
$$
This equation has an immediate physical interpretation: $\int_\Om\rho\,dx$ is the total mass
of the fluid within $\Om$, and $\rho u\cdot n$ is the flux of fluid out of a point on the
boundary of $\Om$. This means that the net change of mass inside $\Om$ equals the influx
of fluid into $\Omega$ via $\pa\Om$. In particular, if the normal component $u\cdot n$
of the velocity on $\pa \Om$ vanishes then there is no change in the mass in $\Om$
and we obtain the law of conservation of mass:
$$
\int_\Om \rho\,dx = \text{ const}.
$$
Furthermore, combining \eqref{cons1} and \eqref{cons2} and re-arranging we get
three further conservation laws:
$$
D_t(\rho u^i) +\sum_{j=1}^3 D_j(\rho u^iu^j + p\delta^j_i) = 0,\quad i=1,2,3.
$$
Again using \eqref{divth} we obtain the laws of conservation of momentum:
$$
\frac{d}{dt}\int_\Om \rho u^i\,dx = - \int_{\pa\Om}(\rho u^i(u\cdot n) + p n_i)\,dS, \quad i=1,2,3.
$$
Here, the first term on the right hand side is the transport of momentum $\rho u^i$ due to the
flow across the surface $\pa\Om$, and the second term is the net change in momentum due to the
pressure across $\pa\Om$. Thus $X_j=\rho u^iu^j + p\delta^j_i$ represents the components of
the momentum flux.
\et
\section{Noether's theorem}\index{Noether's theorem}
In this section we prove a fundamental result on the connection between variational symmetries
and conservation laws, established by Emmy Noether in 1918. It is of central importance in 
mathematical physics as the source of many fundamental laws of nature.
\bt\label{noetherth} Let $G$ be a local one-parameter group of variational symmetries of the variational problem
$\mcl[u]=\int L(x,\un)\,dx$. Let
\begin{equation}\label{vfform2}
v = \sum_{i=1}^p \xi^i(x,u)\pd{}{x^i} + \sum_{\alpha=1}^q \phi_\alpha(x,u)\pd{}{u^\alpha}
\end{equation}
be the infinitesimal generator of $G$, and as in \eqref{chardef} let
$$
Q_\alpha(x,u\heins) := \phi_\alpha(x,u) - \sum_{i=1}^p \xi^i(x,u)u^\alpha_i, \quad \alpha=1,\dots,q.
$$
be the characteristic of $v$. 
Then there exists a $p$-tuple 
$$
F(x,u^{(m)})=(F_1(x,u^{(m)}),\dots,F_p(x,u^{(m)}))
$$ 
such that 
\begin{equation}\label{noether}
\Div F = Q\cdot E(L) = \sumq Q_\al\cdot E_\al(L).
\end{equation} 
Thus $\Div F=0$ is a conservation law 
for the Euler--Lagrange equations $E(L)=0$.
\et
\pr Since $v$ is the infinitesimal generator of a variational symmetry, by \eqref{varcrit}
and \eqref{charprol} we have
\begin{equation*}
\begin{split}
0 &= \prol v(L) + L\,\Div\xi = \prol v_Q(L) +\sum_{i=1}^p\xi^iD_iL + L\,\Div\xi\\
&= \prol v_Q(L) +\Div(L\xi),
\end{split}
\end{equation*}
where $L\xi:=(L\xi^1,\dots,L\xi^p)$. We now `integrate by parts' the first term,
using \eqref{charparts}:
\begin{equation*}
\begin{split}
\prol v_Q(L) &= \sum_{\al,J}D_JQ_\al \pd{L}{u^\al_J} = \sum_{\al,J} Q_\al\cdot (-D)_J 
\pd{L}{u^\al_J} +\Div A\\
&= \sumq Q_\al E_\al(L) + \Div A,
\end{split}
\end{equation*}
where $A=(A_1,\dots,A_p)$ is some $p$-tuple of functions depending on $Q$, $L$ and their
derivatives, but whose concrete form is not important for us. It follows that
\begin{equation}\label{noether1}
\prol v_Q(L) = Q\cdot E(L)+\Div A.
\end{equation}
Therefore, 
\begin{equation}\label{noether2}
0=Q\cdot E(L) + \Div(A+L\xi),
\end{equation}
which means that \eqref{noether} holds with $F=-(A+L\xi)$.
\ep
To explicitly calculate $F$ one may in principle follow the calculation of the previous
proof. In applications, in most cases one encounters first order variational problems.
For these, an explicit formula is given in the following result:
\bc Let $\mcl[u]=\int L(x,u\heins)\,dx$ be a first order variational problem, and let
$v$ as in \eqref{vfform2} be a variational symmetry. Then
\begin{equation}\label{expnoether}
F_i = \sumq \phi_\al \pd{L}{u^\al_i} + \xi^iL -\sumq\sump \xi^j u^\al_j\pd{L}{u^\al_i},
\quad i=1,\dots,p
\end{equation}
form the components of a conservation law $\Div F = 0$ for the corresponding 
Euler--Lagrange equations $E(L)=0$.
\et
\pr Since $L$ depends only on derivatives up to order one, \eqref{charparts} gives
$$
\prolo v_Q(L) = \sumq \Big( Q_\al\pd{L}{u^\al} + \sum_{i=1}^p D_iQ_\al\pd{L}{u^\al_i}\Big).
$$
Now setting $A_i:=\sumq Q_\al \pd{L}{u^\al_i}$ we have
$$
\Div A = \sum_{\al,i} \Big( D_iQ_\al\pd{L}{u^\al_i} + \sum_{\al,i} Q_\al D_i\pd{L}{u^\al_i}\Big).
$$
On the other hand, \eqref{eulerop} implies
$$
Q\cdot E(L) = \sumq Q_\al \pd{L}{u^\al} -\sumq Q_\al D_i \pd{L}{u^\al_i},
$$
so altogether we obtain
$$
\prolo v_Q(L) = Q\cdot E(L)+\Div A,
$$
verifying \eqref{noether1}.
\ep
\bex\label{partmech} Perhaps the most prominent application of the Noether theorem is to the 
mechanics of particles. Consider a system of $n$ particles moving in $\R^3$
in some force field given by a potential. The system is described by the
position of its particles, where $\bfx^\alpha = (x^\al,y^\al,z^\al)$ is the
position of the $\al$-th particle, and $\bfx=(\bfx^1,\dots,\bfx^n)$ is a
vector containing the information about all positions of the particles.
Assuming that the $\al$-th particle has mass $m_\al$, the total kinetic
energy of the system is
$$
K(\dot\bfx) = \frac{1}{2}\sum_{\al=1}^n m_\al |\dot \bfx^\al|^2.
$$
Also, we assume that the force field is determined by a potential $U=U(\bfx,t)$.
Newton's equation of motion then gives
$$
m_\al\ddot \bfx^\al = -\nabla_\al U = -(U_{x^\al},U_{y^\al},U_{z^\al}), \quad \al=1,\dots,n.
$$
The important point to note here is that these equations are the Euler--Lagrange equations of 
the Lagrange function (or action integral)
\begin{equation}\label{lagrangian}
\mcl[\bfx]=\int_{-\infty}^\infty (K-U)\,dt.
\end{equation}
Indeed,
$$
E_{x^\al}(L) = \pd{L}{x^\al} - D_t\pd{L}{\dot x^\al} = -U_{x^\al} - m_\al \ddot x^\al,
$$
and analogously for $y^\al$ and $z^\al$.

By \eqref{varcrit}, a vector field
$$
v = \tau(t,\bfx)\pd{}{t} + \sum_\al {\boldsymbol\xi}^\al(t,\bfx)\cdot \pd{}{\bfx^\al} \equiv
\tau\pd{}{t} +\sum_\al \Big(\xi^\al \pd{}{x^\al} + \eta^\al \pd{}{y^\al} + \zeta^\al \pd{}{z^\al} \Big)
$$
generates a variational symmetry of the Lagrangian \eqref{lagrangian} if and only if
\begin{equation}\label{newtoncrit}
\prolo v(K-U) + (K-U)D_t\tau = 0 \quad \forall (t,\bfx).
\end{equation}
Given a variational symmetry, by Noether's theorem \ref{noether} we obtain a corresponding
conservation law (or first integral) with single component (since $p=1$)
\begin{equation}\label{newcons}
\begin{split}
F &= \sum_{\al=1}^n m_\al {\boldsymbol\xi}^\al\cdot \dot\bfx^\al + \tau (K-U) -
\sum_{\al=1}^n \tau \dot \bfx^\al\cdot (m_\al \dot\bfx^\al)\\
&= \sum_{\al=1}^n m_\al {\boldsymbol\xi}^\al\cdot \dot\bfx^\al -\tau (K+U)
= \sum_{\al=1}^n m_\al {\boldsymbol\xi}^\al\cdot \dot\bfx^\al -\tau E,
\end{split}
\end{equation}
where $E=K+U$ is the total energy. Since $D_tF=0$ it follows that $F$ has to be constant
for any solution of Newton's equations of motion.
We now analyze some examples of variational symmetries of \eqref{lagrangian} that
lead to conservation laws of physical interest.

\begin{itemize}
\item $v=\pa_t$: In this case, $\prolo v = v$, so \eqref{newtoncrit} holds if and only
$\pa_t U = 0$, i.e., if and only if $U$ does not explicitly depend on the time $t$.
Since $\tau=1$ and ${\boldsymbol \xi}=0$, the conserved quantity given in \eqref{newcons} then 
is the total energy $E$. We obtain that {\em invariance of a physical system under time translations
implies conservation of energy}.\index{conservation!of energy}
\item $v=\sum_{\al=1}^n \bfa\cdot\pd{}{\bfx^\al}$: This is the generator of the translation $\bfx\mapsto 
\bfx+\bfa$, so all particles are simultaneously translated in the same fixed direction $\bfa\in \R^3$.
Also in this case, $\prolo v = v$, so \eqref{newtoncrit} holds if and only $v(U)=0$, i.e., if and only
if the potential is translationally invariant in the direction $\bfa$. As $\tau=0$ and ${\boldsymbol \xi}^\al
=\bfa$ for all $\al$, the corresponding conserved
quantity from \eqref{newcons} is the momentum
$$
\sum_{\al=1}^n m_\al \bfa\cdot \dot\bfx^\al = \text{ const}.
$$
In particular, if $U$ is invariant under all translations then the total momentum $\sum_{\al=1}^n m_\al 
\dot\bfx^\al$ is conserved.\index{conservation!of momentum}
\item $v = \sum_{\al=1}^n \big(x^\al\pd{}{y^\al} - y^\al\pd{}{x^\al}\big)$: This is the generator
of a simultaneous rotation of all the masses in the system about some fixed axis, in this case
the $z$-axis. By \eqref{phiform} we have
$$
\prolo v = v +\sum_{\al=1}^n \Big(\dot x^\al\pd{}{\dot y^\al} - \dot y^\al\pd{}{\dot x^\al}\Big).
$$
Note that $\prolo v(K)=0$ (since $K$ is invariant under rotations of velocities), and since
$\tau=0$, \eqref{newtoncrit} shows that $v$ generates a variational symmetry if and only if
$v(U)=0$, i.e., if and only if $U$ is invariant under rotations around the $z$-axis. In this 
case, the conserved quantity corresponding to this symmetry via \eqref{newcons} is
the angular momentum around the $z$-axis
$$
\sum_{\al=1}^n m_\al (x^\al \dot y^\al - y^\al \dot x^\al) = \text{ const.}
$$
Thus {\em rotational invariance implies conservation of angular momentum}.\index{conservation!of angular momentum}
\end{itemize}
In particular, if we assume that the particles only interact through their mutual gravitational
(or electrostatic, \dots) attraction, then the potential energy is of the form
$$
U(t,\bfx) = \sum_{\al\not=\be} \gamma_{\al\be} |\bfx^\al-\bfx^\be|^{-1}.
$$
In this case, all of the above assumptions are satisfied, so we obtain 
conservation of energy, momentum, and angular momentum.
\et
To conclude this chapter we introduce a straightforward generalization of Noether's theorem that
also leads to conservation laws of physical interest.
\bd Let $\mcl[u]=\int L(x,\un)\,dx$ be a variational problem. A vector field $v$ on $M\sse X\times U$
is called an {\em infinitesimal divergence symmetry}\index{infinitesimal divergence symmetry} of $\mcl$
if there exists a $p$-tuple $B(x,u^{(m)}) = (B_1,\dots,B_p)$ of functions of $x$, $u$ and the derivatives
of $u$ such that
\begin{equation}\label{divsym}
\prol v(L) + L\cdot \Div\xi = \Div B
\end{equation}
for all $(x,u)\in M$.
\et
The point to note is that the Noether theorem \ref{noetherth} remains valid if in the hypothesis 
we replace `variational symmetry' by `divergence symmetry'. Indeed, the only thing that has to be changed
in the proof is that we have to incorporate the new term $\Div B$ into \eqref{noether2}, which now becomes
$$
Q\cdot E(L) + \Div(A+L\xi) = \Div B,
$$
but we still obtain a conservation law of the form \eqref{noether}, namely $F=B-A-L \xi$.
In particular, in the case of a variational problem of first order, the conserved quantity
corresponding to \eqref{expnoether} now becomes
\begin{equation}\label{expnoether2}
F_i = \sumq \phi_\al \pd{L}{u^\al_i} + \xi^iL -\sumq\sump \xi^j u^\al_j\pd{L}{u^\al_i} - B_i.
\quad i=1,\dots,p
\end{equation}
\bex Returning to the setup of \ref{partmech}, consider a Galilean boost
$$
(t,\bfx^\al)\mapsto (t,\bfx^\al + \eps t\bfa),
$$
for some $\bfa\in\R^3$. By\eqref{phiform}, the infinitesimal generator 
$v=\sum_{\al=1}^n t\bfa \pa_{\bfx^\al}$ then has prolongation
$$
\prolo v = \sum_{\al=1}^n \Big(t\bfa \pd{}{\bfx^\al} + \bfa \pd{}{\dot\bfx^\al}\Big),
$$
so
$$
\prolo v(L) = \prolo(K-U) = \sum_{\al=1}^n m_\al \bfa\cdot \dot\bfx^\al - t \sum_{\al=1}^n \bfa\cdot \nabla_\al U.
$$
Since $\tau = 0$, \eqref{newtoncrit} demands that this expression should vanish identically, but this is
never the case for $\bfa\not=0$. However, since
$$
\sum_{\al=1}^n m_\al \bfa\cdot \dot\bfx^\al = D_t\Big(\sum_{\al=1}^n m_\al \bfa\cdot\bfx^\al \Big),
$$
$v$ is an infinitesimal divergence symmetry if $\bfa\cdot \nabla_\al U=0$, i.e., if $U$ is translationally
invariant in the direction $\bfa$. The conserved quantity according to \eqref{expnoether2} now reads
$$
t \sum_{\al=1}^n m_\al \bfa\cdot\dot\bfx^\al -
\sum_{\al=1}^n m_\al \bfa\cdot\bfx^\al.
$$
Here the first sum, when divided by the total mass $\sum_\al m_\al$ is the position of the center
of mass of the system in the direction $\bfa$, and the second is the momentum in that direction.
It follows that if $U$ is translationally invariant in a given direction then on the one hand,
as shown in \ref{partmech}, the momentum in that direction is constant. In addition it now 
follows that the center of mass in that direction is a linear function of $t$:
$$
\text{position of center of mass =} \ \text{initial position} + t(\text{momentum})/\text{mass}.
$$
In particular, if $U$ is invariant under all translations, then the center of mass of any such system
moves linearly in one fixed direction.
\et
\section{Trivial conservation laws and characteristics}
There are certain types of conservation laws that do not yield any useful information on the
system under consideration and hence are called trivial.\index{conservation law!trivial} 
A conservation law $F=(F_1,\dots,F_p)$ can be trivial for one of two reasons: triviality of the
first kind holds if $F$ itself vanishes on every solution of the system. This kind of 
triviality\index{trivial conservation law!first kind} 
is usually easy to eliminate by solving the system itself and its prolongations
for certain of the variables $u^\al_J$ and then substituting for these variables wherever they
occur. For example, in the case of an evolution equation $u_t=P(x,\un)$, any time derivative
of $u$, e.g., $u_{tt},\,u_{xt}$, etc.\ can be expressed in terms of $x$, $u$ and spatial derivatives
of $u$. 
\bex\label{heatc} Consider the system of first order evolution equations
$$
u_t=v_x,\quad v_t=u_x,
$$
which is equivalent to the one-dimensional wave equation $u_{tt}=u_{xx}$. Then
$$
D_t\Big(\frac{1}{2}u_t^2+\frac{1}{2}u_x^2\Big) - D_x(u_tu_x) = u_t(u_{tt}-u_{xx}) = 0
$$
is a conservation law for this system. As explained above, we can replace the density and
the flux of this conservation law by ones that depend only on spatial derivatives. This gives
$$
D_t\Big(\frac{1}{2}u_x^2+\frac{1}{2}v_x^2\Big) - D_x(u_xv_x) = 0.
$$
The two conservation laws differ by the trivial conservation law
$$
D_t\Big(\frac{1}{2}u_t^2-\frac{1}{2}v_x^2\Big) + D_x(v_xu_x-u_tu_x) =  0,
$$
for which both density and flux vanish on any solution of the system. In the same way,
for any conservation law of an evolution equation there is, up to addition of a trivial
conservation law of the first kind, one where density and flux depend only on 
spatial derivatives.
\et
A conservation law is called trivial of the second kind\index{trivial conservation law!second kind} 
if the total divergence $\Div F$
in fact vanishes on {\em all} smooth functions $f$, whether or not they solve the equation.
An example of this kind of triviality is the relation
\begin{equation}\label{triv}
D_xu_y - D_yu_x=0,
\end{equation}
which obviously holds for any smooth function $u=f(x,y)$. Hence \eqref{triv} is a trivial 
conservation law of the second kind for any partial differential equation for functions
of $x$ and $y$. The underlying $p$-tuples $(F_1,\dots,F_p)$ of conservation laws of the second 
kind are also called null divergences.\index{null divergence}

There is in fact a complete characterization of null divergences that is a direct analogue
of the characterization of the kernel of the usual divergence operator via the Poincar\'e
lemma. Since the proof of this result, building on the so-called variational complex,
is quite involved, we only state the result and refer to \cite[Sec.\ 5.4]{O} for a proof.
\bt\label{tcc} Let $F=(F_1,\dots,F_p)$ be a $p$-tuple of smooth functions depending on $x=(x^1,\dots,x^p)$,
$(u^1,\dots,u^q)$ and derivatives of $u$, and defined on all of $X\times U\hn$.
Then the following are equivalent:
\begin{itemize}
\item[(i)] $\Div F \equiv 0$.
\item[(ii)] There exist smooth functions $G_{jk}$, $j,k=1,\dots,p$, depending on $x$, $u$,
and derivatives of $u$, such that $G_{jk}=-G_{kj}$ for all $j$, $k$, and 
\begin{equation}\label{tcce}
F_j = \sum_{k=1}^p D_k G_{jk}, \quad j=1,\dots,p
\end{equation}
\end{itemize}
\et
For example, if $p=3$, then 
$$
\Div F = D_1F_1+D_2F_2+D_3F_3\equiv 0
$$
if and only if $F$ is a `total curl': $F=\mathrm{Curl}(G)$, i.e.,
$$
F_1 = D_2G_3 - D_3G_2,\quad F_2 = D_3G_1 - D_1G_3,\quad F_3 = D_1G_2 - D_2G_1
$$
(where we have identified $G_{12}=-G_{21}$ from the theorem with $G_3$, etc.).

We define a trivial conservation law to be any conservation law that is a linear 
combination of trivial conservation laws of the first and second kind. Thus
$\Div F = 0$ is a trivial conservation law if and only if there exist functions
$G_{ij}$ as in  \ref{tcc} (ii) such that \eqref{tcce} holds for all solutions
of the system. Two conservation laws are called {\em equivalent} if they
differ only by a trivial conservation law.\index{conservation law!trivial} 
The interesting objects in the study
of conservation laws therefore are equivalence classes in this sense.

For the following considerations on characteristics of conservation laws we
need some preparations.
\bd Let $P_\nu(x,\un)=0$, $\nu=1,\dots,l$, be a system of differential equations,
with $P:M\hn\to \R^l$. The $k$-th prolongation of this system is the $(n+k)$-th order 
system of differential equations
$$
P\hk(x,u^{(n+k)}) = 0
$$
obtained by differentiating in all possible ways $k$ times. Thus the new system consists
of the 
$\begin{pmatrix}
p+k-1\\ k	
\end{pmatrix}\cdot l
$ equations 
$$
D_J P_\nu(x,u^{(n+k)}) = 0,
$$
where $\nu=1,\dots,l$ and $0\le\sharp J \le k$.
\et
\bex If $P$ is the heat equation $u_t=u_{xx}$, then the first prolongation $P\heins$
is the system 
$$
u_t=u_{xx},\quad u_{xt} = u_{xxx}, \quad u_{tt}=u_{xxt},
$$
and the second one contains, in addition, the equations
$$
u_{xxt}=u_{xxxx},\quad u_{xtt} = u_{xxxt}, \quad u_{ttt}=u_{xxtt}.
$$
\et
We then clearly have:
\blem\label{sowieso} If $u=f(x)$ is a solution of the system $P(x,\un)=0$, then it is also
a solution of every prolongation $P\hk(x,u^{(n+k)})=0$, $k=0,1,2,\dots$.
\et
\bd\label{conspr} A system of differential equations is called {\em totally nondegenerate}\index{nondegenerate!totally}
if it and all its prolongations are nondegenerate (i.e., of maximal rank and locally solvable).
\et
\blem\label{consprlem} Let $P_\nu(x,\un)=0$, $\nu=1,\dots,l$, be a totally nondegenerate system of differential
equations. Let $Q=Q(x,u\hm)$ be a smooth function. Then the following are equivalent:
\begin{itemize}
\item[(i)] $Q$ vanishes for all solutions $u=f(x)$ of the system.
\item[(ii)] There exist differential operators $\D_\nu = \sum_J Q^J_\nu(x,\um)D_J$,
$\nu=1,\dots,l$, such that 
$$
Q =\sum_{\nu=1}^l \D_\nu P_\nu
$$ 
for all $(x,\um)$.
\end{itemize}
\et
\pr (ii)$\Rightarrow$(i) is obvious.

(i)$\Rightarrow$(ii): 
Clearly we can assume that $m\ge n$. By \ref{sowieso}, (i) is equivalent to $Q$ vanishing 
on all solutions of the prolonged system to $P^{(m-n)}$, which is locally solvable by assumption. 
Thus it follows that $Q$ vanishes on the zero set of $P^{(m-n)}$. Since, moreover, $P^{(m-n)}$ is nondegenerate, 
by \ref{hadamard} there exist
smooth functions $Q^J_\nu$ ($\sharp J\le m-n$, $1\le \nu\le l$) such that
$$
Q(x,\um) = \sum_J \sum_\nu Q^J_\nu(x,\um)
D_JP_\nu(x,\um)
$$
\ep
Let now $\Div F=0$ be a conservation law of a totally nondegenerate system
$$
P(x,\un)=0
$$ 
of differential equations. Then by \ref{consprlem}, $\Div F$ vanishes on every solution of
the system if and only if there exist smooth functions $Q^J_\nu(x,u\hm)$ such that
\begin{equation}\label{divfg}
\Div F = \sum_{\nu,J} Q^J_\nu D_JP_\nu.
\end{equation}
We now note that each term in \eqref{divfg} can be `integrated by parts': for example,
if $1\le j\le p$ then
$$
Q^j_\nu D_j P_\nu = D_j(Q^j_\nu P_\nu) - D_j(Q^j_\nu)P_\nu.
$$
Proceeding in this way, we can cast \eqref{divfg} in the form
$$
\Div F = \Div G + \sum_{\nu=1}^l Q_\nu P_\nu \equiv \Div G + Q\cdot P,
$$
where $Q=(Q_1,\dots,Q_l)$ has entries
\begin{equation}\label{qe}
Q_\nu = \sum_J (-D)_J Q^J_\nu,
\end{equation}
and all we need to know about $G=(G_1,\dots,G_p)$ is that it depends linearly on the 
components $P_\nu$ of $P$ and their total derivatives. This means that $G$ defines a 
trivial conservation law of the first kind. Consequently, replacing $F$ by $F-G$
we obtain an equivalent conservation law of the form
\begin{equation}\label{cc}
\Div F = Q\cdot P.
\end{equation}
This is called the {\em characteristic form}\index{conservation law!characteristic form} 
of the conservation law \eqref{divfg}, and $Q=(Q_1,\dots,Q_l)$ is called the characteristic
of the conservation law.\index{conservation law!characteristic of}

Unless $l=1$, the $Q_\nu$ from \eqref{cc} are in general not uniquely determined. If
$Q$ and $\tilde Q$ are two $l$-tuples both satisfying \eqref{cc} for the same $F$,
then $(Q-\tilde Q)\cdot P = 0$. Since $P$ is nondegenerate, by \ref{hadex} it follows
from this that $Q-\tilde Q$ vanishes on all solutions. Based on this observation, we
call a characteristic {\em trivial}\index{characteristic!trivial} 
if it vanishes on all solutions of the system,
and we call two characteristics equivalent if they differ by a trivial one. In this
terminology, characteristics are in general only determined up to equivalence.
\bex To obtain the characteristic form of the conservation law for the heat equation
given in \ref{heatc}, we re-write it in the form \eqref{divfg}:
$$
D_t(\frac{1}{2}u_t^2+\frac{1}{2}u_x^2) - D_x(u_tu_x) = u_tD_t(u_t-v_x) + u_tD_x(v_t-u_x).
$$
Then \eqref{qe} shows that the characteristic is given by
$$
Q = (-D_t(u_t),-D_x(u_t)) = (-u_{tt},-u_{xt}),
$$
and an equivalent conservation law in characteristic form can be found through integration
by parts:
$$
D_t\Big(\frac{1}{2}u_t^2+\frac{1}{2}u_x^2\Big) + D_x(-u_tu_x) \cong -u_{tt}(u_t-v_x) - u_{xt}(v_t-u_x).
$$
\et
One can show that for systems of differential equations that are totally nondegenerate
and normal (i.e., possessing a noncharacteristic direction at every point) the 
two notions of equivalence we have introduced above actually coincide:
\bt Let $P(x,\un)=0$ be a normal and totally nondegenerate system of differential equations.
Let the $p$-tuples $F$ and $\tilde F$ determine conservation laws with characteristics
$Q$ and $\tilde Q$, respectively. Then $F$ and $\tilde F$ are equivalent as conservation
laws if and only if $Q$ and $\tilde Q$ are equivalent as characteristics.
\et 
For a proof we refer to \cite[Sec.\ 4.3]{O}.

Finally, we may apply the terminology developed in this chapter to re-formulate Noether's 
theorem \ref{noetherth} more precisely:
\bt\label{noetherth2} Let $G$ be a local one-parameter group of variational symmetries of the variational problem
$\mcl[u]=\int L(x,\un)\,dx$. Let
\begin{equation*}\label{vfform3}
v = \sum_{i=1}^p \xi^i(x,u)\pd{}{x^i} + \sum_{\alpha=1}^q \phi_\alpha(x,u)\pd{}{u^\alpha}
\end{equation*}
be the infinitesimal generator of $G$, and as in \eqref{chardef} let
$$
Q_\alpha(x,u\heins) := \phi_\alpha(x,u) - \sum_{i=1}^p \xi^i(x,u)u^\alpha_i, \quad \alpha=1,\dots,q.
$$
be the characteristic of $v$. 
Then $Q=(Q_1,\dots,Q_q)$ is also the characteristic of a conservation law
for the corresponding Euler--Lagrange equations $E(L)=0$, i.e.,  
there exists a $p$-tuple 
$$
F(x,u^{(m)})=(F_1(x,u^{(m)}),\dots,F_p(x,u^{(m)}))
$$ 
such that 
\begin{equation*}\label{noether3}
\Div F = Q\cdot E(L) = \sumq Q_\al\cdot E_\al(L)
\end{equation*} 
is a conservation law in characteristic form 
for the Euler--Lagrange equations 
$$
E(L)=0.
$$
\et

\newpage
\addcontentsline{toc}{part}{Index}

\printindex 

\end{document}